\pdfoutput=1
\documentclass[a4paper,11pt]{article}
\usepackage[english]{babel}
\usepackage{amsfonts}
\usepackage{amsthm}
\usepackage{amsmath}
\usepackage{amssymb}
\usepackage{enumerate}
\usepackage{mathtools}
\usepackage{tikz}
\usepackage{stmaryrd}
\expandafter\def\csname opt@stmaryrd.sty\endcsname
{only,shortleftarrow,shortrightarrow}
\usepackage{tikz-cd}
\usepackage[all]{xy}
\usepackage{multirow}
\usepackage{graphicx}
\usepackage{amsmath}
\usepackage{adjustbox}
\usepackage{graphicx}
\usepackage{extpfeil}
\usepackage{comment}
\usepackage{mathabx}
\usepackage{float}
\usepackage[labelformat=empty]{caption}
\usepackage[toc]{appendix}
\usepackage[left=3.2cm,top=3cm,right=3.2cm,bottom=3cm]{geometry}
\usepackage{resizegather}
\usepackage[activate={true, nocompatibility}, final, tracking=true, kerning=true, spacing=true, factor=1100, stretch=10, shrink=10]{microtype}
\usepackage{dsfont}
\usetikzlibrary{matrix}
\usepackage{hyperref}
\hypersetup{
colorlinks=true,
linkcolor=blue,
filecolor=magenta,      
urlcolor=cyan,
citecolor=[RGB]{0,204,0},
}
\usepackage[nameinlink]{cleveref}
\usepackage[normalem]{ulem}
\usepackage{sidecap}
\usepackage{mathrsfs}
\usepackage{mathscinet}

\usepackage{color} %CLASHES with hyperref
\definecolor{darkgreen}{rgb}{0.0, 0.7, 0.0}
\newenvironment{??}{\noindent \color{darkgreen}{\bf ???:} \footnotesize}{}
\definecolor{cyan}{cmyk}{1,0,0,0}

\newcommand{\cb}{\color{blue}}

\newcommand{\cvio}{\color{violet}}
\newcommand{\bdg}{\begin{dg}}

\theoremstyle{plain}
\newtheorem{thm}{Theorem}[section]
\newtheorem{coroll}[thm]{Corollary}
\newtheorem{defn}[thm]{Definition}
\newtheorem{lemma}[thm]{Lemma}
\newtheorem*{claim}{Claim}

\newtheorem{example}[thm]{Example}
\newtheorem{prop}[thm]{Proposition}

\newtheorem{remark}[thm]{Remark}

\newtheorem*{thm*}{Main Theorem}
\newtheorem{notn}[thm]{Notation}
\newtheorem{context}[thm]{Context}
\newtheorem{introthm}{Theorem}

\tikzset{
  symbol/.style={
    draw=none,
    every to/.append style={
      edge node={node [sloped, allow upside down, auto=false]{$#1$}}}
  }
}
\microtypecontext{spacing=nonfrench}

\DeclareMathOperator{\Sym}{Sym}

\newcommand{\fp}{\mathfrak{p}}
\newcommand{\fg}{\mathfrak{g}}
\newcommand{\ft}{\mathfrak{t}}
\newcommand{\fc}{\mathfrak{c}}

\DeclareMathOperator{\colim}{colim}

\DeclareMathOperator{\rank}{rank}

\DeclareMathOperator{\ad}{ad}

\newcommand{\Stab}{\on{Stab}}
\newcommand{\act}{\on{act}}
\newcommand{\Lie}{\on{Lie}}

\newcommand{\codim}{\mathrm{codim}}
\newcommand{\schP}{P}
\newcommand{\stP}{\mathscr{P}}

\newcommand{\DD}{\mathfrak{D}}

\newcommand{\cO}{\mathcal{O}}

\newcommand{\cB}{\mathcal{B}}
\newcommand{\cC}{\mathcal{C}}
\newcommand{\cD}{\mathcal{D}}
\newcommand{\cE}{\mathcal{E}}
\newcommand{\cF}{\mathcal{F}}
\newcommand{\cG}{\mathcal{G}}

\newcommand{\cK}{\mathcal{K}}
\newcommand{\cH}{\mathcal{H}}
\newcommand{\cM}{\mathcal{M}}

\newcommand{\cS}{\mathcal{S}}

\newcommand{\cZ}{\mathcal{Z}}
\newcommand{\cL}{\mathcal{L}}

\newcommand{\bC}{\mathbb{C}}
\newcommand{\bG}{\mathbb{G}}

\newcommand{\bQ}{\mathbb{Q}}

\newcommand{\bZ}{\mathbb{Z}}
\newcommand{\bT}{\mathbb{T}}

\newcommand{\scrX}{\mathscr{X}}
\newcommand{\scrM}{\mathscr{M}}
\newcommand{\scrY}{\mathscr{Y}}
\newcommand{\scrZ}{\mathscr{Z}}

\newcommand{\on}{\operatorname}

\newcommand{\Aut}{ \on{Aut} }

\newcommand{\elli}{\on{ell}}
\newcommand{\reg}{ \on{reg}}

\newcommand{\Hom}{ \on{Hom}}
\newcommand{\Map}{ \on{Hom}}

\renewcommand{\ker}{ \on{ker}}

\newcommand{\Spec}{\on{Spec}}

\newcommand{\uSpec}{\underline{\on{Spec}}}

\newcommand{\Higgs}{\mathscr{H}}
\newcommand{\Higgsrig}{\mathcal{H}}
\newcommand{\MHiggs}{M}
\newcommand{\nilp}{\on{Nilp}}
\newcommand{\Bun}{ \mathscr{B}un} 
\newcommand{\Pic}{\on{Pic}}
\newcommand{\GL}{\on{GL}}

\newcommand{\hitchinbase}{A}

\newcommand{\ic}[1]{\mathcal{IC}_{#1}}
\newcommand{\st}[1]{\mathfrak{#1}}
\usepackage{color} %CLASHES with hyperref
\definecolor{darkgreen}{rgb}{0.0, 0.7, 0.0}

\begin{document}
\title{\textbf{Hitchin fibrations are Ng\^{o} fibrations}}
% \title{\textbf{The Decomposition Theorem
% \\ for logarithmic $G$-Higgs bundles, I:\\
% Ng\^o fibration structure}}
\author{Mark Andrea de Cataldo, Roberto Fringuelli, \\ Andres Fernandez Herrero and  Mirko Mauri}
\date{}
\maketitle
\begin{abstract}
We study the geometry of the Hitchin fibration for $\cL$-valued $G$-Higgs bundles over a smooth projective curve of genus $g$, where $G$ is a reductive group and $\cL$ is a suitably positive line bundle. 
We show that the Hitchin fibration admits the structure of a 
weak Abelian fibration. 
In the case when the line bundle $\cL$ is a twist of the canonical bundle of the curve by a (possibly empty) reduced effective divisor, we prove a cohomological bound and $\delta$-regularity of the weak Abelian fibration.
\end{abstract}

\tableofcontents

\begin{section}{Introduction}

The rich geometric properties of the moduli spaces $M_{G,\cL}$  of $G$-Higgs bundles over a curve for a reductive group $G$  twisted by a line bundle $\cL$
 have been studied from different points of view: representation theory, Teichm\"{u}ller theory, mirror symmetry, and the Langlands programs. 
 The moduli spaces 
 $M_{G,\cL}$  played  a crucial tool in Ng\^o's proof of the Fundamental Lemma \cite{ngo-lemme-fondamental}, and also in the subsequent work of Chaudouard--Laumon 
 on the Weighted Fundamental Lemma \cite{cl-lemme-pondere-I, cl-lemme-pondere-II}.
 Furthermore, the (intersection) cohomology of the Higgs moduli space underlies interesting phenomena such as topological mirror symmetry \cite{hausel-thaddeus-04, GWZ20, Endo, Mauri2021}, degree independence \cite{maulik-shen-independence}, and the $P=W$ conjecture \cite{dCHM2012, dCMS, PW_proof1, HMMS, maulik2023perversefiltrationsfouriertransforms}. %\andres{add references}{\Md I will add [dHM], [dMS], [Ms], [HMMS], [MSQ]}. 
 Currently, except for the case when $G$ is of type $A_n$, there is a rather limited understanding of these moduli spaces and related phenomena.

The goal of this article is to place a first stepping stone towards a better understanding of the cohomology of the Higgs moduli spaces $M_{G,\cL}$ by proving certain strong structural properties of the  Hitchin fibration. To this end, we introduce the notion of Ng\^o fibration 
(cf. \Cref{defn: ngo fibration intro}), which is a strengthening of Ng\^o's notion of weak abelian fibration.
Ng\^o fibrations are defined as to enjoy several symmetries 
which have been considered previously for the Hitchin morphism  in 
\cite{ngo-lemme-fondamental},  \cite{chaudouard-laumon, cl-lemme-pondere-I, cl-lemme-pondere-II},  \cite{de-cataldo-support-sln} and \cite{maulik-shen-independence}.

Before introducing our result, let us fix some terminology and notation.

Let $C$ be a smooth geometrically connected projective curve over a field $k$. We fix a line bundle $\cL$ on $C$ and a split connected reductive group $G$. Recall that an $\cL$-valued $G$-Higgs bundle on $C$ is a pair $(E, \phi)$ consisting of a $G$-bundle $E$ on $C$ and a global section $\phi$   of the $\cL$-twisted adjoint vector bundle $\ad(E) \otimes \cL$ on $C$.
 Canonical Higgs bundles  correspond to the case
 $\cL=\omega_C$ (canonical bundle); logarithmic Higgs bundles  correspond to the case $\cL= \omega_C(D)$ where $D$ is a reduced effective divisor on $C$. The Higgs moduli space $M_{G,\cL}$ is the adequate moduli space of $\cL$-valued $G$-Higgs bundles subject to the standard semistability condition. The connected components $M^d_{G, \cL} \subset M_{G,\cL}$ are quasiprojective varieties and they are indexed by their degree, which is an element $d \in \pi_1(G)$ of the  fundamental group $\pi_1(G)$  of the reductive group $G$.  %the elements $d \in \pi_1(G)$ of the fundamental group $\pi_1(G)$.

A very useful tool to approach the geometry of the Higgs moduli space is the affinization morphism $h \colon M^{d}_{G, \cL} \to A$, which is known as the Hitchin fibration. It plays a prominent role in the study of the cohomology of $M_{G, \cL}^d$ through the study of the BBDG Decomposition Theorem \cite{bbdg} for the direct image $Rh_*IC_{M_{G, \cL}^d}$ of the %($\ell$-adic) 
($\overline{\mathbb{Q}}_ \ell$-adic) intersection complex on $M_{G,\cL}^d$. For Ng\^o's purposes in \cite[\S 7.8]{ngo-lemme-fondamental}, it is sufficient to understand $Rh_*IC_{M_{G, \cL}^d}$ over a certain open set of the Hitchin base $A$ over which $M_{G,\cL}^d$ has at worst quotient singularities, the singularities of the Hitchin fibration are better understood, and the Hitchin fibration is automatically an Ng\^o fibration. 
%; see \cite[\S 7.8]{ngo-lemme-fondamental}.}

However, for the purposes of understanding the global cohomology of $M_{G,\cL}^d$, it is necessary to also study the more singular fibers and understand the Decomposition Theorem over the whole Hitchin base. 
We remark that this has been accomplished in the case of $G= \GL_n$ by Chaudouard--Laumon \cite{chaudouard-laumon} (degree coprime to the rank) and by Maulik--Shen \cite{maulik-shen-independence} (arbitrary degree). Indeed, they proved that when the degree of $\cL$ satisfies $\deg(\cL) > \deg(\omega_C)$, where $\omega_C$ is the canonical bundle, then all of the simple direct summands of $Rh_*IC_{M_{G, \cL}^d}$ have full support. This is an essential ingredient  in some of the phenomena  for $G=\GL_n$ discussed above, such as the degree independence \cite{maulik-shen-independence} and  two of the proofs of $P=W$ \cite{PW_proof1,maulik2023perversefiltrationsfouriertransforms}.

 The main result of this paper is \begin{itemize}
    \item \Cref{thm: thmA}: the Hitchin fibration for canonical and for 
logarithmic Higgs bundles over a curve over a field of characteristic zero  is an Ng\^o fibration. In particular,
 when the ground field is also algebraically closed, the direct image $Rh_*IC_{M_{G, \cL}^d}$ decomposes in Ng\^{o} strings.

\end{itemize}
Theorems B through F below
establish the necessary preliminary geometric and arithmetic results.
We note that for these results, there are varying assumptions of the ground field and on the line bundle $\cL$.
\begin{itemize}
\item \Cref{thmB}: properties of the stack of Higgs bundles and of the stacky Hitchin morphism and its fibers. 
\item \Cref{thmC}: properties of the stacky connected components, of the stacky stable and semistable loci, of the adequate moduli spaces, of the Hitchin morphism and  of its fibers.
\item \Cref{ThmD}:
the Hitchin morphism is a weak abelian fibration. 
\item \Cref{thm: thmE}:
the direct image of the intersection complex via the good moduli space
morphism for a quotient stack admits the intersection complex of the good moduli space as a direct summand (this is key to the  cohomological bound
for the Hitchin fibration). 
\item \Cref{thmF}: the $\delta$-regularity property for the group of symmetries of the Hitchin fibration. 
\end{itemize}

\subsubsection*{The Decomposition Theorem for Ng\^o fibrations}\label{sub: ngo fibration}

In his proof of the Fundamental Lemma \cite{ngo-lemme-fondamental}, Ng\^o introduced the notion of weak abelian fibration $(X, B, P)$ (cf. \Cref{defn: weak Abelian fibration})  and proved the Support Inequality Theorem for such fibrations. 
Ng\^o also considered the
$\delta$-regularity property for the group scheme $P$.

Let us discuss briefly the Support Inequality and the $\delta$-regularity property
in the context of Ng\^o $\delta$-regular weak abelian fibrations.
The identity components $P_b^o$ of the geometric fibers $P_b$ of the smooth group scheme $P/B$ 
are extensions of abelian varieties $A_b$ by affine connected groups schemes $R_b$. The dimension $\delta (b)$ of this affine part  gives rise to  an upper semicontinuos function $\delta$ on the base $B$. For every closed integral subvariety $Z$, the generic value of $\delta$ on $Z$ is denoted by $\delta (Z)$. The Ng\^o Support inequality predicts that if $Z$ is a support for a simple shifted perverse summand in the Decomposition Theorem for a weak abelian fibration, then ${\rm codim}_B(Z) \leq \delta (Z)$. On the other hand,
the $\delta$-regularity property states that, for every integral subvariety $Z$ in $B$, we have the opposite inequality ${\rm codim}_B(Z) \geq \delta (Z)$.
The combination of the Ng\^o Support Inequality and the $\delta$-regularity property %taken 
%together 
ensures that the dimension of every support of any simple direct summand in the Decomposition Theorem for a $\delta$-regular weak abelian fibration is determined by the structural properties of the group scheme (not of the total space of the fibration).  In turn, this  gives rise to the notion of Ng\^o strings (cf. \cite{dCHM2021, dCRS2021}) as the direct summands, supported on a given support, in the Decomposition Theorem for the $\delta$-regular weak abelian fibration. 
In particular,  $\delta$-regularity ensures that each support
has an explicit codimension and contributes precisely one Ng\^o string with an explicit cohomological shift. 
As a result, the intersection cohomology of the total space $X$ is governed by the group scheme $P$ and the highest direct image $R^{top}h_*IC_{X}$; see \Cref{prop: Ngodecompositiontheorem}.

Ng\^o's original notion of  weak abelian fibrations requires
the total space $X$  of the fibration to be rationally smooth (i.e. the intersection complex to be the shifted constant sheaf). Ng\^o proved that the Hitchin fibration is a weak abelian fibration when restricted over a suitable open subset of the Hitchin base. In general, the Higgs moduli space is not rationally smooth over the whole Hitchin base, not even in the case when $G=\GL_n$.   Moreover, it is not known in which general context $\delta$-regularity holds, neither for weak abelian fibration, nor for Hitchin fibrations ($\delta$-regularity holds in the canonical and logarithmic cases, see \Cref{thm: thmA}.2). Ng\^o works on  suitable open subsets of auxiliary moduli spaces where, by definition, the $\delta$-regularity property holds, and then he shows that this suffices for the purposes of proving the Fundamental Lemma. 

In this paper, we follow \cite{maulik-shen-independence} and we define a Ng\^o fibration to be a weak abelian fibration,  whose group scheme is $\delta$-regular over the whole base,  and whose total space, while not necessarily rationally smooth, satisfies a certain cohomological bound, needed to prove the analogue of the Ng\^o Support Inequality. 
The definition on Ng\^o fibration  arises by axiomatizing some of the relevant properties of the Hitchin fibration, following \cite{ngo-lemme-fondamental} and \cite{maulik-shen-independence}.

\begin{defn}[Ng\^o fibration] 
\label{defn: ngo fibration intro} 
Let $B$ be a scheme of finite type over $k$.  Let $h \colon X \to B$ be a quasiprojective scheme, proper over $B$, of relative dimension $n$, equipped with an action of a smooth commutative group scheme $P \to B$ of finite type over $B$.
Then the triple $(X, P, B)$ is an Ng\^o fibration if
\begin{itemize}
    \item (weak Abelian) the triple $(X, P, B)$ is a weak Abelian fibration (see \Cref{defn: weak Abelian fibration});
    \item ($\delta$-regularity) the group scheme $P \to B$ is $\delta$-regular (see \Cref{defn: delta regularity});
    \item (cohomological bound) for the standard  truncation functor $\tau_{> *}(-)$, we have
    \[\tau_{> 2n} Rh_* (IC_X[-\dim(X)])=0,\]
    where $IC_X$ is the perverse intersection complex on $X$.
\end{itemize}
\end{defn}

 The notion of Ng\^o fibration is  crafted so that it yields a rather precise form of the Decomposition Theorem (see \Cref{prop: Ngodecompositiontheorem}). We introduce some necessary notation.
Let $(X,P,B)$ be a Ng\^o fibration.
Denote by $I$ the set of supports $Z_{\sigma} \subseteq B$ of the simple direct summands of the direct image $R h_*IC_{X}$. Up to taking  finite inseparable extensions of the field of rational functions of $Z_\sigma$,  there exist open dense subsets $Z^{\times}_{\sigma}$ such that the restriction $P_{\sigma}$ of the group scheme $P$ admits a Chevalley d\'evissage (cf.\ \Cref{T:sCd-kclosed})
\[0 
\to R_{\sigma} \to P_{\sigma} \to A_{\sigma} \to 0,
\]
where $R_{\sigma} \to Z^{\times}_{\sigma}$ is an affine commutative smooth group scheme with geometrically connected fibers of relative dimension $\delta_{\sigma}:= \dim(R_\sigma/Z^{\times}_\sigma)$ and $g_{\sigma} \colon A_{\sigma} \to Z^{\times}_{\sigma}$ is an Abelian scheme of relative dimension  $\delta^{\mathrm{ab}}_{\sigma}:= \dim(A_\sigma/Z^{\times}_\sigma)$. We introduce the lisse sheaves  $\Lambda^{i}_{\sigma} \coloneqq R^i g_{\sigma *} \overline{\mathbb{Q}}_ {\ell} \, \simeq 
\bigwedge^i ( \Lambda^1_\sigma)$ on $Z^\times_\sigma$.
In the following proposition, the Tate shift can be ignored in the constructible case.

%Denote by $I$ the set of supports $Z_{\sigma} \subseteq A$ of the simple direct summands of the dirct image $R h_*IC_{\MHiggs_{G,\omega_C(D)}^{d}}$. There exist open dense subsets $Z^{\times}_{\sigma}$ such that the restriction $P_{\sigma}$ of the group scheme $P$ admits a Chevalley d\'evissage (cf.\ \Cref{T:sCd-kclosed})
%\[0 
%\to R_{\sigma} \to P_{\sigma} \to A_{\sigma} \to 0,
%\]
%where $R_{\sigma} \to Z^{\times}_{\sigma}$ is an affine commutative smooth group scheme with geometrically connected fibers of relative dimension $\delta_{\sigma}:= \dim(R_\sigma/Z^{\times}_\sigma)$ and $g_{\sigma} \colon A_{\sigma} \to Z^{\times}_{\sigma}$ is an Abelian scheme of relative dimension  $\delta^{\mathrm{ab}}_{\sigma}:= \dim(A_\sigma/Z^{\times}_\sigma)$. Set $\Lambda^{i}_{\sigma} \coloneqq R^i g_{\sigma *} \overline{\mathbb{Q}}_ {\ell} \, \simeq 
%\bigwedge^i ( \Lambda^1_\sigma)$.

%{\Md See his paper, p. 113, bottom. The Ngo strings were defined and made explicit by Jochen, Luca and I; Jochen told me that obviously that is what Ngo intended; I am not  so sure; what is clear, is that the proof of the support theorem makes them evident.}

\begin{prop}[Decomposition Theorem for Ng\^o fibrations]\label{prop: Ngodecompositiontheorem}
Let $(X,P,B)$ be a Ng\^o fibration as in \Cref{defn: ngo fibration intro} over an algebraically closed field $k$.
There exists an isomorphism in the constructible bounded derived category $D^b_c(B, \overline{\mathbb{Q}}_ \ell)$ 
%in $D^bMHM_{\text{alg}}(A)$ (resp. $D^b(A, \overline{\mathbb{Q}}_ \ell)$) 
\begin{equation}\label{eq:decompositiontheoreme}
   R h_*IC_{X} \simeq \bigoplus_{\sigma \in I} \left\{ \bigoplus^{2 \delta^{\mathrm{ab}}_{\sigma}}_{i=0} IC_{Z_{\sigma}} (\Lambda^{i}_{\sigma} \otimes L^d_{\sigma})[-i+\delta_{\sigma}^{\rm ab}]  (-\delta_\sigma) \right\}, 
\end{equation}
where $L^d_{\sigma}$ is a  lisse sheaf on an open and dense subscheme of $Z_{\sigma}$, and $\codim(Z_{\sigma}) = \delta_{\sigma}$.  

Suppose in addition that the field $k$ is the complex numbers $\mathbb{C}$.  
There is an isomorpism analogous to (\ref{eq:decompositiontheoreme}) of pure objects in the bounded derived category $D^b (MHM_{\text{alg}}(B))$ of algebraic mixed Hodge modules over the Hitchin base $A$, where $L^d_{\sigma}$ is a polarizable variation of pure Hodge structures of weight zero.
\end{prop}

\noindent The summand in curly brackets labeled by $\sigma$  in \eqref{eq:decompositiontheoreme} is called the Ng\^o string associated with the support $Z_\sigma$.

\begin{proof}
The proof (in the constructible case, as well as in the case of mixed Hodge modules) runs parallel to the proof of
\cite[Thm. A.0.3]{dCRS2021}, esp.\ of equation (111) in loc.\ cit.  The only difference is that in loc\ cit., weak abelian fibrations are assumed to have total space a rational homology manifold, so that the Ng\^o Support Theorem applies; on the other hand, 
our definition of Ng\^o fibration includes the cohomological bound, in the presence of which the conclusion of the Ng\^o Support Theorem remains valid if we replace the constant sheaf on the total space with its intersection complex; see \cite[\S0.2 and \S1]{maulik-shen-independence} (what is called here cohomological bound is called relative dimension bound in loc.\ cit.)
\end{proof}

\subsubsection*{Main result on Hitchin fibrations}

The main theorem of this paper is the following.

\begin{introthm}[The Hitchin fibration is a Ng\^o fibration]\label{thm: thmA} 
    Let $C$ be a smooth projective geometrically connected curve over a field $k$ of characteristic $0$. Let $D \subset C$ be a (possibly trivial) reduced effective Cartier divisor with $\deg(\omega_C(D))>0$. Fix a connected split reductive group $G$, and a degree $d \in \pi_1(G)$. Let $\MHiggs_{G,\omega_C(D)}^{d}$ be the moduli space of $G$-Higgs bundles of degree $d$ with simple poles at $D$. Then  
   \begin{enumerate}[(1)]
    \item The Hitchin fibration $h: \MHiggs_{G,\omega_C(D)}^{d} 
        \to A$ (\Cref{notn: higgs moduli space and hitchin fibration}) 
        is equipped with the action of a quasiprojective group scheme 
        $\schP^{\circ} \to A$ as in \Cref{coroll: action on the 
        semistable moduli space}.
        
        \item The triple $ (\MHiggs_{G,\omega_C(D)}^{d}, P^{\circ}, A)$ is a Ng\^o fibration. %{\cdg with respect to the %($\overline{\mathbb Q}_\ell$-adic) intersection complex $IC_{\MHiggs_{G,\omega_C(D)}^{d}}$}.
        
        \item  Assume in addition that the ground  field $k$ is algebraically closed. The direct image $Rh_* IC_{\MHiggs_{G,\omega_C(D)}^d}$ of the intersection complex splits into a direct sum of Ng\^o strings as in 
        \Cref{prop: Ngodecompositiontheorem}.
        
    \end{enumerate}
    % Then, the Hitchin fibration $h: \MHiggs_{G,\omega_C(D)}^{d} \to A$ (\Cref{notn: higgs moduli space and hitchin fibration}) for the moduli of degree $d$ $G$-Higgs bundles with simple poles at $D$ equipped with the action of the group scheme $\schP^{\circ} \to A$ of the Hitchin fibration (as in \Cref{coroll: action on the semistable moduli space}) is an Ng\^o fibration with respect to the %($\overline{\mathbb Q}_\ell$-adic) 
    % intersection complex $IC_{\MHiggs_{G,\omega_C(D)}^{d}}$.
\end{introthm}
\begin{proof}
The first two parts of the theorem follows by combining several of our results in the main body of the text (\Cref{thm: weak Abelian fibration structure}, \Cref{thm: ngo inequality supports}, \Cref{cor:relativedimcanonical} and \Cref{thm:deltaregularmero}). In view of \Cref{prop: Ngodecompositiontheorem}, the third part follows from the second. 
\end{proof}

\begin{remark} Here, we  discuss the remaining cases where $D \subset C$ is a (possibly trivial) reduced effective Cartier divisor with $\deg(\omega_C(D))\leq 0$. Denote by $g$ be the genus of $C$.
    \begin{itemize}
        \item If $g=0$ and $\deg(D) = 0 \text{ or } 1$, then $\MHiggs_{G,\omega_C(D)}^{d}$ is either empty or a point.
        \item If $g=0$ and $\deg(D) = 2$, then either $\MHiggs^d_{G,\omega_C(D)}$ is empty, or the Hitchin morphism  $h \colon \MHiggs^d_{G,\omega_C(D)} \to A$ is isomorphic to the identity $A \to A$.
        \item Let $g=1$ and $D=0$ and fix $d \in \pi_1(G)$.  By \cite[Prop. 4.1]{FGPN2019}, there exists an integer $r(d) \leq {\rm rank}(G)$ and a commutative square
        \[
\xymatrix{
C^{r(d)} \times H^0(C, \omega_C)^{r(d)} \ar[r] \ar[d]^-{p} & \MHiggs^d_{G,\omega_C} \ar[d]^-{h}\\
H^0(C, \omega_C)^{r (d)} \ar[r]& A, 
} 
        \]
        where the horizontal arrows are finite, the top horizontal arrow is surjective, and the left vertical arrow is the smooth projection onto the second factor.  
        % Up to a finite modification, the Hitchin fibration $h$ can be identified with the smooth projection $p \colon (C \times H^0(C, \omega_C))^{r} \to H^0(C, \omega_C)^{r}$ with $r$ at most $\rank(G)$. 
        Note that: the Hitchin fibration is not surjective onto $A$ when $r (d) < \rank(G)$; it is an open question whether $\MHiggs^d_{G,\omega_C}$ is normal, unless $G=\GL_{n}$ or $\mathrm{SL}_n$ (see \cite[last paragraph of \S 1]{FGPN2019}). In particular, it is unclear whether $(\MHiggs_{G,\omega_C}^{d}, P^{\circ}, h(\MHiggs_{G,\omega_C}^{d}))$ is an Ng\^o fibration. 
        
        On the positive side, $h$ is genercally a fibration in abelian varieties  onto its image, and all the simple direct summands of $Rh_* IC_{\MHiggs_{G,\omega_C}^d}$ are fully supported on $h(\MHiggs_{G,\omega_C}^d)$, as this holds for the morphism $p$. Therefore, the Decomposition Theorem for the Hitchin fibration takes the form \eqref{eq:decompositiontheoreme} in \Cref{prop: Ngodecompositiontheorem} in this context as well.
    \end{itemize}
\end{remark}

\begin{remark}\label{precise dt}
In the forthcoming paper \cite{DT_paper_2}, by building on the results of this paper, we determine explicitly the Ng\^o strings  in the Decomposition Theorem for the Hitchin fibration in the logarithmic case, where the reduced effective divisor $D>0$ is not empty:  we describe the supports and the coefficients in terms of the theory of reductive algebraic groups.
In future work, we plan to continue this line of investigation for the canonical Higgs bundles, i.e. the case when $D=\emptyset$. 
At present, we can prove \Cref{{prop:intro}}: (in addition to being a Ng\^o fibration) the Hitchin morphism is a (singular) Lagrangian fibration in the sense of \Cref{defn:symplecticvariety}.
\end{remark}

\subsubsection*{Comparison with the known special case when $G=\GL_n$}
\Cref{thm: thmA} is known in the case when $G=\GL_n$; see e.g. \cite{maulik-shen-independence}.
This paper is dedicated to the proof for an arbitrary split connected reductive group $G$.  There are several differences with the special case when $G = \GL_n$:
\begin{itemize}
\item When $G = \GL_n$, there is  an alternative description of the Hitchin fibration as a relative compactified Jacobian for the universal family of spectral curves, and the relative Jacobian acts naturally on it as a group of symmetries. %In the case of an arbitrary reductive group, while there is the notion of cameral curve, the study of the group of symmetries becomes more intricate.

 In the case of an arbitrary reductive group, while there is the notion of the universal family of cameral curves, it is no longer true that their relative Jacobian acts naturally on the Hitchin fibration. Hence, the study of the group of symmetries becomes more intricate.
 %\sout{it is no longer true that their Jacobian acts naturally on the Hitchin fibration. Hence, the study of the group of symmetries becomes more intricate.}
 
 \item When $G = \GL_n$, there is always a smooth irreducible component of the Higgs moduli space (the ones with coprime rank and degree). The existence of such a component is a key ingredient in the known proofs of $\delta$-regularity. 
 For other groups $G$, such a smooth component may fail to exist, thus raising the question of whether $\delta$-regularity holds for arbitrary reductive groups.

% \Roberto{For the moduli space of  $G$-bundles over a complex curve of genus at least two and $\Lie(G)$ contains a simple ideal not of type $A_n$ (i.e. $G^{ad}$ contains a simple normal subgroup not isomorphic to $PGLn$ for some $n$), any connected component contains strictly semistable bundles (see \cite[Prop. 7.8]{ramanathan-stable}). If you take any of those strictly semistable bundles with zero Higgs field, you should get a singular point on $\MHiggs_{G,\cL}$ (I guess).}\Mirko{This is a good remark for the paper on singularities}
\item In \cite{DT_paper_2}, we show that, in the case of an arbitrary reductive group $G$, the shape of the Decomposition Theorem for the direct image $Rh_* IC_{M^d_{G, \omega_C(D)}}$ may depend on the element
$d \in \pi_1(G)$. This is in stark contrast with the case of $G= \mathrm {GL}_n$, where one has independence of $d \in \pi_1(\GL_n)=\mathbb Z$ (cf. \cite{maulik-shen-independence}).
\end{itemize}

\subsubsection*{Foundational results on the Hitchin fibration}
In \Cref{section: stack of higgs bundles} of this paper, we give a complete treatement of several general structural properties of the Higgs moduli space and of the Hitchin fibration. Our general setup for \Cref{section: stack of higgs bundles} is as follows. We fix a split connected reductive group $G$ over a field $k$, a smooth projective geometrically connected curve $C$ of genus $g$ over $k$, and a line bundle $\cL$ on $C$. We fix the choice of a maximal split torus $T$. We denote by $Z_G$ the center subgroup scheme of $G$. We denote by $\Higgs_{G,\cL}$ the moduli stack of $\cL$-valued $G$-Higgs bundles on $C$.

\begin{introthm}[\Cref{P:Higgs-prop} + \Cref{P:Higgs-integral} + \Cref{P:Higgs-normal}]\label{thmB}
    Suppose that $\deg(\cL) \geq  \max\{2g-2,1\}$, and that the characteristic of $k$ is either zero or strictly larger than the height of the adjoint representation of $G$. Then, the following hold:
    \begin{enumerate}[(1)]
        \item The connected components of the stack $\Higgs_{G,\cL}$ are indexed by the elements of the fundamental group $\pi_1(G)$. 
        
        \item For any given $d \in \pi_1(G)$, the corresponding connected component $\Higgs_{G,\cL}^d\subset \Higgs_{G,\cL}$ is a geometrically integral and geometrically normal local complete intersection stack of dimension     
        $$   \dim(\Higgs_{G,\cL}^{d})=\begin{cases}
                \dim(G)(2g-2)+\dim (Z_G),&\text{if }\cL=\omega_C;\\
                \dim(G)\deg(\cL),&\text{otherwise}.
            \end{cases}
            $$

        \item The Hitchin morphism $h:\Higgs_{G,\cL}^d \to A$ is a syntomic morphism, i.e. a flat morphism with local complete intersection fibers  which are equidimensional of dimension 
        \[
        \frac{1}{2} \dim(G)\deg(\cL) - \frac{1}{2}\dim(T) (\deg(\cL) - 2g+2).
        \]
    \end{enumerate}
\end{introthm}

In \Cref{thmC}, we also establish  properties of the moduli space of semistable $G$-Higgs bundles. 
We expect these properties to be   well-known to experts. On the other hand,  since a complete treatment seems to be lacking in the literature,  we have included  it  in the hope that it will be useful to  the community working on the moduli of Higgs bundles.

\begin{introthm}[\Cref{thm: adequate moduli spaces for higgs} + \Cref{cor: higgs moduli space properties} + \Cref{prop: fibers hitchin pure dimensional} + \Cref{prop: smoothness of the semistable hitchin stack}]\label{thmC}
    Suppose that the characteristic of $k$ is either zero, or strictly larger than the height of the adjoint representation of $G$.  
    For every $d \in \pi_1(G)$, the locus of semistable $G$-Higgs bundles $\Higgs_{G,\cL}^{d,ss} \subset \Higgs_{G,\cL}^d$ is an open and dense substack which admits a quasiprojective adequate moduli space $ \MHiggs_{G,\cL}^d$. %$\Higgs_{G,\cL}^{d,ss} \to \MHiggs_{G,\cL}^d$. 
    Furthermore, if in addition $\deg(\cL)\geq \max\{2g-2,1\}$, then the following hold:
    \begin{enumerate}[(1)]
        \item The locus $\Higgs_{G,\cL}^{d,s}$ of stable $G$-Higgs bundles is open and dense.
        \item The moduli space $\MHiggs_{G,\cL}^d$ is a geometrically irreducible and geometrically normal quasiprojective variety of dimension              $$
            \dim(\MHiggs_{G,\cL}^{d})=\begin{cases}
                \dim(G)(2g-2)+2\dim(Z_G),&\text{if }\cL=\omega_C;\\
                \dim(G)\deg(\cL)+ \dim(Z_G),&\text{otherwise}.
            \end{cases}
            $$
        \item The induced Hitchin fibration $h: \MHiggs_{G,\cL}^d \to A$ is proper and has equidimensional fibers of dimension
        \begin{equation}\label{dim fibers Hitchin fibration}
        \frac{1}{2} \dim(G)\deg(\cL) - \frac{1}{2}\dim(T) (\deg(\cL) - 2g+2)
        +\dim(Z_G).
        \end{equation}
        If the characteristic of $k$ is zero, then $h$ is flat.
        \item If furthermore $\deg(\cL)>2g-2$, then the stack $\Higgs_{G,\cL}^{d,ss}$ is smooth.
    \end{enumerate}
\end{introthm}

\subsubsection*{Weak abelian fibration structure}

In \Cref{section: group scheme of symmetries of the hitchin fibration}, by keeping the same setup and notation as in \Cref{section: stack of higgs bundles}, we establish the weak Abelian fibration structure of the moduli space of Higgs bundles, which we discuss next. 

It is expected that any Lagrangian fibration $f \colon X \to B$ carries an action of a group scheme $P \to B$ extending the Hamiltonian automorphisms on the smooth locus of $f$, promoting the triple $(X, P, B)$ to a weak Abelian fibration. When $X$ and $B$ are smooth varieties, this is proved in increasing degree of generality by Arnold--Liouville \cite[\S 49]{Arnold1989} for smooth fibrations (not necessarily holomorphic), by Arinkin--Fedorov \cite[Thm. 2]{AR2016} for algebraic fibrations with integral fibers, by Markushevich \cite[Thm. 2.1]{Markushevich1996} and Kim \cite[Cor. 4.28]{YK2024} for holomorphic (resp.\ algebraic) fibrations with no nowhere-reduced fibers. Note that, in the last case, $X$ can be taken to be a singular symplectic variety; see \cite[Prop. 5.7]{Sacca2024}. 

The Hitchin fibrations $h: \MHiggs_{G,\cL}^d \to A$ are expected to be families of Lagrangian fibrations, i.e., there exists a linear map $ev_D \colon A \to A_D$ such that for any closed point $r \in \mathrm{Im}(ev_D)$ the restriction $h \colon h^{-1}(A(r)) \to A(r) \coloneqq ev_D ^{-1}(r)$ is a Lagrangian fibration with respect to a standard symplectic form; cf. \cite[Main Thm. and Cor. 8.10]{Markmann94} at least for $G=\GL_n$ and the same is expected for arbitrary reductive groups. 
%{\Md Mirko, this next comment of yours is not about fixing this sentence, but about what would need to be done in the future? Right? \Mirko{Right, and an explanation for us why I write only it is expected but I cannot say with certainty. Actually it should come from the argument using parabolic that I sketched in a old email.}}
%\Mirko{For this statement we should prove that $h^{-1}(A(r))$ has rational singularities. Markmann uses a weaker notion of Lagrangian fibration though.} 
Unfortunately, the results mentioned in the previous paragraph cannot be applied in this context, since the Hitchin fibration $h$ may have multiple fibers, e.g., \cite[Prop. (34) and (35)]{Thaddeus89} (except when $\cL$ admits a square root and $d=0$, in which case the existence of the Kostant--Hitchin section \cite[\S 4.2.4]{ngo-lemme-fondamental} forces all fibers of the Hitchin fibration to be reduced at some point).
%\Mirko{Implicitly our proof reduces to this case, so to Markushevich and Yoonjoo, but our group scheme is explicit and we have a prescribed action of it even when we have non-reduced fibers. I wonder what is the class of Lagrangian fibrations we have covered; all the definitions I could think seems a bit ad hoc. Maybe it is pointless since we expect it in any case} 
However, the group stack $\stP\to A$ of symmetries for the Hitchin fibration from the stack of Higgs bundles $\Higgs_{G,\cL}$ is constructed by Ng\^{o} on the whole Hitchin base, regardless of the presence of non-reduced fibers, using the modular description of $\Higgs_{G,\cL}$. We refine this result as follows.

 Suppose that the characteristic of the ground field $k$ does not divide the order of the Weyl group $W$ of $G$ and that $\deg(\cL)>0$. Let $P^{\circ} \to A$ be the rigidification (by the center of $G$) of the open substack of $\stP$ with geometrically connected fibers; see \Cref{P:centr-group}. Then we show that the group scheme $P^{\circ} \to A$  acts fiberwise on the moduli space by preserving each connected component $\MHiggs_{G,\cL}^d$, and that it is a quasiprojective group scheme,  thus endowing the triple  $(\MHiggs_{G,\cL}^d, P^o,A)$ with the structure of weak Abelian fibration; see \Cref{P:centr-group}. 

% More precisely, in \Cref{section: group scheme of symmetries of the hitchin fibration} we keep the same setup and notation with $C$, $\cL$ and $G$ as in \Cref{section: stack of higgs bundles}. Suppose in addition that the characteristic of the ground field $k$ does not divide the order of the Weyl group $W$ of $G$ and that $\deg(\cL)>0$. %Then, in \Cref{P:centr-group} we define the group scheme $\schP^o \to A$ over $A$ of symmetries of the Hitchin fibration  over $A$. This is a smooth commutative group scheme with connected fibers of dimension as in \Cref{dim fibers Hitchin fibration}. 
% The main result in \Cref{section: group scheme of symmetries of the hitchin fibration} is the following.
\begin{introthm}[Weak abelian fibration structure: \Cref{thm: weak Abelian fibration structure}]\label{ThmD}
    Suppose that $\deg(\cL)\geq \max\{2g-2,1\}$ and that the characteristic of $k$ is either zero or strictly larger than the height of the adjoint representation of $G$. The group scheme $\schP^o \to A$ is quasiprojective and acts fiberwise on the moduli space $\MHiggs_{G,\cL}^d \to A$. The triple $(\MHiggs_{G,\cL}^d, \schP^o, A)$ is a weak Abelian fibration in the sense of \Cref{defn: weak Abelian fibration}.
\end{introthm}

% Note further that the locus of non-critical points for the Hitchin fibration is a torsor under the group scheme REF, which fails in general for arbitrary Lagrangian fibration; see \cite[Prop. 1.7]{YK2024} and \cite[Example 5.6]{Sacca2024}.

\subsubsection*{ Cohomological bound: the  case
$\deg (\cL) > \max \{ 2g-2 ,  0\}$}

In order to obtain the cohomological bound appearing in the definition of Ng\^o fibration (\Cref{defn: ngo fibration intro}), Maulik and Shen \cite{maulik-shen-independence} use some ideas of Meinhardt \cite{Meinhardt} (cf. also \cite[\S 3.3]{maulik-shen-independence} and the references therein) that cannot be naturally generalized to the case of reductive groups. %{\Md I suggest to change the following sentnece to the gree one that follows.} Recently, Kinjo \cite{kinjo_decomposition_theorem} has used the recently developed theory of mixed Hodge modules on stacks to prove a version of the Decomposition Theorem for good moduli space morphisms. 
Recently, Kinjo \cite{kinjo_decomposition_theorem}  proved the purity of the direct image of the intersection complex via a good moduli space morphism.

In this paper, we have concurrently obtained a weaker version of
Kinjo's result, i.e. that the intersection complex of the good moduli space is a direct summand of the pushforward of the constant sheaf from the stack; the proof of this fact is fairly elementary and relies on Kirwan's desingularization technique.  As it is pointed out following \Cref{thm: thmE}, this implies the desired cohomological bound.

We work with various versions of the constructible derived categories of $\overline{\mathbb{Q}}_ \ell$-adic sheaves on stacks of finite type and of the derived categories of algebraic mixed Hodge modules;
see \Cref{section: intersection cohomology of gms}
and \cite{tubach_mhm_paper, kinjo_decomposition_theorem}, respectively.

\begin{introthm}[\Cref{thm: intersection cohomology of the stack vs the moduli space}]\label{thm: thmE} Let $k$ be an algebraically closed ground field, and let $\ell$ be a prime distinct from the characteristic of $k$. Let $\scrM=[X/G]$ be a quotient stack of a finite type algebraic space $X$ over $k$ quotiented by an affine group $k$-scheme $G$. Suppose that $\scrM$ admits a properly stable good moduli space $\pi: \scrM \to M$ in the sense of \cite[Defn. 2.5]{edidin-rydh-canonical-reduction}. Then, the derived pushforward $R\pi_*(IC_{\scrM})$ of the $\overline{\mathbb{Q}}_{\ell}$-adic intersection cohomology complex
of $\scrM$  admits the intersection complex $IC_{M}$ as a direct summand, i.e. there is an isomorphism  $R\pi_*(IC_{\scrM}) \simeq IC_M \oplus E$ in the bounded below constructible derived category
$D^+_c(M, \overline{\mathbb{Q}}_ \ell)$ for some complex $E$. If the field $k = \mathbb C$, then the analogous statement holds  in the bounded below derived category of algebraic mixed Hodge modules $D^+ (MHM_{\rm alg} (M))$.
\end{introthm}

We use this theorem, together  with the smoothness of the stack of semistable $\cL$-valued $G$-Higgs bundles when  $\deg (\cL) > \max \{ 2g-2 ,  0\}$,  and obtain the required cohomological bound in \Cref{thm: ngo inequality supports}. We remark that,  at present, for that particular result,  we need to impose the assumption that the characteristic of the ground field $k$ is $0$.

\subsubsection*{Cohomological bound: the canonical case}

The canonical  case $\cL = \omega_C$ (i.e. the logarithmic case with $D=\emptyset$) is of particular interest, among all possible choices of line bundle $\cL$. For instance, the stack $\Higgs_{G,\omega_C}$ is isomorphic to the cotangent bundle of the stack of $G$-principal  bundles. Furthermore, via the non-abelian Hodge correspondence,
the moduli space $\MHiggs^0_{G,\omega_C}$, with $k=\bC$, is homeomorphic to the character variety parametrizing semisimple $G$-representations of the fundamental group of $C$, .

When $\cL = \omega_C$, the stack of semistable $\cL$-valued $G$-Higgs bundles is not necessarily smooth, so the proof of the cohomological bound in \Cref{thm: ngo inequality supports} no longer applies. In \Cref{cor:relativedimcanonical}, we recover the bound in this context by observing that the Hitchin fibration is Lagrangian in the sense of \Cref{defn:symplecticvariety}, and then by  invoking \cite[Prop. 3.3 and 3.6]{MM2024}. 

\begin{prop}[\Cref{prop:LagrangianfibrationHitchin}]\label{prop:intro}
 Suppose that $\text{char}(k)=0$. For any given $d \in \pi_1(G)$, $\MHiggs_{G,\omega_C}^{d}$ is a symplectic variety, and the Hitchin fibration $h:\MHiggs_{G,\omega_C}^{d} \to A$ is a Lagrangian fibration.
\end{prop}

\Cref{prop:intro} essentially amounts to establishing the rationality of the singularities of the moduli space $\MHiggs_{G, \omega_C}$. We refer to \Cref{rmk:literature} for a comparison with the known results in the literature.

\subsubsection*{$\delta$-regularity}
In \Cref{section: delta regularity}, we employ stack-theoretic techniques to prove the $\delta$-regularity of the group scheme of symmetries of the Hitchin fibration.

 A weak Abelian Lagrangian fibration $(X, P, A)$ with smooth total space $X$ is $\delta$-regular: this is sketched in \cite[\S 2]{Ngo2011} and proved in \cite[Prop. 8.9]{AR2016} or \cite[Prop. 2.3.2]{dCRS2021};  see also \cite[Proof of Thm. 4.7]{YK2024}. The idea is to exploit the non-degeneracy of the symplectic form on the tangent space of a point fixed by the action of the maximal affine subgroup of the group scheme $P_a$ for all $a \in A$. If $X$ is singular, then the fixed point may lie in the singular locus, and in general the Poisson bracket does not induce a non-degenerate alternating form on its tangent space, thus undermining the n\"{a}ive extension of the proof of $\delta$-regularity to the singular setting. At the moment, there is no general proof of the $\delta$-regularity of a weak abelian Lagrangian fibration with total space a singular symplectic variety.

For $G=\GL_n$, the previous approach works since $\MHiggs^d_{\GL_n,\cL}$ is smooth for $\gcd(n,d)=1$, or, alternatively, since $\MHiggs^d_{\GL_n,\cL}$ is a smooth family of compactified Jacobians and one can invoke Severi's inequalities; see \cite[paragraph following Thm. 2]{FGvS1999} and also \cite[(41)]{de-cataldo-support-sln} and \cite[Lemma 4.1]{maulik-shen-independence}.

In the reductive case, we cannot rely on the existence of a smooth connected component $\MHiggs^d_{G,\cL}$. However, we can replace the singular moduli space with the stack $\Higgs^{ss}_{G,\cL}$, and exploit the existence of an alternating form on its tangent space, defined via Serre's duality; cf.\ the bivector $\Theta$ in \Cref{square:delta}. 

\begin{introthm}[$\delta$-regularity: \Cref{thm:deltaregularmero}]\label{thmF}
    Suppose that the characteristic of the ground field $k$ is $0$, and let $D \subset C$ be a (possibly empty) reduced effective divisor such that $\deg(\omega_C(D))>0$. Then the group scheme $\schP^\circ \to A$ of symmetries of the Hitchin fibration for the reductive group $G$ and the line bundle $\cL= \omega_C(D)$ (as in \Cref{P:centr-group}(ii)) is $\delta$-regular.
\end{introthm}

The alternating form mentioned above is part of the datum of the standard shifted-symplectic structure on $\Higgs_{G,\omega_C}$ (see \cite[Cor. 2.6 (2)]{shifted_symplectic_paper}), which generalizes to a shifted-Poisson structure for $\Higgs_{G, \omega_C(D)}$ in the case when $D>0$. Hence, we can consider the $\delta$-regularity of the group scheme of symmetries of the Hitchin fibration as a geometric manifestation in classical algebraic geometry of the existence of this shifted-Poisson structure. As we remarked above, although the shifted-Poisson structure induces a Poisson structure on the moduli space $\MHiggs_{G,\omega_C(D)}$, the Poisson bivector fails to be non-degenerate  
at the singular points of the closure of the symplectic leaves of the moduli space, and it is unclear how it could serve  our purpose. It seems then essential to work on the stack $\Higgs^{ss}_{G,\cL}$.

\subsubsection*{Generalization to isotrivial group schemes}
In \Cref{section: isotrivial group schemes}, we generalize some of our foundational results on the Higgs moduli space by replacing $G$ with an isotrivial reductive group scheme $\cG$ over the curve. We establish the weak Abelian fibration structure in \Cref{thm: weak Abelian fibration structure twisted higgs} when $\deg(\cL)>2g-2$, and explain the proof of the $\delta$-regularity of the corresponding group scheme of symmetries of the Hitchin fibration in \Cref{remark: delta regularity for twisted higgs}. This level of generality arises naturally even if one is only interested in the Hitchin fibration for reductive groups over the ground field. Indeed, the results on isotrivial group schemes proved in this paper will be needed in our subsequent work on the Decomposition Theorem for the Hitchin fibration \cite{DT_paper_2} for reductive groups $G$ over $k$.

\subsubsection*{Acknowledgements} We would like to thank Giuseppe Ancona, Tom\'as G\'omez, Jochen Heinloth, Tasuki Kinjo, Luca Migliorini, Davesh Maulik, Junliang Shen and Siqing Zhang for useful mathematical discussions.

MAdC has been supported by  NSF Grants DMS 1901975 and DMS 2200492. RF was supported by the University of Roma «La Sapienza» and University of Pisa. The research was partially supported by MUR Excellence Department Project awarded to the Department of Mathematics of the University of Pisa. AFH was supported by Columbia University and University of Pennsylvania.
MM was supported by the Institute of Science and Technology Austria, the Hausdorff Institute of Mathematics in Bonn, \'{E}cole Polytechnique, and Université Paris Cité and Sorbonne Université, CNRS, IMJ-PRG, F-75013 Paris, France. This project has received funding from the European Union's Horizon 2020 research and innovation programme
under the Marie Sk{\l}odowska-Curie grant agreement No. 101034413, and from PEPS «Jeunes chercheurs et jeunes chercheuses» de l’Insmi 2024.

\begin{subsection}{Notation and setup}\label{subs: not setup}

We work over a (not necessarily algebraically closed) fixed ground field $k$. In this paper, unless otherwise stated, all stacks, algebraic spaces and schemes, and the morphisms between them, are over the field $k$; moreover, all stacks are algebraic. 

We  fix an affine algebraic group $G$ over $k$, which, unless otherwise stated, is assumed to be connected split reductive.  We denote
  the Lie algebra of $G$ by $\mathfrak{g}$, and the center subgroup scheme of $G$
by $Z_{G}$.  
The neutral component $Z_{G}^{o}$ of the center is the unique maximal central torus contained in $G$. We  
 fix a maximal split torus $T$ in $G,$  denote its Lie algebra by $\mathfrak{t}$, and denote by $\Phi= \Phi_G \subset  X^*(T) \coloneqq \mathrm{Hom} (T, {\mathbb G}_m)$ the set of roots of $G$ with respect 
to $T$.  We denote by $W:=N_G(T)/T$ the Weyl group of $G$, where $N_G(T)$ is the normalizer of the torus $T$ in $G$. 

\begin{context}
In this paper, it is always assumed that $\text{char}(k)$ does not divide the order $|W|$ of the Weyl group. 
\end{context}

Let $X_*(T) \coloneqq (X^*(T))^{\vee} = \mathrm{Hom} ({\mathbb G}_m, T)$ be the cocharacter lattice of the maximal split torus $T$, and let $X_{coroots}= \langle \Phi^{\vee}\rangle \subseteq X_*(T)$ be the sublattice generated by the coroots. %$\Phi^{\vee}$. 

\begin{defn}[{\cite[Def.\ 5.4]{hoffmann-connected-components}}]
   The fundamental group of the reductive group $G$ is
   defined by setting $\pi_1(G) \coloneqq X_*(T)/X_{coroots}$. %; see \cite[Def. 5.4]{hoffmann-connected-components}. 
\end{defn}
If $k= \mathbb{C}$, then $\pi_1(G)$ coincides with the topological fundamental group of $G$ as a complex Lie group.

A choice of a Borel subgroup $B \supset T$ determines a set of positive roots $\Phi^+$ inside the character lattice $X^*(T)$.  Let $2 \rho^{\vee} = \sum_{\alpha \in \Phi^{+}} \alpha^{\vee}$ denote the sum of the positive coroots. Given character $\chi \in X^*(T)$, its height is defined by $\text{ht}(\chi) = \langle 2 \rho^{\vee}, \chi\rangle$.

\begin{defn}[{\cite[Defn.\ 4.4]{deligne-balaji-parameswaran-complete-red}}]
   Let $V$ be a linear representation of $G$. Let $V = \oplus_{\chi} V_{\chi}$ denote the decomposition into $T$-weight spaces. Set
   \[ \text{ht}(V) := \text{max} \left\{ \text{ht}(\chi) \; \mid \; \text{$\chi$ is a character with $V_{\chi} \neq 0$} \right\}\]
   We say that the representation $V$ is of low height if either $\text{char}(k)=0$ or $\text{char}(k) > \text{ht}(V)$.
\end{defn}

\begin{remark}
In this paper, we are mostly interested in the height of the adjoint representation $\fg$ of $G$. Recall that an almost simple group $G$ is a semisimple (hence reductive) non-commutative algebraic group over $k$ such that every proper normal subgroup is finite;  %or equivalently, such that its adjoint group is simple; 
see \cite[Def. 19.7]{milne-algebraic-groups}. If $G$ is almost simple,  then we have $\text{ht}(\fg) = 2h(G) -2 = 2\dim(G)/ \rank(G)$, where $h(G)$ denotes the Coxeter number of $G$ \cite[pg. 213, (5.2.6)]{Serre2003-2004}. If $G$ is a torus, then $\text{ht}(\fg) =0$ by convention. For an arbitrary reductive group $G$, it follows from definition that the height $\text{ht}(\fg)$ is the maximum among $\text{ht}(\fg_i) = 2h(G_i) -2$ as we range over all simple quotients $G \twoheadrightarrow G_i$.
\end{remark}

 By $G$-bundle, or $G$-torsor, we mean a principal $G$-bundle.
Let $Y$ be a scheme, and let $E$ be $G$-bundle on $Y$. For any affine scheme $X$ over $k$ equipped with a left $G$-action, we denote by $E(X) \to Y$ the associated fiber bundle with fiber $X$, defined as the fppf quotient $E(X) = (X\times E)/G$.
The quotient $E(X)$ is represented by a relatively affine scheme over $Y$. In particular, when $X =V$ is a linear representation of $G$, then $E(V)$ is the associated vector bundle on $Y$ with fiber $V$.

\begin{remark}
    The height of the adjoint representation is important for us because it yields a bound on the characteristic of the ground field that ensures that the semistability of a $G$-bundle $E$ coincides with the semistability of its adjoint vector bundle $\ad(E) := E(\mathfrak{g})$ (\Cref{lemma: semistability under adjoint representation}). This is used, for example, in our proof of the existence of the moduli space (\Cref{thm: adequate moduli spaces for higgs}) and the smoothness of the stack of semistable objects in \Cref{prop: smoothness of the semistable hitchin stack}.
\end{remark}

We fix a smooth projective geometrically connected curve $C$ over $k$. We denote by $g$ the genus of $C$. We fix a line bundle $\mathcal{L}$ on $C$. 

\end{subsection}
\end{section}

\begin{section}{Lie algebras}\label{S:LieAlg}
 In this section, we collect some facts about the adjoint action $\mathrm{ad}:G\to \GL_n(\fg)$ of a reductive group $G$ on its Lie algebra $\fg$. By definition of the normalizer $N_G(T)$ of $T$ in $G$,
 and by the fact that $T$ is Abelian, the adjoint action restricted to $N_G(T)$ 
 %preserves the Lie algebra $\ft$. Furthermore, the torus $T$ acts trivially on $\ft$ (because it is Abelian). Hence, the adjoint action $\mathrm{ad}$ restricted to $N_G(T)$ 
 induces the adjoint action of the Weyl group $W:=N_G(T)/T$ on the Lie algebra $\ft$.

\subsection{The Chevalley Restriction Theorem}
\begin{thm}\label{T:konstant}
Let $G$ be a split connected reductive group over $k$.  We have that:
	\begin{enumerate}[(i)]
		\item\label{T:konstant1} The inclusion $\ft\subset\fg$ of Lie algebras induces a ring isomorphism $k[\fg]^G\cong k[\ft]^W$, where $G$ acts on $\fg$ by the adjoint representation, and $W$ acts on $\ft$ by restricting the adjoint representation of $G$.
  
  \item\label{T:konstant2} The  $k$-algebra $k[\fg]^G$, resp. $k[\ft]^W$, is freely generated by $r= rank(G):= \dim (T)$ homogeneous   polynomials $p_1,\ldots,p_r$ in $k[\fg]$, resp. in $k[\ft]$, of degrees $e_1,\ldots,e_r$. The un-ordered collection of such degrees does not depend on the choice of such generators.
		\item\label{T:konstant3} The morphism of affine varieties
        \[\chi: \fg \to \mathfrak{c}:=\Spec(k[\ft]^W), \;\; \; \;x \mapsto (p_1(x),\ldots,p_r(x))\]
		% $$
		% \begin{array}{cccl}
		% 	\chi:&\fg&\longrightarrow &\mathfrak{c}:=\Spec(k[\ft]^W)\\
		% 	&x&\longmapsto& (p_1(x),\ldots,p_r(x))
		% \end{array}
		% $$
		is $\mathbb G_m$-equivariant, where $\mathbb G_m$ acts on $\fg$ with weight $1$ and it acts on $\fc$ as follows:  $t \cdot (u_1,\ldots,u_r):=(t^{e_1}u_1,\ldots,t^{e_r}u_r)$, where $t\in\bG_m(S)$ and $(u_1,\ldots,u_r)\in \fc(S)$, for any $k$-scheme $S$.
		\item\label{T:konstant4} The morphism of affine varieties
        \[\pi: \ft \to \mathfrak{c}, \; \; \; \; t \mapsto (p_1(t),\ldots,p_r(t))\]
		% \begin{equation*}
		% 	\begin{array}{cccl}
		% 		\pi:&\ft&\longrightarrow &\mathfrak{c}\\
		% 		&t&\longmapsto& (p_1(t),\ldots,p_r(t))
		% 	\end{array}
		% \end{equation*}
		is a Galois ramified over with Galois group $W$. 
The discriminant (or branch divisor) $\DD$
  on $\mathfrak{c}$ is given by the zeroes of the 
  $W$-invariant function
		$
		\prod_{\alpha\in\Phi}d\alpha 
		$ 
  viewed as a function on $\fc,$ where $d\alpha: \ft \to k$ is the derivative of the root $\alpha: T \to \mathbb{G}_m.$
		\item\label{T:konstant5} For any simple 
  root $\alpha$, we fix a non-zero vector 
  $x_{\alpha}\in\fg_{\alpha}$  in the corresponding root space, and we set 
  $x_+:=\sum_{\alpha\text{ simple }}x_{\alpha}  \, \in \mathfrak{g}$. Then, 
  there exists a unique $\mathfrak sl_2$-triple $\{x_-,h,x_+\}$ in $\fg$, with $h\in\ft$.  
		\item\label{T:kostant5bis}  Notation as in   (\ref{T:konstant5}). 
  Let $\fg^{x_+}  \subset \mathfrak{g}$ be the Lie algebra of the centralizer 
  of $x_+$  in $\mathfrak g$. The affine subspace $x_-
  +\fg^{x_+}\subset\fg$ is contained in the open  and dense subset 
  $\fg^{\mathrm{reg}}$ of regular elements of $\mathfrak{g}$. Furthermore, the restriction of 
  $\chi: x_-+\fg^{x_+}\to \fc$ is an isomorphism. Its 
  inverse $\kappa:\fc\to \fg^{\mathrm{reg}}$ is called  a \emph{Kostant section}.
	\end{enumerate}
\end{thm}

\begin{proof} See \cite{Ko63}. See also \cite[Thm. 2.1]{NG06}.
\end{proof}

\begin{remark}\label{R:h} Here, we collect some facts about the degrees $e_1,\ldots,e_r$ introduced in \Cref{T:konstant}.

	\begin{enumerate}[(i)]
        \item\label{R:h1} $\sum_{i=1}^r e_i=|\Phi|/2+r=\dim(B)$, where $B$ is a Borel subgroup of $G$.
		\item $G$ is semisimple if and only if all the degrees $e_1,\ldots,e_r$ are greater than or equal to $2$.
		\item  $G$ is Abelian, and hence $G=T$,  if and only if all the degrees $e_i$ are equal to $1$.  In this case, the natural morphism $\pi_G:\ft \to \fc$ is the identity.
	\end{enumerate}	
\end{remark}

\begin{remark}\label{R:decom-pi} In general, a reductive group $G$ admits a canonical isogeny $ Z_G^o\times D_G\to G$, where $ Z_G^o\subset G$ is the neutral component of the center $ Z_G$ and $ D_G :=[G,G] \subset G$ is the derived/commutator subgroup, which is a semisimple group (see \cite[Prop. 19.21(c)]{milne-algebraic-groups}). We have the commutative diagram of affine schemes
	\begin{equation}
		\xymatrix{
			\ft_G\ar[rr]^{\pi_G}&\quad&\fc_G\\
			\ft_{ Z_G}\oplus \ft_{ D_G}\ar[u]^\cong \ar[rr]^{(\mathrm{Id},\pi_{ D_G})}&\quad& \ar[u]^\cong\ft_{ Z_G}\oplus \fc_{ D_G}
			}
	\end{equation}
	where the horizontal arrows are the morphisms in \Cref{T:konstant}$(\ref{T:konstant4})$, the left-hand-side vertical arrow is induced by the isogeny $ Z_G^o\times D_G\to G$  and the right-hand-side vertical arrow is obtained by taking the quotient by $W$ of the left-hand-side.
\end{remark}

\subsection{Local regular centralizer group schemes}\label{SS:local-centr-gr}

In this subsection, we assume that $G$ is a split connected reductive group over $k$. The scheme 
$$
I:=\{(x,g)\in \fg\times G\, \vert \,  \mathrm{Ad}_g(x)=x\}\subset \fg\times G
$$
is an affine not necessarily flat group scheme over the Lie algebra $\fg$. The pull-back $J:=\kappa^*I$ to $\fc$ along the Kostant-section (see \Cref{T:konstant}\eqref{T:kostant5bis}), is a  smooth affine group scheme over $\fc$, called the group scheme of regular centralizers.

Let $\pi_*(\ft \times T)$ denote the Weil restriction of the $\ft$-torus $\ft \times  T$ via the finite flat morphism $\pi: \ft \to \fc$. Recall that $\pi_*(\ft \times T)$ is represented by a smooth affine group scheme over $\fc$ \cite[Prop. A.5.2]{pseudo-reductive-groups}. On the other hand, we denote by $J^1:= (\pi_*(\ft\times T))^W$ the closed subfunctor of $W$-invariant points of $\pi_*(\ft\times T)$, which is a $\fc$-smooth closed subgroup scheme of $\pi_*(\ft \times T)$ by \cite[Prop. 4.3]{edixhoven-neron-models}. Equivalently, $J^1$ is the functor which associates to any morphism $S\to \fc$ the set of $W$-equivariant morphisms $S\times_{\fc}\ft\to T$. We denote by $J^0$  the open subgroup scheme on $J^1$ of the neutral components of the fibers \cite[\href{https://stacks.math.columbia.edu/tag/055R}
    {Tag 055R}]{stacks-project}.

The next propositions collect the main facts about these group schemes.

\begin{prop}\label{P:chain-j-c}There exists a natural chain of morphisms 
of group schemes over $\fc$ $$J^0\hookrightarrow J\hookrightarrow J^1\hookrightarrow \pi_*(\ft\times T),$$ where the first two arrows 
are open embeddings and the last one is a closed embedding, and we have that
    \begin{enumerate}[(i)]
    \item\label{P:chain-j-c--1} The group scheme $\pi_*(\ft\times T)$ is a smooth affine commutative group scheme over $\fc$ of relative dimension 
$|W|\cdot\dim(T)$ and with geometrically connected fibers.
\item\label{P:chain-j-c--2} The group schemes $J^0$, $J$ and $J^1$ are smooth commutative group schemes over $\fc$ of relative dimension $\dim(T)$. Furthermore, $J$ and $J^1$ are affine over $\fc$.
\item\label{P:chain-j-c--3} There exists a canonical homomorphism $\chi^*J\to I$ of group schemes over $\fg$, which is a $G$-equivariant isomorphism when restricted to the open locus $\fg^{\mathrm{reg}} \subset \fg$ of regular elements.   
\item\label{P:chain-j-c--4} The image of $J$ in $J^1$ is the subfunctor  whose $S$-points
        are  the $W$-equivariant morphisms 
	$
	f:S\times_{\fc }\ft\to T
	$
	satisfying the following additional property: for every geometric point $x$ of $S\times_{\fc} \ft$  and for every root $\alpha$
 such that the  hyperplane-root reflection $s_{\alpha}$ fixes $x$,  we have that $\alpha(f(x))=1$.
 \item\label{P:chain-j-c--5} The image of $J^0$ in $J$ is the open subgroup scheme of the neutral components of the fibers.
    \end{enumerate}
\end{prop}

\begin{proof} Point \eqref{P:chain-j-c--1}. The group is smooth and affine because it is the Weil restriction of a smooth and affine group scheme \cite[Prop. A.5.2]{pseudo-reductive-groups}, see also \cite[\S 2.4]{ngo-lemme-fondamental}. The statement about the fibers follows from \cite[Prop. A.5.9]{pseudo-reductive-groups}. The relative dimension can be checked by computing the dimension of any  (geometric) fiber. 
Let $c$ be a geometric point in the complement of the branch divisor $\mathfrak D$ of the map $\pi:\ft\to\fc$, and let $\kappa(c)$ be the (algebraically closed) residue field of $c$. Then the fiber of the group scheme $\pi_*(\ft\times T)\to \fc$ over $c$ is isomorphic to $T^{|\pi^{-1}(c)|}_{\kappa(c)} \cong T^{|W|}_{\kappa(c)}$ (because the quotient map $\pi$ is \'etale over $c$). In particular, $\dim(\pi_*(\ft\times T)_c)=|W|\cdot\dim(T)$.

Point \eqref{P:chain-j-c--2}.  We recall that $J^1$ is a closed subgroup scheme of $\pi_*(\ft\times T)$ and that $J$ is affine. By Point \eqref{P:chain-j-c--1}, we have that the group schemes $J^1$ is affine as well.  Furthermore, the smoothness of these groups follows from the smoothness of $J^1$, which has been proved in \cite[Lem. 2.4.1]{ngo-lemme-fondamental}. 

Point \eqref{P:chain-j-c--3}. \cite[Lem. 2.1.1]{ngo-lemme-fondamental} says that there exists a unique smooth commutative group scheme $J'$ over $\fc$  satisfying the conclusion of Point \eqref{P:chain-j-c--3}.  By the definition of a Kostant section $\kappa$, one may check (cf. 
\cite[Paragraph after Lemma 2.1.1 page 30]{ngo-lemme-fondamental}) that $J:=\kappa^*I$ satisfies the desired property, so that so $J'=J$ and the conclusion follows.

Points \eqref{P:chain-j-c--4} and \eqref{P:chain-j-c--5}. See \cite[Lem. 2.4.6, Prop. 2.4.7 and Cor. 2.4.8]{ngo-lemme-fondamental}.  
\end{proof}

The next proposition expands on  the differences between the groups $J^0$, $J$ and $J^1$.
\begin{prop}\label{P:diff-j}We have that
\begin{enumerate}[(i)]
\item\label{P:diff-j-0} 
 The natural morphism
$\pi_0(J)\hookrightarrow\pi_0(J^1)$ of \'etale group algebraic spaces of finite type over $\fc$ is an 
inclusion.
    \item\label{P:diff-j-1}  
    The natural closed embedding $Z_G\times\fc\hookrightarrow J$ of group schemes over $\fc$ induces a surjective morphism of \'etale group algebraic spaces of finite type over $\fc$
    $$
    \pi_0(Z_G)\times\fc\twoheadrightarrow \pi_0(J)=J/J^0.
    $$
    \item\label{P:diff-j-2} 
The open embedding $J\hookrightarrow J^1$ of group schemes over $\fc$ introduced in \Cref{P:chain-j-c} induces a natural isomorphism of \'etale group algebraic spaces of finite type over $\fc$
    $$
    J^1/J\cong \pi_0(J^1)/\pi_0(J)\cong \pi_*\left(\bigoplus_{\alpha\in \Phi}(\iota_{\alpha})_*(\ker(\alpha^\vee)\times \mathfrak t_\alpha)\right)
    $$
    where $\iota_{\alpha}:\mathfrak t_\alpha\hookrightarrow\ft$ is the linear subspace of those elements $t$ such that $d\alpha(t)=0$ and $\alpha^\vee:\mathbb G_m\to T_G$ is the co-root corresponding to the root $\alpha$.
\end{enumerate}
\end{prop}

\begin{proof}Point \eqref{P:diff-j-0}. Since $J$ and $J^1$ are smooth, commutative and of finite type over  $\fc$, their $\pi_0$ are represented by \'etale group algebraic spaces of finite type over $\fc$ (see \cite[Lem. 6.5]{ancona-etal-relative-motive}). Moreover, by \Cref{P:chain-j-c}\eqref{P:chain-j-c--5}, the group scheme $J^0$ is the open subgroup scheme of the neutral components of the fibers of both groups scheme $J$ and $J^1$. Hence, the inclusion follows.

Point \eqref{P:diff-j-1}. When  $\pi_0(Z_G)=\{1\}$, see \cite[Prop. 2.3.1]{ngo-lemme-fondamental}. In general, we have an exact sequence of group schemes over $\fc=\fc_G$
$$
1\to Z_G\times \fc_G\xrightarrow{j} J_G\to a^*J_{G^{\mathrm{ad}}}\to 1,
$$
where $a:\fc_G\to \fc_{G^{\mathrm{ad}}}$ is the natural morphism induced by the morphism of Lie algebras $\ft_G\to \ft_{G^{\mathrm{ad}}}=\ft_G/\ft_{Z_G}$. Since the center of $G^{\mathrm{ad}}$ is trivial (in particular, it is connected), the latter group in the above sequence has connected fibers, i.e.  $\pi_0(a^*J_{G^{\mathrm{ad}}})$ is the trivial group scheme over $\fc_G$. Hence, the morphism $\pi_0(j):\pi_0(Z_G)\times\fc\to \pi_0(J_G)$ is
 surjective by the right-exactness of $\pi_0(-)$. 

%\Roberto{....., i.e. $\pi_0(a^*J_{G^{\mathrm{ad}}})$ is the trivial group scheme over $\fc_G$. In particular, we get the following commutative diagram of commutative group schemes over $\fc$
%\[   \xymatrix{
%1\ar[r]&Z_G\times\fc_G\ar[r]^j&J_G\ar[r]&a^*J_{G^{\ad}}\ar[r]&1 \\
%1\ar[r]&N\ar@{^{(}->}[u]\ar[r]&J^0_G\ar@{^{(}->}[u]\ar[r]&a^*J_{G^{\ad}}\ar@{=}[u]\ar[r]&1    
%    }
%\]
% where the rows are exact, and $N$ contains the subgroup scheme $Z_G^o\times\fc$ (because $Z_G^o\times\fc$ is contained in $J_G^0$). By Snake Lemma, the morphism $j$ induces a surjective morphism of commutative group algebraic spaces of finite type over $\fc$$$
% j':Z_G\times\fc_G\twoheadrightarrow (Z_G\times\fc)/N\cong J_G/J^0_G\cong\pi_0(J_G).$$
%  Since $N$ contains $Z_G^o\times\fc$, the morphism $j'$ factors through $\pi_0(j):\pi_0(Z_G\times\fc)=\pi_0(Z_G)\times\fc\to \pi_0(J_G)$. In particular, $\pi_0(j)$ is surjective.
% }

 Point \eqref{P:diff-j-2}. By \Cref{P:chain-j-c}, we get the following exact sequence of group schemes over $\fc$
 $$
 1\to J\to J^1\xrightarrow{\sum_{\alpha\in\Phi}r_\alpha}\pi_*\left(\bigoplus_{\alpha\in\Phi}(\iota_{\alpha})_*(\alpha(T^{s_\alpha})\times \ft_\alpha)\right)\to 1,
 $$
 where $s_\alpha$ is the reflection with respect to the root $\alpha$, $T^{s_{\alpha}}$ is the subgroup of the torus $T$ of the $s_{\alpha}$-invariant elements and $\alpha(T^{s_{\alpha}})\subset\mathbb G_m$ is the image of $T^{s_\alpha}$ along the root $\alpha:T\to \mathbb G_m$. For each $\alpha$, $r_\alpha$ denotes the morphism which sends a $W$-equivariant morphism $f:S\times_\fc \ft\to T\in J^1(S)$ to 
 \begin{equation}\label{E:diff-j-1}
 S\times_\fc\ft_\alpha\to T^{s_\alpha}\xrightarrow{\alpha} \alpha(T^{s_\alpha}) \in (\pi\circ \iota_\alpha)_*\left(\alpha(T^{s_\alpha})\times \ft_\alpha\right)(S),
 \end{equation}
 where the first arrow is the restriction of $f$ on $S\times_\fc\ft_\alpha$. Note that the image of this restriction is contained in $T^{s_\alpha}$ because the points in $\ft_\alpha$ are $s_\alpha$-invariants and $f$ is $W$-equivariant. Hence, the composition \eqref{E:diff-j-1} is well-defined. 
 
 In order to  prove the statement, it is enough to show that, for each root $\alpha\in\Phi$, the composition $T^{s_\alpha}\xrightarrow{\alpha} \mathbb G_m\xrightarrow{\alpha^\vee}T$ induces an equality of group schemes over $k$
  \begin{equation}\label{E:diff-j-2}
 \alpha(T^{s_\alpha})=\ker\{\alpha^\vee\}.
 \end{equation}
 After replacing $G$ with $C_G(\ft_\alpha)$, we may assume that $G$ is a reductive group with roots $\Phi=\{\alpha,-\alpha\}$ and co-roots $\Phi^\vee=\{\alpha^\vee,-\alpha^\vee\}$. In particular, $G=H\times M$, where $H=\mathrm{SL}_2,\GL_2,\mathrm{PGL}_2$ and $M$ is an algebraic torus. Let $T_H$ be a maximal torus of $H$. We now prove the equality \eqref{E:diff-j-2} by treating the cases $\mathrm{SL}_2,\GL_2,\mathrm{PGL}_2$ separately.
 
 If $H=\mathrm{SL}_2$, we have $T_H=\mathbb G_m$, $T^{\mathrm{s_\alpha}}=\mu_2\times M$, $\alpha(h,m)=h^2$ and $\alpha^\vee(t)=(t,1)$. Hence, $\alpha(T^{s_\alpha})=\ker\{\alpha^\vee\}=\{1\}$.

If $H=\mathrm{GL}_2$, we have $T_H=\mathbb G_m^2$, $T^{\mathrm{s_\alpha}}=Z_{\mathrm{GL}_2}\times M$, $\alpha(h_1,h_2,m)=h_1(h_2)^{-1}$ and $\alpha^\vee(t)=(t,t^{-1},1)$. Hence, $\alpha(T^{s_\alpha})=\ker\{\alpha^\vee\}=\{1\}$.

If $H=\mathrm{PGL}_2$, we have $T_H=\mathbb G_m$, $T^{\mathrm{s_\alpha}}=\mu_2\times M$, $\alpha(h,m)=h$ and $\alpha^\vee(t)=(t^2,1)$. Hence, $\alpha(T^{s_\alpha})=\ker\{\alpha^\vee\}=\mu_2$.
\end{proof}

From \Cref{P:diff-j}, we get the following result.
\begin{coroll}\label{C:diff-j} The cokernel of the natural embedding $J\hookrightarrow J^1$ is a $2$-torsion sheaf. Moreover,
    \begin{enumerate}[(i)]
        \item\label{C:diff-j1} The natural open embedding $J^0\hookrightarrow J$ is an isomorphism if and only if the center $Z_G$ of $G$ is connected.
        \item\label{C:diff-j2} The natural open embedding $J\hookrightarrow J^1$ is an isomorphism  if and only if all the coroots $\alpha^\vee:\bG_m\to T$ of $G$ are injective (e.g. $\pi_1(G)$ is torsion-free).
    \end{enumerate}
\end{coroll}

\begin{proof}The first assertion follows from the last part of the proof of \Cref{P:diff-j}, where we have showed that $\ker\{\alpha^\vee\}$ is either trivial or $\mu_2$. Point (i):  the implication $\Leftarrow$ follows from \Cref{P:diff-j}(ii); the implication $\Rightarrow$ follows from the fact that the fiber of $J$ over the zero-element $[0]\in\fc$ is equal to $Z_G\times U$, where $U$ is a connected unipotent group (see \cite[Thm. 5.9(b)]{springer-arithm}). Point (ii) follows directly from \Cref{P:diff-j}(iii).
\end{proof}

\end{section}
\begin{section}{The stack of \texorpdfstring{$G$}{G}-Higgs bundles} \label{section: stack of higgs bundles}

\subsection{The stack of \texorpdfstring{$\cL$}{L}-valued \texorpdfstring{$G$}{G}-Higgs bundles} \label{subsection: stack of Higgs bundles}

The goal of this section is to introduce,
for a line bundle $\mathcal L$ on the curve $C$,
the algebraic stack $\Higgs_{G,\mathcal L}$ of 
$\mathcal L$-valued $G$-Higgs bundles on  $C$,
establish some of its basic properties and 
prove that, when $\mathcal L$ is sufficiently
positive, this stack is a normal local complete intersection and  a disjoint union of irreducible components parametrized by the fundamental group $\pi_1(G)$.

For any $G$-bundle $E$ over $C$, we denote by $\ad(E):=E\times^G\fg$ the adjoint bundle over $C$, attached to $E$ via the adjoint representation. The vector bundle  $\ad(E)$  has   rank  $\dim(G)$ and, being self-dual \cite[Lem. 4.2.3]{riche-kostant-section},  it has degree zero.

\begin{defn}Let $\mathcal L$ be a line bundle over $C$. 
An $\mathcal L$-valued $G$-Higgs bundle (or $G$-Higgs 
bundle, when it is clear from the context) is a pair 
$(E,\theta)$, where $E$ is a $G$-bundle and $\theta\in 
H^0(C,\ad(E)\otimes\cL)$, the Higgs field, is a global section of the 
$\mathcal L$-valued adjoint bundle $\ad(E)\otimes\cL$.
\end{defn}

We denote by $\Bun_G$ the moduli stack of $G$-bundles on $C$ and by $\Higgs_{G,\cL}$ the moduli stack of $\mathcal L$-valued $G$-Higgs bundles on $C$. Both of these stacks are algebraic stacks with affine diagonal and locally of finite type over $k$, because of their interpretation as mapping stacks $\text{Map}(C, BG)$ and $\text{Sect}_C([\mathfrak{g} \otimes \cL /G])$, respectively; see \cite[Thm. 1.3]{hall-rydh-tannakahom}. There is the natural  morphism of stacks given by forgetting the Higgs field $\theta$
\begin{equation}\label{E:Higgs-Bun}
\pi: \Higgs_{G,\mathcal L}\to \Bun_G, \; \; \;(E,\theta)\mapsto E.
\end{equation}
The geometric fibers of this morphism are the corresponding  vector spaces of Higgs fields.

The next proposition collects some facts about the moduli stack of $G$-bundles.

\begin{prop}\label{P:Bung-prop} The moduli stack $\Bun_G$ is a smooth algebraic stack over $k$ with affine diagonal of dimension $\dim(\Bun_G)=\dim(G)(g-1)$. The connected  components of $\Bun_G$ are irreducible
and they are naturally  labeled by the elements of the fundamental group $\pi_1(G)$ of $G$.
\end{prop}

\begin{proof} The smoothness follows by \cite[Prop. 4.4.1]{behrend-thesis}. The assertion about the dimension follows from \cite[Cor. 8.1.9]{behrend-thesis}.  The smoothness implies that the  connected components are irreducible. 
The assertion about the connected components is in \cite[Thm. 5.8]{hoffmann-connected-components}.
\end{proof}

% These irreducible connected components of BunG are not %quasi 
%compact in general, e.g $GL_n$; they are for $G$ a torus;
% the open of slope semistables in each irreducible %connected component is quasi compact. Similarly for Higgs.

For any $d\in\pi_1(G)$, we denote by $\Bun_G^d$ the corresponding connected component of $\Bun_G$. Note that,
in general, they are not quasicompact. We say that a $G$-bundle $E$ has degree $d \in \pi_1(G)$ if it lies in $\Bun_G^d$. 
Using the forgetful map \eqref{E:Higgs-Bun}, we get a decomposition as a disjoint union of open and closed substacks
\begin{equation}\label{higgs conn comp} \Higgs_{G,\mathcal{L}} = \bigsqcup_{d\, \in \pi_1(G)} \Higgs_{G,\mathcal{L}}^d,
\end{equation}
where $\Higgs_{G,\mathcal{L}}^d$ denotes the stack of $G$-Higgs bundles whose underlying $G$-bundle has degree $d$. Later, we show that, if $\cL$ is positive enough and $char(k)$ is zero, or large enough, then the open and closed substacks $\Higgs_{G,\mathcal{L}}^d$  are  the connected components of $\Higgs_{G,\mathcal{L}}$ (see \Cref{P:Higgs-prop}). 

\begin{example}
If $G=\GL_n,$ then $d \in \pi_1(\GL_n)={\mathbb Z}$
 corresponds to the usual degree of the underlying vector bundle. If $G=SL_n,$ then 
$\pi_1(SL_n)=1$ 
and the underlying vector bundles have degree zero.
\end{example}

\begin{lemma}\label{L:Higgs-lci-criterion}For any degree $d\in\pi_1(G)$, we have $\dim(G)\deg(\cL)\leq \dim \Higgs_{G,\mathcal L}^d$. Moreover, if the equality holds, then $\Higgs_{G,\mathcal L}^d$ is a local complete intersection stack of pure dimension $\dim(G)\deg(\cL)$.
\end{lemma}

  \begin{proof}
  We may assume that the ground field $k$ is algebraically closed.  Let $\cE$ be the universal $G$-bundle over $p:C\times \Bun_G\to \Bun_G$. Let $U\to\Bun_G$ be a smooth cover such that $U$ is a disjoint union of affine schemes. The complex $Rp_*(\mathrm{ad}(\cE_{U})\otimes\cL)$ is perfect and its formation commutes with base change $U' \to U$; see \cite[\href{https://stacks.math.columbia.edu/tag/0A1H}{Tag 0A1H}]{stacks-project}. We may represent $Rp_*(\mathrm{ad}(\cE_{U})\otimes\cL)$ by a two-term complex $\partial:\cF^0\to\cF^1$ of locally free sheaves; see \cite[Section 5, Lem.\ 1]{mumford-abelian}. In particular, we have $\ker(\partial_{U'})=p_*(\mathrm{ad}(\cE_{U'})\otimes\cL)$ for any $U' \to U$. By Riemann--Roch, we have 
\begin{align}\label{E:rk-f0-f1}
    \rank(\cF^0)- \rank(\cF^1) & = \chi(\ad(E)\otimes \cL) \nonumber \\
    & = \deg(\ad(E)\otimes \cL)+\rank(\ad(E)\otimes \cL)(1-g) \nonumber\\
    &= \dim(G) \cdot (\deg(\cL)+1-g),
\end{align}
where $E$ is a $G$-bundle corresponding to any point of a connected component of $U$.

  %  Let  $\partial:\cF^0\to\cF^1$ be a morphism of vector bundles of finite rank over each connected component of $U$, such that, for any arbitrary base change $U'\to U$, 
  %   we have $\ker(\partial_{U'})=p_*(\mathrm{ad}(\cE_{U'})\otimes\cL)$ and $\mathrm{coker}(\partial_{U'})=R^1p_*(\mathrm{ad}(\cE_{U'})\otimes\cL)$. 
  %   Note that for the vector bundle $\ad(E)\otimes \cL$  the rank is $\dim(G)$ and the degree  is $\dim(G) \cdot \deg(\cL).$
  %   By Riemann--Roch, we have
  % \begin{equation}\label{E:rk-f0-f1}
		% 		\rank(\cF^0)=\rank(\cF^1)+\dim(G) \cdot (\deg(\cL)+1-g).
		% 	\end{equation}

			By construction, the inverse image  (cf. \eqref{E:Higgs-Bun}) 
 $\pi^{-1}(U):=U\times_{\Bun_G}\Higgs_{G,\cL}$ is a closed substack of the vector bundle $\mathrm{Tot}(\cF^0)$ over $U,$ and it is locally defined by $\rank(\cF^1)$ equations. In particular, we get the following bound on the dimension
			\begin{align*}
				\dim(\pi^{-1}(U))&\geq \dim (\mathrm{Tot}(\cF^0))-\rank(\cF^1)=\dim(U)+\rank(\cF^0)-\rank(\cF^1)=\\
    &=\dim (U)+\dim(G)(\deg(\cL)+1-g).
			\end{align*}
   The above inequality implies that the dimension of each irreducible component of $\Higgs_{G,\cL}$ must be greater than  or equal to $\dim(\Bun_G)+\dim(G)(\deg(\cL)+1-g)=\dim(G)\deg(\cL)$. In particular, if the equality holds, the natural morphism $\pi^{-1}(U)\hookrightarrow \mathrm{Tot}(\cF^0)$ is a regular embedding into a vector bundle over a smooth base. It follows that $\pi^{-1}(U)$   is a local complete intersection of pure dimension $\dim(U)+\dim(G)(\deg(\cL)+1-g)$, so that
    $\Higgs_{G,\mathcal L}$ is a local complete intersection stack of pure dimension $\dim(G)\deg(\cL)$. 
   \end{proof}
   To end this subsection we recall the following notion, which is used in some forthcoming proofs.

\begin{defn}[Regular Higgs bundle] \label{defn: regular higgs bundle}
An $\mathcal L$-valued $G$-Higgs bundle $(E,\theta)$ is regular if, for every geometric
point $p\to C$, the kernel of the linear morphism of vector spaces
$$
[\theta_p,-]:\ad(E_p)\to \ad(E_p)\otimes\cL_p
$$
has dimension $\dim(T)$. 
\end{defn}

There is an open substack $\Higgs_{G,\cL}^{reg} \subset \Higgs_{G,\cL}$ consisting of regular Higgs bundles.

\begin{remark}
It can be proven that the intersection $\Higgs_{G,\cL}^{reg} \cap \Higgs_{G,\cL}^d$ is nonempty for every $d \in \pi_1(G)$ if the general cameral curve is smooth,
e.g.\ if  $\deg(\cL)>g+\frac{1}{2}$ by \Cref{P:dh-sets1}. We will not make use of this fact in the paper.
\end{remark}

\begin{subsection}{The Hitchin fibration} \label{subsection: hitchin fibration}
We fix a minimal set of homogeneous generators 
of $k[\fg]^G=k[\ft]^W$, and thus identify once and for all the $\mathbb G_m$-scheme $\fc$ with a vector space endowed with a contracting $\mathbb G_m$-action with strictly positive weights $(e_1, \ldots, e_r)$; see \Cref{T:konstant}\eqref{T:konstant3}.

 \begin{defn} \label{defn: hitchin base}
     The Hitchin base is the vector space $\hitchinbase = \hitchinbase_{G, \cL}$ of sections of the vector bundle over the curve $C$ given by $\fc_\cL:=\fc\times^{\mathbb{G}_m}\cL^\times$, where $\cL^\times:=\cL\setminus\{\operatorname{zero-section}\}$ is the $\mathbb G_m$-bundle attached to $\cL$.
 \end{defn}
 
 We have  $\mathbb{G}_m$-equivariant isomorphisms of vector bundles over $C$, resp. of vector spaces
\begin{equation*}\label{E:dec-eq}
\fc_\cL\cong\cL^{e_1}\oplus\ldots\oplus\cL^{e_r},
\quad\text{ resp. }\quad A\cong H^0(C,\cL^{e_1})\oplus\ldots\oplus H^0(C,\cL^{e_r}).
\end{equation*}

From the above decomposition, we deduce that
\begin{align}\label{E:dimA}
    \dim(A)=&\sum_{i=1}^{r}e_i\deg(\cL)+r(1-g)+H^1(C,\fc_{\cL})=\\
    \nonumber=&\dim(B)\deg(\cL)+\dim(T)(1-g)+H^1(C,\fc_{\cL})=\\
    %\nonumber=&\left(\dim(G)(g-1)+
    %\frac{\dim(G)+\dim(T)}{2}(\deg(\cL)-2g+2)\right)+H^1(C,\fc_{\cL}) = \\
    \nonumber=& 
    \left(
    \frac{1}{2} \dim(G)\deg(\cL) + \frac{1}{2} \dim(T) \left(\deg(\cL) -2g+2  \right) 
    \right)
    +H^1(C,\fc_{\cL}),
\end{align}
where the first equality follows by Riemann-Roch, the second one by \Cref{R:h}\eqref{R:h1}, and the third one by using the equality $\dim(G)=2\dim(B)-\dim(T)$.

\begin{defn}
    The quotient morphism $[\fg \otimes \cL/G] \to \fc_{\cL}$ of $C$-stacks induces a morphism between stacks of sections \[h:\Higgs_{G,\cL}=\text{Sect}_C([\mathfrak{g} \otimes \cL /G]) \to \text{Sect}_C(\fc_{\cL})=A,\] which is called Hitchin morphism. 
Given an $S$-point $a\in A(S)$, we call Hitchin fiber over $a$ the inverse image $\Higgs_{G,\cL}(a):=h^{-1}(a)$.
\end{defn}

\begin{remark}
    Alternatively, the Hitchin morphism can be written equationally as:
\begin{equation*}\label{E:hit-map}
		h:\Higgs_{G,\cL} \to A, \; \; \;(E,\theta) \mapsto \left(p_1(\theta),\ldots,p_r(\theta)\right),
\end{equation*}
where, with a small abuse of notation, the maps \[p_i \colon H^0(C, ad(E) \otimes \cL) \to H^0(C, \cL^{e_i})\]
are induced by the invariant homogeneous polynomials $p_i$ defined in \Cref{T:konstant}.\eqref{T:konstant2}.  
\end{remark}

The global nilpotent cone is defined as the Hitchin fiber $\Higgs_{G,\mathcal L}(0)$ over the zero-section $0\in A(k)$. Following \cite{BD}, we provide an upper bound for the dimension of the global nilpotent cone, when $\deg(\cL)$ is big enough.

		\begin{lemma}\label{P:nilp-cone-dim} Suppose that $\deg (\mathcal L)\geq 2g-2$, and that $char(k)$ is either zero or strictly larger than the height of the adjoint representation. Then, the dimension of the global nilpotent cone $\Higgs_{G,\mathcal L}(0)$ is at most $\frac{1}{2} \dim(G)\deg(\cL) - \frac{1}{2} \dim(T) \left(\deg(\cL) -2g+2  \right)$. 
		\end{lemma}
		
		\begin{proof} 
When $G$ is semisimple, $\cL=\omega_C$ and $char(k)=0$, a proof can be found in \cite[2.10.3, Lemma 1 and Remark at p. 51]{BD}. The idea of the proof in loc. cit. is the following: the authors provide a stratification of the global nilpotent cone, where each stratum is a smooth stack. Then they prove that the dimension of each of the stratum is at most $\dim(G) (g-1)$ (i.e. the bound in our statement for $\cL=\omega_C$). In particular, the same bound must hold for the dimension of the global nilpotent cone.
The proof with exactly the same strategy works with our hypotheses. Due to our more general assumptions, the computation of the bound of the dimension of these strata is slightly different from the one in loc. cit. For this reason, we decided to provide a self-contained proof.
  
First, we fix some notation. We denote by 
$\mathcal N$ the nilpotent cone of $\fg$, i.e.\ the $G$-
stable closed subset in $\fg$ of $\mathrm{ad}$-nilpotent 
elements in $\fg$, by $O$ a $G$-orbit in $\mathcal N$ and 
by $\overline{O}$ its closure. Similarly, we denote by 
$\mathcal N_{\mathcal L}:=\mathcal N\times^{\mathbb 
G_m}\mathcal L^\times$ the $\mathcal L$-valued nilpotent cone 
and by $O_{\mathcal L}:=O\times^{\mathbb G_m}\mathcal L^\times$ 
the $\mathcal L$-valued $G$-orbit in $\mathcal N_{\mathcal 
L}$ corresponding to the untwisted $G$-orbit 
$O\subset\mathcal N$. 
Note that that the closure of the twisted orbit is the twist of the closure of the orbit, i.e.
$\overline{O}_\cL=\overline{O}\times^{\mathbb G_m}\cL^\times$.
We also set $\nilp:=\Higgs_{G,\mathcal L}(0)$, which is an algebraic stack.

			The strategy of the proof goes as follows:
			\begin{enumerate}[(1)]
				\item\label{i:strat} we define a stratification $\{\nilp_O\}_{O\subset\mathcal N\text{ orbit}}$ by locally closed substacks of $\nilp$ labeled by the $G$-orbits of $\mathcal N$;
				\item\label{i:strat-cover} for any $G$-orbit $O\subset\mathcal N$, we define a stack $Y_O$ together with a surjective representable morphism $Y_{O}\to\nilp_O$ onto the strata;
				\item\label{i:bound-Y} we show that $Y_O$ is a smooth stack;
                \item\label{i:dim-Y} we show that $\dim(Y_O)\leq \frac{1}{2} \dim(G)\deg(\cL) - \frac{1}{2} \dim(T) \left(\deg(\cL) -2g+2  \right)$.
			\end{enumerate}

			Point \eqref{i:strat}. By construction, the stack $\nilp$ is isomorphic to  moduli stack of sections $\mathrm{Sect}_C([\mathcal N_{\mathcal L}/G])$ from the curve $C$ to the quotient stack $[\mathcal N_{\mathcal L}/G]$. For any $G$-orbit $O\subset \mathcal N$ in the nilpotent cone, we define the substack
			$$
			\nilp_{O}:=\mathrm{Sect}_C\left([\overline{O}_{\mathcal L}/G]\supset [O_{\mathcal L}/G]\right)\subset \nilp=\mathrm{Sect}_C([\mathcal N_{\mathcal L}/G])
			$$ 
		that given a $k$-scheme $V$	parametrizes sections $\eta:C\times V\to [\overline{O}_{\mathcal L}/G]$ such that there exists an open subset $U\subset C\times V$, which is fiberwise dense over $V$ and its image $\eta(U)$ is contained in $[O_{\mathcal L}/G]$. 
			
			%Fix $k$-scheme $V$.
            We remark that the locus in $V$ where a map $C\times V\to [\fg_\cL/G]$ factors along $[\overline{O}_\cL/G]$ is closed in $V$. On the other hand, the locus where a map $C\times V\to [\overline{O}_\cL/G]$ factors generically along the open orbit $[O_\cL/G]$ is open in $V$. In particular, we get that $\nilp_O$ is a locally closed substack of $\nilp$ and is thus  an algebraic stack. Moreover, the collection of these substacks defines a partition of the global nilpotent cone 
			$$
			\bigsqcup_{O\subset \mathcal N\text{ orbit}}\nilp_{O}=\nilp.
			$$

			Point \eqref{i:strat-cover}. Fix a $G$-orbit $O\subset\mathcal N$ and an element $e$ in the orbit $O$. By assumption on $char(k)$, we may apply \cite[Prop. 5.9 at page 54]{jantzen-nilp}: there exists a one-parameter subgroup $\lambda:\mathbb G_m\to G$ such that, if we set
			$$
			\mathfrak g=\bigoplus_{i\in\mathbb Z}\fg^i, \text{ where } \mathfrak g^i=\{x\in\mathfrak g\vert \operatorname{ad}(\lambda(t))(x)=t^ix\}, \quad \text{and} \;
    \fp^j:=\bigoplus_{i\geq j}\fg^i,
			$$
			then the following facts hold true:
			\begin{enumerate}[(a)]
				\item $\fp^0$ is a parabolic subalgebra of $\mathfrak g$;
				\item\label{i:stab} the parabolic subgroup $P$ with Lie algebra $\fp^0$ contains the stabilizer $Z_G(e)$ of $e$;
				\item the parabolic subgroup $P$  preserves the sublinear spaces $\fp^j$ for any $j$;
                \item\label{i:[e]+2} $e\in\fg^2$; in particular, $[e,\fg^i]\subset \fg^{i+2}$, for any $i\in\bZ$;
				\item\label{i:[e]-surj} $P \cdot e=\mathfrak p^2\cap O$, $\overline{P \cdot e}=\fp^2$,  and $[e,\fp^0]=\fp^2$.
			\end{enumerate}

Consider the inclusion of orbits $P \cdot e\subset G \cdot e$. By point (e), after passing to the closure, we get a $P$-equivariant closed embedding $\fp^2\subset \overline{\cO}$. After twisting by the line bundle $\cL$, we get a closed embedding of quotient stacks $i:[\fp^2_\cL/P]\hookrightarrow [\overline{O}_\cL/P]$. On the other hand, the morphism of quotient stacks $p:[\overline{O}_\cL/P]\to[\overline{O}_\cL/G]$ is proper and representable (more precisely, it is a $G/P$-fibration). In particular, the composition

			\begin{equation}\label{E:p-O}
				[\fp^2_\cL/P]\xrightarrow{i} [\overline{O}_\cL/P]\xrightarrow{p}[\overline{O}_\cL/G]
			\end{equation}
			 is proper and representable. Furthermore, it is an isomorphism over the open orbit $[\fp^2_\cL\cap O_\cL/P]\cong [O_\cL/G]$. We then define the moduli stack 
			$$
			Y_{O}:=\mathrm{Sect}_C\left([\fp^2_{\cL}/P]\supset [\fp^2_\cL\cap O_{\mathcal L}/P]\right)\subset \mathrm{Sect}_C([\fp_\cL^2/P])
			$$
			of sections $\eta:C\times V\to [\fp_\cL^2/P]$ such that there exists an open 
   subset $U\subset C\times V$, which is fiberwise dense over $V$ and its image $\eta(U)$ is 
   contained in $[\fp^2_\cL\cap O_\cL/P])$. The morphism \eqref{E:p-O} induces a morphism of 
   mapping stacks $F:Y_O\to \nilp_O$. It remains to show that $F$ is surjective and 
   representable. Since the morphism $\eqref{E:p-O}$ is representable, also the morphism $F$ 
   is representable. We now show the surjectivity. Fix a $k$-point $\eta:C\to 
   [\overline{O}_\cL/G]$ in $\nilp_O$. By definition, there exists an open subset $U\subset 
   C$ such that the map $U\hookrightarrow C\xrightarrow{\eta} [\overline{O}_\cL/G]$ factor 
   through $[\fp_\cL^2/P]$. Using the properness of the representable morphism 
   \eqref{E:p-O}, we extend the lift over the entire curve. Hence, $F$ is surjective.

Point \eqref{i:bound-Y}. For any $P$-bundles $E$ and $P$-representation $V$ (e.g. $\mathfrak p^i$, $\fg$, $\mathfrak g^i=\mathfrak p^i/\mathfrak p^{i+1}$, $\fg/\mathfrak p^i$), we set $E(V):=E\times^P V$. The datum of a section $\eta:C\to [\mathfrak p^2_\cL/P]$ is equivalent to the datum of a pair $(E,\theta)$, where $E$ is a $P$-bundle and $\theta\in H^0(C,E(\mathfrak p^2))$. 

The condition on the section $\eta$ of being in $Y_O$, i.e. meeting the open orbit $[\fp_\cL^2\cap O_\cL/P]$ generically, translates on the pair $(E,\theta)$ by asking that there exists an open subscheme $U\subset C$ such that, for any closed point $p\in U$, we have $\theta_p=e$, for a suitable trivialisation of the fiber $E(\mathfrak p^2)|_p$.

By \cite{halpernleistner2019mapping}, the deformation theory of a $k$-pair $(E,\theta)$ in the mapping stack $Y_O$ is controlled by the hypercohomology of the 
complex 
  \begin{equation*}
				F^\bullet:=\left[	E(\fp^0)\xrightarrow{[\theta,-]} E(\fp^2)\otimes\mathcal L\right],
			\end{equation*}
        where the left term $E\times^P\fp^0$ is concentrated in cohomological degree $0$ and the right term is concentrated in cohomological degree $1$; see also \cite{biswas-ramanan-infinitesimal} and \S\ref{subsection: tangent complex of the Higgs stack}. 
			
   Since the obstruction space is contained in $\mathbb H^2(F^\bullet)$, it is enough to show that $\mathbb H^2(F^\bullet)=0$. This group is the cokernel of the morphism
			\begin{equation}\label{E:smooth-cmplx}
	H^1(C,E(\fp^0))\xrightarrow{H^1([\theta,-])} H^1(C,E(\fp^2)\otimes\mathcal L),
			\end{equation}
			so it is enough to show that the sheaf $\mathrm{coker}([\theta,-])$ is supported on a finite 
   subset of $C$. 
   By the definition of $Y_O$, there exists an open subset $U\subset C$ such that the 
   fiber of the complex $F^\bullet$ over a closed point $p\in U$ is equal to $[e,-]:\fp^0\to \fp^2$, which is surjective by \eqref{i:[e]-surj}. Hence, the support of the sheaf 
   $\mathrm{coker}([\theta,-])$ is contained in the finite subset $C\setminus U$.

Point \eqref{i:dim-Y}. Before presenting the computation of the upper bound for $\dim(Y_O)$, we make the following observation. Take a $G$-invariant non-degenerate symmetric bilinear form $b(-,-)$ on $\mathfrak g$; see \cite[Lem. 4.2.3]{riche-kostant-section}. It induces an isomorphism of $G$-modules $\fg\cong \fg^{\vee}$, between the adjoint representation and its dual. This isomorphism sends the eigenspace $\fg^{-i}$ to $(\fg^i)^\vee$. In particular, the bilinear form $b(-,-)$ induces an isomorphism of $P$-representations
			$(\fp^i)^\vee\cong \fg/\fp^{-i+1}$. Hence, if $E$ is a $P$-bundle, we have a canonical isomorphism of vector bundles
			\begin{equation}\label{E:dual-gi}
				E(\fp^i)^\vee\cong E((\fp^i)^\vee)\cong  E(\fg/\fp^{-i+1}).
			\end{equation}
			We are now ready for proving the dimension bound. We have the following equalities
			\begin{align}\label{E:dimY1}
				\dim_{(E,\theta)} (Y_O)=&
    -\chi(F^\bullet)
    =\chi(E(\fp^2)\otimes\mathcal L)-\chi (E(\fp^0))=\\
				\nonumber=	&-\chi(E(\fg/\fp^{-1})\otimes\mathcal L^\vee\otimes \omega_C)-\chi (E(\fp^0))=\\
				\nonumber=&
    -\left[\chi(E(\fg)\otimes\mathcal L^\vee\otimes \omega_C)\right] +
    \left[ \chi(E(\fp^{-1})\otimes\mathcal L^\vee\otimes \omega_C) \right] - \left[\chi( E(\fp^0))\right],
			\end{align}
			where: the second row is obtained via Serre duality and \eqref{E:dual-gi};  the third row is obtained by using the compatibility of  Euler characteristics with exact sequences of sheaves. 
   
			By using Riemann-Roch on the last row of \eqref{E:dimY1} and the fact that $\dim(\fp^0)=\dim(P)$, we get the following equalities  
			\begin{align}\label{E:dimY2}
				\nonumber\dim_{(E,\theta)} (Y_O)=&
    \left[ \dim(G)(g-1) + \dim(G)(\deg(\cL)+2-2g) \right] + \left[\chi(E(\fp^{-1})\otimes\mathcal L^\vee\otimes \omega_C) \right] \\
					& \nonumber -\left[\chi (E(\fp^0)\otimes\mathcal L^\vee\otimes\omega_C)+\dim(P)(\deg(\cL)+2-2g) \right] =\\
				\nonumber=&\dim(G)(g-1) +(\dim(G)-\dim(P))(\deg \mathcal L+2-2g)+\\
&+\chi((E(\fp^{-1}/\fp^0)\otimes\mathcal L^\vee\otimes \omega_C).
\end{align}
\noindent \textbf{CLAIM:}
    The summand $\chi((E(\fp^{-1}/\fp^0)\otimes\mathcal L^\vee\otimes \omega_C)$ in the last row \eqref{E:dimY2} is a non-positive integer. To prove the CLAIM, we first show that there exists an exact sequence of the form
\begin{equation}\label{eq:exactsequencetheta}
   0 \xrightarrow{}
			E(\fp^{-1}/\fp^0)\xrightarrow{[\theta,-]} E
            (\fp^1/\fp^2)\otimes \mathcal L
   \xrightarrow{} N:= \text{coker}([\theta,-]) \xrightarrow{} 0,
			\end{equation}
   with $N$ 
   supported on a finite subset of the curve $C.$ We observe that the morphism $[\theta,-]$ is well-defined by Point \eqref{i:[e]+2}. By arguing as in the proof of Point \eqref{i:bound-Y},
   %\eqref{E:smooth-cmplx}, 
   we deduce that, for any closed point $p$ in $U$, the morphism $[\theta,-]$ is fiberwise isomorphic to $[e,-]:\mathfrak g^{-1}\to \mathfrak g^{1}$. Since $[\theta,-]$ is a morphism of vector bundles over a smooth connected curve, to obtain the exact sequence \eqref{eq:exactsequencetheta}, it is enough to show 
   that the morphism $[e,-]:\fg^{-1}\to \fg^1$ is an 
   isomorphism of $k$-linear spaces. 
   
   First, observe that 
   the morphism $[e,-]$ is injective since the kernel is equal 
   to $\mathrm{Lie}(Z_G(e))\cap \fg^{-1}$, which is trivial 
   because $\mathrm{Lie}(Z_G(e))$ lies in the sublinear 
   space of positive weights $\fp^0$ (see Point 
   \eqref{i:stab}). On the other hand, since $\fg^{-1}$ and 
   $\fg^{1}$ are dual to each other (via any fixed $G$-invariant nondegenerate bilinear form on $\fg$ \cite[Lem. 4.2.3]{riche-kostant-section}). Since $\dim \mathfrak g^{-1}=\dim\mathfrak g$, the morphism $[e,-]$ is an isomorphism of $k$-linear spaces. 
   
   We then get the following inequalities
  concerning the skyscraper cokernel $N$
			\begin{align}\label{E:dimY3}
				0\geq& -\chi(N)=\chi(E(\fp^{-1}/\fp^{0}))-\chi(E(\fp^{ 1}/\fp^{ 2})\otimes \mathcal L)=\\
				\nonumber=&\chi(E(\fp^{-1}/\fp^{ 0}))+\chi(E(\fp^{ -1}/\fp^{ 0})\otimes\mathcal L^\vee\otimes\omega)=\\
				\nonumber=&2\chi(E(\fp^{-1}/\fp^{ 0})\otimes\mathcal L^\vee\otimes\omega)+\dim(\fp^{-1}/\fp^{ 0})(\deg(\cL)+2-2g)\geq\\
				\nonumber\geq& 2\chi(E(\fp^{ -1}/\fp^{ 0})\otimes\mathcal L^\vee\otimes\omega),
			\end{align}
   where, for the last inequality, we have used that $\deg(\cL)\geq 2g-2$.
   The \textbf{CLAIM} is proved.
			 From \eqref{E:dimY2} and \eqref{E:dimY3}, we deduce the following bound
			\begin{align}\label{E:dimY4}
				\nonumber\dim_{(E,\theta)} (Y_O)\leq &\dim(G)(g-1) +(\dim(G)-\dim(P))(\deg \mathcal L+2-2g)\leq \\
				\nonumber\leq &\dim(G)(g-1) +(\dim(G)-\dim(B))(\deg \mathcal L+2-2g)=\\
				= &\frac{1}{2} \dim(G)\deg(\cL) - \frac{1}{2} \dim(T) \left(\deg(\cL) +2-2g  \right).
			\end{align}
		This concludes the proof of \Cref{P:nilp-cone-dim}.
		\end{proof}

\begin{coroll}\label{C:bound-hitchin-fibers} Suppose that $\deg(\cL)\geq 2g-2$, and that $char(k)$ is either zero or strictly larger than the height of the adjoint representation. Then, the fibers of the Hitchin morphism $\Higgs_{G,\cL}\to A$ have dimension at most $\frac{1}{2} \dim(G)\deg(\cL) - \frac{1}{2} \dim(T) \left(\deg(\cL) -2g+2  \right)$.    
\end{coroll}

\begin{proof}The multiplicative group $\mathbb{G}_m$ acts on $\Higgs_{G,\cL}$ by scaling the Higgs field,  i.e. $t\cdot (E,\theta):=(E,t\theta)$. The Hitchin base $A=\oplus_{i=1}^r H^0(C, \cL^{e_i})$ also admits a $\mathbb{G}_m$-action: if $s\in  H^0(C,\cL^{e_i})\subset A$, then $t\cdot \sigma=t^{e_i}\sigma$. The Hitchin morphism is $\mathbb{G}_m$-equivariant with respect to these actions. 
The zero-section $0\in A(k)$ is the only fixed $k$-point for the $\mathbb{G}_m$-action and it is contained in the closure of every other $\mathbb{G}_m$-orbit in $A$. 
In view of the upper-semicontinuity of the dimension of the fibers \cite[Prop. 2.6]{osserman-dim-stacks}, in order to verify the desired  upper bound,  it is enough to verify it  on the global nilpotent cone $\Higgs_{G,\cL}(0):=h^{-1}(0)$. This  verification is the content of \Cref{P:nilp-cone-dim}.
\end{proof}

\begin{prop}\label{P:Higgs-prop}Suppose that $\deg(\cL)\geq\max\{2g-2,1\}$, and that $char(k)$ is either zero or strictly larger than the height of the adjoint representation. Then, for all $d \in \pi_1(G)$ the following hold: 
        \begin{enumerate}[(i)]
            \item\label{P:Higgs-prop-1} The stack $\Higgs_{G,\cL}^{d}$ is a local complete intersection stack of pure dimension
            $$
            \dim(\Higgs_{G,\cL}^{d})=\begin{cases}
                \dim(G)(2g-2)+\dim (Z_G),&\text{if }\cL=\omega_C;\\
                \dim(G)\deg(\cL),&\text{otherwise}.
            \end{cases}
            $$
            \item\label{P:Higgs-prop-2} The Hitchin morphism $\Higgs_{G,\cL}^{d} \to A$ is a flat morphism with fibers of pure dimension $\frac{1}{2} \dim(G)\deg(\cL) - \frac{1}{2} \dim(T) \left(\deg(\cL) -2g+2  \right)$.
        \end{enumerate}

		\end{prop}
  		\begin{proof} We divide the proof in three parts depending on the structure of the group $G$.
        
       \noindent $\bullet$ Case $G$ semisimple (i.e. $\dim(Z_G)=0$). By Serre duality, we have 
        \[
        H^1(C,\fc_\cL)\cong H^0(C,\fc_{\cL^\vee}\otimes\omega_C)^{\vee}\cong \bigoplus_{i=1}^rH^0(C,\cL^{-e_i}\otimes\omega_C)^{\vee}.
        \]
        Since we assumed that $G$ is semisimple, we have $e_i\geq 2$, for all $i$ (see \Cref{R:h}). By assumptions on $\cL$, all the line bundles $\cL^{-e_i}\otimes\omega_C$ have strictly negative degree. Hence, we have that $H^1(C,\fc_{\cL})=0$. Let $\dim(h):= \dim(\Higgs^d_{G,\cL})-\dim(h(\Higgs^d_{G,\cL}))$ be the relative dimension of the Hitchin morphim. By \Cref{L:Higgs-lci-criterion} and \eqref{E:dimA}, we have
\begin{align}\label{E:ineq-dim(h)}
		\dim(h)&\geq \dim(\Higgs_{G,\cL})-\dim(A)\geq \dim(G)\deg(\cL)-\dim(A)\geq \\
        \nonumber&\geq \frac{1}{2} \dim(G)\deg(\cL) - \frac{1}{2} \dim(T) \left(\deg(\cL) -2(g-1)\right).
\end{align}
Recall that $\dim(h)$ is a lower bound for the  dimension of the irreducible components of the  fibers of the Hitchin morphism $h.$  Then by \Cref{C:bound-hitchin-fibers} and \eqref{E:ineq-dim(h)}, the fibers of $h$ have pure dimension $\frac{1}{2} \dim(G)\deg(\cL) - \frac{1}{2} \dim(T) \left(\deg(\cL) -2(g-1)\right)$. Furthermore, all the inequalities in $\eqref{E:ineq-dim(h)}$ are equalities. In particular, from the second equation in \eqref{E:ineq-dim(h)}, we deduce that $\dim(\Higgs_{G,\cL})=\dim(G)\deg(\cL)$ and, hence, the moduli stack  $\Higgs_{G,\cL}$ is a local complete intersection stack of pure dimension $\dim(G)\deg(\cL)$ by \Cref{L:Higgs-lci-criterion}. In remains to show that $h$ is flat. The source of the  morphism $h$ is equidimensional Cohen-Macaulay and its target is irreducible and regular. Moreover, the fibers of $h$ have pure dimension $dim(h)=\dim(\Higgs_{G,\cL})-\dim(A)$. Hence, $h$ is flat by Miracle Flatness. This concludes the proof of the proposition  in the case when $G$ is semisimple.

\noindent $\bullet$ Case $G=T$ torus (i.e. $G=Z_G$). There is a natural isomorphism of stacks $\Higgs^d_{T,\cL}\cong \Bun^d_T\times H^0(C,\ft_{\cL})$ and the Hitchin morphism is the projection onto the second factor. By \Cref{P:Bung-prop}, the stack $\Higgs_{T,\cL}$ is smooth and the Hitchin morphism is flat of relative dimension $\dim(\Bun_T)=\dim(T) (g-1) \color{black}=\frac{1}{2} \dim(G)\deg(\cL) - \frac{1}{2} \dim(T) \left(\deg(\cL) -2g+2\right)$ (because $G=T$). The statement about the dimension of $\Higgs_{T,\cL}$ follows from the fact that $\dim(H^0(C,\ft_{\cL}))$ is $\dim(T) g$ if $\cL=\omega_C$, and $\dim(T)(1-g + \deg(\cL))$, otherwise.

\noindent $\bullet$ Case $G$ arbitrary reductive group. The morphism of stacks $\Higgs_{G,\cL}^d\to\Higgs^e_{G/Z_G^{\circ},\cL}$, introduced in \Cref{R:G to G/Z Higgs}, is smooth and its geometric fibers are non-canonically isomorphic to the moduli stack $\Higgs^0_{Z_G^{\circ},\cL}$ of $Z_G^{\circ}$-Higgs bundles of degree $0\in\pi_1(Z_G^{\circ})$.

Since $G/Z_G^o$ is semisimple, $Z_G^o$ is a torus, and $\pi$ is smooth, we conclude from the previous cases that $\Higgs_{G,\cL}^d$ is a local complete intersection stack of pure dimension $\dim(\Higgs_{G,\cL}^d)=\dim(\Higgs^e_{G/Z^o_G,\cL})+\dim(\Higgs_{Z_G^o,\cL})$. A direct computation shows that $\Higgs_{G,\cL}^d$ has dimension as stated in Point (i). Furthermore, Serre duality implies that
\[\dim H^1(C,\fc_{\cL})= \begin{cases}
    \dim(Z_G)= \#\{i\,|\,e_i=1\}, & \text{if $\cL=\omega_C$}\\
0 , & \text{otherwise.}
\end{cases}\]
From this, one may show that that $\dim(h)$ satisfies the inequality \eqref{E:ineq-dim(h)}. By arguing as in the semisimple case, we obtain Point (ii).
\end{proof}

Let $\mathscr{Z}^d_{G,\cL} \subset  \Higgs^d_{G,\cL}(0)$ denote the closed substack, with its underlying reduced structure, of $\cL$-valued $G$-Higgs bundles $(E,0)$ of degree $d\in\pi_1(G)$ with zero Higgs field. The forgetful morphism $\mathscr Z_{G,\cL}^d\to\Bun_G^d:(E,0)\mapsto E$ induces an isomorphism of stacks $\mathscr Z_{G,\cL}^d\cong \Bun_G^d$. Hence, $\mathscr{Z}^d_{G,\cL}$ is an irreducible substack of $\Higgs_{G,\cL}^d(0)$ of dimension $\dim(G)(g-1)$. 

If $\deg(\cL)=2g-2$ and $g\geq 2$, the stack $\mathscr Z^d_{G,\cL}$ is an irreducible component of $\Higgs_{G,\cL}^d(0)$ (cf. \Cref{P:Higgs-prop}). As the next result shows, if $\cL$ has negative degree, $\mathscr Z^d_{G,\cL}$ is a distinguished irreducible component of the whole stack $\Higgs_{G,\cL}^d$ of $G$-Higgs bundles.

\begin{lemma}\label{L:Higgs stack l negative}Suppose that $char(k)$ is either zero or strictly larger than the height of the adjoint representation, and assume that one of the following conditions hold
\begin{enumerate}[(i)]
\item\label{L:Higgs stack l negative 1} $G$ reductive, $\deg(\cL)\leq -\max\{2g-2,1\}+(2g-2)=\min\{0,2g-3\}$, and $\cL\neq\cO_C$;
\item\label{L:Higgs stack l negative 2} $G$ semisimple, $g\geq 2$ and $\cL=\cO_C$.  
\end{enumerate}
Then, for all $d\in\pi_1(G)$, we have $\dim(\Higgs_{G,\cL}^d)\leq\dim(G)(g-1)$. Moreover, the substack $\mathscr{Z}^d_{G,\cL}\subset \Higgs^d_{G,\cL}$ of $\cL$-valued $G$-Higgs bundles with zero Higgs field is the unique irreducible component of $\Higgs_{G,\cL}^d$ of maximal dimension $\dim(G)(g-1)$.
\end{lemma}

\begin{proof}Under the assumptions (ii) and $char(k)=0$, a proof can be found in \cite[\S 2.10.5]{BD}. The strategy of the proof in \cite{BD} is to reduce the statement to a certain bound for the dimension of the fibers of the Hitchin morphism $\Higgs_{G,\cL}\to A_{G,\cL}$. We adopt the same strategy, although the way we compute the bound of the dimension of the global nilpotent cone is different from  the one in $\cite{BD}$.

The morphism of stacks $\Higgs_{G,\cL}^d\to\Higgs^e_{G/Z^o_G,\cL}$, introduced in \Cref{R:G to G/Z Higgs}, is smooth and its geometric fibers are non-canonically isomorphic to the moduli stack $\Higgs^0_{Z_G^o,\cL}\cong\Bun_{Z_G^0}\times H^0(C,\Lie(Z_G)\otimes_k\cL)$ of $Z_G^o$-Higgs bundles of degree $0\in\pi_1(Z_G^o)$. By assumptions (i) (resp. (ii)), we have that $H^0(C,\cL)=0$ (resp. $\dim(Z_G)=0$). In both cases, the morphism $\Higgs_{G,\cL}^d\to\Higgs^e_{G/Z^o_G,\cL}$ is smooth with geometrically integral fibers of dimension $\dim(Z_G)(g-1)$. Hence, without loss of generality, we may assume $G$ semisimple.

By assumption, we have $\deg(\cL)\leq 0$. Hence, if we take the usual splitting of the Hitchin base $A_{G,\cL}=\bigoplus_{i=1}^rH^0(C,\cL^{e_i})$, then we get
\begin{equation}\label{E:hitchin base L negative}
    H^0(C,\cL^{e_i})=\begin{cases}
    H^0(C,\cO_C)=k,&\text{ if }\cL^{e_i}=\cO_C,\\
    0,& \text{ otherwise}.\end{cases}
\end{equation}
In particular, there is a natural closed embedding of affine spaces $\iota:A_{G,\cL}=\mathrm{Sect}_C(\fc_{\cL})\hookrightarrow A_{G,\cO_C}=\mathrm{Sect}_C(\fc_{\cO_C})=\fc$. We prove the statement by studying the dimension of the fibers of the composition 
\begin{align}\label{E:composition hitchin trivial base}
    \Higgs^d_{G,\cL}\xrightarrow{h}A_{G,\cL}\xhookrightarrow{\iota} A_{G,\cO_C}=\fc.
\end{align}
Let $L(G)$ be the (finite) set of Levi subgroups $L$ in $G$ containing the maximal torus $T$. We then have a partition 
$
\ft=\bigsqcup_{L\in L(G)}\ft^{L-\reg}
$ of the Lie algebra $\ft$,
where $\ft^{L-\reg}$ is the locally closed subscheme in $\ft$ of those points $t$ such that $L$ is the centralizer of $t$ in $G$ (i.e. $C_G(t)=L$) or, equivalently, $\Lie(L)$ is the kernel of the adjoint morphism $\ad(t):\fg\to\fg$. By definition, $\ft^{L-\reg}$ is an open subscheme of the Lie subalgebra $\Lie(Z_L)\subset \ft$ of the center of $L$. In particular, $\dim(\ft^{L-\reg})=\dim(Z_L)$.
Similarly, we have a partition on $\fc$
\begin{equation}\label{E:stratification fc}
    \fc=\bigsqcup_{L\in L(G)}\fc^{L-\reg}, \quad \text{ with }\dim(\fc^{L-\reg})=\dim(Z_L)
\end{equation}
where $\fc^{L-\reg}:=\pi(\ft^{L-\reg})$ and $\pi:\ft\to \ft/W=\fc$ is the quotient morphism. 
Using the composition \eqref{E:composition hitchin trivial base} and the partition \eqref{E:stratification fc}, we deduce that the statement follows from
\begin{align}
   \label{E:inequality 1}& \dim(\Higgs_{G,\cL}(a))<\dim(G)(g-1)-\dim(Z_{L}),\text{ for } a\in\iota^{-1}(\fc^{L-\reg}\setminus\{0\})\cap  h(\Higgs_{G,\cL});\\
    \label{E:inequality 2}&\dim(\Higgs_{G,\cL}(0)\setminus\mathscr Z_{G,\cL} )<\dim(G)(g-1).
\end{align}

Proof of \eqref{E:inequality 2}. Let $\mathscr X$ be an irreducible component of $\Higgs^d_{G,\cL}(0)$. Observe that Points \eqref{i:strat}, \eqref{i:strat-cover} and \eqref{i:bound-Y} in the proof of \Cref{P:nilp-cone-dim} hold for any line bundle $\cL$: the assumption $\deg(\cL)\geq 2g-2$ is used only in the fourth row of \ref{E:dimY3} and in \ref{E:dimY4} in Point \eqref{i:dim-Y}. Then, there exists a $G$-orbit $O:=G \cdot e$ in the nilpotent cone of $\fg$ such that the image of the morphism $Y_O\to \Higgs_{G,\cL}(0)$ (cf. Point \eqref{i:strat-cover} in the proof of \Cref{P:nilp-cone-dim}) contains an open and dense substack of $\mathscr X$. Set $\delta:=\deg(\cL)-(2g-2)$. With the notations of the proof of \Cref{P:nilp-cone-dim}, there exists a pair $(E,\theta)$ in $Y_O$ such that $\dim(\mathscr X)\leq\dim_{(E,\theta)}(Y_O)$. Moreover, we have the following inequalities
\begin{align}\label{E:bound-X}
     \nonumber\dim(\mathscr X)\leq&\dim(G)(g-1)+(\dim(G)-\dim(P))\delta+\chi(E(\fp^{-1}/\fp^0)\otimes\mathcal L^\vee\otimes \omega_C)\leq\\
    \nonumber\leq &\dim(G)(g-1)+(\dim(G)-\dim(P))\delta-\dim(\mathfrak p^{-1}/\mathfrak p^0)\frac{\delta}{2}=\\
=&\dim(G)(g-1)+(\dim(G)-\dim(P))\frac{\delta}{2}+\dim(\mathfrak p^2)\frac{\delta}{2},
\end{align}
where the first inequality follows from \eqref{E:dimY2}, the second one follows from the first three rows in \eqref{E:dimY3}, the third one by $\Lie(P)=\mathfrak p^0$ (cf. Point \eqref{i:stab} in proof of \Cref{P:nilp-cone-dim}) and $\mathfrak p^{2}\cong\mathfrak g/\mathfrak p^{-1}$ (cf. \eqref{E:dual-gi}). 
By assumptions on $\cL$, we have $\delta=\deg(\cL)-(2g-2)\leq-\max\{2g-2,1\}<0$. 
Hence, the second and third item in the last row of \eqref{E:bound-X} are non-positive integers. In particular, $\dim(\mathscr X)\leq\dim(G)(g-1)$. The equality $\dim(\mathscr X)=\dim(G)(g-1)$ can occur if and only if $P=G$ and $\mathfrak p^2=\{0\}$. The latter conditions hold if and only if $O$ is the closed $G$-orbit $\{0\}\in\fg$ (cf. Point \eqref{i:[e]-surj} in proof of \Cref{P:nilp-cone-dim}). In other words, the generic $G$-Higgs bundle in $\mathscr X$ has zero Higgs field and, so, $\mathscr X=\mathscr Z^d_{G,\cL}$.

Proof of \eqref{E:inequality 1}. First, observe that, if $\cL$ has stricly negative degree, then the Hitchin base $A_{G,\cL}=\{0\}$ contains only the zero-section (cf. \eqref{E:hitchin base L negative}). Hence, without loss of generality we may assume $\deg(\cL)=0$, $g\geq 2$ and $G$ semisimple. 

Since the Hitchin base $A_{G,\cL}$ is a subscheme of $\fc$ \eqref{E:composition hitchin trivial base} and $\Higgs_{G,\cL}$ is defined as the stack of sections $\mathrm{Sect}_{C}([\fg_\cL/G])$, we have that  $\Higgs_{G,\cL}(a)=\mathrm{Sect}_C([\chi^{-1}(a)_\cL/G])$, where $\chi^{-1}(a)$ is the fiber of the morphism $\chi:\fg\to \fg/G=\fc$ and $\chi^{-1}(a)_\cL$ denotes its $\cL$-twist. Fix an element $t\in\ft$ in the fiber $\pi^{-1}(a)$. Let $L$ be the centralizer of $t$ and let $\mathcal N_L$ be the nilpotent cone of the Lie algebra of $L$. By \cite[Rmk pg. 84]{jantzen-nilp}, we have the following isomorphisms of stacks
\begin{align}\label{E:reduction to L}
    [\chi^{-1}(a)/G]\cong [G(t+\mathcal N_L)/G]\cong [\mathcal N_L/L].
\end{align}
 The $\cL$-twisted versions of the isomorphisms in \eqref{E:reduction to L} induce an isomorphism of stacks
\begin{align*}
    \Higgs_{G,\cL}(a)=\mathrm{Sect}_C([\chi^{-1}(a)_\cL/G])\cong \mathrm{Sect}_C([(\mathcal N_L)_\cL/L])=\Higgs_{L,\cL}(0).
\end{align*}
Even if $L$ is a reductive and non-semisimple group, the same proof of \eqref{E:inequality 1} applies with no changes. Thus, we obtain $\dim(\Higgs_{G,\cL}(a))=\dim(\Higgs_{L,\cL}(0))=\dim(L)(g-1)=\dim(G)(g-1)-\dim(G/L)(g-1)\leq \dim(G)(g-1)-\dim(G/L)$ (in the last inequality we have used $g\geq 2$). 

In order to conclude the proof,  it is enough to show $\dim(G/L)\geq 2\dim(Z_L)>\dim(Z_L)$. Since the center of a Levi subgroup has always positive dimension, the second strict inequality holds. It remains to show the first inequality. For any reductive group $H$ with maximal torus $T_H$, we have $\dim(H)=\dim(T_H)+2|\Phi^+_H|$, where $|\Phi^+_H|$ is the number of positive roots of $H$. Since a set positive roots of $G$ contains a set of positive roots of $L$, we get that $\frac{1}{2}\dim(G/L)=|\Phi^+_G\setminus\Phi^+_L|\geq |\Delta_G\setminus \Delta_L|$, where $\Delta_G$, resp. $\Delta_L$, is a set of simple roots of $G$, resp. $L$. Since the center of any reductive group is the intersection of the kernels of the simple roots, we have $\frac{1}{2}\dim(G/L)\geq |\Delta_G|-|\Delta_L|=[\dim(T)-\dim(Z_G)]-[\dim(T)-\dim(Z_L)]=\dim(Z_L)-\dim(Z_G)=\dim(Z_L)$, where in the last equality we have used that $G$ is semisimple. This concludes the proof of \eqref{E:inequality 1}. 
\end{proof}

\begin{coroll}\label{P:hj-cod}Suppose that $char(k)$ is either zero or strictly larger than the height of the adjoint representation, and assume that one of the following conditions hold
\begin{enumerate}[(i)]
\item\label{P:hj-cod 1} $G$ reductive, $\deg(\cL)\geq \max\{2g-2,1\}$ and $\cL\neq \omega_C$;
\item\label{P:hj-cod 2} $G$ semisimple, $g\geq 2$ and $\cL=\omega_C$.  
\end{enumerate}
Let $j \in {\mathbb Z}^{\geq 0}.$ Then the codimension in $\Bun_G$ of the locally closed substack 
			\begin{align*}
				\cZ_j:=&\{E\in\Bun_G \; \vert \; \dim H^0(\operatorname{ad}(E)\otimes\cL^\vee\otimes\omega_C)=j\}\\
				=&\{E\in\Bun_G \; \vert \; \dim H^1(\operatorname{ad}(E)\otimes\cL)= j \}
			\end{align*}
			is zero if $j=0$, and strictly larger than $j$ otherwise.
		\end{coroll}
\begin{proof}Take the forgetful morphism $\pi:\Higgs_{G,\cL^\vee\otimes\omega_C}\to \Bun_G$. The statement is equivalent to show $\dim(\pi^{-1}(\cZ_0))=\dim(\Bun_G)$ and $\dim(\pi^{-1}(\cZ_j)<\dim(\Bun_G)$, for any $j> 0$. By assumptions on $\cL$ and the equality $\dim(\Bun_G)=\dim(G)(g-1)$, the latter statement is a direct consequence of \Cref{L:Higgs stack l negative}.
\end{proof}

\Cref{P:Higgs-integral}
fills a gap in the literature pointed out in \cite[\S 2.2.2]{Dalakov}.

\begin{prop}\label{P:Higgs-integral}
Suppose that $\deg(\cL)\geq\max\{2g-2,1\}$, and that $char(k)$ is either zero or strictly larger than the height of the adjoint representation. Then for all $d \in \pi_1(G)$, the moduli stack $\Higgs_{G,\mathcal L}^d$ is geometrically integral.
\end{prop}

  \begin{proof}We may assume that the ground field $k$ is algebraically closed. Since the morphism of stacks $\Higgs_{G,\cL}^d\to\Higgs^e_{G/Z^o_G,\cL}$ (see \Cref{R:G to G/Z Higgs}) is smooth with geometrically connected fibers, it is enough to show the statement for $G$ semisimple.  Consider the natural forgetful morphism $\pi:\Higgs_{G,\mathcal L}\to\Bun_G$. Let $\cZ_j\subset\Bun_G$ be the locally closed substacks defined in \Cref{P:hj-cod}.
By Riemann-Roch, the restriction $\pi^{-1}(\mathcal \cZ_j)\to \cZ_j$ is the total space of a vector  bundle, and
			\begin{align}\label{E:bound-hj}
				\dim(\pi^{-1}( \cZ_j))= &\dim(G) (\deg(\cL)+1-g)+j+\dim(\cZ_j)=\\
				\nonumber=&\dim(G) \deg(\cL)+j-\mathrm{codim}_{\Bun_G}(\cZ_j)\leq \dim(G) \deg(\cL),
			\end{align}
which is an equality if and only if $j=0$. We claim that $\pi^{-1}(\cZ_0)$ is dense in $\Higgs_{G,\cL}$. Note that, if the claim holds true, then the irreducibility of $\Bun_G^d$, and hence  of $\cZ_0\cap \Bun_G^d$, implies the irreducibility of $\Higgs_{G,\cL}^d$. By way of  contradiction, assume that  $\pi^{-1}(\cZ_0)$ is not dense in $\Higgs_{G,\cL}$. Then, there would exist an irreducible component $Z$ of $\Higgs_{G,\cL}$ such that $\dim(Z)\leq \dim(\pi^{-1}(\cZ_j))$, for some $j>0$. Since $j\neq 0$, the inequality \eqref{E:bound-hj} is strict and so we would get $\dim(Z)<\dim(\Higgs_{G,\cL})$, which contradicts the pure dimensionality of $\Higgs_{G,\cL}$ (see \Cref{P:Higgs-prop}\eqref{P:Higgs-prop-1}).
   \end{proof}

\begin{prop}\label{P:Higgs-normal}
Suppose that $\deg(\cL)\geq\max\{2g-2,1\}$, and that $char(k)$ is either zero, or strictly larger than the height of the adjoint representation. Then, for all $d \in \pi_1(G)$, the moduli stack $\Higgs_{G,\mathcal L}^d$ is normal.
\end{prop}

\begin{proof}Since $\Higgs_{G,\cL}$ is a locally complete intersetion stack (see \Cref{P:Higgs-prop}\eqref{P:Higgs-normal 1}), it is enough to show that the singular locus of $\Higgs_{G,\cL}$ has codimension at least two. 

Let $h^\heartsuit:\Higgs_{G,\cL}^{\heartsuit}\to A^\heartsuit$ be the restriction of the Hitchin morphism over the open subscheme $A^\heartsuit$ (cf. \Cref{defn: a-heart}). By \cite[Cor. 5.29(iii)]{DT_paper_2}%\Cref{C:codimension complement elliptic}\eqref{C:codimension complement elliptic 3}%
, the complement of $A^\heartsuit$ in $A$ has codimension at least two.  Since the Hitchin morphism is flat (see \Cref{P:Higgs-prop}\eqref{P:Higgs-normal 2}), it is enough to show that the smooth locus of $\Higgs_{G,\mathcal L}^\heartsuit$ is fiberwise dense along $h^{\heartsuit}$. 

If $\deg(\cL)>2g-2$, the whole stack $\Higgs_{G,\cL}^{\heartsuit}$ is smooth (see \cite[Thm 4.14.1]{ngo-lemme-fondamental}) and, so, we have the statement. We now focus on $\deg(\cL)=2g-2$. Let $\Higgs_{G,\cL}^{\heartsuit,\reg}$ be the open substack of regular $G$-Higgs bundles in $\Higgs_{G,\cL}^\heartsuit$. The statement is a consequence of the following two facts
\begin{enumerate}[(1)]
\item\label{P:Higgs-normal 1} the open substack $\Higgs_{G,\cL}^{\heartsuit,\reg}$ is smooth.
   \item\label{P:Higgs normal 3} there exists an open and dense subscheme $A^{\mathrm{ell}}\subset A^{\heartsuit}$ such that $\Higgs_{G,\cL}\times_AA^{\mathrm{ell}}$ is a smooth stack.
    \item\label{P:Higgs-normal 2} the open substack $\Higgs_{G,\cL}^{\heartsuit,\reg}$ is fiberwise dense along $h^\heartsuit$.
\end{enumerate}
Since $\cL$ has even degree, the regular Hitchin fibration 
$\Higgs^{\reg}_{G,\cL}\to A$ admits a section $\kappa:A\to \Higgs_{G,\cL}^{\reg}$ 
(cf. \cite[\S 2.2.3]{ngo-lemme-fondamental} and \Cref{notn: kostant section}). Then, the statement \eqref{P:Higgs-normal 1} 
follows from \Cref{prop: torsor regular higgs under group stack} and 
\Cref{P:gcgrp-algebraic}. By the assumptions on $\cL$, the elliptic locus $A^{\elli}$ (see \cite[Def. 4.3]{DT_paper_2})
%\Cref{defn: elliptic locus}%
is an open and dense subscheme of $A^\heartsuit$ 
(cf. \cite[Cor. 5.29(i)]{DT_paper_2}).
%\Cref{C:codimension complement elliptic}\eqref{C:codimension complement elliptic 1}% 
Hence, the statement \eqref{P:Higgs normal 3} follows from \cite[Prop. 5.31]{DT_paper_2}.
%\Cref{P:elliptic locus is smooth}%

The statement \eqref{P:Higgs-normal 2} has been shown 
for $\deg(\cL)>2g-2$ in \cite[Prop. 4.16.1]{ngo-lemme-fondamental}. The same argument can be used in the case  $\deg(\cL)=2g-2$. Here a sketch of the proof. Let $\Higgs_{G,\cL}(a)$ be an Hitchin fiber over a geometric point $a$ in $A^\heartsuit$. Since the regular Hitchin fibration admits a section, the Product Formula as in \cite[Prop. 4.15.1]{ngo-lemme-fondamental} holds for the Hitchin fiber $\Higgs_{G,\cL}(a)$ (cf. \cite[\S 4.15]{ngo-lemme-fondamental} and \cite[Thm. 4.6]{NG06}). By the Product Formula and \cite[Prop. 3.7.1]{ngo-lemme-fondamental}, the locus of non-regular $G$-Higgs bundles has dimension strictly smaller than $\dim(\Higgs_{G,\cL}(a))$. Since the Hitchin fibers are equidimensional (cf. \Cref{P:Higgs-prop}\eqref{P:Higgs-normal 2}), each irreducible component of $\Higgs_{G,\cL}(a)$ intersects non-trivially the locus of regular $G$-Higgs bundles, as claimed in \eqref{P:Higgs-normal 2}.
\end{proof}
\end{subsection}

\begin{subsection}{The semistable locus}\label{subs: sslocus}

 We recall the notion of semistability for $G$-Higgs bundle. 

Let $F$ be any field over $k$. Let us consider a $G$-Higgs bundle on $(E, \theta)$ on $C_{F}$. Let $P \subset G$ be a parabolic subgroup and let $\mathfrak{p} = \text{Lie}(P)$ denote its Lie algebra. For any reduction from $G$ to $P$ of structure group $E_{P} \subset E$,  we get an induced inclusion of associated vector bundles (see \S\ref{subs: not setup})
\[
\text{ad}(E_{P})=E_{P}(\mathfrak{p}) \hookrightarrow E(\mathfrak{g}) = \text{ad}(E).
\]
\begin{defn}
We say that the parabolic reduction $E_{P}$ is compatible with $\theta$ if the global section $\theta \in H^0(C, \text{ad}(E)\otimes \mathcal{L})$ is contained in the subbundle $\text{ad}(E_{P})\otimes \mathcal{L}$.
\end{defn}

Recall that parabolic subgroups $P$ admit 
Levi decompositions $P = L \ltimes U$, 
where $L$ is a split Levi subgroup, and 
$U$ is the unipotent radical of $P$. The 
conjugation action of $L$ induces a 
representation of $L$ on the Lie algebra 
$\mathfrak{u} = \text{Lie}(U)$. When we 
restrict to the maximal central torus 
$Z_{L}^{o} \subset L$, the 
representation decomposes %breaks 
into one-dimensional weight spaces: $\mathfrak{u} = 
\oplus_{\chi \in X^*(Z^{o}_L)} {\mathfrak u}_\chi$. We say that a 
cocharacter $\lambda: \mathbb{G}_m \to 
Z_{L}^{o}$ is $P$-dominant if all 
$\lambda$-weights $\langle \lambda, \chi \rangle$ of $\mathfrak{u}$ are 
positive. The adjoint action 
of $Z_{L}^{o}$ on $\mathfrak{g}$ 
yields the weight-space decomposition
$\mathfrak{g}= \bigoplus_\chi \mathfrak{g}_{\chi}.$ 
We define the bilinear trace pairing  
$\text{tr}_{\mathfrak{g}}$ 
on the cocharacters of $Z_{L}^{o}$
by setting $\text{tr}_{\mathfrak{g}}(\lambda, \mu):= \sum_{\chi} \dim{(\mathfrak{g}_{\chi})} \lambda(\chi) \mu (\chi)$. The trace pairing can also be viewed as 
a morphism
$\text{tr}_{\mathfrak{g}}: X_*(Z_{L}^{o}) \to  X^*(Z_{L}^{o})$
from cocharacters to characters.

\begin{defn}
    We say that a character $\chi: P \to \mathbb{G}_m$ is $P$-dominant if its restriction to $Z_{L}^{o}$ is of the form $\text{tr}_{\mathfrak{g}}(\lambda)$ for some $P$-dominant cocharacter $\lambda$. This notion does not depend on the choice of a Levi subgroup $L \subset P$.
\end{defn}

\begin{remark}\label{remark:trivial on center}
Note that the neutral component $Z_G^o$ of the center $Z_G$ acts trivially on the adjoint representation $\mathfrak{g}$. It follows that any dominant character $\chi \in X^*(P)$
is automatically trivial on the identity component $Z_G^{o}$
of the center of $G.$
\end{remark}

The following definition of semi-stability is a re-wording of the definition \cite{heinloth-hilbertmumford, halpernleistner2018structure} of $\Theta$-semistability with respect to the determinant of cohomology of the universal adjoint bundle \cite[\S4.1]{herrero2023meromorphic}. It agrees with the classical notion of semistability originally due to Ramanathan \cite{ramanathan-stable} in the case of $G$-bundles without Higgs field. The semistability of the $G$-Higgs bundle $(E, \theta)$ is equivalent to the semistability of the $GL(\mathfrak{g})$-Higgs bundle $(\ad(E), \ad(\theta))$ as explained in \Cref{lemma: semistability under adjoint representation}. 

\begin{defn}[(Semi)stable $G$-Higgs bundles]\label{D:semistability higgs bundles}
Let $F$ be an algebraically closed field over $k$. Let $p=(E, \theta)$ be a $G$-Higgs bundle on $C_{F}$. We say that $p$ is semistable if for all proper parabolic subgroups $P \subsetneq G$, all parabolic reductions $E_{P} \subset E$ compatible with $\theta$, and all $P$-dominant characters $\chi: P \to \mathbb{G}_m$, we have the inequality $\text{deg}(E_{P}(\chi)) \leq 0$ for the line bundle $E_P(\chi)$ associated to $E_P$ via the character $\chi$ of $P$. We say that the $G$-bundle is stable if these inequalities are strict.
\end{defn}

We remark that $(E, \theta)$ is semistable (resp. stable)  if and only if 
it is so after base-change to any algebraically closed field $K \supset F$.

\begin{remark} \label{remark: action of torsors for the center}
    There is a natural action of the commutative group stack $\Bun_{Z_G^o}$ on $\Higgs_{G, \cL}$ given by
    \[ \Bun_{Z_G^o} \times \Higgs_{G,\cL} \to \Higgs_{G,\cL}, \; \; \; (F, (E,\theta)) \mapsto (F\times^{Z_G^o}E, \theta),\]
    where we are using  the natural identification $\rm{ad}(F\times^{Z_G^o}E) = \rm{ad}(E)$, stemming from the fact that  $Z_G^o$ acts trivially on the adjoint representation. This action is compatible with parabolic reductions, and, in view of  \Cref{remark:trivial on center},  the action preserves the locus of semistable 
    (resp. stable) geometric points of $\Higgs_{G,\cL}$.
\end{remark}

If $\text{char}(k) =0,$ 
or if $\text{char}(k)$ is strictly larger than the height 
the adjoint representation, then the 
geometric points of $\Higgs_{G,\mathcal{L}}$ that are 
semistable are exactly the geometric points of an open 
substack of $\Higgs_{G,\mathcal{L}}$, which we denote by 
$\Higgs_{G,\mathcal{L}}^{ss} \subset 
\Higgs_{G,\mathcal{L}}$ (\Cref{thm: adequate moduli spaces 
for higgs}). We let $\Higgs_{G,\mathcal{L}}^{d, ss} \subset 
\Higgs_{G,\mathcal{L}}^d$ be the corresponding open 
substack of semistable $G$-Higgs bundles of degree $d$. 
By \Cref{thm: adequate moduli spaces for higgs}, the stack $\Higgs_{G,\mathcal{L}}^{d,ss}$ is of finite 
type over $k$ and admits a quasi-projective adequate moduli 
space.  This is furthermore a good moduli space if the characteristic of $k$ is $0$. 

\begin{notn} \label{notn: higgs moduli space and hitchin fibration}
    We denote the quasi-projective adequate moduli space by 
$\MHiggs_{G,\mathcal{L}}^d$. We call the induced morphism $h: \MHiggs_{G,\mathcal{L}}^d \to A$ the Hitchin morphism (or Hitchin fibration).
\end{notn}

\begin{prop} \label{prop: smoothness of the semistable hitchin stack}
Suppose that $\deg(\cL)>2g-2$ and that either $\text{char}(k) =0$, or $\text{char}(k)$ is strictly larger than the height of the adjoint representation $\ad: G \to \GL(\mathfrak{g})$. Then the stack $\Higgs_{G,\mathcal{L}}^{d,ss}$ is smooth over $k$.
\end{prop}
\begin{proof}
  We can assume without loss of generality that $k$ is algebraically closed. 
  Let $p=(E, \theta)$ be a semistable Higgs bundle in $\Higgs_{G,\mathcal{L}}^{d,ss}(k)$. We 
  shall show 
  that $p$ is a smooth point of the stack. Using the Lie algebra structure on the vector 
  bundle $\ad(E)$, we can form the following complex of vector bundles on $C$
    \[ \mathscr{K}^\bullet_{(E, \theta)}:= \left[ \ad(E) \xrightarrow{[\theta, -]} \ad(E) \otimes \cL\right],
    \]
    where the first term $\ad(E)$ sits in cohomological degree $-1$ and the second term sits in cohomological degree $0$. According to  the deformation theory in \cite[\S2, \S3]{biswas-ramanan-infinitesimal} (see also \S\ref{subsection: tangent complex of the Higgs stack}), in order to show that $p$ is a smooth point of the stack it suffices to show that the hypercohomology group $\mathbb{H}^1(C, \mathscr{K}^\bullet_{(E, \theta)})$ vanishes. By Serre duality, $\mathbb{H}^1(C, \mathscr{K}^\bullet_{(E, \theta)})$ is dual to the kernel $K$ of the morphism
    \[ 
    H^0([\theta \otimes id_{\omega_C \otimes \cL^{-1}}, -]): H^0(\ad(E) \otimes \omega_C \otimes \cL^{-1}) \to H^0(\ad(E) \otimes \omega_C).
    \]
    Fix an element $\psi \in K$, we need to show that $\psi=0$. We may view $\psi$ as a homomorphism $\psi: \cO_{C} \to \ad(E) \otimes \omega_C \otimes \cL^{-1}$ such that $[\theta\otimes id_{\omega_C \otimes \cL^{-1}}, \psi] =0$, i.e., $\psi: (\cO_{C}, 0) \to (\ad(E) \otimes \omega_C \otimes \cL^{-1}, \ad(\theta)\otimes Id)$ is a morphism of $GL$-Higgs bundles with value in $\mathcal{L}$.

    The pair $(\cO_{C}, 0)$ is a $\cL$-valued Higgs bundle of rank 1 with the $0$ Higgs field. Notice that $(\cO_C, 0)$ is automatically semistable, and its slope as a vector bundle is $0$. On the other hand $(\ad(E)  \otimes  \omega_C \otimes \cL^{-1}, \ad(\theta)\otimes Id)$ is  a 
    $GL$-Higgs bundle, which is obtained from the 
    adjoint $GL$-Higgs bundle $(\ad(E), 
    \ad(\theta))$ by tensoring by $\omega_C 
    \otimes \cL^{-1}$. By \Cref{lemma: semistability under 
    adjoint representation}, the Higgs bundle $(\ad(E), 
    \ad(\theta))$ is semistable, and so it 
    follows that the twist $(\ad(E)  \otimes  \omega_C \otimes \cL^{-1}, \ad(\theta)\otimes Id)$ is 
    also semistable. The slope of $\ad(E) \otimes \omega_C 
    \otimes \cL^{-1}$ is strictly negative by the 
    assumption that $\deg(\cL) > 2g-2$. 
    By the semistability of the Higgs bundles and the fact that the slope of the source of $\psi$ is strictly larger than the slope of the target, it follows that $\psi$ must be $0$.
\end{proof}

\begin{lemma} \label{lemma: stable higgs is saturated} Suppose that $char(k)=0,$ or that $\text{char}(k)$ is strictly larger than the height of the adjoint representation. The stable geometric points of $\Higgs_{G, \cL}^{d,ss}$ form an open substack  $\Higgs_{G, \cL}^{d,s} \subset \Higgs_{G, \cL}^{d,ss}$ which is saturated with respect to the moduli space morphism $f:\Higgs_{G, \cL}^{d,ss} \to \MHiggs_{G,\cL}^d$. The fibers of the restriction $f: \Higgs_{G, \cL}^{d,s} \to f(\Higgs_{G, \cL}^{d,s})$ each contain a single geometric point.
\end{lemma}
\begin{proof}
    The openness of the stable locus follows from the same argument as in \cite[Lem.A.1]{herrero2023semistable}, which is stated in the special case $\cL = \omega_C$. The rest of the lemma is proven in the special case $\cL=\omega_C$ in \cite[Cor. 4.7]{herrero2023semistable}, and again the same arguments apply without changes for an arbitrary line bundle $\cL$.
\end{proof}

\begin{lemma} \label{lemma: dimension automorphisms stable higgs}
    Suppose that $char(k)=0,$ or that $\text{char}(k)$ is strictly larger than the height of the adjoint representation. Let $x$ be a geometric point of $\Higgs_{G, \cL}^{d,s}$ with algebraic group of automorphisms $\Aut(x)$. Then we have equality of dimensions as schemes $\dim(\Aut(x)) = \dim(Z_G)$.
\end{lemma}
\begin{proof}
After extending the ground field, we may assume without loss of generality that $x$ is a $k$-point. By \Cref{lemma: stable higgs is saturated}, $x$ is a closed point of the stack $\Higgs_{G,\cL}^{d,ss}$. Since $\Higgs_{G,\cL}^{d,ss}$ admits an adequate moduli space, it follows that the neutral component of the reduced subgroup scheme of $\Aut(x)$ is reductive (see \cite[p. 3 (A)]{herrero_automorphism}). 

Suppose for the sake of contradiction that $\dim(\Aut(x)) > \dim(Z_G)$. The reductivity explained above implies that there exist some homomorphism $\lambda: \mathbb{G}_m \to \Aut(x)$ that doesn't factor through $Z_G \subset \Aut(x)$.
% Take the quotient of that reductive by Z_G; it is reductive of positve dim; so it has a cocharacter (take max torus); lift it back, possibly by taking a power
We shall use the alternative descriptions of stability in \cite[Lem. 4.9]{herrero2023semistable}, which we recall now.

Set $\Theta_k := [\mathbb{A}^1_k /\mathbb{G}_m]$. A morphism $\Theta_k \to \Higgs_{G,\cL}^d$ is called 
noncentral if the composition to the rigidification $\Theta_k \to \Higgs_{G,\cL}^d \to \Higgs_{G,\cL}^d \fatslash Z_G$ does not factor through a 
morphism $\Spec(k) \to \Higgs_{G,\cL}^d \fatslash Z_G$. In \cite[Lem. 4.9]{herrero2023semistable} it is explained 
that there is a line bundle $\cL_{det}$ on $\Higgs_{G,\cL}^d$ such that a $k$-point $x \in \Higgs_{G,\cL}^d(k)$ is 
stable if and only if for all noncentral $f:\Theta_k \to \Higgs_{G,\cL}^d$ with $f(1) \cong x$, the $\mathbb{G}_m$-weight $wt(f^*(\cL_{\det})|_0)$ of the zero fiber of the pullback $f^*(\cL_{\det})$ satisfies $wt(f^*(\cL_{\det})|_0) < 0$.

The morphism $\lambda: \mathbb{G}_m \to \Aut(x)$ corresponds  to a morphism $B\mathbb{G}_m \to \Higgs_{G,\cL}^d$ by descent. Taking the $\mathbb{G}_m$-quotient of the structure morphism $\mathbb{A}^1_k \to \Spec(k)$, we obtain a map $\pi: \Theta_k = [\mathbb{A}^1_k/\mathbb{G}_m] \to B\mathbb{G}_m$. The composition $f_{\lambda}: \Theta_k \xrightarrow{\pi} B\mathbb{G}_m \to \Higgs_{G,\cL}^d$
 satisfies $f_{\lambda}(1) \cong x$, and it is noncentral because $\lambda$ does not factor through $Z_G$. Similarly, using the inverse  $\lambda^{-1}: \mathbb{G}_m \xrightarrow{(-)^{-1}} \mathbb{G}_m \xrightarrow{\lambda} \Aut(x)$, we get another noncentral morphism $f_{\lambda^{-1}}$ with $f_{\lambda^{-1}}(1) \cong x$. By stability of $x$, we must have $wt(f_{\lambda}^*(\cL_{\det})|_0) < 0$ and $wt(f_{\lambda^{-1}}^*(\cL_{\det})|_0)<0$. But by construction  $wt(f_{\lambda}^*(\cL_{\det})|_0) = - wt(f_{\lambda^{-1}}^*(\cL_{\det})|_0)$. We have reached the desired contradiction.
\end{proof}

\begin{lemma}[Jordan-H\"older reductions] \label{lemma: jordan-holder}
    Suppose that the characteristic of $k$ is zero, or that it is strictly larger than the height of the adjoint representation. Let $x= (E,\theta)$ be a geometric point of $\Higgs_{G, \cL}^{d,ss}$ that is polystable (i.e., it is closed in its fiber under the adequate moduli morphism $\Higgs_{G, \cL}^{d,ss} \to \MHiggs_{G, \cL}^{d,ss}$). Then there exists a unique Levi subgroup $L\subset G$ such that $(E,\theta)$ is obtained from a stable $L$-Higgs bundle by extension of structure group from $L$ to $G.$ In particular, if $x$ is not stable then $L \subsetneq G$.
\end{lemma}
\begin{proof}
    After base-change, we may assume that $k$ is algebraically closed and that $x$ is defined over $k$. If the characteristic of $k$ is zero, then this lemma follows from \cite[Thm. 3.18]{otero-jordan-holder}. We note that the arguments in \cite{otero-jordan-holder} may also be applied under the assumption that the characteristic of $k$ is strictly larger than the height of the adjoint representation.
\end{proof}

	\begin{prop} \label{cor: higgs moduli space properties}
 Suppose that $\deg (\mathcal L)\geq \max\{2g-2, 1\}$, and that $char(k)=0,$ or $\text{char}(k)$ is strictly larger than the height of the adjoint representation. Then the following statements hold.
		\begin{enumerate}[(i)]
  \item The stable locus $\Higgs_{G,\cL}^{d,s}$ is nonempty.
  
			\item The moduli space $\MHiggs_{G,\mathcal L}^d$ is a geometrically normal and geometrically irreducible quasiprojective scheme of dimension
             $$
            \dim(\Higgs_{G,\cL}^{d})=\begin{cases}
                \dim(G)(2g-2)+2\dim(Z_G),&\text{if }\cL=\omega_C;\\
                \dim(G)\deg(\cL)+ \dim(Z_G),&\text{otherwise}.
            \end{cases}
            $$
		\end{enumerate}
	\end{prop}
 \begin{proof}
      All the statements can be checked after passing to an algebraic closure of $k$, since the formation of adequate moduli spaces commutes with flat base change \cite[Prop. 5.2.9(1)]{alper_adequate}. Hence, we assume without loss of generality that $k$ is algebraically closed. The quasi-projectivity of $\MHiggs_{G,\cL}^d$ is proved in \Cref{thm: adequate moduli spaces for higgs}.

      \noindent \textit{Proof of (i):} By \cite[Prop. 5.28(i) and Cor. 5.29(i)]{DT_paper_2},
      %\Cref{P:dh-sets2} and \Cref{C:codimension complement elliptic}\eqref{C:codimension complement elliptic 1}
      there exists a nonempty open set $A^{\mathrm{ell}} \subset A$ such that $h^{-1}(A^{\mathrm{ell}})$ is contained in the stable locus $\Higgs_{G,\cL}^{d,s}$. Since the morphism $h \colon \Higgs_{G,\cL}^{d,ss} \to A$ is flat, thus open, (\Cref{P:Higgs-prop}), it follows that $h^{-1}(A^{\mathrm{ell}})$ is nonempty.

      \noindent \textit{Proof of (ii):}
    The stack $\Higgs_{G, \cL}^{d,ss}$ is irreducible (\Cref{P:Higgs-integral}) and normal (\Cref{P:Higgs-normal}).
    Its adequate moduli space $\MHiggs_{G, \cL}^d$ is thus  normal and irreducible \cite[Prop. 5.4.1]{alper_adequate}. From \Cref{lemma: stable higgs is saturated} we know thatthe compositum $f: \Higgs_{G, \cL}^{d,s} \hookrightarrow \Higgs_{G,\cL}^{d,ss} \to \MHiggs_{G, \cL}^d$ is a universal homeomorphism onto an open subscheme $U \subset\MHiggs_{G, \cL}^d$. By \cite[\href{https://stacks.math.columbia.edu/tag/06RC}
    {Tag 06RC}]{stacks-project}, there is an open dense substack $\mathscr{W} \subset \Higgs_{G, \cL}^{d,s}$ that is a gerbe over an algebraic space  $\mathscr{W} \to W$, where $W \subset U$ is open. For 
    every stable Higgs bundle $x$ in $\Higgs_{G, \cL}^{d,s}$ we have $\dim(\Aut(x)) 
    = \dim(Z_G)$ by \Cref{lemma: dimension automorphisms stable higgs}. By 
    \Cref{P:Higgs-integral} and \Cref{P:Higgs-prop}, $\mathscr{W} \subset \Higgs_{G,\cL}^{d,ss}$ is irreducible 
    of dimension
   $$
            \dim(\Higgs_{G,\cL}^{d})=\begin{cases}
                \dim(G)(2g-2)+\dim (Z_G),&\text{if }\cL=\omega_C;\\
                \dim(G)\deg(\cL),&\text{otherwise}.
            \end{cases}
            $$
            Since the dimension of all stabilizers of $\mathscr{W}$ are equal to $\dim(Z_G)$, it follows that $W$  has the desired dimension. Hence, the same holds for the irreducible scheme $\MHiggs_{G,\cL}^{d} \supset W$.
 \end{proof}

	\begin{prop} \label{prop: fibers hitchin pure dimensional}
 Suppose that $\deg (\mathcal L)\geq \max\{2g-2, 1\}$, and that $char(k)=0$ or $\text{char}(k)$ is strictly larger than the height of the adjoint representation. 
 \begin{enumerate}[(i)]
     \item The Hitchin morphism $h:\MHiggs_{G,\cL}^{d} \to A$ is a proper
     surjective morphism with fibers of pure dimension \[\frac{1}{2}\dim(G) \deg(\cL) -\frac{1}{2} \dim(T) (\deg(\cL)-2g+2)+\dim(Z_G).\]
     \item If the characteristic of $k$ is zero, then $h$ is flat.
 \end{enumerate} 
	\end{prop}
	\begin{proof}
	    
	\textit{Proof of (i):} The assumptions on the characteristic imply that Behrend's conjecture holds for $G$ (cf. \cite[Cor. 6.12]{biswas-holla-hnreduction}). Therefore, the properness follows from \cite[Cor. 6.21+Rem. 6.22]{alper2019existence}.
 
 Since the morphism $ \Higgs_{G,\cL}^{d,ss} \to A$ is flat, every point of $\Higgs_{G,\cL}^{d,ss}$ admits a generalization that lies in the generic fiber. Since the adequate moduli space morphism $\Higgs_{G,\cL}^{d,ss} \to \MHiggs_{G,\cL}^{d}$ is surjective, it follows that every point in $\MHiggs_{G,\cL}^d$ admits a generalization that lies in the generic fiber of $h: \MHiggs_{G,\cL}^d \to A$. By \Cref{cor: higgs moduli space properties} and the formula \eqref{E:dimA} for the dimension of $A$, the generic fiber of $h: \MHiggs_{G,\cL}^d \to A$ is irreducible and has the desired dimension $\frac{1}{2}\dim(G) \deg(\cL) -\frac{1}{2} \dim(T) (\deg(\cL)-2g+2)+\dim(Z_G)$. By the upper semicontinuity of fiber dimension, it follows that the fibers of $h$ have dimension at least $\frac{1}{2}\dim(G) \deg(\cL) - \frac{1}{2} \dim(T) (\deg(\cL)-2g+2)+\dim(Z_G)$ at every point. To conclude, it suffices to show that the fibers of $h$ have dimension at most $\frac{1}{2}\dim(G) \deg(\cL) - \frac{1}{2} \dim(T) (\deg(\cL)-2g+2)+\dim(Z_G)$.

    We first look at the stable locus. We use an argument similar to the proof of \Cref{cor: higgs moduli space properties}. Consider the open substack $\Higgs_{G, \cL}^{d,s} \subset \Higgs_{G, \cL}^{d,ss}$ of stable Higgs bundles as in \Cref{lemma: stable higgs is saturated}. For every stable Higgs bundle $x$ in $\Higgs_{G, \cL}^{d,s}$ we have $\dim(\Aut(x)) = \dim(Z_G)$ by \Cref{lemma: dimension automorphisms stable higgs}. From \Cref{lemma: stable higgs is saturated} we know that $f: \Higgs_{G, \cL}^{d,s} \hookrightarrow \Higgs_{G,\cL}^{d,ss} \to \MHiggs_{G, \cL}^d$ is a universal homeomorphism onto an open subscheme $U \subset\MHiggs_{G, \cL}^d$. Recall that the fibers of $h: \Higgs_{G, \cL}^{d,s} \to A$ are equidimensional of dimension $\frac{1}{2}\dim(G) \deg(\cL) - \frac{1}{2} \dim(T) (\deg(\cL)-2g+2)$ (\Cref{P:Higgs-prop}). Let $x$ be a geometric point of $A$, and consider the reduced closed substack $\scrX$ of the fiber $\left(\Higgs_{G, \cL}^{d,s}\right)_x$. By \cite[Lem 5.2.11]{alper_adequate}, the stack $\scrX$ admits an adequate moduli space morphism $\scrX \to X$ which is a homeomorphism, and there is an induced morphism to the fiber $X \to U_x$ that is a universal homeomorphism. By \cite[\href{https://stacks.math.columbia.edu/tag/06RC}
    {Tag 06RC}]{stacks-project}, it follows that there is an open dense substack $\mathscr{W} \subset \scrX$ that is a gerbe over an algebraic space $\mathscr{W} \to W$, where $W \subset X \to U_x$ is dense. Since $\mathscr{W} \subset \scrX$ is equidimensional of dimension $\frac{1}{2}\dim(G) \deg(\cL) - \frac{1}{2} \dim(T) (\deg(\cL)-2g+2)$ and the dimension of all stabilizers of $\mathscr{W}$ are equal to $\dim(Z_G)$, it follows that $W$ is equidimensional of the desired dimension $\frac{1}{2}\dim(G) \deg(\cL) - \frac{1}{2} \dim(T) (\deg(\cL)-2g+2)+\dim(Z_G)$. Since $W \subset X \to U_x$ maps homeomorphically onto an open dense subset, we have shown that the fibers of $h:U \to A$ are equidimensional of the desired dimension.
    
    In order  to conclude that the fibers of the morphism $h: \MHiggs_{G,\cL}^d 
    \to A$ have dimension at most 
    $\frac{1}{2}\dim(G) 
    \deg(\cL) - \frac{1}{2} \dim(T) 
    (\deg(\cL)-2g+2)+\dim(Z_G),$ we use a similar argument 
    as in the proof of \cite[Thm. II.6(i)]{faltings-projective-conn} to cover them by Hitchin fibers for Levi subgroups. For any Levi subgroup $L \subsetneq G$, we consider their Hitchin fibrations $h: \MHiggs_{L,\cL}^{d',ss} \to A_L$ with stable loci $U_L \subset \MHiggs_{L,\cL}^{d',ss}$. Just as in \cite[Thm. II.6(i)]{faltings-projective-conn}, the existence of Jordan-H\"older reductions (\Cref{lemma: jordan-holder}) implies that for every geometric point $x$ of $A$, every irreducible component of the fiber $(\MHiggs_{G,\cL}^{d})_x$ is dominated by some open subscheme of a Hitchin fiber $(U_L)_y$ of the stable locus $U_L \subset \MHiggs_{L,\cL}^{d',ss}$ for some Levi subgroup $L\subsetneq G$. By our dimension computation for the dimension of $h$-fibers over the stable locus, we have 
    \begin{align*}
        \dim((U_L)_y) &=  \frac{1}{2}\dim(L) \deg(\cL) - \frac{1}{2} \dim(T) (\deg(\cL)-2g+2)+\dim(Z_L) \leq \\
        &\leq\frac{1}{2}\dim(G) \deg(\cL) - \frac{1}{2} \dim(T) (\deg(\cL)-2g+2)+\dim(Z_G)
    \end{align*}
    Hence, every $h$-fiber of $\MHiggs_{G, \cL}^d$ has dimension at most $\frac{1}{2}\dim(G) \deg(\cL) - \frac{1}{2} \dim(T) (\deg(\cL)-2g+2)+\dim(Z_G)$, as desired.

       \medskip

    \noindent \textit{Proof of (ii):}  The statement about $h$-flatness in characteristic $0$ follows from the flatness of $h: \Higgs_{G, \cL}^{d,ss} \to A$ (\Cref{P:Higgs-prop}) by \cite[Thm. 4.16(ix)]{alper-good-moduli}, because in characteristic zero the adequate moduli space $\MHiggs_{G, \cL}^d$ is a good moduli space.
	\end{proof}

   \begin{notn}[Kostant section] \label{notn: kostant section}
     Given the choice of a square root of the line bundle $\cL$ in $\Pic(C)$ and choices of $x_{\alpha}^+$ as in \Cref{T:konstant}\eqref{T:konstant5}, Ng\^o 
     defines a section of the Hitchin morphism $\kappa: A \to \Higgs_{G,\cL}^0$ known as 
     a Kostant section \cite[\S2.2.3]{ngo-lemme-fondamental}. By construction, the 
     image of $\kappa$ lies on the open substack of regular $G$-Higgs bundles as in \Cref{defn: regular higgs bundle}.
 \end{notn}

 \begin{lemma} \label{lemma: kostant section is stable}
   Suppose that $\text{deg}(\cL)>0$, and that $\text{char}(k)$ is  either zero, or is strictly larger than the height of the adjoint representation. If $\mathcal{L}$ admits a choice of square root, then the Kostant section $\kappa: A \to \Higgs_{G, \cL}^0$ lands in the locus of stable Higgs bundles.
 \end{lemma}
 \begin{proof}
     This is proven in \cite[App. A]{herrero2023semistable}. Note that the proof in \cite{herrero2023semistable} is written for the special case when $\mathcal{L} = \omega_{C}$, but it applies to an arbitrary line bundle $\mathcal{L}$ 
     admitting a square root as long as $\text{deg}(\mathcal{L})>0$.
 \end{proof}

We end this section by explaining the behavior of the Higgs moduli space under finite \'etale covers. This is to be used in Subsections \ref{subsection: group of symmetries of hitchin untwisted} and \ref{subsection: action on the semistable locus untwisted}.

\begin{notn} \label{notn: pullback higgs under cover}
    Let $n: \widetilde{C} \to C$ be a finite \'etale Galois cover of degree $N$. Pullback under $n$ induces a morphism of stacks of Higgs bundles
\[n_{\Higgs}: \Higgs_{G,\cL} \to \Higgs_{G,n^*(\cL)}, \; \; \;  (E, \theta) \mapsto (n^*E, n^*\theta).\]
It follows from the definition of degree of a $G$-bundle in \cite{hoffmann-connected-components} that the pullback morphism $n_{\Higgs}$ sends a connected component $\Higgs_{G,\cL}^{d}$ into $\Higgs_{G,n^*(\cL)}^{N \cdot d}$.
\end{notn}

\begin{lemma} \label{lemma: pullback morphism on higgs moduli stacks}
     With setup as in \Cref{notn: pullback higgs under cover} above, the morphism $n_{\Higgs}: \Higgs_{G,\cL} \to \Higgs_{G,n^*(\cL)}$ is affine and respects the loci of semistable geometric points. In particular, if $\text{char}(k)$ is either zero or stricly larger than the height of the adjoint representation, then it induces an affine morphism of moduli spaces $n_M: \MHiggs_{G,\cL}^d \to \MHiggs_{G,n^*(\cL)}^{N \cdot d}$.
\end{lemma}
\begin{proof}
    Let $\Gamma$ denote the Galois group of the cover $n: \widetilde{C} \to C$. There is an action of the group $\Gamma$ on the stack $\Higgs_{G,n^*(\cL)}^{N \cdot d}$ given by pulling back the $G$-bundle and Higgs field:
    \[ \Gamma \times \Higgs_{G,n^*(\cL)}^{N \cdot d} \to \Higgs_{G,n^*(\cL)}^{N \cdot d}, \; \; \; \; (\gamma, (E, \theta)) \mapsto (\gamma^*E, \gamma^*\theta).\]
    \'Etale descent for the cover $n: \widetilde{C} \to C$ implies that the pullback morphism $n_{\Higgs}: \Higgs_{G,\cL}^d \to \Higgs_{G,n^*(\cL)}^{N \cdot d}$ exhibits $\Higgs_{G,\cL}^d$ as the stack of fixed points $\left(\Higgs_{G,n^*(\cL)}^{N \cdot d}\right)^{\Gamma} \to \Higgs_{G,n^*(\cL)}^{N \cdot d}$ as in \cite[Def. 2.3]{romagny-actions-stacks} (cf. the proof in \cite[Prop. 5.4.2.4]{lurie-gaitsgory-tamagawabook}). Since  $\Higgs_{G,n^*(\cL)}^{N \cdot d}$ has affine diagonal, it follows  that $n_{\Higgs}: \Higgs_{G,\cL}^d \to \Higgs_{G,n^*(\cL)}^{N \cdot d}$ is affine \cite[Prop. 3.7]{romagny-actions-stacks}. Furthermore, by the uniqueness of Harder-Narasimhan parabolic reductions \cite[\S 4]{herrero2023meromorphic}, it follows as in the proof of \cite[Prop. 4.3]{gauged_theta_stratifications} that the pullback morphism preserves semistability, and so it induces a morphism of adequate moduli spaces $n_M: \MHiggs_{G,\cL}^d \to \MHiggs_{G,n^*(\cL)}^{N \cdot d}$. Since the morphism on stacks is affine, if follows that the morphism of adequate moduli spaces is also affine \cite[Prop. 5.2.11]{alper_adequate}.
\end{proof}
\end{subsection}
\end{section}

\begin{section}{Group scheme of symmetries of the Hitchin fibration} \label{section: group scheme of symmetries of the hitchin fibration}

\subsection{Cameral Covers}

\begin{defn}[\cite{donagi-spectral}]
The cameral curve $\cC_a$ associated with a $k$-point $a$ of the Hitchin base $A$
is the curve defined by the following Cartesian diagram
{\cb (cf. \Cref{defn: hitchin base})}
$$	\begin{tikzcd}
        \cC_a\arrow[d,"\pi_a"]\arrow[r]& \ft_\cL\arrow[d,"\pi"]\\
        C\arrow[r,"a"]&\fc_\cL.
        \end{tikzcd}
$$
The cameral cover associated with $a$ is the projection morphism $\pi_a:\cC_a\to C$.
\end{defn}

	In the paper, we use the same symbol $\pi$ for the morphism $\ft\to \fc$ in \Cref{T:konstant} (\ref{T:konstant4}) and its $\mathcal{L}$-valued version $\pi:\ft_{\mathcal{L}}\to \fc_{\mathcal{L}}$ in the above diagram.  We call universal cameral cover the left-hand side vertical arrow in the Cartesian diagram
$$	
\begin{tikzcd}
        \cC\arrow[d,"\pi"]\arrow[r]& \ft_\cL\arrow[d,"\pi"]\\
        C\times A\arrow[r,"ev"]&\fc_\cL,
\end{tikzcd}
$$	
where the horizontal arrow $ev$ is defined by the assignment 
$(p,a)\mapsto a(p)$. We call universal cameral curve the composition 
$\eta:=(pr_2\circ \pi):\cC\to C\times A\to A$ of the 
universal cameral cover with the projection onto the 
Hitchin basis. Given a  morphism $f:U\to A$, we denote by 
$\eta_U:\cC_U\to U$, the pull-back of the universal cameral 
curve along $f$. Similarly, 
we denote by $\pi_U:\cC_U\to C\times U$ the pull-back of the 
universal cameral cover along $\mathrm{Id}_C\times f:C\times 
U\to C\times A$.\\

% \begin{remark}Let $G=\GL_n$ and let $T$ be the maximal torus given by the diagonal matrices.
% For $a\in A_{\GL_n}$, the spectral curve $\mathcal S_a$ corresponding to $a$ is canonically isomorphic to the quotient of the cameral curve $\cC_a$ by the subgroup $S_{n-1}$ of the Weyl group $W=S_n$ of those elements which fix the first entry of the diagonal matrices in $T$, \Mirko{It is important that we choose this particular subgroup in its conjugacy class. Up to conjugation, there is a unique  subgroup $S_{n-1}$ in $S_n$} i.e. $\mathcal S_a\cong \cC_a/S_{n-1}$ (see \cite[Section 9]{donagi-gaitsgory}). 
% For more details on the topic, we refer the reader to \cite{donagi-spectral}.
% \end{remark}

The discriminant divisor $\DD_{\mathcal{L}}:=\DD\times^{\mathbb G_m}\mathcal{L}^{\times}$ is a hypersurface in the vector bundle $\fc_{\mathcal{L}}$. 

%The following lemma shows that as soon as $\cL$ is ample, every section meets the discriminant divisor.

\begin{lemma}\label{L:a-meet-d}Assume $\deg(\cL)>0$.  Every section $a\in A(\overline{k})$ meets every irreducible component of the discriminant divisor $\DD_\cL$.
\end{lemma}

\begin{proof} We may assume that $k$ is algebraically closed.

 By construction, every irreducible component of $\DD_\cL=\DD\times^{\bG_m}\cL^\times$ is of the form $\DD'_\cL:=\DD'\times^{\bG_m}\cL^\times$, for a unique  irreducible component   $\DD'$  of $\DD\subset \fc.$ We fix such an irreducible component $\DD'.$  
 
%Let $\DD^1,\ldots, \DD^n$ be the irreducible components of $\DD$. Then, the union $\bigcup_{i=1}^n\DD^i_\cL$, where $\DD^i_{\cL}:=\DD^i\times^{\bG_m}\cL$, is contained in $\DD_\cL$. We claim that the inclusion is an equality. This can be checked easily, locally over $C$, where the line bundle trivialises. Furthermore, $\DD^i_\cL$ is irreducible because it is the $\bG_m$-quotient of $\DD'\times\cL$, which is irreducible. So, $\DD^1_\cL,\ldots,\DD^n_\cL$ are the irreducible components of $\DD_\cL$. We fix a such irreducible component $\DD'_\cL$.}. 
We have that $\DD'=\Spec(k[\ft]^W/(F))\subset \fc=\Spec(k[\ft]^W)$ for some $W$-invariant polynomial $F$ in $k[\ft]$. Since the polynomial describing the discriminant $\DD$ is the product of linear homogeneous polynomials in $k[\ft]$ (see \Cref{T:konstant}\eqref{T:konstant4}), we  have that the same holds for $F,$ so that $F$ is homogeneous in $k[\ft]$. In particular, the divisor $\DD'_\cL$ may be described as the inverse image of the zero-section along the morphism
$$
F:\fc_\cL\to\cL^{N},
$$
where $N$ is the total degree of $F$ as a homogeneous polynomial in $k[t].$

 Thus, any section $a\in A(k)$ produces a section $F(a)\in H^0(C,\cL^N)$ of the line bundle $\cL^N$. By construction, $a$ meets $\DD'_\cL$ if and only if $F(a)$ vanishes at least at one point of $C$. The latter statement is always true since we assumed the line bundle $\cL$ to be ample.
\end{proof}

% {\Md Remove this sentence?}
% The next result deals with the problem of connectedness of a cameral curve. 

	\begin{prop}\label{P:cameral-curve-conn} Assume $\deg(\cL)>0$. The homomorphism of structure sheaves $\cO_A\to \eta_*\cO_{\cC}$ of the universal cameral curve $\eta:\cC\to A$  is universal isomorphism. In particular, every cameral curve is geometrically connected.
	\end{prop}

\begin{proof} Without loss of generality, we may assume that the ground field $k$ algebraically closed. It is enough to show that for every section $a\in A(k)$, the global regular functions of the corresponding cameral curve are only the constants, i.e. $H^0(\cC_a,\cO_{\cC_a})=k$. Indeed, then the surjectivity hypothesis in \cite[Cor. 8.3.11(a)(i)]{fga-explained} is automatically satisfied at every closed point of $A$ for $q=0$ and $F = \cO_{\cC}$.

% \Mark{Roberto, reference?}.

% \Roberto{I do not have a reference. I wrote an argument. Hope it is clear

% In order to prove the statement, it is enough to show that for every geometric point $a$ in $A$, the induced morphism of $k(a)$-algebras 
% \begin{equation}\label{E:k(a)}
%     k(a)\xrightarrow{\cong} H^0(C_a,\cO_{\cC_a})
% \end{equation} is an isomorphism. Indeed, let $p:U\to A$ be an arbitrary morphism. We need to check that the morphism of structure sheaves
% $$
% \eta_U^\#:\cO_U\to (\eta_U)_*\cO_{\cC_U}
% $$
% is an isomorphism. Now assume that \eqref{E:k(a)} holds for any geometric point. This implies that the function
% $$
% U\ni u\mapsto \dim_{k(u)}H^0(C_u,\cO_{C_u})\in \bZ
% $$
% is the constant function with value $1$. Then, the pushforward $(\eta_U)*\cO_(\cC\times_AU)$ is a line bundle over $U$ and for any Zariski point $u\in U$ the canonical morphism
% \begin{equation}\label{E:commutes}
% (\eta_U)_*\cO_{\cC_U}\otimes_Uk(u)\to H^0(C_u,\cO_{C_u})
% \end{equation}
% is bijective, see \cite[Cor. 2 at pp. 50-51]{mumford-Abelian}. In particular, the morphism $\eta_U^\#$ is a morphism of line bundles over $U$. Hence, we can check the bijectivity at level of fibers over Zariski points. Using the isomorphism \eqref{E:commutes}, we see that $\eta_U^\#$ restricted over a fiber, with respect to a Zariski point $u$ in $U$, is exactly the morphism of $k(u)$-algebras
% $$
% k(u)\to H^(C_u,\cO_{C_u})
% $$
% induced by the curve $C_u\to u$. We can check the bijectivity of the above morphism on the algebraic closure of $k(u)$ and this latter condition is precisely the claim \eqref{E:k(a)}.
% }
	
Let us write $a = (a_1, a_2, \ldots, a_r) \in \bigoplus_{i =1}^r H^0(C, \cL^{e_i})$. Fix a set $\{p_1,\ldots,p_r\}$ of homogeneous polynomials in $k[\mathfrak t]$, which freely generates $k[\mathfrak t]$ as a $k$-algebra (see \Cref{T:konstant}\eqref{T:konstant2}). For $i=1,\ldots,r$, set $e_i:=\deg(f_i)$. 

Let $\pi:\mathbb{P}:=\mathbb{P}(\ft_\cL\oplus\cO_C)=\text{Proj}(\Sym^{\bullet}((\ft_\cL\oplus\cO_C)^\vee))\to C$ be the usual  compactification of the vector bundle $\ft_\cL\to C$. By construction, a cameral curve is a closed subscheme of $\mathbb{P}$.  Since $\mathbb P$ is an integral projective scheme over $k$, it is enough to show that the restriction $H^0(\mathbb P,\cO_{\mathbb{P}})\to H^0(\cC_a,\cO_{\cC_a})$ is surjective. If we denote by $\mathcal{I}$ the ideal sheaf of $\cC_a$ in $\mathbb P$, we only need to show that the first cohomology group $H^1(\mathbb P,\mathcal I)$ vanishes. By definition, the cameral cover $\cC_a$ is defined by the system of equations
\begin{equation*}
	\begin{cases*}
		P_1:=p_1(T_1,\ldots,T_r)-a_1Y^{e_1}=0,\\
		\vdots\\
		P_r:=p_r(T_1,\ldots,T_r)-a_rY^{e_r}=0,
	\end{cases*}
\end{equation*}
where $T_1,\ldots, T_r$ are the coordinate sections in $H^0(\mathbb P,\pi^*\cL\otimes\cO_{\mathbb P}(1))$ of the vector bundle $\ft_\cL$ and $Y$ is the infinity section of $\mathbb P$ in $H^0(\mathbb P,\cO_{\mathbb P}(1))$. In other words, the curve is a complete intersection of $r$ ample divisors $D_1,\ldots,D_r$, where $D_i=div(P_i)$ and $P_i\in H^0(\mathbb P,\pi^*\cL^{e_i}\otimes\cO_{\mathbb{P}}(e_i))$. In particular, if we denote by $\cE$ the vector bundle dual to the direct sum $\bigoplus_{i=1}^r\pi^*\cL^{e_i}\otimes\cO_{\mathbb{P}}(e_i)$, the Koszul complex 
$$
0\to\bigwedge^r\cE\to \cdots\to \bigwedge^2 \cE\to \cE \to \cO_{\mathbb{P}}\to \cO_{C_a}\to 0,
$$
induced by the global sections $P_1,\ldots,P_r$ is exact. By the hypercohomology spectral sequence, in order to show that $H^1(\mathcal I)=0,$ it suffices to show $H^i(\bigwedge^i\cE)=0$ for all $0 \leq i \leq r$.

By construction, these wedge products in the Koszul complex split as direct sum of line bundles of the form $\pi^*\cL^d\otimes\cO_{\mathbb{P}}(d)$, with $d$ strictly negative integer. Hence, it is enough to check the vanishing of the cohomology groups on these line bundles in degrees $0 \leq i \leq r$. We actually prove that
\begin{equation*}
H^i(\mathbb P, \pi^*\cL^d\otimes\cO_{\mathbb{P}}(d))=\begin{cases*}
	H^1(C,R^r\pi_*(\pi^*\cL^d\otimes\cO_\mathbb{P}(d)))& if $i=r+1$;\\
	0& otherwise.
\end{cases*}
\end{equation*}
We first compute the derived push-forward along $\pi$. By the projection formula and  \cite[\href{https://stacks.math.columbia.edu/tag/01XX}{Tag 01XX}]{stacks-project}, we get the following equalities (for $d$ strictly negative):
\begin{equation}\label{E:rel-coh-pb}
R^q\pi_*(\pi^*\cL^d\otimes\cO_\mathbb{P}(d))=
\begin{cases*}
	\cL^{d+r}\otimes \Sym^{-r-1-d}(\ft_\cL \oplus\cO_C),&$q=r$ and $d+r+1\leq 0$;\\
	0& otherwise.
\end{cases*}
\end{equation}
Note that the symmetric product $\Sym^{-r-1-d}(\ft_{\cL}\oplus\cO_C)$ splits as a direct  sum of line bundles. Each of these line bundles is isomorphic to $\cL^\nu$, with $0\leq \nu \leq -d-r-1 $.
Thus, when $q=r$ and $d+r+1\leq0$, the $r$-th derived push-forward splits as direct sum of line bundles of the form
$
\cL^{d+r+\nu}$, where $d+r+\nu \leq -1 $.
Since $\deg(\cL)>0$ (by hypothesis), these line bundles have no global sections. Putting all together, we arrive at the conclusion that 
\begin{equation}\label{E:no-global-section}
H^0(C,R^q\pi_*(\pi^*\cL^d\otimes\cO_\mathbb{P}(d))=0, \text{for any } q\in\bZ\text{ and }d<0.
\end{equation}
Consider now the Leray-spectral sequence
$$
H^p(C,R^q\pi_*(\pi^*\cL^d\otimes\cO_\mathbb{P}(d)))\Rightarrow H^{p+q}(\mathbb P, \pi^*\cL^d\otimes\cO_\mathbb{P}(d)).
$$
By \eqref{E:no-global-section} and the fact the $C$ is one-dimensional, we deduce that 
$
H^{i}(\mathbb P, \pi^*\cL^d\otimes\cO_\mathbb{P}(d))=H^1(C,R^{i-1}(\pi_*(\pi^*\cL^d\otimes\cO_\mathbb{P}(d)))
$.
By \eqref{E:rel-coh-pb}, the latter group is trivial if $0\leq i\leq r$. This concludes the proof.
\end{proof}

\subsection{The loci of reduced and smooth cameral covers \texorpdfstring{$A^{\heartsuit}$}{A heart} and \texorpdfstring{$A^{\diamondsuit}$}{A diamond}.}

Let $A$ be the Hitchin base as defined in \Cref{defn: hitchin base}. We define two subsets of sections $a\in A$ depending on the position of the curve $a(C)$ with respect to the discriminant divisor $\DD_{\mathcal{L}}$. 

\begin{defn} \label{defn: a-heart}
The locus $A^{\heartsuit}$ is the space of sections in $A$ not contained in $\DD_{\mathcal{L}}$, i.e.
    %Let $A^{\heartsuit}$ be the subset of $A$ whose sections are not contained in $\DD_{\mathcal{L}}$, i.e.
\begin{equation*}\label{E:ag-hearth}
	A^{\heartsuit}:=\{a\in A\;| \;a(C)\nsubseteq \DD_{\mathcal{L}}\}.
\end{equation*}
\end{defn}

\begin{defn} \label{defn: a-diamond}
The locus $A^{\diamondsuit}$ is the space of sections in $A$ cutting the divisor $\DD_{\mathcal{L}}$ transversally, i.e.
    %Let $A^{\diamondsuit}$ be the subset of $A^{\heartsuit}$, whose sections cut the divisor $\DD_{\mathcal{L}}$ transversally, i.e. 
\begin{equation*}
	A^{\diamondsuit}:=\{a\in A\;|\; a(C) \text{ intersects   }\DD_{\mathcal{L}}\text{ on the smooth locus, with multiplicity }1 \}.
\end{equation*}
\end{defn}

The next proposition summarizes the main facts about these sets.

\begin{prop}\label{P:dh-sets1} There is a chain of open embeddings $A^{\diamondsuit}\subset A^{\heartsuit}\subset A$, and the following facts hold.
\begin{enumerate}[(i)]
    \item\label{P:dh-sets-h} A geometric point $a$ of $A$ lies in $A^{\heartsuit}$ if only if the corresponding cameral curve $\cC_a$ is reduced.
    \item\label{P:dh-sets-d} A geometric point $a$ of $A$ lies in $A^{\diamondsuit}$ if and only if the corresponding cameral curve $\cC_a$ is smooth.
\end{enumerate}
Moreover, if $\cL^{m}$ is very ample for any $m\geq 2$ (e.g. $2\deg(\cL)\geq 2g+1$), then the sets $A^{\diamondsuit}$ and $A^{\heartsuit}$ are non-empty.
\end{prop}

\begin{proof} The first inclusion of open sets follows by the definition. 
If $a$ lies in $A^{\heartsuit}$, then the cameral curve $\cC_a$ is generically reduced because $\chi \colon \mathfrak{t} \to \mathfrak{c}$ is \'{e}tale away from $\mathfrak{D}$, so reduced everywhere since $\chi$ is flat; see \cite[Lem. 4.5.1]{ngo-lemme-fondamental}. On the other hand, if $a(C)$ is contained in the discriminant divisor $\mathfrak D_{\cL}$, then the generic fiber of the corresponding cameral cover $\cC_a\to C$ is non-reduced. Hence, the cameral curve $\cC_a$ is non-reduced. This gives Point \eqref{P:dh-sets-h}, while 
Point \eqref{P:dh-sets-d} is the content of \cite[Lem. 4.7.3]{ngo-lemme-fondamental}. %Moreover, the last assertion follows from the proof of \cite[Prop. 4.7.1]{ngo-lemme-fondamental}.

It remains to show the last assertion. When $\deg(\cL)>2g$, this is proven in \cite[Prop. 4.7.1]{ngo-lemme-fondamental}. We now briefly explain how to extend the result to $\cL^m$ very ample for any $m\geq 2$. Let $G\twoheadrightarrow G^{\mathrm{ad}}:=G/Z_G$ be 
the natural quotient morphism of reductive groups, and let $\mathrm{ad}:\fc_{G}\to \fc_{G^{\mathrm{ad}}}$ be the 
corresponding morphism of affine spaces. By definition of discriminant divisor, we have an equality $\mathrm{ad}^{-1}(\mathfrak D_{G^{\mathrm{ad}}})=\mathfrak D_G$ of closed 
subschemes of $\fc_G$. In particular, we have an equality of open subschemes 
$A_G^\diamondsuit=\mathrm{ad}^{-1}(A^\diamondsuit_{G^{\mathrm{ad}}})$ in $A_G$, where $\mathrm{ad}$ denotes 
the induced morphism $A_{G}\to A_{G^{\mathrm{ad}}}$ of Hitchin bases. %Hence, the last assertion follows from the proof of \cite[Prop. 4.7.1]{ngo-lemme-fondamental} for $G^{\mathrm{ad}}$. 
Hence,  without loss of generality, we may assume $G=G^{\mathrm{ad}}$, so in particular, we may assume $G$ 
 is semisimple. By \Cref{R:h}, we have a splitting $\fc_{\cL}\cong \cL^{e_1}\oplus\cdots\oplus \cL^{e_r}$, where 
 $e_i\geq 2$. By our assumptions on   $\cL$, the sections of $\fc_{\cL}$ separate infinitesimal points. Using this fact, the assertion follows by arguing as in the proof of \cite[Prop. 4.7.1]{ngo-lemme-fondamental}.
\end{proof}

\subsection{The group scheme of symmetries of the Hitchin fibration} \label{subsection: group of symmetries of hitchin untwisted}
In this subsection, we assume $char(k)\nmid |W|$.

Consider the group schemes $J^0$, $J$ and $J^1$ over $\fc$ introduced in Subsection \ref{SS:local-centr-gr}. They come equipped with compatible $\bG_m$-actions, as we now explain.

\begin{defn}\label{D:I equivariant} The group $\bG_m$ acts on the group $\fg$-scheme $I=:\{(x,g)\in\fg\times G|\,\mathrm{Ad}_g(x)=x\}$ as follows: for any test scheme $S$, we set
$$
    \lambda.(x,g):=(\lambda x,g),\text{ where }\lambda\in\bG_m(S)\text{ and } (x,g)\in I(S).
    $$
The structural morphism $I\to \fg$ is $\bG_m$-equivariant, where the action on the target is given by scalar multiplication
\end{defn}

\begin{prop}\label{P:J equivariant} There exists a canonical $\bG_m$-action on $J$ (resp. $J^0$) such that 
\begin{enumerate}[(i)]
    \item\label{P:J equivariant 1} the structural morphism $J\to \fc$ (resp. $J^0\to \fc$) is $\bG_m$-equivariant, where the action on the target is the one described in \Cref{T:konstant}\eqref{T:konstant3}
    \item\label{P:J equivariant 2} the canonical morphism $\chi^*J\to I$ (resp. $\chi^*J^0\to I$) of group schemes over $\fg$ (see \Cref{P:chain-j-c}\eqref{P:chain-j-c--3}) is $\bG_m$-equivariant, where the equivariant structure on $I$ is defined as in \Cref{D:I equivariant} and the equivariant structure on $\chi^*J$ is given by pulling-back the one of Point \eqref{P:J equivariant 1}, via the $\bG_m$-equivariant morphism $\chi:\fg\to \fc$ (see \Cref{T:konstant}\eqref{T:konstant3}), 
\end{enumerate}
    \end{prop}

\begin{proof}The assertion for $J$ follows from \cite[Prop. 3.3]{NG06}. Since $\bG_m$ is connected and $J^0$ is open subgroup scheme on $J$ of the neutral components of the fibers (see \Cref{P:chain-j-c}\eqref{P:chain-j-c--5}), the $\bG_m$-equivariant structure on $J$ restricts to a $\bG_m$-equivariant structure on $J^0$.
\end{proof}

\begin{prop}\label{P:equivariant embeddings}The morphisms 
of group schemes over $\fc$ $$J^0\hookrightarrow J\hookrightarrow J^1\hookrightarrow \pi_*(\ft\times T),$$ introduced in \Cref{P:chain-j-c} are $\bG_m$-equivariant, where
\begin{enumerate}[(i)]
\item\label{P:equivariant embeddings 1} the $\bG_m$-equivariant structure on $J^0\to\fc$ and $J\to\fc$ is defined as in \Cref{P:J equivariant}.
    \item\label{P:equivariant embeddings 2} the $\bG_m$-equivariant structure on $J^1\to\fc$ (resp. $\pi_*(\ft\times T)\to\fc$) is defined as follows:
    $$
    \lambda.(S\times_\fc\ft\xrightarrow{f} T):=\big(S\times_\fc \ft\xrightarrow{(\mathrm{Id}_S,\lambda^{-1})}S\times_\fc\ft\xrightarrow{f} T\big),
    $$
where $S$ is a test scheme, $\lambda\in\bG_m(S)$ and $f\in J^1(S)=\pi_*(\ft\times T)^W(S)$ (resp. $f\in \pi_*(\ft\times T)(S)$).
\end{enumerate}
\end{prop}

\begin{proof}From the proof of \Cref{P:J equivariant}, we get that the embedding $J^0\hookrightarrow J^1$ is $\bG_m$-equivariant. On the other hand, it follows directly from the definitions in Point \emph{(ii)} %\eqref{P:equivariant embeddings 2}
that the embedding $J^1=\pi_*(\ft\times T)^W\hookrightarrow\pi_*(\ft\times T)$ is $\bG_m$-equivariant. 
It remains to show that the embedding $J\hookrightarrow J^1$ is $\bG_m$-equivariant. Since $J^1=\pi_*(\ft\times T)^W\hookrightarrow\pi_*(\ft\times T)$ is a closed embedding (see \Cref{P:chain-j-c}), it is enough to show that the composition $J\hookrightarrow J^1\hookrightarrow \pi_*(\ft\times T)$ is $\bG_m$-equivariant. By adjunction, the latter statement is equivalent to showing that the morphism $j:\pi^*J\to \ft\times T$ is $\bG_m$-equivariant, where the equivariant structure on the source is given by pulling-back the equivariant structure on $J$ and the one on the target is given by scalar multiplication on $\ft$ and by the trivial action on $T$. The morphism $j$ has been constructed in the proof of \cite[Prop. 2.4.2]{ngo-lemme-fondamental}. It follows from the construction that the morphism $j$ is $\mathbb G_m$-equivariant.
\end{proof}

By \Cref{P:equivariant embeddings}, we may thus perform a twist by $\cL$ and obtain three group schemes $J^0_\cL$, $J_\cL$ and $J_\cL^1$ over $\fc_\cL$. 

\begin{notn}
    By pulling back these group schemes along the evaluation map $ev:C \times A\to \fc_\cL:(p,a)\mapsto a(p)$, we obtain three group schemes $J_A^0$, $J_A$ and $J_A^1$
    over $C\times A$.
They are related by a chain of open embeddings
$$
J_A^0\hookrightarrow J_A\hookrightarrow J_A^1
$$
of group schemes over $C\times A$, given by pulling-back the chain of inclusions in \Cref{P:chain-j-c}. For any morphism of schemes $f:U\to A$, we denote by $J_U^0$, $J_U$ and $J_U^1$ the group schemes over $C\times U$, obtained by pull-back of $J_A^0$, $J_A$ and $J_A^1$ along ${\rm id}_C \times f$, respectively.
\end{notn}

\begin{prop}\label{P:chain-j-cam}There exists a natural chain of homomorphisms of group schemes over $C \times A$ $$J_A^0\hookrightarrow J_A\hookrightarrow J^1_A\hookrightarrow \pi_*(\cC\times T),$$ where the first two arrows are open embeddings and the last one is a closed embedding, such that
   \begin{enumerate}[(i)]
       \item\label{P:chain-j-cam1} the image of $J_A^1$ in $\pi_*(\cC\times T)$ is the subfunctor of morphisms
	$$
	\Big(S\times_{(C \times A)}\cC\to T\big)\in \Hom(S\times_{(C \times A)}\cC,T)=\pi_*(\cC\times T)(S)
	$$
 that are $W$-equivariant;
    \item\label{P:chain-j-cam2} the image of $J_A$ in $J^1_A$ is the subfunctor of $W$-equivariant morphisms 
    $$
	\Big(S\times_{(C \times A)}\cC\to T\big)\in \Hom_W(S\times_{(C \times A)}\cC,T)=J_A^1(S)
	$$
	such that for any geometric point $x\in S\times_{(C\times A)}\cC$ fixed by the hyperplane-root reflection $s_{\alpha}$, for some root $\alpha:T\to \bG_m$, we have $\alpha(f(x))=1$;
 \item\label{P:chain-j-cam3} the image of $J_A^0$ in $J_A$ is the open subgroup scheme of the neutral components of the fibers.
 \end{enumerate}
\end{prop}

\begin{proof}It follows directly from \Cref{P:chain-j-c}
\end{proof}

\begin{prop}\label{P:galois-descr-j-cam}Then we have that
\begin{enumerate}[(i)]
\item\label{P:galois-descr-j-cam1} The group scheme $\pi_*(\cC\times T)$ is a smooth affine commutative group scheme over $C \times A$ of relative dimension $|W|\cdot\dim(T),$ with geometrically connected fibers.
\item\label{P:galois-descr-j-cam2} The group schemes $J_A^0$, $J_A$ and $J^1_A$ are smooth affine commutative group schemes over $C\times A$ of relative dimension $\dim(T)$.
\item\label{P:galois-descr-j-cam3} For any $a\in A(\overline{k})$, we have the following isomorphisms $$\mathrm{Lie}(J^0_a)\cong\mathrm{Lie}(J_a)\cong\mathrm{Lie}(J_a^1)\cong \fc_{\cL}^\vee\otimes\cL\cong \cL^{1-e_1}\oplus\cdots\oplus\cL^{1-e_r}$$
of vector bundles over $C$.
\end{enumerate}
\end{prop}

\begin{proof}The first two points follow directly from \Cref{P:chain-j-c}. The last one follows from \cite[Prop. 4.13.2]{ngo-lemme-fondamental}.
\end{proof}

\begin{defn} \label{defn: stack of symmetries of the hichin fibration}
The stack of symmetries of the Hitchin fibration, denoted by  $\stP \to A$, is the relative moduli stack of $J_A$-torsors for the proper morphism $C \times A \to A$. 
\end{defn}

Since the group scheme $J_A$ is commutative, the moduli stack $\stP$ comes naturally equipped with a commutative group-stack structure in sense of \cite[XVIII 1.4]{sga4}. Given two $J_A$-torsors $E$ and $F$, we denote their sum by $E\otimes F$. For any morphism of schemes $U\to A$, we set $\stP_U:=\stP\times_{A}U$.

\begin{prop} \label{prop: torsor regular higgs under group stack} 
The stack of symmetries of the Hitchin fibration $\stP\to A$ acts on $\Higgs_{G,\cL}\to A$. Moreover, the open substack $\Higgs_{G, \cL}^{reg} \to A$ of regular $G$-Higgs bundles (see \Cref{defn: regular higgs bundle}) is a $\stP$-pseudo-torsor over $A$. 
\end{prop}

\begin{proof} The morphism $\chi_\cL^*J\to I$ of \Cref{P:chain-j-c}\eqref{P:chain-j-c--3} induces for any $S$-point $(E,\theta)\in\Higgs_{G,\cL}(S)$ over $a:=h(E,\theta)\in A(S)$ a morphism of group schemes $J_a\to Aut_{C\times S}(E,\theta)$ over $S$. Hence, we may twist the $G$-Higgs bundle by any $J_a$-torsor (see also \cite[\S3, 4.3]{ngo-lemme-fondamental}). This defines an action morphism $\stP \times_{A} \Higgs_{G,\cL} \to \Higgs_{G,\cL}$. By \cite[Prop. 2.2.1]{ngo-lemme-fondamental} (see also proof of \cite[Prop. 4.3.3]{ngo-lemme-fondamental}), %\Roberto{Pseudo-torsor=simply-transitive action. This is guaranteed by the fact that $[\fg^{reg}/G]\to \fc$ is a $J$-gerbe. The assumption on $D'$ is needed for the existence of the Kostant-Hitchin section that tells you that it is a torsor. 
% Since the natural morphism $[\fg^{reg}/G]\to \fc$ is a $J$-gerbe (cf. \cite[Prop. 2.2.1]{ngo-lemme-fondamental}), the $\stP$-action on $\Higgs_{G,\cL}^{reg}=\mathrm{Sec}(C,[\fg^{reg}_\cL/G])$ is simply-transitive, i.e. $\Higgs_{G,\cL}^{reg}$ is $\stP$-pseudo torsor (cf. proofs of \cite[Prop. 4.3.3]{ngo-lemme-fondamental} and \cite[Prop. 4.3]{NG06}).}
the open substack $\Higgs_{G,\cL}^{reg}$ becomes a $\stP$-pseudo-torsor under this action.
\end{proof}

\begin{remark}\label{rmk: BJ action}
Let $\pi \colon \mathscr{X} \to S$ be an algebraic stack over a scheme $S$, and let $I_{\mathscr{X}} \to \mathscr{X}$ be the inertia stack. Let $G \to S$ be an fppf group algebraic space
endowed with a homomorphism $\pi^*G \to I_{\mathscr{X}}$ of relative group algebraic spaces over $\mathscr{X}$. There is a natural morphism $BG\times_S \mathscr{X} \to \mathscr{X}$ and, if $G \to S$ is commutative,
then $BG$ acts on
$\mathscr{X}$ over $S$. %If further $S \to C$ is fibred over $C$, then the (not necessarily algebraic) stack of sections $\text{Sect}_C(BG)$ acts on the (not necessarily algebraic) stack of sections $\text{Sect}_C(\mathscr{X})$. 
We may apply this to the case when $S= \fc_{\cL}$, the stack is $\mathscr{X}= [\fg \times \cL/G]$ and $G = J_{\cL}$, thus obtaining an action of $BJ_{\cL}$ on $[\fg \times \cL/G]$. Taking stacks of sections over $C$, we recover the action described in \Cref{prop: torsor regular higgs under group stack}.
\end{remark}

\begin{prop}\label{P:gcgrp-algebraic}The moduli stack $\stP$ is a smooth commutative group stack over $A$, with affine diagonal, and of relative dimension $\frac{1}{2}\dim(G) \deg(\cL) - \frac{1}{2} \dim(T) (\deg(\cL)-2g+2)$.
\end{prop}

\begin{proof}
The algebraicity and smoothness follow from \cite[Prop. 1.1]{heinloth-uniformization}, and the affineness of the diagonal is also shown in the proof of loc. cit.
% \underline{Algebraic and locally of finite type.} Let $\mathcal BJ_A\to C \times A$ be the classifying stack of $J_A$-torsors. By definition, we have a natural isomorphism of $A$-stacks
% $$
% \stP\cong (pr_{A})_*(\mathcal BJ_A)
% $$
% where the right-hand side is the Weil restriction of the stack $\mathcal BJ_A$ along the projection $pr_{A}: C \times A\to A$, i.e. for any $A$-scheme $U$, the groupoid $(pr_{A})_*(\mathcal BJ_A)(U)$ \Mark{is the agreed-upon notation for the $B$ in  $BJ$?} denotes the groupoid of maps $\Hom_{C \times A}(C\times U, BJ_A)$. By applying \cite[Thm. 1.3]{hall-rydh-tannakahom}, with $S=A$, $Z=C \times A$ and $X=\mathcal BJ_A$, the Weil restriction (and so $\stP$) is an algebraic stack locally of finite type over $A$.\\
% \underline{Smoothness.} Fix a geometric point $a$ in $A$. Since $J_a$ is a smooth group scheme over the curve $C$, we have $H^2(C,\mathrm{Lie}(J_a))=0$. Hence, every deformation of a $J_a$-torsor is unobstructed  so that $\stP_a$ is smooth. Hence, $\stP$ is smooth over $A$ \Mark{why? don't we need some flatness?}\andres{I agree, you need to show the formal criterion for smoothness also for Artin rings with nontrivial morphsims to $A$ (this also follows from the same cohomological computation, but it is not enough to say that $\stP_a$ is smooth. I would just cite \cite[Prop. 1.1]{heinloth-uniformization}}.\\
We are left to check the dimension. Fix a geometric point $a$ in $A$. We have the following equalities
\begin{align*}
\dim\stP_a=&-\chi(\mathrm{Lie}(J_a))=-\chi(\oplus_{i=1}^r\cL^{1-e_i})=\left(\sum_{i=1}^r(e_i-1)\right)\deg(\cL)+r(g-1)=\\
=&(\dim(B)-\dim(T))\deg(\cL)+\dim(T)(g-1)=\\
=&\frac{1}{2}\dim(G) \deg(\cL) - \frac{1}{2} \dim(T) (\deg(\cL)-2g+2).
\end{align*}
The first equality follows from the smoothness of the stack, the second one from  \Cref{P:galois-descr-j-cam}\eqref{P:galois-descr-j-cam3}, the third one from Riemann--Roch. The fourth equality  follows from the three equalities $\sum e_i=|\Phi|/2+r$ (see Remark \ref{R:h}),  $|\Phi|=2(\dim(B)-\dim(T)),$ and $r=\dim(T)$. The fifth and last  equality follows easily from the identity $\dim(G)=2\dim(B)-\dim(T)$.
\end{proof}

\begin{prop}\label{P:centr-group}
Assume $\deg(\cL)>0$ and $g(C)>0$. Then, the following statements hold:
	\begin{enumerate}[(i)]
		\item\label{P:centr-group2} The moduli stack $\stP\to A$ admits a good moduli space $\schP\to A$, which is a quasi-separated smooth 
 commutative group algebraic  space of relative  dimension $\frac{1}{2}\dim(G) \deg(\cL) - \frac{1}{2} \dim(T) (\deg(\cL)-2g+2)+\dim(Z_G)$. Furthermore, the canonical morphism
  $\stP\to\schP$ is the rigidification of $\stP$ with respect to the flat group scheme $Z_G\times A\to A$ (in the sense of \cite[Definition 5.1.9]{acv-twisted-bundles-covers}). In particular, it is a  $Z_G$-gerbe.
		\item \label{P:centr-group3} 
  There is a unique open group substack $\stP^o \subset \stP$ with geometrically connected fibers over $A$. Its $Z_G$-rigidification $\schP^o := \stP^o \fatslash Z_G$ is an open subgroup space $\schP^o \subset \schP$ with geometrically connected fibers over $A$. We have that $\stP^o \to \schP^o$ is a good moduli space. Furthermore, if $\deg(\cL)>g+\frac{1}{2}$, then the generic fiber of $\schP^o \to A$ is an Abelian variety.

  \item \label{P:centr-group4} For any $d \in \pi_1(G)$, the action of $\stP^o$ on $\Higgs_{G,\cL}$ preserves the open and closed substack $\Higgs_{G,\cL}^d \subset \Higgs_{G,\cL}$.
	\end{enumerate}
\end{prop}

\begin{proof}
\underline{Point \eqref{P:centr-group2}} Fix a $k$-scheme $U$ and a $U$-point in $A(U)$. The automorphism group scheme of a $J_U$-torsor over $p:C\times U\to U$ is equal to the Weil restriction $p_*(J_a)$, which is equal to the group scheme $Z_G\times U\to U$ (see \Cref{P:Weil-Restr} below). We denote by $\schP:=\stP\fatslash(Z_G\times A)$ the rigidification of the stack $\stP\to A$ along the flat group scheme $Z_G\times A\to A$. By the above discussion, the automorphism group scheme of any object in $\schP$ is trivial. Hence, $\schP$ is an algebraic space.

By \cite{acv-twisted-bundles-covers}, there exists a canonical map
$\nu:\stP\to \schP$, which is a $Z_G$-gerbe. In particular, \'etale locally over $\schP$, 
it is equal to $\mathcal B(Z_G)\times V= \mathcal B(Z_G\times V)\to V$. It follows that 
$\nu$ is a good moduli space morphism, because this property can be checked flat locally 
on $\schP$ (see \cite[Prop. 4.7]{alper-good-moduli}) and $\mathcal 
B(Z_G\times V)\to V$ is a good moduli space (because $Z_G$ is diagonalizable, see 
\cite[I 5.3.3]{sga3} and \cite[Example 12.4]{alper-good-moduli}). It also follows from 
the local description of $\nu$ that $\stP \to \schP$ flat. Hence the rigidification 
$\schP$ is smooth over $A$, because $\stP$ is smooth and smoothness can be checked 
flat locally. Moreover, by the universal property of good moduli spaces \cite[Thm. 6.6.]
{alper-good-moduli}, the commutative groups stack structure on $\stP$ descends to show 
that $\schP$ is a commutative group algebraic space over $A$.

 \underline{Point \eqref{P:centr-group3}}. 
 The existence of the neutral component $\stP^o$ follows from \cite[Lem. 3.1]{herrero2023semistable}, and the fact that the rigidification $\schP^o \subset \schP$ is an open subgroup space with geometrically connected $A$-fibers follows immediately from the properties of $\stP^o \subset \stP$. The fact that the generic fiber of $\schP^o$ is an Abelian variety follows from the fact that $P^o\to A$ restricted to the open subset $A^\diamondsuit\subset A$ is an Abelian scheme (see \cite[Prop. 4.7.7]{ngo-lemme-fondamental}) and that $A^{\diamondsuit}$ is non-empty if $2\deg(\cL)>2g+1$ (see \Cref{P:dh-sets1}).

 \underline{Point \eqref{P:centr-group4}} Since the fibers of $\stP^o \to A$ are geometrically connected, it follows immediately that the action of $\stP^o$ preserves any open and closed substack of $ \Higgs_{G,\cL}$.

\end{proof}

\begin{prop}\label{P:Weil-Restr}Assume that $\deg(\cL)>0$. The Weil restriction $(pr_A)_*(J_A)$ is isomorphic to the constant group scheme $Z_G\times A$ over $A$. On the other hand, the Weil restriction $(pr_A)_*(J^1_A)$ is isomorphic to the constant group scheme $T^W\times A$ over $A$, where $T^W$ is the subgroup of $W$-invariant elements of the torus $T$. 
\end{prop}
\begin{proof}Fix a morphism $U\to A$. By definition of the local regular centralizer group scheme $J$ over $\fc$ (see Subsection \ref{SS:local-centr-gr} and \Cref{P:chain-j-cam}), we have a chain of inclusions
	$$
	Z_G\times (C\times U)\hookrightarrow J_U\hookrightarrow J_U^1
	$$
 of group schemes over $C\times U$. Taking the Weil restriction 
 along the projection $pr_U:C\times U\to U$, we get an inclusion of group schemes over $U$
	\begin{equation}\label{E:chain.pf}
	(pr_U)_*(Z_G\times C\times U)=Z_G\times U\hookrightarrow (pr_U)_*J_U\hookrightarrow (pr_U)_*J_U^1
	\end{equation}
Let  $\eta$ denote the composition $\cC_U\xrightarrow{\pi_U} C\times U\xrightarrow{pr_U} U.$ By \Cref{P:galois-descr-j-cam}\eqref{P:galois-descr-j-cam2}, the group scheme $(pr_U)_*J_U^1$ is equal to the subgroup scheme of $W$-invariant elements in the group scheme  
 $\eta_*(\cC_U\times T)\to U$. By  \Cref{P:cameral-curve-conn}, we have $\eta_*(C_U\times T)=U\times T$ and, so, $(pr_U)_*J_U^1=U\times T^{W}$. This proves the second assertion.

 We now focus on the first assertion. By what we have proved before, we have a chain of inclusions
 $$
 Z_G\times U\hookrightarrow (pr_U)_*J_U\hookrightarrow T^W\times U.
 $$
 Since we are assuming $\deg(\cL)>0$,  \Cref{L:a-meet-d} implies that, for any geometric point $u$ in $U$, the corresponding section $u:C\to \fc_\cL$ meets all the irreducible components of $\DD_\cL\subset \fc_\cL$. Note that the equation of $\DD$ is given by the product of all the roots (see \Cref{T:konstant}\eqref{T:konstant4}). Hence, the locus in $\cC_U^{\alpha}\subset \cC_U$ of fixed points, with respect to the action of the hyperplane-root reflection $s_\alpha$, is non-empty, for any root $\alpha$. Combining this fact with the description of the image of $J_U$ in $J_U^1$ (see \Cref{P:chain-j-cam}\eqref{P:chain-j-cam3}), we get the following equality
 $$
 (pr_U)_*J_U\cong \bigcap_{\alpha\in\Phi}\ker\left(T^W\hookrightarrow T\xrightarrow{\alpha}\bG_m\right)\times U.
 $$
 Since the roots are the non-trivial weights of the adjoint representation $\mathrm{ad}:G\to \GL(\fg)$ and the kernel of the representation is equal to the center $Z_G$, we must have that the right-hand side in the isomorphism above is equal to $Z_G\times U$. This concludes the proof.
\end{proof}

Let $n: \widetilde{C} \to C$ be a finite \'etale Galois cover and let 
$n_\Higgs$ be the morphism induced by pullback at the level of stacks of Higgs bundles; see  
\Cref{notn: pullback higgs under cover} and \Cref{lemma: pullback morphism on higgs moduli stacks}.
 Let $n_A: A_C \to A_{\widetilde{C}}$ be the morphism induced by $n$
by pullback on the  Hitchin bases for $\cL$ on $C$, and for  $n^*(\cL)$ on $\widetilde{C}$ respectively. 
The morphisms $n_\Higgs$, $n_A$ and the Hitchin morphisms for $\cL$ and $\widetilde{\cL}$ Higgs bundles give rise to the commutative square diagram
    \[
\begin{tikzcd}
  \Higgs_{G,\cL} \ar[r, "n_{\Higgs}"] \ar[d, "h"] & \Higgs_{G,n^*(\cL)} \ar[d, "h"] \\ A_C \ar[r, "n_A"] & A_{\widetilde{C}}.
\end{tikzcd}
\]
Let $\stP_{C} \to A_{C}$ (resp. $\stP_{\widetilde{C}} \to A_{\widetilde{C}}$) denote the corresponding commutative group stack of torsors for $J_{A_{C}}$ (resp. $J_{A_{\widetilde{C}}}$). Recall that  the stacks of torsors $\stP_{C}$ (resp. 
    $\stP_{\widetilde{C}}$) acts naturally on the corresponding Hitchin fibration 
    $\Higgs_{G,\cL} \to A_C$ (resp. $\Higgs_{G,n^*(\cL)}\to A_{\widetilde{C}}$); see
    \Cref{prop: torsor regular higgs under group stack}. 
    As we see in the next proposition, the morphism $n$
    induces  by pullback a morphism $n_\stP: \stP_{C} \to n_A^*(\stP_{\widetilde{C}}):= \stP_{\widetilde{C}} \times_{A_{\widetilde{C}}} A_C$
    of commutative group stacks  over $\hitchinbase$.

\begin{prop} \label{prop: stacky pullback morphism group scheme}
The morphism  $n$ induces by pullback a natural affine morphism $n_\stP: \stP_{C} \to n_A^*(\stP_{\widetilde{C}})$. By the universal property of fiber products, this induces a morphism $n_{\stP}: \stP_{C} \to \stP_{\widetilde{C}}$ over $n_A$, which fits into the following
commutative diagram, where the horizontal morphisms are the respective actions
    \begin{equation} \label{diagram: stacky cover action}
        \begin{tikzcd}[ampersand replacement=\&]
         \stP_{C} \times_{A_{C}} \Higgs_{G,\cL} \ar[d,"n_{\stP}\times n_\Higgs"] \ar[r,"{\rm act}"]\& \Higgs_{G,\cL} \arrow[d,"n_\Higgs"]\\
        \stP_{\widetilde{C}} \times_{A_{\widetilde{C}}} \Higgs_{G,n^*(\cL)} \ar[r,"{\rm act}"] \& \Higgs_{G,n^*(\cL)}.
        \end{tikzcd}
    \end{equation}	
\end{prop}
\begin{proof}
We have the commutative affine 
    smooth group scheme $({\rm id}_{\widetilde{C}} \times n_A)^*(J_{A_{\widetilde{C}}})$ over $\widetilde{C} \times A_C$.  
    Note that $n_A^*(\stP_{\widetilde{C}})$ is the stack over   
    $A_C$  of torsors  for the 
    group scheme $({\rm id}_{\widetilde{C}}\times n_A)^*(J_{A_{\widetilde{C}}})$, i.e. $n_A^*(\stP_{\widetilde{C}})(S\to A_C)$
    is the groupoid of said torsors over $\widetilde{C}\times S.$
    It follows from the 
    definition of global regular centralizer group that 
    $({\rm id}_{\widetilde{C}}\times n_A)^*(J_{A_{\widetilde{C}}})$ is 
    canonically isomorphic to the pullback $(n\times 
    {\rm id}_{A_C})^*(J_{A_C})$ via the morphism 
    $n \times {{\rm id}_{A_C}}: 
    \widetilde{C} \times A_{C} \to C \times A_{C}$.  We 
    conclude that $n_A^*\stP_{\widetilde{C}}$ can also be described as the 
    stack of torsors for  the pullback group scheme 
    $(n \times {{\rm id}_{A_C}})^*(J_{A_C})$ under the finite \'etale 
    Galois cover $n \times {{\rm id}_{A_C}}: \widetilde{C} \times 
    A_{C} \to C\times A_{C}$. Pulling back $J_{A_C}$-torsors via this cover yields a morphism
    \[ n_\stP: \stP_C \to n_A^*(\stP_{\widetilde{C}}), \; 
    \; \; F \mapsto (n\times {\rm id})^*(F).\]
    The compatibility of this morphism with the action 
    on $\Higgs_{G,\cL}$ 
    %\sout{Higgs} \Mirko{shall we replace with $\Higgs_{G,\cL}$?}\Roberto{Just for being sure: before there was "$\mathrm{Higgs}$" instead of "$\Higgs$"?}\Mirko{exactly. If correct let's remove \sout{Higgs} and comments} 
    as in diagram \eqref{diagram: stacky cover action} follows directly from the definitions.

We are only left to check that $n_\stP$ is an affine morphism. Let $\Gamma$ denote the Galois group $\Gamma$ of the cover $n: \widetilde{C} \to C$ acts on the stack $n_A^*(\stP_{\widetilde{C}})$ of $(n\times {\rm id}_A)^*(J_{A_C})$-torsors by pulling back via the isomorphisms $\gamma \times {\rm id}: \widetilde{C} \times A_C \xrightarrow{\sim} \widetilde{C} \times A_C$,
    $\gamma \in \Gamma$, as follows:
    \[ \Gamma \times n_A^*(\stP_{\widetilde{C}}) \to 
    n_A^*(\stP_{\widetilde{C}}),\; \; \; F \mapsto (\gamma \times {\rm id})^*(F).\]
    \'Etale descent for the Galois cover $n\times {\rm id}: \widetilde{C} \times A_C \to \widetilde{C} \times A_C$ implies that the pullback morphism $n_\stP: \stP_C \to n_A^*(\stP_{\widetilde{C}})$ exhibits $\stP_C$ as the stack of $\Gamma$-fixed points $\left(n_A^*(\stP_{\widetilde{C}})\right)^{\Gamma} \to n_A^*(\stP_{\widetilde{C}})$ as in \cite[Defn. 2.3]{romagny-actions-stacks} (cf. the proof in \cite[Prop. 5.4.2.4]{lurie-gaitsgory-tamagawabook}). Since the stack $n_A^*(\stP_{\widetilde{C}})$ has affine diagonal (by \Cref{P:centr-group}), it follows that the morphism $n_\stP: \stP_C \to n_A^*(\stP_{\widetilde{C}})$ is affine \cite[Prop. 3.7]{romagny-actions-stacks}.
\end{proof}

\begin{subsection}{Action on the semistable locus} \label{subsection: action on the semistable locus untwisted}

For the next definition, we use the quotient stack $\Theta =  \Theta_k:= [\mathbb{A}^1_k /\mathbb{G}_m]$. The stack $\Theta$ has exactly two $k$-points, the open point $1$ (corresponding to the orbit of $1 \in \mathbb{A}^1_k(k)$) and the closed point $0$ (corresponding to the origin of $\mathbb{A}^1_k$). Note that the closed point $0$ has automorphism group $\mathbb{G}_m$.

\begin{defn}[\cite{heinloth-hilbertmumford, halpernleistner2018structure}] \label{defn: theta semistability}
Let $\scrM$ be an algebraic stack, and fix a line bundle $\cK$ on $\scrM$. A geometric point $x: \Spec(F) \to \scrM$ defined over some algebraically closed field $F \supset k$ is called $\Theta$-semistable with respect to $\cK$ if for all morphisms $f: \Theta_F \to \scrM$ such that $f(1) \cong x$ (a.k.a. filtrations of $x$), the $\mathbb{G}_m$-weight of the $0$-fiber of the pullback $f^*(\cK)$ satisfies $wt(f^*(\cK)|_{0}) \leq 0$.
\end{defn}

\begin{prop} \label{prop: general prop action of connected stack semistable locus}
Let $\scrM \to B$ be an algebraic stack locally of finite presentation over a Noetherian scheme $B$. Let $\scrX \to B$ be an $B$-group algebraic stack with geometrically connected $B$-fibers, equipped with an action $\scrX \times_{B} \scrM \to \scrM$. Suppose that the locus of $\Theta$-semistable geometric points with respect to a given line bundle $\mathcal{K} \in \Pic(\scrM)$ forms an open substack $\scrM^{ss} \subset \scrM$. Then, the action of $\scrX$ preserves $\scrM^{ss}$, i.e. there is a factorization
\[
\begin{tikzcd}
  \scrX \times_{B} \scrM^{ss} \ar[r, dashed] \ar[d, symbol= \hookrightarrow] & \scrM^{ss} \ar[d, symbol = \hookrightarrow] \\ \scrX \times_{B} \scrM \ar[r] & \scrM.
\end{tikzcd}
\]
\end{prop}
\begin{proof}
The factorization through an open substack $\scrM^{ss} \subset \scrM$ can be checked at the level of geometric points. Let $x \to \scrM$ be a geometric point with image $b$ in $B$. After base-change, we may assume that $b=B = \Spec(k)$ is the spectrum of an algebraically closed field. It is enough to show that for all $k$-points $g$ of $\scrX$, we have that $x$ is semistable if and only if $g \cdot x$ is semistable.

The action morphism $g\cdot(-): \scrM \to \scrM$ sends a filtration $f: \Theta_k \to \scrM$ of $x$ to a filtration $g \cdot f$ of $g \cdot x$. By the definition of $\Theta$-semistability, it suffices to show that for all filtrations $f: \Theta_k \to \scrM$ of $x$, the $\mathbb{G}_m$-weight $wt(f^*(\mathcal{K})|_0)$ is the same as the $\mathbb{G}_m$-weight $wt((g\cdot f)^*(\mathcal{K})|_0)$. In order to see this, consider the action morphism
\[ 
\widetilde{f}: \scrX  \times  \Theta_k
\xrightarrow{{\rm id} \times f} \scrX \times \scrM \to \scrM.\]
We have that $\widetilde{f}|_{1 \times \Theta_k} =f$ and that $\widetilde{f}|_{g \times   \Theta_k } = g \cdot f$. If we restrict under the morphism $\scrX 
 \times  0 \to  \scrX \times \Theta_k$, then we obtain a map $\widetilde{f}_0: \scrX \to \scrM$. The pullback $\widetilde{f}_0^*(\mathcal{K})$ is a $\mathbb{G}_m$-equivariant line bundle on the connected stack $\scrX$. Since the $\mathbb{G}_m$-weight of a line bundle is locally constant, we conclude that the $\widetilde{f}_0^*(\mathcal{K})$ agrees for all points of $\scrX$. By restricting to $1$ and $g$ we obtain our desired equality $wt(f^*(\mathcal{K})|_0)=wt((g\cdot f)^*(\mathcal{K})|_0)$.
\end{proof}

\begin{coroll} \label{coroll: action on the semistable moduli space}
Suppose that $\text{char}(k)$ is either zero, or strictly larger than the height of the adjoint representation. Then, the action of $\stP^o$ on $\Higgs_{G,\mathcal{L}}^{d}$ (see \Cref{P:centr-group}(iii)) preserves the semistable locus $\Higgs_{G, \cL}^{d,ss} \subset \Higgs_{G, \cL}^{d}$,  so that it descends to an action of $\schP^{o}$ on $\MHiggs_{G,\mathcal{L}}^{d}$ over $A$.
\end{coroll}
\begin{proof}
By \cite[\S4.1]{herrero2023meromorphic}, the locus of semistable Higgs bundles $\Higgs_{G, \cL}^{d,ss} \subset \Higgs_{G, \cL}^{d}$ is the $\Theta$-semistable locus for a line bundle on $\Higgs_{G, \cL}^{d}$. Hence, the corollary follows by applying \Cref{prop: general prop action of connected stack semistable locus} to the action of the commutative group stack $\stP^o \to A$ on $\Higgs_{G, \cL}^{d}$.
\end{proof}

\begin{lemma} \label{lemma: po preserves stable locus}
    The action of $\stP^o$ on $\Higgs_{G, \cL}^d$ preserves the stable locus $\Higgs_{G, \cL}^{d,s}$.
\end{lemma}
\begin{proof}
    This follows from the proof of \cite[Lem. 4.8]{herrero2023semistable}.
\end{proof}

For the next lemma, we note that for every scheme $S$ and any morphism $a: S \to A$, the base-change $\left(\Higgs_{G, \cL}^{ss}\right)_a$ admits an adequate moduli space $\left(\Higgs_{G, \cL}^{ss}\right)_a \to M_a$ by \cite[5.2.9(3)]{alper_adequate}.

\begin{lemma} \label{lemma: saturated orbit of group scheme}
    Let $S$ be a scheme, and let $x: S \to \Higgs_{G, \cL}^{ss}$ be an $S$-point that corresponds to a regular $G$-Higgs bundle (as in \Cref{defn: regular higgs bundle}). Let $a \in A(S)$ be the image in the Hitchin base. Then, the following hold.
    
    \begin{enumerate}[(1)]
    \item The orbit morphism $x \times_{S} \stP^o_{a} \to \left(\Higgs_{G, \cL}^{ss}\right)_a$ is an open immersion into the Hitchin fiber $\left(\Higgs_{G, \cL}^{ss}\right)_a$. 
    
    \item If $x$ lands in the locus of stable Higgs bundles, then the open image $\stP^o_{a} \cdot x \subset \left(\Higgs_{G, \cL}^{ss}\right)_a$ of the orbit morphism is saturated with respect to the adequate moduli space morphism $\left(\Higgs_{G, \cL}^{ss}\right)_a \to M_a$. Furthermore, $\stP^o_{h(x)} \cdot x$ is a $Z_G$-gerbe over an open subscheme of $M_a$.
    \end{enumerate}
\end{lemma}
\begin{proof}
    The fact that the orbit morphism is an open immersion follows from \Cref{prop: torsor regular higgs under group stack}. Since $\stP^o_{a}$ is a $Z_G$-gerbe over its adequate moduli space {\cb by} \Cref{P:centr-group}, the same holds for the isomorphic open image $\stP^o_{a} \cdot x$ under the orbit morphism. 
    
    Suppose now that $x$ lands in the stable locus. By \Cref{lemma: po preserves stable locus}, it follows that the open substack $\stP^o_{a} \cdot x$ lies in the stable locus. Hence, $\stP^o_{a} \cdot x$ is saturated with respect to the composition  $\left(\Higgs_{G, \cL}^{ss}\right)_a \to M_a \to \left(\MHiggs_{G,\cL}\right)_a$ by \Cref{lemma: stable higgs is saturated}. Since the morphism $M_a \to \left(\MHiggs_{G,\cL}\right)_a$ is a universal homeomorphism \cite[Prop. 5.2.9(3)]{alper_adequate}, it follows that the open substack $\stP^o_{a} \cdot x$ is saturated with respect to the morphism $\left(\Higgs_{G, \cL}^{ss}\right)_a \to M_a$. The adequate moduli space of $\stP^o_{a} \cdot x$ is its open image in $M_a$, and we have already established that it is a $Z_G$-gerbe over it.
\end{proof}

The rest of this section is dedicated to showing that
the group scheme $\schP^o$ over $\hitchinbase$ is quasiprojective and it acts with affine stabilizers on $\MHiggs_{G, \cL}^d$. 

First, we prove the affineness of stabilizers in a special case.

\begin{lemma} \label{lemma: affine stabilizers special case}
    Suppose that the line bundle $\cL$ admits a square root in $\Pic(C)$, that $\text{deg}(\cL)>0$, and that $\text{char}(k)$ is either zero, or strictly larger than the height of the adjoint representation. Set $d=0 \in \pi_1(G)$. Then for any geometric point $x \to \MHiggs_{G, \cL}^0$ with image $h(x) \to A$, the stabilizer $\Stab_{\schP^o_{h(x)}}(x)$ is an affine group scheme.
\end{lemma}
\begin{proof}
    By \cite[Lem. 6.1]{migliorini-shende} and the connectedness of the fibers of the Hitchin morphism \Cref{cor: higgs moduli space properties}, it suffices to show that there is a section $A \to \MHiggs_{G, \cL}^0$ such that the stabilizer of every point in the image is affine. We take the composition $s: A \xrightarrow{\kappa} \Higgs_{G,\cL}^{0,ss} \to \MHiggs_{G, \cL}^0$, where $\kappa$ is the Kostant section. By the stability of the Kostant section (\Cref{lemma: kostant section is stable}), it follows that we can apply \Cref{lemma: saturated orbit of group scheme} to obtain a saturated open substack $\mathcal{U} \subset \Higgs_{G,\cL}^{0,ss}$ containing the image of $\kappa$. Furthermore, by \Cref{lemma: saturated orbit of group scheme}, $\mathcal{U}$ is a gerbe over an open subscheme $U \subset \MHiggs_{G, \cL}^0$ isomorphic to $\schP^o$ with its free right action. We conclude that all the stabilizers of $\schP^o$ of points in the image of $\kappa$ are trivial, hence affine.
\end{proof}

To generalize to the case when $\cL$ does not admit a square root, we shall use the following.

\begin{lemma} \label{lemma: pullback group scheme galois cover}
Suppose that $\text{char}(k)$ is either zero, or strictly larger than the height of the adjoint representation, and that $\text{deg}(\cL)>0$. Let $n: \widetilde{C} \to C$ be a finite \'etale Galois cover of some degree $N$. Then the pullback morphism in \Cref{prop: stacky pullback morphism group scheme} induces an affine morphism of group algebraic spaces $n_{\schP}: \schP^o_{C} \to n_A^*(\schP^o_{\widetilde{C}})$. The induced  pullback morphism $n_{\MHiggs}: \MHiggs_{G,\cL}^d \to \MHiggs_{G,n^*(\cL)}^{N \cdot d}$ (as in \Cref{lemma: pullback morphism on higgs moduli stacks}) is $\schP^o$-equivariant.
\end{lemma}
\begin{proof}
    The morphism $n_{\stP}: \stP_{C} \to n_A^*(\stP_{\widetilde{C}})$ is affine (\Cref{prop: stacky pullback morphism group scheme}), and so from \cite[Lem. 4.14]{alper-good-moduli} it follows that it induces an affine morphism of the good moduli spaces $n_{\schP}: \schP_{C} \to n_A^*(\schP_{\widetilde{C}})$ granted by \Cref{P:centr-group}. By passing to relative neutral components, we obtain an affine morphism $n_{\schP}: \schP^o_{C} \to n_A^*(\schP^o_{\widetilde{C}})$. The equivariance of the morphism of Higgs moduli spaces $n_{\MHiggs}: \MHiggs_{G,\cL}^d \to \MHiggs_{G,n^*(\cL)}^{N \cdot d}$ is a direct consequence of diagram \eqref{diagram: stacky cover action} after passing to neutral components and semistable loci, and taking adequate moduli spaces.
\end{proof}

As a corollary of the proof of \Cref{lemma: affine stabilizers special case}, we obtain the quasiprojectivity of $\schP^o$.
\begin{coroll} \label{coroll: quasiprojectivity of algebraic space}
    Suppose that $\text{deg}(\cL)>0$, and that $\text{char}(k)$ is either zero, or strictly larger than the height of the adjoint representation. Then the smooth commutative group algebraic space $\schP^o \to A$ is a quasiprojective group scheme.
\end{coroll}
\begin{proof}
    If $\cL$ admits a square root, then the quasiprojectivity follows from the proof of \Cref{lemma: affine stabilizers special case}, which exhibits $\schP^o$ as an open subscheme of the quasiprojective (cf. \Cref{cor: higgs moduli space properties}) scheme $\MHiggs_{G, \cL}^0$.

    If $\cL$ does not admit a square root, then there is a finite Galois cover $n:\widetilde{C} \to C$ such that the pullback $n^*(\cL)$ admits a square root. Then the group algebraic space $\schP_{\widetilde{C}}^o$ for the curve $\widetilde{C}$ is quasiprojective. Let $n_A: A_{C} \to A_{\widetilde{C}}$ denote the pullback morphism on Hitchin bases. By \Cref{lemma: pullback group scheme galois cover}, the pullback morphism on group algebraic spaces $n_{\schP}: \schP_C^o \to n_A^*(\schP_{\widetilde{C}}^o)$ is affine, and hence it follows that $\schP_C^o$ is also quasiprojective.
\end{proof}

\begin{lemma} \label{lemma: divisibility of connected components}
Suppose that $k$ is algebraically closed and that $\text{char}(k)$ is either zero, or strictly larger than the height of the adjoint representation. There exists an integer $H$ depending only on $G$ such that for all $d \in \pi_1(G)$ and all $n \in \mathbb{Z}$, there is a $\schP^o$-equivariant isomorphism $\MHiggs_{G,\cL}^{n H \cdot d }\cong \MHiggs_{G,\cL}^{0}$.
\end{lemma}
\begin{proof}
    Recall the natural action of the commutative group stack $\Bun_{Z_G^o}$ on $\Higgs_{G, \cL}$ on $\Higgs_{G,\cL}$ described in \Cref{remark: action of torsors for the center}. The inclusion $Z_G^o \times C \times A \hookrightarrow J_{A}$ induces a morphism of commutative group stacks $\Bun_{Z_G^o} \times A \to \stP$ over $A$. It follows from definitions that the natural induced action of $\Bun_{Z_G^o} \times A$ agrees with the one obtained by using the morphism to $\stP$ and acting via $\stP$. Since $\stP$ is a commutative group stack, if follows that for all $k$-points $x: \Spec(k) \to \Bun_{Z_G^o}$, the induced isomorphism $x \cdot(-): \Higgs_{G, \cL} \to \Higgs_{G,\cL}$
    is $\stP$-equivariant. Since $x \cdot (-)$ preserves the semistable locus (\Cref{remark: action of torsors for the center}), this induces an isomorphism at the level of moduli spaces $x \cdot(-): \MHiggs_{G, \cL} \to \MHiggs_{G,\cL}$
    which is $\schP^o$-equivariant.

    Recall that the connected components of $\Bun_{Z_G^o}$ are indexed by $\pi_1(Z_G^o)$. The inclusion $Z_G^o \hookrightarrow G$ induces a group homomorphism 
    \[ \iota: \pi_1(Z_G^o) = X_*(Z_G^o) \to X_*(T) \to X_*(T)/X_{coroots} = \pi_1(G)\]
    with finite cokernel. Hence, there exists an integer $H$ such that for all $d \in \pi_1(G)$, the multiple $d \cdot H$ is in the image of $\iota$. If we choose a preimage $d' \in \pi_1(Z_G^o)$ of $nH\cdot d$ and a $k$-point $x: \Spec(k) \to \Bun_{Z_G^o}^{d'}$, then the morphism $x\cdot (-)$ restricts to a $\schP^o$-equivariant isomorphism  of connected components $x \cdot(-): \MHiggs_{G,\cL}^{0} \to \MHiggs_{G,\cL}^{nH\cdot d}$.
\end{proof}

By using the same covering trick as in \Cref{coroll: quasiprojectivity of algebraic space}, we get the following generalization of \Cref{lemma: affine stabilizers special case}.
\begin{prop} \label{prop: afineness of stabilizers general}
    Suppose that $\text{deg}(\cL)>0$, and that $\text{char}(k)$ is either zero, or strictly larger than the height of the adjoint representation. Let $x$ be a geometric point of $\MHiggs_{G, \cL}^d$. Then the stabilizer $\Stab_{\schP^o_{h(x)}}(x)$ is an affine group scheme.
\end{prop}
\begin{proof}
    Let $n: \widetilde{C} \to C$ be a finite \'etale 
    Galois cover of some degree $N$, inducing a 
    pullback homomorphism of Hitchin bases $n_A: 
    A_{C} \to A_{\widetilde{C}}$. We get affine 
    morphisms $n_{\MHiggs}: \MHiggs_{G,\cL}^{d} \to 
    \MHiggs_{G,n^*(\cL)}^{N \cdot d}$ and $n_{\schP}: 
    \schP_C^o \to n_A^*(\schP_{\widetilde{C}}^o)$ 
    that are compatible with the corresponding action 
    morphisms \Cref{lemma: pullback group scheme 
    galois cover}. By \Cref{lemma: divisibility of 
    connected components}, if $N$ is sufficiently 
    divisible, the moduli space $\MHiggs_{G,n^*
    (\cL)}^{N \cdot d}$ is $\schP_{\widetilde{C}}^o$-
    equivariantly isomorphic to $\MHiggs_{G,n^*
    (\cL)}^{0}$. Furthermore, by possibly enlarging 
    the Galois cover $\widetilde{C}$, we may assume 
    that $n^*(\cL)$ admits a square root. Therefore, 
    by \Cref{lemma: affine stabilizers special case}, 
    it follows that $\schP_{\widetilde{C}}^o$ acts 
    with affine stabilizers on $\MHiggs_{G,n^*(\cL)}^{0} \cong \MHiggs_{G,n^*(\cL)}^{N \cdot d}$. 
    Since the morphism $n_{\schP}: \schP_C^o \to n_A^*(\schP_{\widetilde{C}}^o)$ is affine and $n_{\MHiggs}: \MHiggs_{G,\cL}^{d} \to \MHiggs_{G,n^*(\cL)}^{N \cdot d}$ is equivariant, it follows that the stabilizers of $\schP_C^o$ acting on $\MHiggs_{G,\cL}^{d}$ are also affine, as desired.
\end{proof}
\end{subsection}

\begin{subsection}{Polarizability of the group of symmetries of the Hitchin fibration}

In this subsection, we show that the Tate module of the group scheme $\schP^o\to A$ is polarizable (see \Cref{C:TateModule-Polarizable}).

\begin{defn}
	Let $S$ be a scheme. Let $\pi: G \to S$ be a smooth commutative group scheme
	over $S$. Assume that $\pi$ has connected fibers of dimension $d$. Let $\ell$ be a prime invertible in $S.$ 
	The Tate module
	of $G/S$ is the sheaf on $S$ defined by $T_{\overline{\mathbb{Q}}_\ell}(G)\coloneqq
	R^{2d-1}g_! {\overline{\mathbb{Q}}_\ell}(d).$
\end{defn}

The fiber
at any geometric points $s$ of $S$ is
naturally isomorphic to the Tate module of the fiber $G_s/s,$
i.e. $(T_{\overline{\mathbb{Q}}_\ell} (G))_s= T_{\overline{\mathbb{Q}}_\ell}(G_s).$ Moreover, we have the following two alternative interpretations of the Tate module of $G_s$: (1) $T_{\overline{\mathbb{Q}}_\ell}(G_s)=H_1 (G_s, {\overline{\mathbb{Q}}_\ell});$
(2) $T_{\overline{\mathbb{Q}}_\ell}(G_s) = (\varprojlim_n G_s[l^n])\otimes_{\mathbb Z_\ell}
{\overline{\mathbb{Q}}_\ell}.$

\begin{defn} Let $G$ be a commutative group scheme over a field $k$. We say that $G/k$ admits a Chevalley d\'evissage if there exists an exact sequence of group schemes over $k$
	\begin{equation}\label{E:Weak-Chev-dec}
		1\to G^{\mathrm{aff}}\to G\to G^{\mathrm{ab}}\to 1,
	\end{equation}
	where $G^{\mathrm{ab}}$ is an Abelian variety over $k$ and $G^{\mathrm{aff}}$ is a connected affine group scheme over $k$.
\end{defn}

\begin{thm}\label{T:sCd-kclosed} If $G$ is a connected smooth commutative group scheme over a perfect field $k$. Then, $G$ admits a Chevalley d\'evissage and its formation is compatible with extension of the base field.
\end{thm}

\begin{proof}See \cite[\S 9.2, Thm. 1]{blr-neron}.
\end{proof}

Let $G/S$ be a smooth commutative group scheme with connected fibers over $S$. By \Cref{T:sCd-kclosed}, the fiber $G_s$ over a geometric point $s$ of $S$ admits a Chevalley d\'evissage. In particular, the exact sequence \eqref{E:Weak-Chev-dec} induces an exact sequence of Tate modules
$$
0\to T_{\overline{\mathbb Q}_\ell}(G_s^{\mathrm{aff}})\to T_{\overline{\mathbb Q}_\ell}(G_s)\to T_{\overline{\mathbb Q}_\ell}(G_s^{\mathrm{ab}})\to 0.
$$
We call the left-hand side the affine submodule of $T_{\overline{\mathbb Q}_\ell}(G_s)$ and the right-hand side the proper quotient of $T_{\overline{\mathbb Q}_\ell}(G_s)$. 

\begin{defn}\label{defn:polarizable tate module} Let $G/S$ be a smooth commutative group scheme with connected fibers. We say that its Tate module $T_{\overline{\mathbb Q}_\ell}(G)$ is polarizable if it admits a polarization, i.e. a skew-symmetric bilinear pairing
	$$
	e:T_{\overline{\mathbb Q}_\ell}(G)\otimes_{\overline{\bQ}_\ell}T_{\overline{\mathbb Q}_\ell}(G)\to T_{\overline{\mathbb Q}_\ell}(\bG_{m,S}),
	$$
	such that, for any geometric point $s$ of $S$, the pull-back $e_s:T_{\overline{\mathbb Q}_\ell}(G_s)\otimes_{\overline{\bQ}_\ell}T_{\overline{\mathbb Q}_\ell}(G_s)\to T_{\overline{\mathbb Q}_\ell}(\bG_{m,s})$ over the geometric fiber $G_s/s$ is trivial on the affine submodule of $T_{\overline{\mathbb Q}_\ell}(G_s),$ and the induced pairing on the proper quotient of $T_{\overline{\mathbb Q}_\ell}(G_s)$ is non-degenerate.
\end{defn}

\begin{prop}\label{C:TateModule-Polarizable} Suppose that $\deg(\cL)> 0$, and that the characteristic of $k$ is either zero or strictly larger than the height of the adjoint representation. Then the Tate module of the group scheme $\schP^o\to A$ is polarizable.
\end{prop}
\begin{proof}By \Cref{coroll: quasiprojectivity of algebraic space}, we know that $\pi:\schP^o\to A$ is a smooth quasi-projective commutative group scheme with connected fibers. Therefore, polarizability follows from \cite[Thm. 6.2]{ancona-fratila-polarizability}.
\end{proof}

\begin{remark}
    A previous version of this paper included an explanation of the polarizability of the Tate module for smooth commutative group schemes with connected fibers over a Noetherian normal base scheme, under the assumption that the generic fibers are Abelian varieties. The proof proceeded using standard theory of biextensions, and could be mostly extracted from the contents of \cite[VII, VIII]{SGA7-1}. Since then, the proof of \cite[Thm. 6.2]{ancona-fratila-polarizability} appeared, so we decided to omit this exposition for the sake of brevity. 
\end{remark}

\end{subsection}

\begin{subsection}{Weak Abelian fibration structure on the Higgs moduli space}

\begin{defn}[Weak Abelian fibration {\cite[\S7.1]{ngo-lemme-fondamental}}, {\cite[1.1]{maulik-shen-independence}}] \label{defn: weak Abelian fibration}
    Let $B$ be a scheme of finite type over $k$. Let $P \to B$ be a smooth commutative group scheme of finite type over $B$. Let $h: X \to B$ be quasiprojective scheme proper over $B$ equipped with an action of $P$. We say that the triple $(X,P,B)$ is a weak Abelian fibration of relative dimension $s$ if the following are satisfied:
    \begin{enumerate}[(a)]
        \item The group scheme $P \to B$ has fibers of dimension $s$ and the scheme $X$ is equidimensional of dimension $\dim(X) = s + \dim(B)$.
        
        \item The action of $P$ on $X$ has affine stabilizers.
        \item The Tate module $T_{\overline{\mathbb{Q}}_{\ell}}(P^o)$ of the open subgroup scheme $P^o \subset P$ of fiberwise neutral components is polarizable in the sense of \Cref{defn:polarizable tate module}.
    \end{enumerate}
\end{defn}

\begin{remark}\label{rmk: p pzero ok}
Let $(X,P,B)$ be a weak Abelian fibration
and let $P^o/B$ be the open subgroup scheme
of fiberwise neutral components of $P/B$.
Then $(X,P^o,B)$ is a weak Abelian fibration.
\end{remark}

By using our results so far, we can now state one of our main theorems.
\begin{thm} \label{thm: weak Abelian fibration structure}
    Suppose that the characteristic of $k$ is either zero, or strictly larger than the height of the adjoint representation. Assume that $\deg(\cL) \geq \max\{2g-2,1\}$. Fix $d \in \pi_1(G)$. The action of the smooth commutative group scheme $\schP^o \to A$ on the moduli space $h:\MHiggs_{G, \cL}^d \to A$ (as in \Cref{coroll: action on the semistable moduli space}) makes the triple $(\MHiggs_{G, \cL}^d, \schP^o, A)$ into a weak Abelian fibration of relative dimension $\frac{1}{2}\dim(G) \deg(\cL) - \frac{1}{2} \dim(T) (\deg(\cL)-2g+2)+\dim(Z_G)$.
\end{thm}
\begin{proof}
The Hitchin fibration $h: \MHiggs_{G, \cL}^d \to A$ is proper by \Cref{prop: fibers hitchin pure dimensional} and $\MHiggs_{G,\cL}^d$ is quasiprojective by \Cref{thm: adequate moduli spaces for higgs}. Its dimension is the required one by \Cref{cor: higgs moduli space properties}. The commutative group scheme $\schP^o \to A$ is smooth (\Cref{P:centr-group}\eqref{P:centr-group3}), quasiprojective (\Cref{coroll: quasiprojectivity of algebraic space}), of the required relative dimension (\Cref{P:centr-group}\eqref{P:centr-group2}). By \Cref{prop: afineness of stabilizers general}, the stabilizers of the action of $\schP^o$ are affine. By \Cref{C:TateModule-Polarizable}, the Tate module of $\schP^o \to A$ is polarizable. 
    The conclusion follows.
\end{proof}
\end{subsection}

\end{section}

\begin{section}{Generalization to isotrivial reductive group schemes} \label{section: isotrivial group schemes}

In this section we generalize our results by replacing $G$ with an isotrivial reductive group scheme over $C$. This is achieved by means of a covering trick. Our main references for reductive group schemes are \cite{sga3} and \cite{conrad_reductive}. Recall that a reductive group scheme $\cG \to C$ is a smooth affine group scheme with connected fibers over $C$ such that there exist an  \'etale and surjective morphism $U \to C$, a connected split reductive group $G$ over $k$, and an isomorphism of group $U$-schemes $G\times U \xrightarrow{\sim} \cG\times_{C} U$. The group $G$ is uniquely determined up isomorphism, and we call it the fiber of $\cG$ over $C.$

\begin{defn} \label{defn: isotrivial group scheme} A reductive group scheme $\cG \to C$ with fiber $G$ is called isotrivial if there is a finite \'etale morphism $\widetilde{C} \to C$ and an isomorphism of group $\widetilde{C}$-schemes $G \times \widetilde{C} \xrightarrow{\sim} \cG \times_{C} \widetilde{C}$.
\end{defn}

\begin{context} \label{context: isotrivial section}
    For the rest of this section, we fix 
    an isotrivial reductive group scheme $\cG \to C$ with fiber $G$, and we fix a choice of a trivializing morphism  $n: \widetilde{C} \to C$ 
    for $\cG$ as in \Cref{defn: isotrivial group scheme}. Furthermore, after modifying the cover, we can and do assume without loss of generality that  $n:\widetilde{C} \to C$ is a connected Galois cover, whose Galois group we denote by $\Gamma$.
    All of the results in this section can be checked after extending the base field. After passing to a finite extension of $k$ and modifying the Galois cover, we might assume that $\widetilde{C}$ is furthermore geometrically connected.
\end{context}

\begin{context}
    As usual, we implicitly assume that the characteristic of $k$ does not divide the order $|W|$ of the Weyl group of $G$.
\end{context}

\begin{notn}
    We  denote by $\fg$ the Lie algebra of $\cG$. This is a vector bundle over the curve $C$, and it is equipped with 
    the adjoint action of $\cG$.
\end{notn}

\begin{remark}
Ultimately, for the study of the Decomposition Theorem, we are only  interested in the special case when the isotrivial group scheme $\cG$ is an outer form of $G$. This means that the group scheme $\cG$ is obtained by letting $n: \widetilde{C} \to C$ be a torsor for a subgroup of the group of outer automorphisms $\text{Out}(G)$, and then twisting $G \times C$ by this torsor, via the action of $\text{Out}(G)$ on $G$ induced by the choice of a pinning. These are the twisted group schemes considered in \cite{ngo-lemme-fondamental}, which we use in \cite{DT_paper_2}.

Since the additional property of $\cG$ being an outer form is not needed in any of our arguments in this section, we have chosen to work in the more general setting of an arbitrary isotrivial reductive group scheme $\cG$. 
\end{remark}

\subsection{The stack of \texorpdfstring{$\cL$}{L}-valued \texorpdfstring{$\cG$}{G}-Higgs bundles}

For any $\cG$-bundle $E$ on the curve $C$, we denote by $\ad(E) := E(\fg)$ the adjoint vector bundle associated with the action of $\cG$ on $\fg$.

\begin{defn}
    A ($\cL$-valued) $\cG$-Higgs bundle is a pair $(E, \theta)$ where $E$ is a $\cG$-bundle on $C$ and $\theta \in H^0(C, \ad(E) \otimes \cL)$ is a global section of the $\cL$-valued adjoint bundle $\ad(E) \otimes \cL$.
\end{defn}

We denote by $\Higgs_{\cG,\cL}: = \text{Sect}_{C}([(\fg\otimes \cL)/\cG])$ the stack parametrizing $\cG$-Higgs bundles. It is an algebraic stack locally of finite type and with affine diagonal over $k$ \cite[Thm. 1.3]{hall-rydh-tannakahom}.

\begin{notn}
    We denote by $\Higgs_{G,n^*(\cL)}(\widetilde{C})$ the stack of $n^*(\cL)$-valued $G$-Higgs bundles over the curve $\widetilde{C}$.
\end{notn}

There is a pullback morphism
\[ n_{\Higgs}: \Higgs_{\cG, \cL} \to \Higgs_{G,n^*(\cL)}(\widetilde{C}), \; \; \; (E, \theta) \mapsto (n^*(E), n^*(\theta)).\]

Our next goal is to understand this pullback morphism from the point of view of \'etale descent (see \Cref{prop: twisted higgs as fixed points}). Recall that we have a fixed isomorphism of $\widetilde{C}$-group schemes $\psi: n^*(\cG) \xrightarrow{\sim} G_{\widetilde{C}}$. For any automorphism  $\gamma \in \Gamma$, we get an induced automorphism of group schemes over $\widetilde{C}$
\[ \psi_{\gamma}: G_{\widetilde{C}} \xrightarrow{\psi^{-1}} n^*(\cG) \xrightarrow{\sim} (n \circ \gamma)^*(\cG) \xrightarrow{\gamma^*(\psi)} \gamma^*(G_{\widetilde{C}}) = G_{\widetilde{C}}\]
The tuple $(\psi_{\gamma})_{\gamma \in \Gamma}$ yields the descent data needed to descend $G_{\widetilde{C}}$ to $\cG$ under the Galois cover $n: \widetilde{C} \to C$.

    \begin{notn}
        If we are given a $G$-bundle $E$ on $\widetilde{C}$, we can view it as a $G_{\widetilde{C}}$-torsor on $\widetilde{C}$. We set $\gamma \cdot E$ to be the $G_{\widetilde{C}}$-torsor on $\widetilde{C}$ given by 
    \[ \gamma \cdot E := \gamma^*(E) \times^{G_{\widetilde{C}},\,  \psi_{\gamma}} G_{\widetilde{C}}.\]
    We view $\gamma \cdot E$ as a $G$-bundle on 
    $\widetilde{C}$.
    \end{notn}

    Let $\ad: G_{\widetilde{C}} \to \GL(\text{Lie}(G))_{\widetilde{C}}$ denote the adjoint representation of the $\widetilde{C}$-group scheme $G_{\widetilde{C}}$. Since this representation descends to the adjoint representation of $\cG$, it follows by descent that for every $\gamma$ there is a canonical isomorphism $\ad  \xrightarrow{\sim} \ad \circ \psi_{\gamma}$. In other words, the automorphism $\psi_{\gamma}$ is compatible with the adjoint representation. Hence, for any $G_{\widetilde{C}}$-bundle $F$, we get an isomorphism of adjoint bundles $(\psi_{\gamma})_*: \ad(F) \xrightarrow{\sim} \ad(F \times^{\psi_{\gamma}} G_{\widetilde{C}})$.

    \begin{notn}
        Given $G$-bundle $E$ on $\widetilde{C}$ and a Higgs field $\theta \in H^0(\widetilde{C}, \ad(E) \otimes n^*(\cL))$, we set $\gamma \cdot \theta \in H^0(\widetilde{C}, \ad(\gamma \cdot E)\otimes n^*(\cL))$ to be the Higgs field on $\gamma \cdot E$ given by $\gamma \cdot \theta = (\psi_{\gamma})_*(\theta)$. (Here we are using the canonical identification $\gamma^* n^*(\cL) = n^*(\cL)$).
    \end{notn}

\begin{defn} \label{defn: action of gamma on higgs}
    We define an action of the group $\Gamma$ on the stack $\Higgs_{G, \cL}(\widetilde{C})$ given by
    \[ \Gamma \times \Higgs_{G, \cL}(\widetilde{C}), \; \; \; (\gamma, (E,\theta)) \mapsto (\gamma \cdot E, \gamma \cdot \theta)\]
\end{defn}

\begin{prop} \label{prop: twisted higgs as fixed points}
    The pullback morphism $n_{\Higgs}: \Higgs_{\cG, \cL} \to \Higgs_{G,n^*(\cL)}$ exhibits $\Higgs_{\cG, \cL}$ as the stack of fixed points $\Higgs_{G,n^*(\cL)}(\widetilde{C})^{\Gamma} \to \Higgs_{G,n^*(\cL)}(\widetilde{C})$ for the action of $\Gamma$ in \Cref{defn: action of gamma on higgs} (in the sense of \cite[Defn. 2.3]{romagny-actions-stacks}). In particular the pullback morphism $n_{\Higgs}: \Higgs_{\cG, \cL} \to \Higgs_{G,n^*(\cL)}$ is affine and of finite type.
\end{prop}
\begin{proof}
The identification of $\Higgs_{\cG, \cL}$ with the stack of fixed points is a direct consequence of \'etale descent for the finite Galois cover $n: \widetilde{C} \to C$. Indeed, the equivariant data encoded in a point of $\Higgs_{G,n^*(\cL)}(\widetilde{C})^{\Gamma}$ (as explained in the proof of \cite[Thm. 3.3]{romagny-actions-stacks}) is in our case exactly the data of a cocycle needed to descend to a $\cG$-Higgs bundle. The affineness of the morphism follows from \cite[Prop. 3.7]{romagny-actions-stacks}, because the stack $\Higgs_{G,n^*(\cL)}(\widetilde{C})$ has affine diagonal.
\end{proof}

\begin{defn}
    For any $d \in \pi_1(G)$, we denote by $\Higgs_{\cG,\cL}^d$ the preimage of the connected component $\Higgs_{G,n^*(\cL)}^d(\widetilde{C})$ under the pullback morphism $n_\Higgs$.
\end{defn}

By definition $\Higgs_{\cG,\cL}^d \subset \Higgs_{\cG,\cL}$ is an open substack, and the pullback morphism induces an affine morphism of finite type $n_{\Higgs}: \Higgs_{\cG,\cL}^d \to \Higgs_{G,n^*(\cL)}^d(\widetilde{C})$.

\subsection{Hitchin fibration}

\begin{defn} 
    We denote by $\fc_{\cL}$ the $C$-affine GIT quotient $\fc_{\cL} = (\fg \otimes \cL)\sslash \cG := \Spec_{C}\left(\left(\Sym^{\bullet}(\fg^{\vee} \otimes \cL^{\vee})\right)^{\cG}\right)$. There is an adequate moduli space morphism
    \[ \chi_{\cL}: [(\fg \otimes \cL)/\cG] \to \fc_{\cL}.\]
\end{defn}

\begin{defn}
    The Hitchin base $A$ is the affine scheme of finite type over $k$ representing the functor of sections $\text{Sect}_{C}(\fc_{\cL})$.
\end{defn}

\begin{defn}
    The Hitchin morphism $h: \Higgs_{\cG,\cL} \to A$ is the morphism of stacks of sections $\text{Sect}_{C}([(\fg\otimes \cL)/\cG]) \to \text{Sect}_{C}(\fc_{\cL})$ induced by the moduli space morphism $\chi_{\cL}$.
\end{defn}

We denote by $\widetilde{A}$ the $n^*(\cL)$-valued Hitchin base for the reductive group $G$ over the curve $\widetilde{C}$. There is a pullback morphism of schemes of sections $n_{A}: A \to \widetilde{A}$ inducing a commutative diagram of Hitchin morphisms
\[
\begin{tikzcd}
  \Higgs_{\cG,\cL} \ar[r, "n_{\Higgs}"] \ar[d, "h"] & \Higgs_{G,n^*(\cL)}(\widetilde{C}) \ar[d, "h"] \\ A \ar[r, "n_A"] & \widetilde{A}.
\end{tikzcd}
\]

For the next lemma, we denote by $\rank(\cG)$ the dimension of any split maximal torus $T$ of the fiber $G$, and we set $\dim(\cG/C)$ to be the relative dimension of $\cG \to C$.

\begin{lemma}\label{L:hitchin base for isotrivial group schemes}
    Suppose that $|\Gamma|$ does not divide the characteristic of $k$. Then, the affine scheme $A$ is connected and smooth over $k$. If $\deg(\cL)>2g-2$, then $\dim(A) = \frac{1}{2}\dim(\cG/C)\deg(\cL) +\frac{1}{2}\rank(\cG)(\deg(\cL)-2g+2)$.
\end{lemma}
\begin{proof}
Recall that the Hitchin base $\widetilde{A}$ classifies sections $a: \widetilde{C} \to \fc_{n^*(\cL)}$, where $\fc_{n^*(\cL)} \to \widetilde{C}$ is the affine GIT quotient $\text{Lie}(G)\sslash G$ twisted by $\cL$. By the compatibility of the formation of the ring of invariants $\left(\Sym^{\bullet}(\fg^{\vee} \otimes \cL^{\vee})\right)^{\cG}$ with flat base-change \cite[Prop. 5.2.9(1)]{alper_adequate}, the identification $n^*(\cG) \cong G_{\widetilde{C}}$ induces an identification $\fc_{\cL}\times_{C} \widetilde{C} \cong \fc_{n^*(\cL)}$. The corresponding descent datum  with respect to the finite Galois cover $n: \widetilde{C} \to C$ amounts to an action of $\Gamma$ on $\fc_{n^*(\cL)}$ such that the morphism $\fc_{n^*(\cL)} \to \widetilde{C}$ is $\Gamma$-equivariant. The actions of $\Gamma$ on $\widetilde{C}$ and $\fc_{n^*(\cL)}$ jointly induce an action of $\Gamma$ on the scheme of sections $\widetilde{A} := \text{Sect}_{\widetilde{C}}(\fc_{n^*(\cL)})$ given by $f \mapsto \gamma \circ f \circ \gamma^{-1}$. \'Etale descent for the cover $n: \widetilde{C} \to C$ implies that the scheme of sections $A$ is the schematic locus of fixed points of the action of $\Gamma$ on $\widetilde{A}$. Note that, by assumption on the characteristic of $k$, the group $\Gamma$ is linearly reductive. Since $\widetilde{A}$ is smooth (in fact it is an affine space), the scheme of fixed points $A = \widetilde{A}^{\Gamma}$ is smooth as well \cite[Prop. 4.3]{edixhoven-neron-models}.

The connectedness of $A$ follows from the fact that it admits a $\mathbb{G}_m$ action such that the limit of every point is the origin $0 \in A(k) \subset \widetilde{A}(k)$.

Finally, to compute the dimension, we may extend the ground field to assume $k$ is algebraically closed, and it suffices to compute the dimension of the tangent space at any $k$-point $a: \Spec(k) \to A$, which corresponds to a section $s_a: C \to \fc_{\cL}$. We use the description of a tangent obstruction theory in \cite[Prop. 5.1.10]{halpernleistner2019mapping}. Indeed, there is a natural derived enhancement $A \hookrightarrow A^{der}$ of the section scheme $A$, and the restriction $a^*(\mathbb{T}_{A^{der}/k})$ of the tangent complex of $A^{der}$ is given by $R\pi_*(s_a^*T_{\fc_{\cL}/C})$, where $\pi: C \to \Spec(k)$ is the structure morphism and $T_{\fc_{\cL}/C}$ is the relative tangent bundle of the smooth morphism $\fc_{\cL} \to C$. We  show that $H^1(C,(s_a^*T_{\fc_{\cL}/C}))=0$, and then that $\dim(T_{A,a}) = \dim(H^0(a^*(\mathbb{T}_{A^{der}/k})) = \chi(s_a^*T_{\fc_{\cL}/C})$ has the desired dimension by Riemann-Roch.

Since the finite \'etale morphism $n: \widetilde{C} \to C$ has degree coprime to the characteristic of $k$ by assumption, it follows that we have a splitting $n_*n^*(s_a^*T_{\fc_{\cL}/C})) = s_a^*T_{\fc_{\cL}/C} \oplus \cE$ for some vector bundle $\cE$. Hence, in order to show that $H^1(C, (s_a^*T_{\fc_{\cL}/C}))=0$, it suffices to show that $H^0(C, n_*n^*(s_a^*T_{\fc_{\cL}/C}))) = H^1(\widetilde{C}, n^*(s_a^*T_{\fc_{\cL}/C})))=0$. The image of $\widetilde{a} := n_A(a)$ induces a section $s_{\widetilde{a}}: \widetilde{C} \to \fc_{n^*(\cL)}$. Note that the pullback $n^*(s_a^*T_{\fc_{\cL}/C})$ is isomorphic to $s_{\widetilde{a}}^*(T_{\fc_{n^*(\cL)}/\widetilde{C}})$. But we know that $\fc_{n^*(\cL)} \to \widetilde{C}$ is isomorphic to the total space of the vector bundle $\fc_{n^*(\cL)} \cong \bigoplus_{i} n^*(\cL)^{\otimes e_i}$ for some positive integers $e_i>0$ (see Subsection \ref{subsection: hitchin fibration}). Therefore $s_{\widetilde{a}}^*(T_{\fc_{n^*(\cL)}/\widetilde{C}}) \cong \bigoplus_{i} n^*(\cL)^{\otimes e_i}$. By our assumption that $\deg(\cL)>2g-2$, it follows that $\deg(n^*(\cL)) > 2 g(\widetilde{C}) -2$, and hence by Serre duality we have $H^1(\widetilde{C}, n^*(s_a^*T_{\fc_{\cL}/C})))=0 = \bigoplus_i H^1(\widetilde{C}, n^*\cL^{\otimes e_i}) =0$, as desired.

In order to conclude we just need to compute $\dim(H^0(a^*(\mathbb{T}_{A^{der}/k}))$, which is equal to the Euler characteristic $ \chi(s_a^*T_{\fc_{\cL}/C})$ by the vanishing of $H^1$. We may use the multiplicativity of Euler characterisitic under the finite \'etale morphism $n: \widetilde{C} \to C$ to get
\[ \chi(s_a^*T_{\fc_{\cL}/C}) = |\Gamma|^{-1} \chi(n^*s_a^*T_{\fc_{\cL}/C}) = |\Gamma|^{-1} \chi(\bigoplus_i n^*\cL^{\otimes e_i}).\]
The last dimension for the group $G$ and the curve $\widetilde{C}$ has been computed in the equation \eqref{E:dimA} to be $\frac{1}{2} \dim(G)\deg(n^*\cL) + \frac{1}{2} \dim(T) \left(\deg(n^*\cL) -2(g(\widetilde{C}-1)  \right)$. Dividing by $|\Gamma|$ we get our desired dimension $\dim(T_{A,a}) = \chi(s_a^*T_{\fc_{\cL}/C}) = \frac{1}{2}\dim(\cG/C)\deg(\cL)-\frac{1}{2}\rank(\cG)(\deg(\cL)-2g+2)$.
\end{proof}

\begin{prop} \label{prop: valuative criterion stacky hichin twisted morphism}
    Suppose that the characteristic of $k$ is not a torsion prime for $G$ \cite[\S2]{steinberg-torsion-primes} and is very good for $G$ \cite[\S2.9]{herpel-smoothness}. Then the Hitchin morphism $h: \Higgs_{\cG,\cL} \to A$ satisfies the existence part of the valuative criterion for properness for stacks.
\end{prop}
\begin{proof}
   Let $\tau$ denote the generic point of the curve $C$. After extending the ground field $k$, we may assume without loss of generality, that $k$ is algebraically closed. By Steinberg's theorem \cite[Chapitre III, \S 3.4, pp. 139-140, Theorem $1'$ and Remark 1)]{serre-galois}, the reductive group $\cG_{\tau}$ is quasisplit. We note that then the argument in \cite[Lem. 6.20]{alper2019existence} applies in our case. Let us point out the steps where one needs to be careful: 
   \begin{enumerate}[(1)]

   \item Loc. cit. first proves the assertion in the case when the derived subgroup $\cD_{G}$ is simply-connected. We need to consider the base-change $\cG_{R(\eta)}$, where $R(\eta)$ is the local ring at the generic point of the special fiber of $C_R$. Note that $\cG_{R(\eta)}$ is the base-change of the quasisplit reductive $k(\tau)$-group $\cG_{\tau}$. With this in mind, everything in the proof of the simply-connected case in \cite[Lem. 6.20]{alper2019existence} goes through. The only thing that might not be standard is the existence of Kostant sections for $\cG_{R(\eta)}$, which is needed to define $Y$ in the proof. Since the base-change $\cG_{R(\eta)}$ is quasisplit, we may appeal to \cite[Prop. 4.3.2]{bouthier-cesnavicius} for the existence of the Kostant section.
   
    \item The second step in loc. cit. is to reduce to the case when the derived subgroup of $\cG$ is simply-connected. In our case, we can set $\widetilde{\cG} = Z_{\cG}^o \times \widetilde{\cD}_{\cG}$, where $Z_{\cG}^o \subset \cG$ is the maximal central torus and $\widetilde{\cD}_{\cG} \to \cD_{\cG}$ is the simply-connected central cover of the derived subgroup $\cD_{\cG} \subset \cG$, as in \cite[Exerc. 5.6.2(i)]{conrad_reductive}. Then it is still the case that the morphism $\Higgs_{\widetilde{\cG}, \cL} \to \Higgs_{\cG, \cL}$ has local sections as needed in the proof of \cite[Lem. 6.20]{alper2019existence}.
   \end{enumerate}
\end{proof}

\subsection{Semistability and moduli spaces}

Recall that by \cite[\S4.1]{herrero2023meromorphic}, the locus of semistable $G$-Higgs bundles on $\Higgs_{G,n^*(\cL)}(\widetilde{C})$ is the $\Theta$-semistable locus with respect to certain line bundle $\cL_{det}$ on the stack $\Higgs_{G,n^*(\cL)}(\widetilde{C})$.

\begin{notn}
    We also denote by $\cL_{det}$ the pullback of this line bundle under the morphism $n_{\Higgs}: \Higgs_{G,\cL} \to \Higgs_{G,n^*(\cL)}(\widetilde{C})$.
\end{notn}

\begin{lemma} \label{lemma: compatibility of semistability under pullback morphism}
    Suppose that the characteristic of $k$ is zero, or strictly larger than the height of the adjoint representation of $G$. Then a geometric point $x$ of $\Higgs_{\cG,\cL}$ is $\Theta$-semistable with respect to $\cL_{det}$ (\Cref{defn: theta semistability}) if and only if the pullback $n_{\Higgs}(x)$ is a semistable $G$-Higgs bundle on $\widetilde{C}$.
\end{lemma}
\begin{proof}
    This follows by a similar argument as in the proof of \cite[Prop. 4.3]{gauged_theta_stratifications}, by using our interpretation of $\Higgs_{\cG,\cL}$ as a stack of $\Gamma$-fixed points (\Cref{prop: twisted higgs as fixed points}) and the uniqueness of Harder-Narasimhan parabolic reductions established in \cite[\S 4]{herrero2023meromorphic} under the given characteristic assumptions on $k$.
\end{proof}

\begin{prop} \label{prop: moduli spaces for twisted higgs}
    Suppose that the characteristic of $k$ is either zero, or strictly larger than the height of the adjoint representation of $G$. Then we have that:
    \begin{enumerate}[(1)]
        \item The set of semistable geometric points of $\Higgs_{\cG,\cL}$ defines an open substack $\Higgs_{\cG, \cL}^{ss}$.
        \item For any $d \in \pi_1(G)$, the  open substack $\Higgs_{\cG,\cL}^{d,ss} := \Higgs_{\cG,\cL}^{ss} \cap \Higgs_{\cG,\cL}^d$ admits a quasiprojective adequate moduli space $\MHiggs_{\cG,\cL}^d$.

        \item The  morphism $n_{\Higgs}$ restricts to an affine morphism $n_{\Higgs}: \Higgs_{\cG,\cL}^{d,ss} \to \Higgs_{G,n^*(\cL)}^{d,ss}(\widetilde{C})$ which induces an affine morphism of adequate moduli spaces $n_{\MHiggs}: \MHiggs_{\cG,\cL}^d \to \MHiggs_{G,n^*(\cL)}^d$.
    \end{enumerate}
\end{prop}
\begin{proof}
    By \Cref{lemma: compatibility of semistability under pullback morphism}, the locus of semistable geometric points in $\Higgs_{\cG,\cL}$ is the preimage under the pullback morphism $n_{\Higgs}$ of the locus of semistable geometric points of $\Higgs_{G,n^*(\cL)}(\widetilde{C})$. By \Cref{thm: adequate moduli spaces for higgs}, the locus of semistable points in $\Higgs_{G,n^*(\cL)}(\widetilde{C})$ is open, and hence the same holds for $\Higgs_{\cG,\cL}$. This shows part (1).

    Since the pullback morphism $n_{\Higgs}: \Higgs_{\cG, \cL} \to \Higgs_{G,n^*(\cL)}(\widetilde{C})$ is affine, the same holds for the restriction 
    \[n_{\Higgs}: \Higgs_{\cG,\cL}^{d,ss} = (n_{\Higgs})^{-1}\left(\Higgs_{G,n^*(\cL)}^{d,ss}(\widetilde{C})\right) \to \Higgs_{G,n^*(\cL)}^{d,ss}(\widetilde{C}).\]
    Since $\Higgs_{G,n^*(\cL)}^{d,ss}(\widetilde{C})$ admits an adequate moduli space $\MHiggs_{G,n^*(\cL)}^d$ (\Cref{thm: adequate moduli spaces for higgs}), it follows that $\Higgs_{\cG,\cL}^{d,ss}$ admits an adequate moduli space $\MHiggs_{\cG,\cL}^d$ that is affine and of finite type over $\MHiggs_{G,n^*(\cL)}^d$ \cite[Lem. 5.2.11 + Thm. 6.3.3]{alper_adequate}. This shows part (3).

    Finally, since the moduli space $\MHiggs_{G,n^*(\cL)}^d$ is quasiprojective (\Cref{thm: adequate moduli spaces for higgs}), the same holds for the $\MHiggs_{G,n^*(\cL)}^d$-affine scheme $\MHiggs_{\cG,\cL}^d$. This shows part (2) and concludes the proof of the proposition.
\end{proof}

\begin{prop} \label{prop: smoothness of stack of twisted higgs}
    Suppose that $\deg(\cL)>2g-2$, that $|\Gamma| \nmid \text{char}(k)$, and that $\text{char}(k)$ is either zero, or strictly larger than the height of the adjoint representation. Then the open substack of semistable points $\Higgs_{\cG,\cL}^{ss}$ is smooth over $k$.
\end{prop}
\begin{proof}
    It follows from \Cref{prop: twisted higgs as fixed points} and the proof of \Cref{prop: moduli spaces for twisted higgs} that $\Higgs_{\cG,\cL}^{ss}$ is a stack of $\Gamma$-fixed points of an action of $\Gamma$ on $\Higgs_{G,n^*(\cL)}^{ss}(\widetilde{C})$. Since the stack $\Higgs_{G,n^*(\cL)}^{ss}(\widetilde{C})$ is smooth (\Cref{prop: smoothness of the semistable hitchin stack}) and the group $\Gamma$ is linearly reductive (because $|\Gamma| \nmid \text{char}(k)$), it follows that the stack of fixed points $\Higgs_{\cG,\cL}^{ss}$ is also smooth \cite[Thm. 4.3.6]{romagny2022algebraicity}.
\end{proof}

\begin{coroll} \label{coroll: properties twisted isotrivial Higgs moduli space}
    Suppose that $\deg(\cL)>2g-2$, that $|\Gamma| \nmid \text{char}(k)$, and that $\text{char}(k)$ is either zero, or strictly larger than the height of the adjoint representation. Then the moduli space $\MHiggs_{\cG,\cL}^d$ is a geometrically normal quasiprojective scheme.
\end{coroll}
\begin{proof}
    This is a consequence of the smoothness of the stack $\Higgs_{\cG,\cL}^{ss}$ (\Cref{prop: smoothness of stack of twisted higgs}) and \cite[Prop. 5.4.1]{alper_adequate}.
\end{proof}

\subsection{\texorpdfstring{$\Theta$}{Theta}-stratification and semistable reduction}

\begin{prop} \label{prop: theta-stratification twisted higgs}
    Suppose that the characteristic of $k$ is either zero, or strictly larger than the height of the adjoint representation of $G$. Then the semistable locus $\Higgs_{\cG,\cL}^{d,ss}$ is the open stratum $\left(\Higgs_{\cG,\cL}^{d}\right)_{\leq 0}$ of a well-ordered $\Theta$-stratification $\Higgs_{\cG,\cL}^{d} = 
 \bigcup_{\sigma \in \Sigma} \left(\Higgs_{\cG,\cL}^{d}\right)_{\leq \sigma}$ (as in \cite[Defn. 6.1]{alper2019existence}).
\end{prop}
\begin{proof}
    We shall construct a well-ordered $\Theta$-stratification of $\Higgs_{\cG,\cL}$ with lowest open stratum $\Higgs_{\cG, \cL}^{ss}$. Then our desired well-ordered $\Theta$-stratification is obtained by restricting the stratification to the open and closed substack $\Higgs_{\cG,\cL}^{d} \subset \Higgs_{\cG,\cL}$. Under our hypotheses on the characteristic of $k$, such $\Theta$-stratification was constructed in \cite[Thm. 4.22]{herrero2023meromorphic} for the stack $\Higgs_{G,n^*(\cL)}(\widetilde{C})$ (in loc. cit., it is stated in the context of $t$-connections, but setting $t=0$ one recovers the case of Higgs bundles, and the arguments hold verbatim for an arbitrary line bundle such as $n^*(\cL)$). Let us explain how this stratification induces a $\Theta$-stratification on the stack of $\Gamma$-fixed points $\Higgs_{\cG,\cL}$.

   By definition of $\Theta$-stratification, we have an exhaustive union of open substacks $\Higgs_{G,n^*(\cL)}(\widetilde{C}) = \bigcup_{\sigma \in \Sigma} \left(\Higgs_{G,n^*(\cL)}(\widetilde{C})\right)_{\leq \sigma}$ indexed by a well-ordered set $\Sigma$ with minimum $0 \in \Sigma$. For every $\sigma \in \Sigma$ there is an open substack $\cS_{\sigma} \subset \Map(\Theta,\Higgs_{G,n^*(\cL)}(\widetilde{C}))$ 
   %of the mapping stack 
   such that the natural evaluation morphism $ev_1: \Map(\Theta,\Higgs_{G,n^*(\cL)}(\widetilde{C})) \to \Higgs_{G,n^*(\cL)}(\widetilde{C})$ restricts to a locally 
   closed immersion $ev_1: \cS_{\sigma} \hookrightarrow \Higgs_{G,n^*(\cL)}(\widetilde{C})$ (cf. \cite[Defn. 6.1]{alper2019existence}). The image of $|\cS_\sigma|$ in the 
   topological space $\left|\Higgs_{G,n^*(\cL)}(\widetilde{C}))\right|$ is exactly the locally closed 
   stratum $\left|\left(\Higgs_{G,n^*(\cL)}(\widetilde{C})\right)_{\leq \sigma}\right| \; \setminus \; \left|\left(\Higgs_{G,n^*(\cL)}(\widetilde{C})\right)_{< \sigma}\right|$.

 There is an action of $\Gamma$ on the stack $\Higgs_{G,n^*(\cL)}(\widetilde{C})$, which induces an action on the mapping stack $\Map(\Theta,\Higgs_{G,n^*(\cL)}(\widetilde{C}))$ given by post-composition $\gamma \cdot f := \gamma \circ f$. By the universal property of the stack of $\Gamma$-fixed points $\Higgs_{\cG,\cL}$(\cite[Defn. 2.3]{romagny-actions-stacks}, a point in the mapping stack $\Map(\Theta,\Higgs_{\cG,\cL})$ is exactly the data corresponding to a fixed point in $\Map(\Theta,\Higgs_{G,n^*(\cL)}(\widetilde{C}))^{\Gamma}$. In other words, we have a canonical identification of stacks  $\Map(\Theta,\Higgs_{\cG,\cL}) = \Map(\Theta,\Higgs_{G,n^*(\cL)}(\widetilde{C}))^{\Gamma}$. If each stratum $\cS_{\sigma} \subset \Map(\Theta,\Higgs_{G,n^*(\cL)}(\widetilde{C}))$ is preserved by the action of $\Gamma$, then the well-ordered set of open substacks of fixed points $\cS_{\sigma}^{\Gamma} \subset \Map(\Theta,\Higgs_{\cG,\cL})$  yields our desired $\Theta$-stratification of $\Higgs_{\cG,\cL}$. Hence we are reduced to showing that each stratum $\cS_{\sigma} \subset \Map(\Theta,\Higgs_{G,n^*(\cL)}(\widetilde{C}))$ is stable under the action of $\Gamma$.

 In \cite[Thm. 4.22]{herrero2023meromorphic} the $\Theta$-stratification is induced by a numerical invariant $\mu$ on the stack $\Higgs_{G,n^*(\cL)}(\widetilde{C})$ in the sense of \cite{halpernleistner2018structure}. This induces a locally constant function $\mu: \left| \Map(\Theta,\Higgs_{G,n^*(\cL)}(\widetilde{C})) \right| \to \mathbb{R}$, and the $\Theta$-stratification is defined by solving an optimization problem for $\mu$ on $\Map(\Theta,\Higgs_{G,n^*(\cL)}(\widetilde{C}))$ to define  ``HN-filtrations" (we refer the reader to \cite[\S2.4]{gauged_theta_stratifications} and the more complete \cite{halpernleistner2018structure} for more details). In particular, in order to show that each $\cS_{\sigma}$ is $\Gamma$-stable, it suffices to show that the function $\mu: \left| \Map(\Theta,\Higgs_{G,n^*(\cL)}(\widetilde{C})) \right| \to \mathbb{R}$ is preserved by the action of $\Gamma$. This can be checked directly from our definition of the $\Gamma$-action and the definition of the numerical invariant in \cite[Defn. 4.4]{herrero2023meromorphic}.
\end{proof}

\begin{defn}
    For any given $d \in \pi_1(G)$, the Hitchin fibration $h: \MHiggs_{\cG, \cL}^d \to A$ is the morphism induced by $h: \Higgs_{\cG,\cL}^{d,ss} \to A$.
\end{defn} 
\begin{prop} \label{prop: properness of the hitchin fibration}
    Suppose that the characteristic of $k$ is either zero, or larger than the height of the adjoint representation of $G$. Then the Hitchin fibration $h: \MHiggs_{\cG, \cL}^d \to A$ is proper.
\end{prop}
\begin{proof}
    It suffices to check that $h: \MHiggs_{\cG, \cL}^d \to A$ is universally closed, for which we might apply the valuative criterion for properness. Since the adequate moduli space morphism $\Higgs_{G,\cL}^{d,ss} \to \MHiggs_{\cG,\cL}^{d}$ is surjective, it suffices to check that $\Higgs_{G,\cL}^{d,ss} \to A$ satisfies the existence part of the valuative criterion for properness. We note that $\Higgs_{G,\cL}^d \to A$ satisfies the existence part of the valuative criterion for properness by \Cref{prop: valuative criterion stacky hichin twisted morphism}, because our assumption on the characteristic of $k$ implies in particular that it is not a torsion prime of $G$ and that it is very good for $G$ by the discussion in \cite[4.5]{deligne-balaji-parameswaran-complete-red}. To conclude the existence part of the valuative criterion for $\Higgs_{G,\cL}^{d,ss} \to A$ we use \Cref{prop: theta-stratification twisted higgs}, which allows us to apply the semistable reduction theorem \cite[Cor. 6.12]{alper2019existence}.
\end{proof}

\subsection{Group scheme of symmetries of the Hitchin fibration} \label{subsection: group scheme of symmetries twisted higgs}

Similarly, as in Subsection \ref{SS:local-centr-gr}, the scheme
$$
I_\cL:=\{(x,g)\in \fg_\cL\times_C\cG|\,\ad_g(x)=x\}\subset \fg_\cL\times_C\cG
$$
is a group scheme over the Lie algebra bundle $\fg_\cL$. We denote by $\fg^{\mathrm{reg}}_\cL$ the open subset of those elements $x$ in $\fg_\cL$ such that $\dim((I_\cL)_x)$ is minimal. We denote by $\chi_\cL:\fg_\cL\to \fg_\cL \sslash \cG=\fc_\cL$ the quotient morphism of $C$-schemes.

\begin{lemma}\label{L:J-quasisplit} There exists a unique smooth commutative affine group scheme of finite type over $\fc_\cL$ equipped with a $\cG$-equivariant isomorphism $(\chi^*J_\cL)_{\fg^{\mathrm{reg}}_\cL}\to I_{\fg^{\mathrm{reg}}_\cL}$.
Furthermore, the isomorphism extends to a morphism $\chi_\cL^*J_\cL\to I_\cL$ of group schemes over $\fg_\cL$.
\end{lemma}

\begin{proof} See \cite[Thm. 4.2.8 and Rem. 4.2.10]{bouthier-cesnavicius}.
\end{proof} 

\begin{defn}We denote by $J_A:= ev^*(J_{\cL})$ the group scheme over $C\times A$ obtained by pulling-back the group scheme $J_\cL$ along the evaluation map $ev:C\times A\to\fc_\cL$ given by $(x,a)\mapsto a(x)$.
\end{defn}

\begin{defn}
    The stack of symmetries of the Hitchin fibration $\stP \to A$ is the relative moduli stack of $J_{A}$-torsors for the morphism $C\times A \to A$.
\end{defn}

\begin{prop}
The following hold:
\begin{enumerate}[(1)]
    \item The stack $\stP \to A$ is a smooth commutative group stack over $A$ with affine diagonal.
    \item There is a unique open group sub-stack $\stP^o \subset \stP$ with geometrically connected fibers over $A$.
\end{enumerate}
\end{prop}
\begin{proof}
    The argument is the same as in the proof of \Cref{P:centr-group}, and is  omitted.
\end{proof}

Arguing as in \cite[Subsection 4.3]{ngo-lemme-fondamental}, the morphism $\chi_\cL^*J_\cL\to I_\cL$ of \Cref{L:J-quasisplit} induces for any $S$-point $(E,\theta)\in\Higgs_{\cG,\cL}(S)$ over $a:=h(E,\theta)\in A(S)$ a morphism of group schemes $J_a\to Aut_{C\times S}(E,\theta)$ over $S$. Hence, we may twist the $\cG$-Higgs bundle by any $J_a$-torsor. This defines an action morphism $\stP \times_{A} \Higgs_{\cG,\cL} \to \Higgs_{\cG,\cL}$.

For the following proposition, we denote by $\stP_G \to \widetilde{A}$ the commutative group stack of symmetries of the Hitchin fibration $\Higgs_{G,n^*(\cL)}(\widetilde{C}) \to \widetilde{A}$ as in \Cref{defn: stack of symmetries of the hichin fibration}. We denote by $n^*_A(\stP) \to A$ the pullback under the morphism $n_A: A \to \widetilde{A}$.
\begin{lemma} \label{lemma: pullback morphism groups stack twisted}
    There is an affine homomorphism $n_{\stP}: \stP \to n^*_A(\stP_G)$ of commutative group stacks over $A$, which restricts to an affine homomorphim of open substacks $\stP^o \to n_A^*(\stP_G^o)$. Furthermore, it is compatible with respect to the actions on the Higgs moduli stacks as follows 
    \begin{equation}
        \begin{tikzcd}[ampersand replacement=\&]
         \stP \times_{A} \Higgs_{\cG,\cL} \ar[d, "n_{\stP}\times n_{\Higgs}"] \ar[r, "act"]\& \Higgs_{\cG,\cL} \arrow[d,"n_{\Higgs}"]\\
        \stP_{G} \times_{\widetilde{A}} \Higgs_{G,n^*(\cL)}(\widetilde{C}) \ar[r, "act"] \& \Higgs_{G,n^*(\cL)}(\widetilde{C})
        \end{tikzcd}
    \end{equation}	
\end{lemma}
\begin{proof}
    This follows from the same argument as in \Cref{prop: stacky pullback morphism group scheme} by viewing $\stP$ as a stack of $\Gamma$-fixed points of $n_A^*(\stP_G)$ using \'etale descent for the cover $n: \widetilde{C} \to C$.
\end{proof}

Suppose that $\deg(\cL)>0$. Let us denote by $\schP_G^o \to \widetilde{A}$ the quasiprojective group scheme, which is the good moduli space of $\stP^o_G$ granted by \Cref{P:centr-group}.
\begin{prop}
    Suppose that $\deg(\cL)>0$ and the characteristic of $k$ is either zero, or strictly larger than the height of the adjoint representation of $G$. Then, the commutative group stack $\stP^o \to A$ admits a quasiprojective good moduli space $\schP^o \to A$, which is a flat commutative group scheme over $A$. 
\end{prop}
\begin{proof}
    Note that the stack $n_A^*(\stP^o_G) \to A$ admits a quasiprojective good moduli space given by the pullback $n_A^*(\schP^o_G)$ \cite[Prop. 4.7(i)]{alper-good-moduli}. Since there is an affine morphism $n_{\stP}: \stP^o \to n_A^*(\stP^o_G)$ (\Cref{lemma: pullback morphism groups stack twisted}), it follows that $\stP^o$ admits a good moduli space $\schP^o \to A$ which is affine and of finite type over $n_A^*(\schP^o_G)$ \cite[Lem. 4.14]{alper-good-moduli}. In particular, $\schP^o \to A$ is also quasiprojective. The commutative group stack stucture on $\stP^o \to A$ induces a commutative group scheme structure on $\schP^o \to A$ by the universal property of good moduli spaces \cite[Thm. 6.6]{alper-good-moduli} and their compatibility with base-change \cite[Prop. 4.7(i)]{alper-good-moduli}. Flatness of $\schP^o \to A$ holds because the stack $\stP^o \to A$ is flat over $A$ \cite[Thm. 4.16(ix)]{alper-good-moduli}.
\end{proof}

\begin{prop} \label{prop: action on twisted higgs}
    Suppose that $\deg(\cL)>0$, and the characteristic of $k$ is either zero or strictly larger than the height of the adjoint representation of $G$. The action of the commutative group stack $\stP^o \to A$ preserves the open substack $\Higgs_{\cG,\cL}^{d,ss} \subset \Higgs_{\cG, \cL}$. By taking moduli spaces, this induces an action $\schP^o \times \MHiggs_{\cG,\cL}^d \to \MHiggs_{\cG,\cL}^d$ with affine stabilizers.
\end{prop}
\begin{proof}
    The fact that $\stP^o$ preserves $\Higgs_{\cG,\cL}^{d,ss}$ follows from $\Higgs_{\cG,\cL}^{d,ss} = (n_{\Higgs})^{-1}\left( \Higgs_{G,n^*(\cL)}^{d,ss}(\widetilde{C})\right)$ (\Cref{lemma: compatibility of semistability under pullback morphism}), the diagram in \Cref{lemma: pullback morphism groups stack twisted}, and the fact that the action of $\stP_G^o$ preserves the open substack $\Higgs_{G,n^*(\cL)}^{d,ss}(\widetilde{C}) \subset \Higgs_{G,n^*(\cL)}(\widetilde{C})$ (\Cref{coroll: action on the semistable moduli space}).

    We get an induced action on the moduli space, because $\schP^o \to A$ is flat and the formation of adequate moduli spaces commutes with flat base-change \cite[Prop. 5.2.9(1)]{alper_adequate}. The morphism $n_\MHiggs: \MHiggs_{\cG,\cL}^d \to \MHiggs_{G,n^*(\cL)}^d$ is equivariant with respect to the affine morphism of group schemes $n_{\schP}: \schP^o \to n_A^*(\schP^o_G)$. Since $\schP^o_G$ acts with affine stabilizers on $\MHiggs_{G,n^*(\cL)}^d$ (\Cref{prop: afineness of stabilizers general}), the same holds for the action of $\schP^o$ on $\MHiggs_{\cG,\cL}^d$.
\end{proof}

Let $\pi: C \to \Spec(k)$ denote the structure morpism, and let $Z_{\cG} \subset \cG$ be the central group subscheme of $\cG$, which is a multiplicative group over $C$ \cite[Thm. 3.3.4]{conrad_reductive}. We denote by $\pi_*(Z_\cG)$ the Weil restriction of $Z_{\cG}$ under the structure morphism. We note that $\pi_*(Z_{\cG})$ is an affine commutative group scheme of finite type over $k$ (\cite[Thm. 1.3]{hall-rydh-tannakahom} and \cite[Thm. 2.3(ii)]{hall-rydh-hilbert-quot}).

For the next proposition, recall that we denote by $\rank(\cG)$ the dimension of a split maximal torus $T$ of the fiber $G$, and we set $\dim(\cG/C)$ to be the relative dimension of $\cG \to C$.
\begin{prop} \label{prop: gerbe group scheme of symmetries}
    Suppose that $\deg(\cL)>0$, and the characteristic of $k$ is either zero, or strictly larger than the height of the adjoint representation of $G$. Then the good moduli space morphism $\stP^o \to \schP^o$ is a $\pi_*(Z_{\cG})$-gerbe. The quasiprojective commutative group scheme $\schP^o \to A$ is smooth with geometrically connected fibers and relative dimension $\frac{1}{2}\dim(\cG/C)\deg(\cL)-\frac{1}{2}\rank(\cG)(\deg(\cL)-2g+2) + \dim(\pi_*(Z_{\cG}))$ over $A$.
\end{prop}
\begin{proof}
    Note that $Z_{\cG} \times_{C} \fc$ admits a natural 
    closed embbedding into $J$. By twisting by $\cL$ and pulling back under the evaluation morphism $ev: C \times A \to \fc_{\cL}$, we obtain a natural inclusion of $Z_{\cG} \times A$ into $J_A$. This induces a compatible central action of the Weil restriction $\pi_*(Z_{\cG})$ on every point of the stack of $J_A$-torsors $\stP$. We may form the $\pi_*(Z_{\cG})$-rigidification $\stP \fatslash \pi_*(Z_{\cG})$ in the sense of \cite[Defn. 5.1.9]{acv-twisted-bundles-covers}. By the universal property of rigidifications, there is a canonical morphism $\stP \fatslash \pi_*(Z_{\cG}) \to \schP$. In order to conclude that $\stP \fatslash \pi_*(Z_{\cG}) = \schP$, it suffices to show that $\stP \fatslash \pi_*(Z_{\cG})$ is an algebraic space. We are reduced to showing that the automorphism group of every point of $\stP$ is $\pi_*(Z_{\cG})$. For this, we use an alternative description of $Z_{\cG}$.

    Note that the pullback $n^*(Z_{\cG})$ under the Galois cover $n: \widetilde{C} \to C$ is identified with $Z_G \times \widetilde{C}$. The action of $\Gamma$ on $G \times \widetilde{C}$ induces an action on the Weil restriction $(\pi \circ n)_*(Z_G \times \widetilde{C}) = Z_G$, and it follows from \'etale descent for morphisms that the Weil restriction $\pi_*(Z_{\cG})$ is isomorphic to the group scheme of fixed points $(Z_G)^{\Gamma}$. Hence we need to show that the automorphism of every point in $\stP$ is isomorphic to $(Z_G)^{\Gamma}$. But, as remarked in the proof of \Cref{lemma: pullback morphism groups stack twisted}, there is a pullback morphism $n_{\stP}: \stP \to n_A^*(\stP_G)$ that witnesses $\stP$ as a stack of $\Gamma$-fixed points for a natural action of $\Gamma$ on $\stP_G$. By the definition of the stack of fixed points, the group of automorphisms $\Aut(x)$ of any geometric point $x$ of $\stP$ is the group scheme of $\Gamma$-fixed points $\Aut(n_{\stP}(x))^{\Gamma}$ with respect to the induced action of $\Gamma$ on the group of automorphisms of the image $n_{\stP}(x)$ in $n_A^*(\stP_G)$. Under the assumption $\deg(\cL)>0$, we have shown in \Cref{P:centr-group} that the group of automorphisms of any point in $\stP_G$ is $Z_G$. The induced action of $\Gamma$ on $Z_G$ is the one induced by descent as explained above, so it follows that the group of automorphisms of any geometric point $x$ in $\stP$ is $(Z_G)^{\Gamma}$, as desired.

For the computation of the relative dimension, choose a geometric point $a \to A$. Since the fiber $\stP_a$ is a $\pi_*(Z_{\cG})_a$-gerbe over $\schP_a$, we have $\dim(\schP_a) = \dim(\stP_a)+ \dim(\pi_*(Z_{\cG})$. Therefore it suffices to show that $\dim(\stP_a) = \frac{1}{2}\dim(G)\deg(\cL)-\frac{1}{2}\rank(G)(\deg(\cL)-2g+2)$. After base extension, we may assume without loss of generality that $a$ is a $k$-point. Let $J_{a}$ denote the restriction of $J_{A}$ to $C = C \times a \hookrightarrow C \times A$. By standard deformation theory, the tangent complex of the smooth stack of $J_a$-torsors $\stP_a$ restricted to any $k$-point of $\stP_a$ is given by $R\pi_*(\text{Lie}(J_a)[1])$. Therefore $\stP_a$ is equidimensional of dimension equal to minus the Euler characteristic $-\chi(\text{Lie}(J_a)$. Note that $J_a$ is a vector bundle on $C$ of dimension $\rank(\cG)$. By Riemann-Roch, we have
    \begin{align*}
           \dim(P^o/A)= -\chi(\mathrm{Lie}(J_{a}))=&-\deg(\mathrm{Lie}(J_{a}))+\rank(\cG)(g-1)=\\
            =&-|\Gamma|^{-1}\deg(n^*\mathrm{Lie}(J_{a}))+\rank(\cG)(g-1).
     \end{align*}
    The pullback $n^*(J_a)$ under $n: \widetilde{C} \to C$ is the corresponding global centralizer group scheme $J_{n_A(a)}$ on $\widetilde{C}$ for the group reductive $G$ at $n_A(a) \in \widetilde{A}(k)$. We may then use the formula for $\text{Lie}(J_{n_A(a)})$ in \Cref{P:galois-descr-j-cam} to express the degree in terms of the dimension of a Borel subgroup $B \subset G$ to get:      
    \begin{align*}
               \dim(P^o/A)&= -|\Gamma|^{-1}\deg(n^*\mathrm{Lie}(J_{a}))+\rank(\cG)(g-1)=\\
            &=-|\Gamma|^{-1}\dim(B)\deg(n^*\cL)+\rank(\cG)(g-1)=\\
            &=-\dim(B)\deg(\cL)+\rank(\cG)(g-1)=\\
            &=\frac{1}{2}\dim(\cG/C)\deg(\cL)-\frac{1}{2}\rank(\cG)(\deg(\cL)-2(g-1)).
            \end{align*}
\end{proof}

 Note that there is a central inclusion of $\pi_*(Z_{\cG}) \times \Higgs_{\cG, \cL}$ into the inertia stack of $\Higgs_{\cG, \cL}$, because the group scheme $\pi_*(Z_{\cG})$ has a canonical inclusion into the center of the automorphism group of very point of $\Higgs_{\cG, \cL}$.

\begin{defn}
    We set $\Higgsrig_{\cG, \cL}:= \Higgs_{\cG,\cL} \fatslash \pi_*(Z_{\cG})$ to be the $\pi_*(Z_{\cG})$-rigidification of $\Higgs_{\cG,\cL}$ in the sense of \cite[Definition 5.1.9]{acv-twisted-bundles-covers}.
\end{defn}

\begin{coroll} \label{coroll: action group on rigidication of twisted higgs}
    Suppose that $\deg(\cL)>0$, and the characteristic of $k$ is either zero or strictly larger than the height of the adjoint representation of $G$. Then the action of $\stP^o$ on $\Higgs_{\cG,\cL}$ induces an action of the moduli space $\schP^o$ on the rigidification $\Higgsrig_{\cG, \cL}$.
\end{coroll}
\begin{proof}
    We have seen in the proof of \Cref{prop: gerbe group scheme of symmetries} that $\schP^o$ is the $\pi_*(Z_{\cG})$ rigidification under a natural central inclusion of $\pi_*(Z_{\cG})$ inside the inertia of $\stP^o$. It follows from construction that the action of $\stP^o$ on $\Higgs_{\cG, \cL}$ is compatible with the inclusions of $\pi_*(Z_{\cG})$ inside the inertias of the stacks. Therefore, by the universal property of rigidification \cite[Thm. 5.1.5(2)]{acv-twisted-bundles-covers}, we get an induced action of $\schP^o = \stP^o \fatslash \pi_*(Z_{\cG})$ on $\Higgsrig_{\cG,\cL} = \Higgs_{\cG,\cL} \fatslash \pi_*(Z_{\cG})$.
\end{proof}

\subsection{Weak Abelian fibration structure}

\begin{thm} \label{thm: weak Abelian fibration structure twisted higgs}
    Assume that all of the following hold:
    \begin{enumerate}[(1)]
    \item $\text{char}(k)$ is either zero or strictly larger than the height of the adjoint representation of $G$.
    \item $\text{char}(k)$ does not divide $|\Gamma|$.
     
    \item $\deg(\cL) > \max\{2g-2,0\}$
    \end{enumerate}
     Fix $d \in \pi_1(G)$. Then the action of the smooth commutative group scheme $\schP^o \to A$ on the moduli space $h:\MHiggs_{\cG, \cL}^d \to A$ (cf. \Cref{prop: action on twisted higgs}) makes the triple $(\MHiggs_{\cG, \cL}^d, \schP^o, A)$ into a weak Abelian fibration (\Cref{defn: weak Abelian fibration}).
\end{thm}
\begin{proof}
    The scheme $\MHiggs_{\cG,\cL}^d$ is quasiprojective and geometrically normal (and hence equidimensional) by \Cref{coroll: properties twisted isotrivial Higgs moduli space}. By \Cref{prop: properness of the hitchin fibration} the Hitchin morphism $h: \MHiggs_{\cG,\cL}^d \to A$ is proper. By \Cref{prop: gerbe group scheme of symmetries} we have that $\schP^o \to A$ is a quasiprojective smooth commutative group scheme, and therefore it follows that the relative Tate module $T_{\overline{\mathbb Q}_\ell}(\schP^o)$ is polarizable by \cite[Thm. 6.2]{ancona-fratila-polarizability}. Finally, by \Cref{prop: action on twisted higgs} the action of $\schP^o$ on $\MHiggs_{\cG,\cL}^d$ has affine stabilizers. We have therefore checked all the properties in the definition of weak Abelian fibration (\Cref{defn: weak Abelian fibration}).
\end{proof}

\begin{remark}\label{rmk: gamma and char}
    The most interesting case for us will be when $\cG$ is an isotrivial outer form of $G$, which is obtained by twisting by a torsor for a subgroup of the group $\text{Out}(G)$ of outer automorphisms. In this case the condition that $\text{char}(k)$ does not divide $|\Gamma|$ is automatic as long as $\text{char}(k)$ does not divide the order of the group of automorphisms of the Dynkin diagram of $G$. 
\end{remark}
    
\end{section}

\begin{section}{Intersection cohomology of good moduli spaces} \label{section: intersection cohomology of gms}

In this section, all  stacks are algebraic stacks of finite type over an algebraically closed ground field $k$ and all morphisms are of finite type.

  Our reference for perverse 
sheaves on a stack is \cite{laszlo-olsson-perverse-sheaves} and for the 
derived categories 
of complexes of constructible sheaves
on a stack  
is \cite{laszlo-olsson-adic-sheaves}.
We work with the field of coefficients $\overline{\mathbb{Q}}_{\ell}$, where $\ell$ 
is a prime that is invertible in the ground field $k$ and with the derived
categories $D^*_c(-)= D^*_c(-,\overline{\mathbb{Q}}_{\ell} )$ of  constructible complexes of sheaves
with $*=b, +, -, \emptyset$ (bounded, bounded below, bounded above and unbounded, resp.).
By a  lisse (smooth)  sheaf we mean
a lisse 
$\overline{\mathbb{Q}}_{\ell}$-sheaf of finite rank.
We work with both the standard $t$-structure (with full subcategories
$D_c^{*,\leq a}(-), D_c^{*,\geq d}(-))$ 
and with hearts the category of constructible sheaves),
and with the middle-perversity $t$-structure
on $D^*_c(-)$
(with full subcategories
$^p\!D_c^{*,\leq a}(-), {^p\!D_c}^{*,\geq d}(-)$) and with hearts the category of $\overline{\mathbb{Q}}_{\ell}$-adic perverse sheaves).

The main result of this section is \Cref{thm: intersection cohomology of the stack vs the moduli space}, which relates the intersection cohomology of a stack to the one of its good moduli space. This generalizes some results in \cite{kirwan-blowups, woolf-intersection}. A similar result holds in the context of mixed Hodge modules when the ground field is $\mathbb{C}$, see \Cref{rmk: ic and ic mhm}.

\begin{remark}
The recent preprint \cite{kinjo_decomposition_theorem} establishes, over the field of complex numbers,  the purity
of the direct image complex in the context of the derived category of bounded below algebraic mixed Hodge modules on a stack, and can be used to yield a different proof of \Cref{thm: intersection cohomology of the stack vs the moduli space}.
\end{remark}

Before stating and proving \Cref{thm: intersection cohomology of the stack vs the moduli space} , we review and develop some preliminaries on cohomology of stacks.

\begin{subsection}{Cohomology of algebraic stacks} \label{subsection: cohomology of algebraic stacks}

In this section, we work with the bounded derived categories $D^b_c(-)$; this is made possible by the fact that all morphisms of stacks in this section have finite relative inertia.

Let
$J: \scrZ \stackrel{j}\to
\overline{\scrZ} \stackrel{i}\to \scrX$  
be a locally closed embedding of integral stacks factored
as the open immersion $j$ of a smooth $\scrZ$ into
its closure $\overline{\scrZ}$ in $\scrX$,  followed by the closed immersion $i$. For a lisse sheaf $\mathcal{L}$ on $\scrZ$, the intersection complex $IC_{\overline{\scrZ}} (\mathcal L)$ on $\overline{\scrZ}$ is the perverse sheaf on $\overline{\scrZ}$ 
\[IC_{\overline{\scrZ}} (\mathcal L) \coloneqq 
%Ri_* 
j_{!*} (\mathcal{L} [\dim (\scrZ)])\]
defined via  the 
%pushforward via $i$ of  
intermediate extension   functor
$j_{!*}$
%[\dim (\scrZ)])
(cf. \cite[\S6]{laszlo-olsson-perverse-sheaves}).
As is well-known, we obtain the same complex if we start 
with another such pair $(\scrZ',\mathcal{L}')$
such that $\mathcal L$ and $\mathcal L'$ agree on the
non-empty intersection $\scrZ\cap \scrZ'.$
By a slight and common abuse of notation, rooted in the fact that the functor $Ri_*$ is fully faithful, we
use the same notation and terminology for the perverse sheaf obtained by  extension by zero to $\scrX$, i.e.
we blur the distinction between $IC_{\overline{\scrZ}}$
 and $Ri_* IC_{\overline{\scrZ}}.$
The intersection complex $IC_{\scrX}$
is obtained in the special case where the lisse sheaf is the constant sheaf of rank one on
any fixed open  dense and smooth substack of $\scrX$. In general, if $\scrX$ is a possibly non-integral stack, then the intersection complex $IC_{\scrX}$ is the direct sum of the intersection complexes of its irreducible  components, each equipped with the reduced substack structure. 
The simple perverse sheaves on a stack are the intersection complexes 
$IC_{\overline{\scrZ}} (\mathcal L)$ on $\scrX$ 
associated with
irreducible lisse sheaves $\mathcal{L}$ on smooth locally closed irreducible substacks $\mathfrak Z$ of $\scrX$
(cf.  \cite[Thm. 8.2(ii)]{laszlo-olsson-perverse-sheaves}). 
A semisimple perverse sheaf is a finite direct sum of simple perverse sheaves. A semisimple complex in $D(\scrX)$ is one isomorphic to a finite
direct sum of shifted semisimple perverse sheaves.

\begin{notn}
    When seemingly convenient, in essence to avoid using too many shifts, we occasionally  employ the shifted version of the intersection complex (a.k.a. the topologist's intersection complex) for an integral stack $\scrZ$, i.e. $\mathcal{IC}_{\scrZ}(\mathcal L):= IC_{\scrZ} (\mathcal L)
[-\text{dim} (\scrZ)],$ which is a bounded complex with trivial cohomology sheaves in negative degrees.
 With this notation $IC_{{\mathbb A}^1}  = 
(\overline{\mathbb Q}_{\ell})_{{\mathbb A}^1} [1],$ whereas
$\mathcal{IC}_{{\mathbb A}^1}  = 
(\overline{\mathbb Q}_{\ell})_{{\mathbb A}^1}.$
For not necessarily integral stacks, we take as the topologist's intersection complex, the 
direct sum of the topologist's intersection complexes of
the irreducible components with their reduced structure.
\end{notn}

\begin{prop}[Decomposition theorem for proper morphisms of stacks] \label{lemma: decomposition theorem for stacks}
Let $\pi: \scrX \to \scrY$ be a proper morphism between stacks with affine stabilizers and of finite type over an algebraically closed field $k.$ 
Then the push-forward of a semisimple complex is a semisimple complex. In particular, there is a finite direct sum decomposition: 
\[
\theta:
\bigoplus_{i \in I} \mathcal{K}_i[s_i] \stackrel{\simeq}\to 
R\pi_*(IC_{\scrX}), \]
where the summands $\mathcal{K}_i[s_i]$ are
shifts of simple perverse sheaves on $\scrY$.
\end{prop}
\begin{proof}
The result is proven
%The existence of a  decomposition $R\pi_*(IC_{\scrX}) \simeq \bigoplus_{i \in I} \mathcal{K}_i[s_i]$ 
in \cite{shenghao-decomposition-stacks} in the case when the ground field $k$ is the algebraic closure of a finite field. This is generalized to the case when $k=\mathbb C$ in \cite{sun-generic-base-change},
 both for the \'etale and classical topologies.
 Loc. cit. can be generalized in the same way to arbitrary algebraically closed fields $k$ and for the \'etale topology
 by a standard spreading out technique  as in \cite[\S6]{bbdg}, coupled with the Generic Base Change Theorem for Stacks \cite[\S4.1]{sun-generic-base-change}
 (replace $\mathbb C$ with $k$ in loc. cit.).
\end{proof}

\begin{remark}
We note that the requirement in \cite{shenghao-decomposition-stacks} that the morphism $\pi$ has finite relative diagonal is redundant, as it is implied by the properness of $\pi$ and the fact that $\scrX$ and $\scrY$ have affine stabilizers. Indeed, it follows that $\pi$ has affine relative stabilizers. As a consequence of properness it follows that the diagonal of $\pi$ must be quasifinite, because every affine proper scheme over a field is finite. Since a proper quasifinite representable morphism is automatically finite and schematic \cite[\href{https://stacks.math.columbia.edu/tag/03XX}{Tag 03XX}]{stacks-project}, it follows that the relative diagonal is finite and schematic.
 Some results below will be stated for stacks with affine diagonal, a stronger condition that implies that they have affine stabilizers. In particular, we may apply \Cref{prop: decomposition theorem with generic connected fibers} for stacks with affine diagonal without further remarks.
\end{remark}

\begin{remark}\label{remark: basic proeprties of dt}
We collect remarks on the Decomposition Theorem that are useful in the forthcoming sections. They are all standard in the case of schemes which implies, by checking on smooth covers, that they  also hold  for stacks.
%For ease of exposition, we assume the stacks to be integral.

\begin{enumerate}

\item
Given a morphism $f: \mathscr{A} \to \mathscr{B}$ of stacks, we denote the unit morphism
of functors
as
$f^*: {\rm Id} \to Rf_* f^*.
$
In cohomology, it corresponds to the pull-back.
The natural unit morphism
$f^*: (\overline{\mathbb{Q}}_{\ell})_{\mathscr{B}} \to Rf_* (\overline{\mathbb{Q}}_{\ell})_{\mathscr{A}}$
induces the pull-back of locally constant functions
$f^*: (\overline{\mathbb{Q}}_{\ell})_{\mathscr{B}} \to R^0f_* (\overline{\mathbb{Q}}_{\ell})_{\mathscr{A}}.$
In particular, when $f$ is dominant, this latter $f^*$ is injective.

\item
Let $j\colon \mathscr{B}^o \to \mathscr{B} \leftarrow \mathscr{A}\colon i$ be complementary open and closed immersions of stacks. Let $j_{!*}P$ be the intermediate extension to $\mathscr{B}$ of a perverse sheaf
on the open $\mathscr{B}^o.$ Then 
$Rj_* P \in {^p\!D^{\geq 0} (\mathscr{B})}$ and  $i_! i^! j_{!_*} P 
\in {^p\!D^{\geq 1} (\mathscr{B})}$  
(cf. \cite[4.2.4]{bbdg}).
% {\cre Something seems to be off here: If, in addition, we have that $\dim{(\mathscr{A})}
% <\dim{(\mathscr{B})},$ then
% $i_! i^! j_{!_*} P
% \in {D^{\geq -  \dim{(\mathscr{B})}} (\mathscr{B})}$
% .}  
 As a consequence,  we have that 
$i_! i^! j_{!_*} P
\in {D^{\geq -  \dim{(\overline{\mathscr{B}^o})}+2} (\mathscr{B})}.$ 

\item Let $\scrZ  \stackrel{j}\to \overline{\scrZ} \stackrel{i}\to \scrX$, $J:= i \circ j$ and  $\mathcal L$ be as at the beginning of this section. As an immediate consequence of part 2. of this remark, 
we have that the natural unit morphism
$J^*:\mathcal{IC}_{\overline{\scrZ}}(\mathcal L) \to RJ_* \mathcal L [\dim (\scrZ)]$
induces  a natural  isomorphism on the $0$-th cohomology sheaves,  which  restricts to the identity $\mathcal{L} \stackrel{=}\to \mathcal{L}$ on $\scrZ.$

\item
Let $\scrX$ be an integral stack.
There is a natural morphism
$
\epsilon: (\overline{\mathbb{Q}}_\ell)_{\scrX} 
\to \mathcal{IC}_{\scrX},
$
which is characterized by being the identity 
on the smooth locus; to see this
apply the functor $\text{Hom}((\overline{\mathbb{Q}}_\ell)_{\scrX},-)$ to the attaching triangle
for $\ic{\st{X}}$ as in part 2. of this remark.
By using this same technique, and the last statement in part 2., we also see the following.
Let $j=j_{\scrX^{o}}: \scrX^{o}  \to \scrX^{reg} \to \scrX$ be the open immersion  of any open dense substack  of the regular open substack
${\scrX^{reg}}.$
There is a canonical factorization
\[j^*:(\overline{\mathbb{Q}}_\ell)_{\scrX}
\stackrel{\epsilon}\to \mathcal{IC}_{\scrX}
\stackrel{j_{\mathcal{IC}}^*}\to 
R j_* (\overline{\mathbb{Q}}_\ell)_{\scrX^{o}}\]
of the unit morphism $j^*$ for 
$(\overline{\mathbb{Q}}_\ell)_{\scrX}$
through the unit morphism $j_{\mathcal{IC}}^*$ for $\mathcal{IC}_{\scrX}.$
We thus have morphisms of  $0$-th cohomology sheaves:
\begin{equation}\label{iso hzero}
j^*: (\overline{\mathbb{Q}}_\ell)_{\scrX}
\stackrel{\epsilon}\to 
{\mathcal H}^0 (\mathcal{IC}_{\scrX})
\xrightarrow[=]{j^*_{\mathcal IC}}
R^0 j_* (\overline{\mathbb{Q}}_\ell)_{\scrX^o},
\end{equation}
where the second one, $j^*_{\mathcal IC},$ 
is an isomorphism.
%(cf. \cite[Lem. 1]{Durfee95}}).
In this way, we  identifty the morphism of $0$-th cohomology sheaves $\epsilon = j^*.$
This morphism on $0$-th cohomology sheaves is injective;
in general, it is not an isomorphism. 

\item
Let $\scrX$ be an integral stack that is universally unibranch, e.g. when it is normal.
Then the morphism $\epsilon$
induces the canonical 
isomorphism
$
(\overline{\mathbb{Q}}_\ell)_{\scrX}
\stackrel{=}\to \mathcal{H}^0
(
\mathcal{IC}_{\scrX}).$
To see this, we argue as follows.
First, we claim that the natural morphism $j^*$ (\ref{iso hzero}) is an isomorphism. This can be checked at the level of stalks, where we use, having passed to schemes, that by removing a proper closed 
subscheme, the strict Henselianization of a universally unibranch scheme remains 
connected
(cf. \cite[\href{https://stacks.math.columbia.edu/tag/06DM}{06DM}]{stacks-project}).
The desired assertion now follows from the
identification 
$ \mathcal{H}^0
(
\mathcal{IC}_{\scrX}) = R^0 j_* (\overline{\mathbb{Q}}_\ell)_{\scrX^o},$

%% Leave paragraph below in case 
\iffalse
CLAIM: the compositum morphism $j^*: (\overline{\mathbb{Q}}_\ell)_{X}
\to j_* ( \overline{\mathbb{Q}}_\ell )_{U}$ 
is an isomorphism. This can be verified at the level of stalks at the geometric points of $Y.$
The stalk of the morphism  $j^*$  at a geometric  point is the restriction morphism
on the $0$-th cohomology groups from 
 the (spectrum of the) strict henselianization at the point
 to the same scheme with a proper closed subscheme removed (cf. [SGA 4.5, p23, 3rd equation)]).
Since by assumption (cf. \href{https://stacks.math.columbia.edu/tag/06DM}{06DM}) both 
spaces are irreducible,  hence connected, the stalk of $j^*$ is an isomorphism and the 
CLAIM is proved.
\fi

\item A perverse sheaf  $\mathcal P$ on a stack has a closed reduced substack support. If $\mathcal P$ 
is perverse semisimple,
    then it admits (apply the method of \cite[Cor. 5.3.11]{bbdg}) a 
canonical decomposition 
$\iota = \sum_{\mathfrak S} \iota_{\mathfrak S} : \bigoplus_\mathfrak{S}
IC_{\mathfrak S} ({\mathcal L}_{\mathfrak S})
\stackrel{=}\to \mathcal P$
into a finite direct sum of semisimple perverse sheaves with distinct non-empty irreducible closed  supports
${\mathfrak S}$.
This is called the canonical decomposition by supports
of $\mathcal P.$

\item \label{item:theta}
\Cref{lemma: decomposition theorem for stacks}
can be re-stated as the existence of a non-canonical isomorphism 
\[\theta: \bigoplus_{b \in {\mathbb Z}} \mathcal{S}_b[-b] \stackrel{\simeq}\to R\pi_* ( IC_{\scrX} ), \]
where the $\mathcal{S}_b$  are the semisimple perverse cohomology sheaves of the direct image.
Each of the $\mathcal{S}_b$ admits its own canonical decomposition by supports.
It is easy to modify any given isomorphism $\theta$ to a new one, still noncanonical,  so that it induces the identity
on the $\mathcal{S}_b.$
To this end, modify the morphism  $\theta$ by pre-composing it with the sum of the appropriately shifted inverses of the morphisms induced by $\theta$  in perverse cohomology, namely $\sum_b
(({^p\!H^b} (\theta))^{-1} [-b]).$

\item
If $\scrY^o \subseteq \scrY$
is an open and dense substack, we usually denote by $-^o$ 
the objects obtained by base change via this open immersion.
In \Cref{lemma: decomposition theorem for stacks},
after perhaps shrinking $\scrY$ to a suitably small dense open substack
${\scrY}^o,$ with pre-image $\scrX^o,$ we may assume that we have a noncanonical 
isomorphism 
$
\bigoplus_{i\geq 0} \widehat{\mathcal{L}}_i[-i]
\simeq 
R\pi_*(\mathcal{IC}_{\scrX^o}),
$ 
inducing the identity on perverse cohomology sheaves, and where the $\widehat{\mathcal{L}}_i$ are  semisimple lisse sheaves on  $\scrY^o.$ The resulting isomorphism $\widehat{\mathcal{L}}_0
\stackrel{=}\to 
R^0 \pi_* (\mathcal{IC}_{\widehat{\scrX}^o}),$ being obtained by standard $0$-th truncation, is in fact canonical.

\end{enumerate}
\end{remark}

\begin{notn}
By abuse of notation, we also denote by $\epsilon$ and $f^*$ (as in \Cref{remark: basic proeprties of dt}) some of the morphisms naturally associated with 
them, e.g. by push-forward, pull-back, etc.  We also label natural isomorphisms simply by the symbol $``="$.
\end{notn}

\begin{prop}[Existence of a split unit morphism for intersection cohomology] \label{prop: decomposition theorem with generic connected fibers}
    Let $\pi: \scrX \to \scrY$ be a proper and surjective morphism of  stacks 
with affine stabilizers and of finite type over an 
algebraically closed field $k.$ Then there is a commutative square: 
\begin{equation}
\begin{tikzcd}[ampersand replacement=\&] \label{fig: diagram 2}
(\overline{\mathbb{Q}}_{\ell})_{\scrY}
\arrow{r}{\pi^*} 
\arrow{d}{\epsilon} 
\& 
R\pi_* (\overline{\mathbb{Q}}_{\ell})_{\scrX}
\arrow{d}{\epsilon}
\\ 
\mathcal{IC}_{\scrY}
\arrow[dotted]{r}{[\widehat{\pi}^*]} 
\& 
R\pi_* \mathcal{IC}_{\scrX},
\end{tikzcd}
\end{equation}
where: $[\widehat{\pi}^*]$ is a noncanonical split injection; 
the canonical morphisms $\epsilon$ and $\pi^*$ are defined in Remark \ref{remark: basic proeprties of dt}. 
In particular: the induced morphisms on $H^0(\scrY,-)$ send  unit sections to the corresponding unit sections; the restriction of 
$[\hat{\pi}^*]$ to a suitably small open and dense substack of $\scrY$ coincides with the corresponding restriction of $\epsilon \circ \pi^*.$
\end{prop}

\begin{proof}
By working with irreducible components given their reduced structure,
we may assume that $\scrX$ and $\scrY$ are integral.
 The natural morphism $\epsilon: (\overline{\mathbb{Q}}_{\ell})_{\scrX}
\to \mathcal{IC}_{\scrX}$ induces an injection
$\epsilon: (\overline{\mathbb{Q}}_{\ell})_{\scrX} \to  \mathcal{H}^0(\mathcal{IC}_{\scrX}),$
and thus the   injections
$(\overline{\mathbb{Q}}_{\ell})_{\scrY} 
\stackrel{\pi^*}\to R^0 \pi_* (\overline{\mathbb{Q}}_{\ell})_{\scrX}
\stackrel{\epsilon}\to R^0 \pi_*  (\mathcal{IC}_{\scrX}).$
The unit sections correspond to each other
via these morphisms.

Let $\theta$ be an isomorphism as in \Cref{lemma: decomposition theorem for stacks}
and inducing the identity on the semisimple perverse cohomology sheaves. 
%REMOVE???? Any such automatically preserves the canonical decomposition by supports (see \Cref{remark: basic proeprties of dt}.\eqref{item:theta}).

 Let $\scrY^o \subseteq \scrY$ be any smooth and dense open substack with preimage $\scrX^{o} := \pi^{-1}(\scrY^{o})$ such that  
$\theta$ takes the form of an isomorphism 
$\oplus_{\geq 0} \widehat{\mathcal{L}_i} [-i]
\stackrel{=}\to 
R\pi_* (\mathcal IC_{{\scrX}^o}),$
with $\widehat{\mathcal{L}}_i$ semisimple lisse sheaves on $\mathfrak {Y}^o$. 
Note that the resulting isomorphism $\widehat{\mathcal{L}}_0 \stackrel{=}\to \mathcal{H}^0 R^0 \pi_* (\mathcal{IC}_{\scrX^o}),$
being obtained by applying the earliest possible non-trivial truncation, is canonical.
Furthermore, by shrinking $\scrY^o,$ we can assume that the domains and codomains of
the  injections, still denoted by  $\pi^*$ and
$\epsilon,$ obtained by restricting  to $\scrY^o$ the injections 
of sheaves $\pi^*$ and $\epsilon$ above, are
split injections of semisimple lisse sheaves.

 By the functoriality of the intermediate extension functor, 
the split morphism $\epsilon \circ \pi^*: (\overline{\mathbb{Q}}_{\ell})_{{\scrY}^o} 
\to  \widehat{\mathcal{L}}_0$ extends to a split morphism 
$[\epsilon \circ \pi^*]: \mathcal{IC}_{\scrY} \to \mathcal{IC}_{\scrY} ( \widehat{\mathcal{L}}_0);$ this extension is easily seen to be unique, however, we do not need this uniqueness.
%use: attaching triangle; i$^! \in ^pD^{\geq 1}$ and adjnction $(j^*,Rj_*)$ 
This morphism sends the unit section to the unit section.

By considering the canonical decomposition by supports,
there is a canonical split morphism $\iota_{\scrY}: \mathcal{IC}_{\scrY} ( \widehat{\mathcal{L}}_0) \to ({^p\!H^{-\dim ( \scrY )}} (R\pi_* \mathcal{IC}_{\scrX}))[\dim ( \scrY )].$

We set:
\begin{equation}\label{the morphism}
[\widehat{\pi}^*]: = \theta \circ \iota_{\scrY} \circ [\epsilon \circ \pi^*].
\end{equation}
By  construction, all the morphisms in 
diagram (\ref{fig: diagram 2}) send unit sections to unit sections. %\Mark{Is this already enough? In general no; think of $Q \to Q \oplus Q_{pt}$, where we have two different morphisms agreeing away from the point. This does not happen here, but how to say it? I wrote a long-winded argument. Would be nice to have a simpler one.}

Next, we show that the square diagram (\ref{fig: diagram 2}) is commutative. Since the upper left corner is a sheaf, and the other complexes are in $D^{\geq 0}(\scrY),$ the desired commutativity
is equivalent to the commutativity 
of the diagram obtained by taking 
 $0$-th cohomology sheaves: this is because truncation
 $\tau_{\leq 0}$
 is a right adjoint to  the inclusion $\iota_0:
 D^{\leq 0} (\scrY) \to D(\scrY)$,
 so that we have the canonical identification:
\begin{equation}\label{rrr}
\text{Hom} 
((\overline{\mathbb{Q}}_{\ell})_{{\scrY}},  R\pi_* \mathcal{IC}_{\scrX})   = 
\text{Hom} 
((\overline{\mathbb{Q}}_{\ell})_{{\scrY}},  R^0\pi_* \mathcal{IC}_{\scrX}).
\end{equation}
By applying the functor $\text{Hom} ((\overline{\mathbb{Q}}_{\ell})_{{\scrY}},-)$
 to the attaching triangle
 for the open immersion $j: \scrY^o \to \scrY$, with complementary closed immersion $i,$
 we get the exact sequence:
 \begin{equation}\label{ttt}
%\begin{tikzcd}\label{fig: diagram 2}
\text{Hom} ((\overline{\mathbb{Q}}_{\ell})_{{\scrY}}, i_! R\pi_* i^!\mathcal{IC}_{\scrX} ) 
\to
\text{Hom} 
((\overline{\mathbb{Q}}_{\ell})_{{\scrY}},  R^0\pi_* \mathcal{IC}_{\scrX})  
\to
\text{Hom} ((\overline{\mathbb{Q}}_{\ell})_{{\scrY}^0}, \widehat{\mathcal{L}}_0), 
%Rj_* j^* R\pi_* \mathcal{IC}_{\scrX}  ) ,
%\end{tikzcd}
  \end{equation}
where: in  the first term, we have used the base change isomorphism $i^! R\pi_* = R\pi_* i^!;$ 
for the second  term we have used (\ref{rrr});
for the last term
 we have used the adjunction identities for $(j^*, Rj_*),$ as well as the argument leading to (\ref{rrr}).
 
Since the diagram (\ref{fig: diagram 2}) is commutative over $\scrY^o$ by construction,
 in order 
 to prove the  commutativity over $\scrY,$ it remains to show that the first term in
 (\ref{ttt}) vanishes.
 This follows from the fact that $i_!i^! \mathcal{IC}_{\scrX} \in D^{\geq 2} (\scrX)$ (c.f. \Cref{remark: basic proeprties of dt}, part 2). 

Finally, by shrinking $\scrY^o$ further if necessary, we may assume that the restriction morphism $\epsilon: (\overline{\mathbb{Q}}_{\ell})_{\scrY}\to \mathcal{IC}_{\scrY}$ in diagram (\ref{fig: diagram 2}) to $\scrY^o$  is the identity morphism, so that the remaining assertion of the proposition follows.
\end{proof}

\begin{remark}\label{pistar}
If, in \Cref{prop: decomposition theorem with generic connected fibers}, the morphism $\pi$
has the property that there is an open and dense substack $\st{Y}'\subseteq \st{Y}$ such that
$\st{Y}'_{red}$ and $\pi^{-1}(\st{Y}')_{red}$
are smooth (e.g. $\pi$
is birational schematic, or  
generically finite schematic) then there is an open and dense substack of $\scrY$ where  the restrictions of 
$[\widehat{\pi}^*]$
and $\pi^*$ coincide. Indeed the morphisms of type $\epsilon$ are the identity on these open sets.
In this case, we denote $[\widehat{\pi}^*]$
simply by $[\pi^*].$
\end{remark}

\begin{remark}\label{not pistar}
Note that in general it is not true that the restriction of $[\widehat{\pi}^*]$
to some open and dense substack of $\scrY$ coincides with $\pi^*.$ This is true if $\scrX$ is normal, and this explains the choice of notation
$[\widehat{\pi}^*]$ and $\widehat{\mathcal{L}_0}.$ 
Consider the  normalization $\nu: \widehat{\scrX} \to \scrX$ \cite[Defn. 1.18]{vistoli-intersection-theory}. We denote
$\hat{\pi}: \widehat{\scrX} \stackrel{\nu}\to {\scrX} \stackrel{\pi}\to {\scrY}$. 
To see this, we argue as follows.
By using the same method of proof of
\Cref{prop: decomposition theorem with generic connected fibers} and \Cref{remark: basic proeprties of dt}, one can show that: 
$\widehat{\mathcal{L}}_0 = 
R^0 \widehat{\pi}_* (\overline{\mathbb{Q}}_{\ell})_{\widehat{\scrX}^o};$
$\widehat{\pi}^*$ splits canonically;
$[\widehat{\pi}^*]$ extends from ${\scrY}^o$
to ${\scrY}$
the canonically split unit morphism $\widehat{\pi}^* =  \nu^* \circ \pi^*:
(\overline{\mathbb{Q}}_{\ell})_{\widehat{\scrX}^o}
\to 
R^0 \pi_* (\overline{\mathbb{Q}}_{\ell})_{{\scrX}^o}
\to
R^0 \widehat{\pi}_* (\overline{\mathbb{Q}}_{\ell})_{\widehat{\scrX}^o}.$
\end{remark}

The following lemma is probably well-known. We include it for lack of a reference.

\begin{lemma}[Proper base change for stacks]
Let $\pi: \scrX \to \scrY$ be a proper morphism between stacks $\scrX$ and $\scrY$ of finite type over an algebraically closed field $k$. Let $f: \mathscr{W} \to \scrY$ be a morphism from a stack $\mathscr{W}$ of finite type over $k$, inducing a Cartesian diagram 
\[
\begin{tikzcd}
\scrX_{\mathscr{W}} 
\arrow{r}{f} 
\arrow{d}{\pi} 
& 
\scrX
\arrow{d}{\pi}
\\ 
\mathscr{W}
\arrow{r}{f} 
& 
\scrY.
\end{tikzcd}\]
For all constructible complexes of $\overline{\mathbb{Q}}_{\ell}$-sheaves $F \in D(\scrX),$ the induced base-change morphism $b: f^*R\pi_*F \to R\pi_*f^*F$ is an isomorphism.
\end{lemma}
\begin{proof}
Let $W \to \mathscr{W}$ be 
a smooth atlas, with $W$ a scheme. By smooth base-change 
\cite[12.1]{laszlo-olsson-adic-sheaves} and the compatibility 
of the formation of the base-change morphism $b$ with further 
base-change, we can without loss of generality replace 
$\mathscr{W}$ with $W$. Similarly, let $Y \to \scrY$ be 
a smooth atlas with $Y$ a scheme. By base-changing everything 
with $Y,$ and  possibly passing to an \'etale cover of $W 
\times_{\scrY}Y$, we can assume without loss of 
generality that both $\scrY=Y$ and $\mathscr{W} = W$ 
are schemes. 
 
 We need to check that the base change  morphism $b: f^*R\pi_*F \to R\pi_*f^*F$  in the derived category of the scheme $W$ is an isomorphism. 
We can check this at the level of stalks at the geometric points of $W$. 
Since $W$ is of finite type over $k,$ it is enough to check this at the level of stalks of $k$-points
of $W$, because the cohomologies of the complexes are constructible \cite[Prop. 8.2]{laszlo-olsson-adic-sheaves}. For any $x: \Spec(k) \to W$, both $x$  and the composition $f \circ x: \Spec(k) \to Y$ are closed immersions.
In view of  the properness of $\pi,$ we have the natural identification $R\pi_* = R\pi_!$ (cf. \cite[Thm. 1.12]{olsson-fujiwara}), so that, by applying twice the base-change result for closed immersions in \cite[12.3]{laszlo-olsson-adic-sheaves}, we conclude that the $x$-stalk of the base-change morphism $b$ is an isomorphism, as desired.
\end{proof}

In the rest of the section, we prove a refinement (\Cref{lemma: intersection cohomology of 
tame moduli spaces}) of \Cref{prop: decomposition theorem with generic connected fibers} for separated universal homeomorphisms of stacks, e.g.\ the morphism from a tame stack to its coarse moduli space. See \Cref{lemma: intersection cohomology of moduli space generically} for an application.   To this end, we start with a topological preliminary.

Recall the notion of universal homeomorphism of stacks \cite[\href{https://stacks.math.columbia.edu/tag/0CI5}{0CI5}]{stacks-project}.

\begin{lemma}\label{claim:homeo}
    Let $\pi: \scrX \to \scrY$ be a separated universal homemorphism  
between stacks with affine stabilizers and of finite type over an algebraically closed field $k$. Then the unit morphism:
\[
\begin{tikzcd}\label{fig: diagram 2.1}
(\overline{\mathbb{Q}}_{\ell})_{\scrY}
\arrow{r}{\pi^*}[swap]{=} 
& 
R\pi_* (\overline{\mathbb{Q}}_{\ell})_{\scrX}
\end{tikzcd}
\]
is an isomorphism.
\end{lemma}

\begin{remark} \Cref{claim:homeo} is clear for schemes and thus for  algebraic spaces. It is false for arbitrary universal homeomorphism of stacks, e.g.\ the nonseparated morphism $B\mathbb{G}_m \to \Spec(k)$.
% a stack is separated is the diagonal is proper
\end{remark}

\begin{proof}[Proof of \Cref{claim:homeo}]
Note that the hypotheses imply that $\pi$ is proper and quasifinite. In particular, by properness, we have a natrual identification  $R\pi_* = R\pi_!$ by \cite[Thm. 1.12]{olsson-fujiwara}. By smooth base-change \cite[12.1]{laszlo-olsson-adic-sheaves}, the statement can be checked smooth 
locally on $\scrY$. Hence, we can assume without loss of 
generality that $\scrY$ is a scheme of finite type over $k$. Since $\scrY$ is of finite type over the 
algebraically closed field $k$ and $\pi^*$ is a morphism of constructible sheaves, it suffices to check this at the stalk of every $k$-point $y: \Spec(k) \to \scrY$. By the base change isomorphism in \cite[12.3]{laszlo-olsson-adic-sheaves}, there is an identification $(R\pi_*(\overline{\mathbb{Q}}_{\ell})_{\scrX})_y = R\Gamma(\scrX_y, \overline{\mathbb{Q}}_{\ell})$. Hence, we can reduce to the case when $\scrY= \Spec(k)$, and so $\scrX \to \Spec(k)$ is a quasifinite proper stack universally homeomorphic to $\Spec(k)$. We need to show that $R\Gamma(\scrX, \overline{\mathbb{Q}}_{\ell}) = \overline{\mathbb{Q}}_{\ell}$. 
    
    Without loss of generality we can replace $\scrX$ with its reduced closed substack, and assume that $\scrX$ is reduced. We claim that  $\scrX=BG$ for some finite group scheme $G$ over $k$. Indeed, $\scrX$ is a gerbe over a reduced locally Noetherian algebraic space $Z$ over $k$ with $|Z|$ a singleton \cite[\href{https://stacks.math.columbia.edu/tag/06QK}{Tag 06QK}]{stacks-project}. Since $\scrX$ is of finite type over $k$, it has a $k$-point, and so does $Z$. The morphism $\Spec(k) \to Z$ is a monomorphism (because the composition $\Spec(k) \to Z \to \Spec(k)$ is a monomorphism), and therefore by \cite[\href{https://stacks.math.columbia.edu/tag/06QY}{Tag 06QY}]{stacks-project} and \cite[\href{https://stacks.math.columbia.edu/tag/06QW}{Tag 06QW}]{stacks-project}, we have that  $Z= \Spec(k)$. Hence $\scrX \to \Spec(k)$ is a gerbe equipped with a section $\Spec(k) \to \scrX$. This means that we have $\scrX = BG$ for some group scheme $G$ over $k$ \cite[\href{https://stacks.math.columbia.edu/tag/06QG}{Tag 06QG}]{stacks-project}. Since $\scrX$ is quasifinite, we know that $G$ is a finite group scheme over $k$.

So we have reduced to showing that, if $G$  is finite group scheme over $k$, then $R\Gamma(BG, \overline{\mathbb{Q}}_{\ell}) = \overline{\mathbb{Q}}_{\ell}$. Consider the natural morphisms $w: pt \stackrel{u}\to BG \stackrel{v}\to pt,$ where $pt := \Spec(k).$
 By a twice-repeated application of  \Cref{prop: decomposition theorem with generic connected fibers},   
 we have a splitting   $Ru_* (\overline{\mathbb{Q}}_{\ell})_{pt} \simeq (\overline{\mathbb{Q}}_{\ell})_{BG} 
 \oplus \mathcal{E},$ and thus a splitting 
 $(\overline{\mathbb{Q}}_{\ell})_{pt}
 \simeq Rw_* 
 (\overline{\mathbb{Q}}_{\ell})_{pt} 
 \simeq
 (Rv_* (\overline{\mathbb{Q}}_{\ell})_{BG}) 
 \oplus Rv_* \mathcal{E} 
 \simeq 
(\overline{\mathbb{Q}}_{\ell})_{pt} 
 \oplus \mathcal{E}'
\oplus Rv_* \mathcal{E},$ 
which implies the desired conclusion $R\Gamma(BG, 
\overline{\mathbb{Q}}_{\ell})) =\overline{\mathbb{Q}}_{\ell}.$ 
\end{proof}

 \begin{lemma}[Intersection cohomology of universal homeomorphisms] \label{lemma: intersection cohomology of 
tame moduli spaces}
Let $\pi: \scrX \to \scrY$ be a separated universal homemorphism  
between stacks with affine stabilizers and of finite type over an algebraically closed field $k$.
     Then there is a commutative square:
\[
\begin{tikzcd}\label{fig: diagram 2.2}
(\overline{\mathbb{Q}}_{\ell})_{\scrY}
\arrow{r}{\pi^*}[swap]{=} 
\arrow{d}{\epsilon} 
& 
R\pi_* (\overline{\mathbb{Q}}_{\ell})_{\scrX}
\arrow{d}{\epsilon}
\\ 
\mathcal{IC}_{\scrY}
\arrow{r}{[\pi^*]}[swap]{=}
& 
R\pi_* \mathcal{IC}_{\scrX},
\end{tikzcd}\]
where:  the canonical morphisms $\epsilon$ and $\pi^*$ are defined in Remark \ref{remark: basic proeprties of dt};
$\pi^*$ is a
canonical isomorphism (cf. \Cref{claim:homeo});
$[\pi^*]$ is a
canonical isomorphism. In particular: the induced morphisms on $H^0(\scrY,-)$ send  unit sections to the corresponding unit sections; the restriction of 
$[{\pi}^*]$ to a suitably small open and dense substack of $\scrY$ (e.g. the one sitting over the regular locus of its reduced substack) coincides with the corresponding restriction of $ \pi^*.$
 \end{lemma}

\begin{proof} 
%We use again that the hypotheses imply that $\pi$ is proper and quasifinite, and so in particular we have $R\pi_* = R\pi_!$. \Mirko{This is true but why we mention it here?} {\Md
%I am not sure. It seems out of place, perhaps  a remnant of days past. Remove?}
Without loss of generality we can assume that both $\scrX$ and $\scrY$ are reduced. Moreover, we can suppose that $\scrY =Y$ is a scheme, as the statement can be checked smooth locally on $\scrY$ by smooth base-change \cite[12.1]{laszlo-olsson-adic-sheaves} and by the compatibility of $\mathcal{IC}$ with passing to smooth covers.  Hence, we can assume that $\scrY =Y$ is a reduced scheme. Since $\pi$ is a universal homeomorphism, there is a dense open smooth subscheme $Y^o \subset Y$ with smooth preimage $\scrX^o$. \Cref{claim:homeo}  applied to the restriction $\pi: \scrX^o \to Y^o$ implies that that the (shifted) intermediate extension $\mathcal{IC}_{Y}$ of $(\overline{\mathbb{Q}}_{\ell})_{Y^o}= R\pi_*(\overline{\mathbb{Q}}_{\ell})_{\scrX^o}$ is the direct summand of $R\pi_* \mathcal{IC}_{\scrX}$ consisting of the simple components in the decomposition theorem with full support $Y.$ In order to prove \Cref{lemma: intersection cohomology of 
tame moduli spaces}, we need to show that there are no other summands, i.e.\ that $R\pi_*\mathcal{IC}_{\scrX}$ is a (shifted) intermediate extension of a lisse sheaf on a smooth open dense subset in $Y$. 
     
     Under these assumptions, by \cite[Thm. B]{rydh-noetherian-approximmation}, there is a finite and surjective schematic morphism $f: X \to \scr X$ with $X$ a scheme.
     Let $b: X \xrightarrow{f} \scrX \xrightarrow{\pi} Y$   be the composition.
     Since $\pi$ is quasifinite and proper, the same follows for 
     the morphism of schemes $b: X \to Y$, which must then be 
     finite. This means that $b$ is $t$-exact for the perverse 
     $t$-structure, and so it follows that $Rb_*\mathcal{IC}_{X}$ 
     is of the form $\mathcal{IC}_{Y}(\mathcal{L})$ for some 
     semisimple lisse sheaf on a smooth open dense subset in $Y$. 
     By applying \Cref{prop: decomposition theorem with generic 
     connected fibers} to the morphisms $f$ and $\pi$, we 
     conclude that $R\pi_*\mathcal{IC}_{\scrX}$ is a direct 
     summand of $Rb_*\mathcal{IC}_{X} = \mathcal{IC}_{Y}
     (\mathcal{L})$, and so in particular it is also of the form 
     $R\pi_*\mathcal{IC}_{\scrX} = \mathcal{IC}_{Y}
     (\mathcal{L}')$ for some direct summand $\mathcal{L}' 
     \subset \mathcal{L}$, as desired.

\end{proof}

\begin{example} \label{example: cohomology of coarse space}
Let $\scrX$ be a stack of finite type over $k$ with finite inertia. The coarse moduli space morphism $\pi: \scrX \to M$ is a separated universal homeomorphism by \cite[\href{https://stacks.math.columbia.edu/tag/0DUT}{Tag 0DUT}]{stacks-project}. Then \Cref{lemma: intersection cohomology of tame moduli spaces} implies that we have a 
canonical identification 
$[\pi^*]: \mathcal{IC}_{M}  \stackrel{=}\to  R\pi_* {\mathcal IC}_{\scrX}$.
In particular, this applies to the coarse moduli space morphism of any tame stack, i.e.\ stack locally modelled in the \'{e}tale topology
of the coarse moduli space $M$ on the quotient $U/G$ of a
scheme $U$ by the action of a finite flat linearly reductive group scheme; see \cite[\S 3]{abramovich2007tame}.
\end{example}
\end{subsection}
\begin{subsection}{Intersection cohomology as a direct summand of stack cohomology}

 In this section, we work with the bounded below derived categories $D^+_c(-,\overline{\mathbb{Q}}_{\ell})$; we cannot restrict ourselves to working with the bounded $D^b_c(-,\overline{\mathbb{Q}}_{\ell})$, because, in general, pushing forward under good moduli space morphisms does not preserve boundedness.

\begin{defn}[{\cite[Defn. 2.5]{edidin-rydh-canonical-reduction}}]
Let $\scrX$ be a finite type stack over $k$ admitting a good moduli space $\pi: \scrX \to M$. A geometric point $x$ of $\scrX$ is stable if we have $\pi^{-1}(\pi(x)) = \{x\}$ under the induced morphism of topological spaces $\pi: |\scrX| \to |M|$. Furthermore, a stable point $x$ is called properly stable if the dimension of the automorphism group scheme $\Aut(x)$ is zero.

We say that the moduli space $\scrX \to M$ is properly stable if the set of properly stable geometric points is dense.
\end{defn}
If a stack $\scrX$ admits a properly stable moduli space, then
there is a saturated open and dense substack of the source
whose stabilizers are zero-dimensional group schemes (cf. \cite[Prop. 2.6]{edidin-rydh-canonical-reduction}). The largest of these open sets is the locus of properly stable points relative to the good moduli space morphism, and it coincides with the locus of stable points.

We are now ready to prove the main result of this section.
\begin{thm}\label{thm: intersection cohomology of the stack vs the moduli space}
Let $X$ be an algebraic space of finite type over an algebraically closed field $k$, and let $G$ be an affine algebraic group. Suppose that the stack $[X/G]$ admits a good moduli space $M$, and the morphism $\pi: [X/G] \to M$ is properly stable. 
There is a commutative square in $D^+_c(M,\overline{\mathbb{Q}}_{\ell})$:
\[
\begin{tikzcd}\label{fig: diagram 2.3}
(\overline{\mathbb{Q}}_{\ell})_{M}
\arrow{r}{\pi^*} 
\arrow{d}{\epsilon} 
& 
R\pi_* (\overline{\mathbb{Q}}_{\ell})_{[X/G]}
\arrow{d}{\epsilon}
\\ 
\mathcal{IC}_{M}
\arrow[dotted]{r}{[\pi^*]} 
& 
R\pi_* \mathcal{IC}_{[X/G]},
\end{tikzcd}\]
where: $[\pi^*]$ is a noncanonical split morphism; the canonical morphisms $\epsilon$ and $\pi^*$ are defined in \Cref{remark: basic proeprties of dt}. 
%In particular: the induced morphisms on $H^0(M,-)$ send  unit sections to the corresponding unit sections; the restriction of 
%$[\pi^*]$ to suitably small open and dense subspaces of $M$ coincides with the corresponding restriction of $\pi^*.$
\end{thm}

\begin{remark}\label{rmk: ic and ic mhm}
Note that, by the assumption that the moduli space is properly stable, the morphism $\pi:[X/G] \to M$ has relative dimension zero. Therefore, we also have that the un-shifted  $IC_{M}$ is a direct summand of $R\pi_*(IC_{[X/G]})$ in $D^+_c(M,\overline{\mathbb{Q}}_{\ell})$.
 The analogous result holds, with an analogous proof, in the context of the bounded below derived category of algebraic mixed Hodge modules. A reference for the derived category of algebraic mixed Hodge modules on stacks of finite type over the field of complex numbers is \cite{tubach_mhm_paper}.
\end{remark}

Since the morphism $\pi$ is not proper in general, 
the Decomposition Theorem (\Cref{prop: decomposition theorem with generic connected fibers}) does not apply.
For our proof of \Cref{thm: intersection cohomology of the stack vs the moduli space}, we use a similar strategy as in \cite{woolf-intersection}.

We first note that a result stronger than \Cref{thm: intersection cohomology of the stack vs the moduli space} holds over a dense open subset of $M.$
By their very constructions, the forthcoming morphisms $\lambda$ and $\gamma$ coincide respectively  with the restrictions $\pi^*$ and $(\pi^*)^{-1}$ of the morphisms
$[\pi^*]$ and $([\pi^*])^{-1}$ 
to a suitable small open and dense 
subspace of $M.$  %By its  very construction, the forthcoming morphism $\lambda$ restricted to a suitable dense open subset of $M$ agrees with the restrictions of $\pi^*$ and of $[\pi^*]$ to this open set, where they are isomorphisms.  
The forthcoming morphism  $\gamma$ restricted to said open set  is the inverse of the restriction of $\lambda.$

\begin{lemma}[The morphisms $\lambda$ and $\gamma$ on a suitable open subspace of $M$] \label{lemma: intersection cohomology of moduli space generically}
Under the hypotheses of \Cref{thm: intersection cohomology of the stack vs the moduli space}, there is an open dense sub algebraic space $U \subset M$ such that the restriction $\pi|_{U}: \pi^{-1}(U) \to U$ induces the  natural identification 
$ [\pi^*]: 
\mathcal{IC}_U \xrightarrow{=} 
R\pi_*(\mathcal{IC}_{\pi^{-1}(U)})$
provided by \Cref{lemma: intersection cohomology of 
tame moduli spaces} applied to $\pi_{|U}.$
 In particular, the induced morphism on $H^0(M,-)$ sends the   unit section to the corresponding unit section,
 and the restriction of 
$[{\pi}^*]$ to  a suitably small open and dense subspace  of $M$ coincides with the corresponding restriction of $\pi^*.$
\end{lemma}
\begin{proof}
By the assumption that the moduli space morphism $\pi: [X/G] \to M$ is properly stable, in view of  \cite[Prop. 2.6]{edidin-rydh-canonical-reduction}, there exists an open subspace $U \subset M$ of the algebraic space $M$ such $\pi^{-1}(U)$ is a tame stack and $\pi^{-1}(U) \to U$ is its coarse space moduli morphism. The lemma then follows from \Cref{lemma: intersection cohomology of tame moduli spaces}
applied to this restricted morphism (cf. \Cref{example: cohomology of coarse space}). 
\end{proof}

The main result \Cref{thm: intersection cohomology of the stack vs the moduli space} 
follows from the following two lemmas.

\begin{lemma}[The morphism $\lambda$] \label{lemma: morphism lambda}
Under the hypotheses of \Cref{thm: intersection cohomology of the stack vs the moduli space} and the additional assumption that $X$ is integral, there exists a morphism $\lambda: \mathcal{IC}_M \to R\pi_*(\mathcal{IC}_{[X/G]})$
and a commutative square in $D^+_c(M,\overline{\mathbb{Q}}_{\ell})$:
\[
\begin{tikzcd}\label{fig: diagram 2.4}
(\overline{\mathbb{Q}}_{\ell})_{M}
\arrow{r}{\pi^*} 
\arrow{d}{\epsilon} 
& 
R\pi_* (\overline{\mathbb{Q}}_{\ell})_{[X/G]}
\arrow{d}{\epsilon}
\\ 
\mathcal{IC}_{M}
\arrow[dotted]{r}{\lambda} 
& 
R\pi_* \mathcal{IC}_{[X/G]}.
\end{tikzcd}\]
%where the canonical morphisms $\epsilon$ and $\pi^*$ are defined in Remark \ref{remark: basic proeprties of dt}. 
In particular: the induced morphisms on $H^0(M,-)$ send  unit sections to the corresponding unit sections; the restriction of 
$\lambda$ to suitably small open subspace of $M$ coincides with the corresponding restriction of $\pi^*.$
\end{lemma}

The only difference between \Cref{thm: intersection cohomology of the stack vs the moduli space} and \Cref{lemma: morphism lambda} is that in \Cref{lemma: morphism lambda} $\lambda$ is not yet known to be a split injection. To this end, we need the following lemma.

\begin{lemma}[The morphism $\gamma$] \label{lemma: morphism gamma}
Under the hypotheses of \Cref{thm: intersection cohomology of the stack vs the moduli space} and the additional assumption that $X$ is integral, there exists a morphism $\gamma: R\pi_*(\mathcal{IC}_{[X/G]}) \to \mathcal{IC}_M$ in $D^+_c(M,\overline{\mathbb{Q}}_{\ell})$ such that
it sends the unit section to the unit section and  
the restriction of 
$\gamma$ to a suitably small open subspace of $M$ coincides with the corresponding restriction of $(\pi^*)^{-1}$
(cf. \Cref{lemma: intersection cohomology of moduli space generically}).
\end{lemma}

We can now present the proof of the main theorem of the section.
\begin{proof}[Proof of \Cref{thm: intersection cohomology of the stack vs the moduli space}]
Since the intersection cohomology complex  of a scheme (or stack) is the direct sum of 
the intersection cohomology complexes of the reduced irreducible components, we can assume, for the 
sake of proving \Cref{thm: intersection cohomology of the stack vs the moduli space}, 
that $X$ is an integral algebraic space (note that under our properly stable assumption, \cite[Lem. 4.14]{alper-good-moduli} implies that the moduli spaces of the reduced irreducible components of $[X/G]$ are the corresponding reduced irreducible components of $M$).

Let $\lambda$ and $\gamma$ be as in \Cref{lemma: morphism 
lambda} and \Cref{lemma: morphism gamma}. 
We need to show  that $\gamma \circ \lambda: 
\mathcal{IC}_M \to \mathcal{IC}_M$ is the identity morphism. Since the intersection cohomology complex $\mathcal{IC}_M$ 
of the integral algebraic space $M$ is simple 
\cite[Thm. 8.2]{laszlo-olsson-perverse-sheaves}, it suffices to check that 
$\gamma \circ \lambda = {\rm Id}$ over an open  and dense sub algebraic space $U \subset M.$ To see this, we 
use the identifications  
 $\lambda|_{U} = \pi^*|_{U} = (\gamma|_{U})^{-1}$ for a suitably small open dense subspace $U\subseteq M$ as in Lemmas \ref{lemma: morphism lambda} and 
\ref{lemma: morphism gamma}.
\end{proof}

We are left with the task of constructing the morphisms $\lambda$ and $\gamma$ as in Lemmas \ref{lemma: 
morphism lambda} and \ref{lemma: morphism gamma}.
\end{subsection}

\begin{subsection}{Construction of the morphism \texorpdfstring{$\lambda$ in \Cref{lemma: 
morphism lambda} }{lambda}}\label{subsect: construction lambda}

In this section, we work with the bounded below derived category $D^+_c(-)$. We start with the following two general facts, which we intend to apply in the special case when
the forthcoming morphism $r^\#$  is the morphism $\pi$ of \Cref{lemma: morphism lambda}.

\begin{lemma}\label{lemma BK}
Let $r^\#: \scrX^\# \to \scrM$ be a morphism of integral stacks  that factors 
as in the following diagram of morphisms of integral stacks:

\begin{figure}[H]
\centering
\begin{tikzcd}
  & \scrX^\# 
  \arrow[hook]{d}[swap]{j}
  \arrow{rd}[description]{p^\#}
  \arrow{rrd}{r^\#}
  &   &
  \\ 
  r: \hspace{-1.8 cm} & \scrX 
  \arrow{r}{p}
   & \mathscr{B}
   \arrow{r}{q}
   & \scrM,
\end{tikzcd}
\end{figure}
\noindent where $j$ is an open immersion with dense image; $p$ is a proper surjetive  morphism of stacks with affine stabilizers; and $q$ is smooth.

Then there is a commutative diagram: 

\begin{tikzcd}\label{fig: diagram 10}
(r^\#)^*: (\overline{\mathbb{Q}}_{\ell})_{\scrM} 
\arrow{d}{\epsilon}
\arrow{r}{q^*}
&
Rq_* q^* (\overline{\mathbb{Q}}_{\ell})_{\scrM} 
= Rq_* (\overline{\mathbb{Q}}_{\ell})_{\mathscr{B}}
\arrow{d}{\epsilon}
\arrow{r}{p^*}
& 
Rq_* Rp_* (\overline{\mathbb{Q}}_{\ell})_{\scrX}
\arrow{d}{\epsilon} 
\arrow{r}{j^*}
&
R r^\#_*
 (\overline{\mathbb{Q}}_{\ell})_{{\scrX}^\#}
\arrow{d}{\epsilon}
\\ 
{[(r^\#)^*} ]:
\mathcal{IC}_{\scrM}
\arrow{r}{q^*} 
&
Rq_* q^* \mathcal{IC}_{\scrM}
= Rq_* \mathcal{IC}_{\mathscr{B}}
\arrow[dotted]{r}{ [\widehat{p}^*]}
& 
Rq_* Rp_*  \mathcal{IC}_{\scrX}
\arrow{r}{j^*}
&
R r^\#_* \mathcal{IC}_{{\scrX}^\#},
\end{tikzcd}

\noindent
where:
all the solid arrows are canonical; the dotted arrow is a noncanonical split injection obtained as the $Rq_*$ pushforward of any morphism stemming from \Cref{prop: decomposition theorem with generic connected fibers}. In particular, the unit section 
$1_{\scrM}$
 is sent to the corresponding  unit sections as one moves along the 
 diagram obtained by taking $H^0 (\scrM, -).$ If in 
 addition $q$ is a universal homeomorphism, then the restriction 
 of the composition $ j^* \circ [\widehat{p}^*] \circ q^*$ to a suitable 
 small open and dense substack of 
 $\scrM$ coincides with the corresponding restriction of 
 $(r^{\#})^*= j^* \circ p^* \circ q^*.$
\end{lemma}

\begin{proof}
Let $[\widehat{p}^*]: \mathcal{IC}_{\mathscr{B}} \to Rp_* 
\mathcal{IC}_{\scrX}$ be a morphism as in \Cref{prop: 
decomposition theorem with generic connected fibers}.
By pushing forward via $Rq_*,$
pre-composing with $q^*,$ post-composing with $j^*,$
and considering the morphisms of type $\epsilon,$ we obtain the desired 
commutative diagram.  The equality $q^*\mathcal{IC}_{\scrM} = \mathcal{IC}_{\mathscr{B}}$ in the second row follows from the smoothness of $q.$ By construction, unit sections are mapped by unit sections.

If in addition $q$ is a
universal homemorphism, then, by using the notation in the proof
of \Cref{prop: 
decomposition theorem with generic connected fibers}, where the $\mathscr{B}$ here is the $\scrY$ there,
we have  the open and dense $\mathscr{B}^o$, which is automatically $q$-saturated, i.e.  of the form $q^{-1}(\scrM^o)$ with 
$\scrM^o:= q(\mathscr{B}^o).$
The lemma follows by taking $\scrM^o \subseteq \scrM$ to be the required open and dense substack.
\end{proof}

\begin{lemma}[Relative $G$-equivariant completion]\label{lemma: G-completion}
Let $f: X \to [X/G] \to M$
be as in \Cref{thm: intersection cohomology of the stack vs the moduli space} and endow $M$ with the trivial $G$ action.
The $G$-equivariant morphism $f: X  \to M$
admits a factorization $X \stackrel{j}\to Z \stackrel{p}\to M$
into $G$-equivariant morphisms of $G$-algebraic spaces  such that $j$ is an open immersion with dense image and $p$ is proper schematic.
\end{lemma}

\begin{proof}
\iffalse
We let $G$ act trivially on $M$ and we  construct a relative $G$-equivariant compactification for $X \to M.$  By this we mean a 
diagram of $G$-equivariant morphisms of  algebraic spaces with $G$-actions

\begin{figure}[H]
\centering
\begin{tikzcd}
  X \ar[r, symbol= \hookrightarrow, "j"] \ar[dr,"f", labels= below] & Z \ar[d, "p"] \\  & M,
\end{tikzcd}
\end{figure}
{\cvio spacing ???}where $j$ is an open immersion and $\pi$ is proper.
\fi
Since $\pi$ is a good moduli space morphism and $G$ is affine,
the compositum morphism $f: X \to [X/G] \to M,$ being a composition of cohomologically affine morphisms, it  is cohomologically affine, hence affine \cite[Prop. 3.3.]{alper-good-moduli}. 

Recall that we can identify $G$-equivariant (quasi-)coherent sheaves on 
$M$ with (quasi-)coherent sheaves on the Noetherian algebraic stack $[M/G] = M \times BG$, by fppf descent and a standard argument (cf. \cite[Exer.9.H, pg.207]{olsson-algebraic-stacks} when $G$ is smooth).
The pushforward $f_*\mathcal{O}_X$ is a $G$-equivariant 
quasicoherent sheaf on $M$. 
By \cite[Prop. 15.4]{lmb-champ}, and the aforementioned equivalence, 
we can write $f_*\mathcal{O}_X$ as a filtered colimit $f_*\mathcal{O}_X = \colim_{i} 
\mathcal{F}_i$, of 
$G$-equivariant coherent subsheaves  $\mathcal{F}_i 
\subset f_*\mathcal{O}_X$.
Since the algebraic space $M$ is quasicompact, there exists some affine \'etale atlas $U \to 
M$. Since $X$ is of finite type over $M$, for the base-change $f_U: X_U \to U$ there is some index 
$k \gg 0$ such that the natural morphism of $\mathcal{O}_U$-algebras $(\Sym_{\mathcal{O}_M}
(\mathcal{F}_i))|_{U} \to (f_U) _*(\mathcal{O}_{X_U})$ is surjective for every $i \geq k$. It follows by \'etale descent that 
$\Sym_{\mathcal{O}_{M}}(\mathcal{F}_i) \to f_*\mathcal{O}_{X}$ is surjective for every $i \geq k$. 

Fix  $i \geq k$. 
We conclude 
that the inclusion $\mathcal{F}_i \to f_*\mathcal{O}_X$ induces a closed immersion of 
$M$-algebraic spaces $X \hookrightarrow \uSpec_M(\Sym_{\mathcal{O}_{M}}(\mathcal{F}_i))$. We also have an open immersion $\uSpec_M(\Sym_{\mathcal{O}_{M}}(\mathcal{F}_i)) 
\hookrightarrow \text{Proj}_M((\Sym_{\mathcal{O}_{M}}(\mathcal{F}_i))[t]),$  where the 
grading on the algebra  $(\Sym_{\mathcal{O}_{M}}(\mathcal{F}_i))[t]$ is prescribed by 
assigning  degree $1$ to $t$ and $\mathcal{F}_i.$ By composition, we obtain a locally closed 
immersion $X \hookrightarrow \text{Proj}_M(((\Sym_{\mathcal{O}_{M}}(\mathcal{F}_i))[t]))$. 
The target $\text{Proj}_M (\Sym_{\mathcal{O}_{M}}(\mathcal{F}_i))[t]$ acquires a 
$G$-action, induced by the $G$-equivariant structure on $\mathcal{F}_i$ (we set $t$ to be $G$-invariant).  We thus obtain a $G$-equivariant diagram of algebraic spaces:
\begin{equation}\label{fig: diagram 1}
\begin{tikzcd}[ampersand replacement=\&] 
  X \ar[r, symbol= \hookrightarrow, "J"] \ar[dr,"f", labels= below] \& \text{Proj}_M ((\Sym_{\mathcal{O}_{M}}(\mathcal{F}_i))[t]) \ar[d, "P"] \\  \& M,
\end{tikzcd}
\end{equation}
\noindent 
where $J$ is a locally closed immersion and $P$ is schematic
and proper. We conclude the construction of our compactification by setting $Z$ to be the scheme-theoretic closure  of $X \to \text{Proj}_M ((\Sym_{\mathcal{O}_{M}}(\mathcal{F}_i))[t])$, equipped with its induced $G$-action, and by taking $j$ and $p$ to be the resulting restrictions of $J$ and $P.$
\end{proof}

\begin{proof}[Proof of \Cref{lemma: morphism lambda}]
We let $G$ act trivially on $M.$ 
By \Cref{lemma: G-completion},  after having passed to quotient-by-$G$ stacks, we obtain a factorization  
$\pi:[X/G] \stackrel{j}\to [Z/G] \stackrel{p}\to [M/G]=M\times BG \stackrel{\text{pr}_M}\to M$ of the morphism $\pi,$
with $j$ an open immersion, $p$ proper surjective between stacks with affine stabilizers,
and $\text{pr}_M$
smooth.

\Cref{lemma: morphism lambda} follows by applying the general
\Cref{lemma BK} to this situation and by  setting $\lambda:= [(r^\#)^*].$
\end{proof}

\end{subsection}

\begin{subsection}{Construction of the morphism \texorpdfstring{$\gamma$}{gamma}}

In this section, we work  with the bounded below derived
categories $D^+_c(-,\overline{\mathbb{Q}}_{\ell})$.

In order to construct the morphism $\gamma$, we will use the work of Edidin and Rydh \cite{edidin-rydh-canonical-reduction} generalizing the theory of Kirwan blowups to the setting of stacks with good moduli spaces. We start by stating a simplified less general version of the main result in \cite{edidin-rydh-canonical-reduction}.

\begin{remark}
    For the following result, it may help the reader to recall that in a nonsingular and  scheme-theoretic GIT set-up,
the morphisms $\psi_i$ are the blow-ups of the loci where the stabilizers have maximal dimensions; the open immersions detect the semistable points for a suitable linearization; the morphisms
$\pi_i$ are actual GIT quotients; the open sets $M_i^o$ are the intersection of the
smooth part with the image via $\pi_i$ of the 
stable part.
\end{remark}

\begin{prop}[{\cite[Thm. 2.11]{edidin-rydh-canonical-reduction}}] \label{thm: edidin-rydh theorem} Let $\mathfrak{X}$ be an integral algebraic stack with affine stabilizers and of finite type over a field $k$. Suppose that $\mathfrak{X}$ admits a good moduli space $ \pi: \mathfrak{X} \to M$ such that $\pi$ is properly stable. Then, there is a finite sequence of integral stacks $\mathfrak{X}_0, \mathfrak{X}_1, \ldots, \mathfrak{X}_n$, and open substacks $\mathfrak{U}_i \subset \mathfrak{X}_i$ for all $i$, and proper schematic birational morphisms $\psi_i: \mathfrak{X}_i \to \mathfrak{U}_{i-1}$ satisfying the following properties:
\begin{enumerate}[(a)]
\item $\mathfrak{X}_0 = \mathfrak{U}_0 = \mathfrak{X}$.
\item For all $i$, the stack $\mathfrak{U}_i$ admits a good moduli space $\pi_i: \mathfrak{U}_i \to M_i$ with $\pi_i$ properly stable.
\item For all $i$, the morphism $\psi_i: \mathfrak{U}_i \to \mathfrak{U}_{i-1}$ is an isomorphism over the open locus of $\mathfrak{U}_{i-1}$ where $\pi_{i-1}$ is properly stable.
\item For all $i$, the morphism on moduli spaces $\xi_i: M_i \to M_{i-1}$ induced by $\psi_i: \mathfrak{U}_i \to \mathfrak{U}_{i-1}$ is proper and birational.
\item The last moduli space morphism $\pi_n: \mathfrak{U}_n \to M_n$ is a proper tame moduli space morphism.
\end{enumerate}
\end{prop}

Using this result, we shall construct our desired morphism $\gamma$.

\begin{proof}[Proof of \Cref{lemma: morphism gamma}]
We can apply \Cref{thm: edidin-rydh theorem} to the integral stack $[X/G]$ with properly stable good moduli space $\pi: [X/G] \to M$. This way we obtain a sequence of good moduli space morphisms $\pi_n: \mathfrak{U}_n \to M_n$. We shall show that \Cref{lemma: morphism gamma} holds for all $\pi_i: \mathfrak{U}_i \to M_i$ by descending induction on $i.$ In particular, this will show it for $\pi_0 =\pi: \scrX \to M$, as desired. The base case $\pi_n: \mathfrak{U}_n \to M_n$ is a coarse moduli space of a tame stack, and so we can conclude \Cref{lemma: morphism gamma} directly by using \Cref{lemma: intersection cohomology of tame moduli spaces}.

For the induction step, we have a diagram:
\begin{figure}[H]
\centering
\begin{tikzcd}
  \scrX_i \ar[d, "\psi_i"] \ar[r, symbol = \hookleftarrow, "j"]  & \mathfrak{U}_i  \ar[r,"\pi_i"] & M_i \ar[d, "\xi_i"] \\  
  \mathfrak{U}_{i-1} \ar[rr, "\pi_{i-1}"] & & M_{i-1},
\end{tikzcd}
\end{figure}
\noindent
where $j$ is an open immersion. By assumption, we have a morphism $\gamma_i: R(\pi_i)_*IC_{\mathfrak{U}_i} \to IC_{M_i}$ that restricts to the natural identification $\varphi_i$ over an open subspace of the properly stable locus of $M_i$.

By using \Cref{prop: decomposition theorem with generic connected fibers} for the proper birational morphism $\psi_{i} \colon \scrX_i \to \mathfrak{U}_{i-1}$, we get that $IC_{\mathfrak{U}_{i-1}}$ is a direct summand of $R(\psi_i)_* IC_{\scrX_i}$. In particular we get a morphism $IC_{\mathfrak{U}_{i-1}} \to R(\psi_i)_* IC_{\scrX_i}$, which yields the following morphism after applying $R(\pi_{i-1})_*(-)$
\[ R(\pi_{i-1})_*IC_{\mathfrak{U}_{i-1}} \to R(\pi_{i-1} \circ \psi_i)_* IC_{\scrX_i}\]
Composing with the morphism $R(\pi_{i-1} \circ \psi_i)_* IC_{\scrX_i} \to R(\pi_{i-1} \circ \psi_i)_* Rj_* j^* IC_{\scrX_i}$ induced by the counit of the open immersion $j$, we get a morphism
\[ a: R(\pi_{i-1})_*IC_{\mathfrak{U}_{i-1}} \to R(\pi_{i-1} \circ \psi_i)_* IC_{\scrX_i} \to R(\pi_{i-1} \circ \psi_i)_* Rj_*IC_{\mathfrak{U}_i} = R(\xi_i)_*R(\pi_i)_* IC_{\mathfrak{U}_i}\]
By composing with the morphism $R(\xi_i)_*\gamma_i: R(\xi_i)_* R(\pi_i)_*IC_{\mathfrak{U}_i} \to R(\xi_i)_*IC_{M_i}$, we get
\[ R(\pi_{i-1})_*IC_{\mathfrak{U}_{i-1}} \xrightarrow{a} R(\xi_i)_*R(\pi_i)_* IC_{\mathfrak{U}_i} \xrightarrow{R(\xi_i)_*\gamma_i}R(\xi_i)_*IC_{M_i} \]
Finally, by \Cref{prop: decomposition theorem with generic connected fibers} applied to the proper birational morphism $\xi_i: M_i \to M_{i-1}$, we have that $IC_{M_{i-1}}$ is a direct 
summand of $R(\xi_i)_*IC_{M_i}$. In particular, there is a 
projection morphism $R(\xi_i)_*IC_{M_i} \to IC_{M_{i-1}}$. 
Composing with our previous map we obtain our desired morphism
\[\gamma_{i-1}:  R(\pi_{i-1})_*IC_{\mathfrak{U}_{i-1}} 
\xrightarrow{R(\xi_i)_*\gamma_i \circ a} R(\xi_i)_*IC_{M_i} \to 
IC_{M_{i-1}}\]
If we restrict over the open locus where $\pi_{i-1}: 
\mathfrak{U}_{i-1} \to M_{i-1}$ is properly stable, then the restrictions of $\psi_{i}$, $j$ and $\xi_i$ are isomorphisms. It follows by construction that then the restriction of the morphism $\gamma_{i-1}$ is identified with the restriction of the morphism $\gamma_i$. Since by the induction assumption $\gamma_i$ restricts to  the natural identification over an open subspace of the properly stable locus of $M_i$, similarly $\gamma_{i-1}$ satisfies this property, thus concluding the induction step.
\end{proof}
\end{subsection}

\begin{subsection}{Cohomological bound for the Hitchin morphism: the case with {\texorpdfstring{$\deg(\cL)>2g-2$}{deg(L)>2g-2}}}

In this section, we work with $D_c(-)$ without bounding conditions. However, given that we are working with morphisms $\pi$ of finite type, we use the facts that $R\pi_*$ preserves
bounded below complexes, $R\pi_!$ preserves bounded above complexes 
and Verdier duality exchanges bounded below with bounded above and viceversa.

The results developed in this section allow us to prove the necessary cohomological bound in \Cref{thm: thmA} in the case when $\deg(\cL) > 2g-2$. Let us first state more general results for weak Abelian fibrations.
\begin{prop}[Cohomological bound] \label{cohomological bound} Let $(P,h\colon M \to B)$ be a weak Abelian fibration as in \Cref{defn: weak Abelian fibration}. Suppose that $M$ is a properly stable good moduli space of a quotient stack $\scrM$ which is smooth over $k$ and such that every fiber of $\scrM \to M \to B$ is equidimensional of dimension $\dim(h) = \dim(M) - \dim(B)$. Then we have \[\tau_{> 2\dim(h)} R h_*(IC_M[-\dim(M)]) =0.\]
\end{prop}
\begin{proof}
    Our proof follows \cite[Proof of Prop.\ 3.2 and 3.3]{maulik-shen-independence}. Consider the composition $\scrM \xrightarrow{\pi} M \xrightarrow{h} B$. By the splitness of the morphism $[\pi^*]$ in \Cref{thm: intersection cohomology of the stack vs the moduli space}, we have a direct sum decomposition $R\pi_*(IC_{\scrM}) \simeq IC_M \oplus \cE$. By applying Verdier duality, we get
\[ R\pi_!(IC_\scrM) = \mathbb{D} R\pi_*(IC_{\scrM}) \simeq IC_M \oplus \mathbb{D}(\cE).\]
We may now apply $Rh_!$ to the previous equality,  use $R(h \circ \pi)_! = Rh_! \circ R\pi_!$ (functoriality) and $Rh_!=Rh_*$ (properness of $h$) to get
\[ R(h \circ \pi)_! (IC_{\scrM}) \simeq Rh_!(IC_M) \oplus Rh_!(\mathbb{D}(\cE)) =  Rh_*(IC_M) \oplus  Rh_*(\mathbb{D}(\cE)). \]
 We conclude that, in order to show $\tau_{> 2\dim(h)} Rh_*(IC_M[-\dim(M)]) =0$, it suffices to show that $\tau_{> 2\dim(h)} R(h \circ \pi)_!(IC_\scrM[-\dim(M)])  =0$. 

Since $\scrM$ is smooth, we have $IC_{\scrM} = (\overline{\mathbb{Q}}_{\ell})_{\scrM}[\dim(\scrM)]$. On the other hand, the condition about $M$ being a properly stable good moduli space guarantees that $\dim(\scrM) = \dim(M)$. Therefore we have $IC_\scrM[-\dim(M)] = (\overline{\mathbb{Q}}_{\ell})_{\scrM}$, and we need to show that $\tau_{> 2\dim(h)} R(h \circ \pi)_!((\overline{\mathbb{Q}}_{\ell})_{\scrM}) =0$. Since the fibers of $h \circ \pi$ have dimension $\dim(h)$ by assumption, this follows from the proof of \cite[Prop. 3.3]{maulik-shen-independence}. 
\end{proof}

% \begin{prop} \label{prop: cohomological bound general}
%     Let $(P,h:M \to B)$ be a weak Abelian fibration as in \Cref{defn: weak Abelian fibration}. Suppose that $M$ is a properly stable good moduli space of a quotient stack $\scrM$ which is smooth over $k$ and such that every fiber of $\scrM \to M \to B$ is equidimensional of dimension $\dim(h) = \dim(M) - \dim(B)$. Then, any support $Z\subset B$ of the Decomposition Theorem of $h$ satisfies $\codim_B(Z) \leq \delta(Z)$, where $\delta(Z)$ is defined using the group scheme $P \to B$ as in \Cref{defn: delta invariant}.
% \end{prop}
% \begin{proof}
% In view of \cite[Thm. 1.1]{maulik-shen-independence}, this is a consequence of \Cref{cohomological bound}.
% \end{proof}

We now return to our main example of the Hitchin fibration $h: \MHiggs_{G, \cL}^{d} \to A$ (\Cref{notn: higgs moduli space and hitchin fibration}). Recall that under the assumption $\deg(\cL)>\max\{2g-2,0\}$ and $\text{char}(k) =0$, we have shown that the action of the group scheme $\schP^o$ (\Cref{P:centr-group}(ii)) on $\MHiggs^d_{G,\cL}$ makes the triple $(\schP^o, \MHiggs_{G,\cL}^d,A)$ into a weak Abelian fibration (\Cref{thm: weak Abelian fibration structure}).
\begin{thm} \label{thm: ngo inequality supports}
    Suppose that $\text{char}(k)=0$ and that $\deg(\cL) >\max\{2g-2,0\}$. For any given $d \in \pi_1(G)$, the Hitchin fibration $h:\MHiggs_{G,\cL}^{d} \to A$ satisfies
    \[\tau_{> 2\dim(h)} Rh_*(IC_{\MHiggs_{G,\cL}^{d}}[-\dim(\MHiggs_{G,\cL}^{d})]) =0.\]
    %for all $n>2\dim(h)$.
\end{thm}
\begin{proof}
    By \Cref{cohomological bound}, it suffices to show that $\MHiggs_{G, \cL}^{d}$ is a properly stable good moduli space $\pi: \scrM \to M$, where $\scrM$ is a smooth quotient stack such that the fibers of $(h \circ \pi): \scrM \to A$ are equidimensional of dimension $\dim(h)$. We set $\scrM$ to be the rigidification $\Higgs_{G,\cL}^{d,ss}\fatslash Z_G$. The smoothness of $\scrM$ follows from \Cref{prop: smoothness of the semistable hitchin stack}, and the fact that the fibers of $h \circ \pi$ have the expected dimension follows from combining \Cref{P:Higgs-prop} with the formula for $\dim(h)$ provided by \Cref{cor: higgs moduli space properties}. The fact that $\scrM \to M$ is a good moduli space is an immediate consequence of the composition $\Higgs_{G,\cL}^{d,ss} \to \scrM \to M$ being a good moduli space, where $\Higgs_{G,\cL}^{d,ss} \to \scrM$ is a $Z_G$-gerbe by construction. Furthermore, by combining \Cref{lemma: stable higgs is saturated} with \Cref{lemma: dimension automorphisms stable higgs}, it follows that all geometric points of $\scrM$ in the image of the stable locus of $\Higgs_{G,\cL}^{d,ss}$ are properly stable. Since $\scrM$ is irreducible (\Cref{P:Higgs-integral}) and the stable locus is nonempty (by the proof of \Cref{cor: higgs moduli space properties}), it follows that $\scrM \to M$ is a properly stable good moduli space. We are only left to show that $\scrM$ is a quotient stack. The proof of \cite[Prop. 5.4]{herrero_automorphism} implies that $\Higgs_{G,\cL}^{d,ss}$ is a quotient stack of the form $[Y/G]$, where $Y$ is a quasiprojective scheme classifying Higgs bundles with a framing at a fixed $k$-point of $C$. It follows that the rigidification $\scrM$ is isomorphic to the quotient stack $[Y/(G/Z_G)]$.
\end{proof}
\end{subsection}

\begin{subsection}{Cohomological bound for the Hitchin morphism: canonical case}

 The proof of the cohomological bound for meromorphic Higgs bundles \Cref{thm: ngo inequality supports} relies on the smoothness of the stack $\Higgs_{G,\cL}^{d,ss}\fatslash Z_G$, where $\deg(\cL) >\max\{2g-2,0\}$. When $\cL=\omega_C$, the stack $\Higgs_{G,\cL}^{d,ss}\fatslash Z_G$ is singular, but the conclusion of \Cref{thm: ngo inequality supports} still holds, because the Hitchin fibration $h: \MHiggs_{G, \cL}^{d} \to A$ is Lagrangian; see \Cref{cor:relativedimcanonical}.

 \begin{defn}\label{defn:symplecticvariety}\label{def: lagrangian fibration}
     Suppose that $\text{char}(k)=0$. 
     
     \begin{enumerate}[(1)]
         \item A \emph{symplectic variety} $(X, \sigma)$ is a normal variety $X$  with rational singularities whose regular part admits a symplectic form $\sigma$, i.e. a nondegenerate closed 2-form, or equivalently, a normal variety
         whose regular locus admits a symplectic form that extends to a regular 2-form on some (hence any) resolution. 
\item Let $(X, \sigma)$ be a symplectic variety. Then a subvariety $F \subset X$, with $F \not \subset X^{\mathrm{sing}}$, is \emph{Lagrangian} (with respect to the symplectic form $\sigma$) if, for a general closed point $x \in F$, the tangent space $T_{x} F$ is a Lagrangian subspace of $T_{x} X$ with respect to $\sigma$, or equivalently $\sigma |_{X^{\mathrm{reg}} \cap F^{\mathrm{reg}}}=0$ and $\dim F = \frac{1}{2}\dim X$.
\item A \emph{Lagrangian fibration} $f:X \to A$ is a proper surjective morphism with connected fibers between normal varieties, with $X$ symplectic, such that the general fiber of $f$ (or equivalently any irreducible component of any closed fiber) is a Lagrangian subvariety. This implies that $f$ is equidimensional by \cite{Mat1999}.
     \end{enumerate}
 \end{defn}

\begin{prop}\label{prop:LagrangianfibrationHitchin}
 Suppose that $\text{char}(k)=0$ and that $\cL=\omega_C$. For any given $d \in \pi_1(G)$, $\MHiggs_{G,\cL}^{d}$ is a symplectic variety, and the Hitchin fibration $h:\MHiggs_{G,\cL}^{d} \to A$ is a Lagrangian fibration.
\end{prop}

\begin{proof}
The moduli space $\MHiggs_{G,\cL}^{d}$ is a normal variety since the stack $\Higgs_{G,\cL}^{ss, d}$ is normal by \Cref{P:Higgs-normal}. The regular locus of $\MHiggs_{G,\cL}^{d}$ admits a symplectic form extending the standard symplectic form on the open set $T^* N_G^{rs,d} \subseteq \MHiggs_{G,\cL}^{d}$, where $N_G^{rs,d}$ is the moduli space of regularly stable $G$-bundles of degree $d$, see \cite[Thm.\ 4.3]{biswas-ramanan-infinitesimal}. %It remains to prove that 

We show that $M_{G,\cL}^d$ has rational singularities. We divide the proof in several cases depending on the type of derived subgroup $D_G:= [G,G]$ and on the genus $g$. 
\begin{enumerate}[(a)] 
    \item\label{i:higgs rational 2} If $G$ does not admit any quotient onto $\mathrm{PGL}_2$ or $g>2$, then $M_{G,\cL}^d$ has rational singularities.
    \item\label{i:higgs rational 3} If $D_G$ is either $\mathrm{SL}_2$ or $\mathrm{PGL}_2$ and $g=2$, then $M_{G,\cL}^d$ has rational singularities.
    \item\label{i:higgs rational 4} If $D_G$ is simply-connected, then any geometric fiber of the isotrivial (cf. \Cref{L:twisted locally trivial}) morphism $\nu:\MHiggs_{G,\cL}^d\to \MHiggs_{G/D_G,\cL}$ has rational singularities.
    \item\label{i:higgs rational 5} $M_{G,\cL}^d$ has rational singularities for any $G$ and $g$.
\end{enumerate}

% Case \eqref{i:higgs rational 1}. \sout{After fixing a splitting $G\cong\mathbb G_m^{\dim(G)}$, we see that}\Mirko{This have been fixed at the very beginning no?} $\MHiggs_{G,\cL}^d$ is isomorphic to the $\dim(G)$ copies of the cotangent bundle of the Jacobian of the base curve. Hence, the assertion follows.

 Case \eqref{i:higgs rational 2}. %$G$ does not admit any quotient onto $\mathrm{PGL}_2$ or $g>2$. %\color{black} %\color{blue} $D_G$ is almost-simple, different from $SL_2$ or $PGL_2$, or $g>2$. \color{black}
%\footnote{\Mirko{I added almost simple, then I changed again to admit a quotient to $PGL_2$, this is the greatest generality where the codim 4 condition is known and maybe it is better to mention it}} %$D_G\neq SL_2,PGL_2$ or $g>2$. 
By \cite[Thm.\ II.6]{faltings-projective-conn}, the singular locus of $\MHiggs_{G,\cL}^{d}$ has codimension at least $4$. Therefore, the symplectic form on the regular locus of $\MHiggs_{G,\cL}^d$ 
extends to a regular 2-form on some (hence any) resolution by \cite{Flenner1988}. In particular, $\MHiggs_{G,\cL}^d$ has canonical (hence rational) singularities; see \cite[Prop.\ 1.2.(I)]{Reid1980} and \cite[Cor.\ 5.24]{KM1998}. %\color{black} \footnote{\Mirko{The symplectic form does not extend to a volume form. I changed a bit the phrase}}

%By a result of Flenner \cite{Flenner1988}, the symplectic form on the regular locus of $\MHiggs_{G,\cL}^d$ extends to a regular volume form on any resolution, so $\MHiggs_{G,\cL}^d$ has canonical (hence rational) singularities. 

Case \eqref{i:higgs rational 3}. We have that $G\cong H\times T'$, where $T'$ is torus and $H$ is isomorphic to $\GL_2$, $\mathrm{SL}_2$ or $\mathrm{PGL}_2$; see for instance \cite[Thm.\ 20.33]{milne-algebraic-groups}. Since %the induced morphism 
$M_{G,\cL}\cong M_{H,\cL}\times M_{T',\cL}$, % is an isomorphism, 
it suffices to show the statement for each factor. The case $T'$ follows from Case \eqref{i:higgs rational 2} (actually $M_{T',\cL}$ is even smooth), while the case $H=\GL_2,\mathrm{SL}_2,\mathrm{PGL}_2$ follows by direct inspection; see for instance \cite[Cor.\ 8.4]{CB2003}, \cite[Cor.\ 1.2]{BS19} and \cite[Thm.\ 1.10]{Budur2021}. %\cite[Thm.\ 2.10 and 2.11]{Mauri2023}.  

Case \eqref{i:higgs rational 4}. %Fix a $k$-point $(E,\theta)$ in $\MHiggs_{G/D_G,\cL}$. 
Following \Cref{N:twisted higgs}, we denote by $M_{G,\cL}^{(E,\theta)} \coloneqq \nu^{-1}((E,\theta))$ the fiber along $\nu$ of a $k$-point $(E,\theta)$ in $\MHiggs_{G/D_G,\cL}$. By \Cref{L:twisted split}, we may assume that $D_G$ is almost-simple and non-trivial (this reduction uses the simply-connectedness of $D_G$). By \Cref{L:twisted locally trivial}, there exists an \'etale and surjective morphism of varieties \[\MHiggs_{Z_G^0}^0\times M_{G,\cL}^{(E,\theta)}\to \MHiggs_{G,\cL}^d.\] 
Now, $\MHiggs_{Z_G^0}^0$ and $\MHiggs_{G,\cL}^d$ have rational singularities by Cases \eqref{i:higgs rational 2} or \eqref{i:higgs rational 3}. Therefore, $M_{G,\cL}^{(E,\theta)}$ has rational singularities too. 
% Since $Z_G^o$ is a torus, the moduli space $\MHiggs^0_{Z_G^0}$ is smooth, see Case \eqref{i:higgs rational 1}. Hence, the fiber $M_{G,\cL}^{(E,\theta)}$ has rational singularities if and only if $M_{G,\cL}^d$ has rational singularities. Since we assumed $D_G$ almost-simple, the latter statement follows from Cases \eqref{i:higgs rational 2} and \eqref{i:higgs rational 3}.

Case \eqref{i:higgs rational 5}. 
By \Cref{L:isogeny-finite} and \Cref{L:reduction to D(G) almost-simple}, for any $d \in \pi_1(G)$, there exists a morphism $\widehat{G} \to G$ of reductive groups, and a finite surjective morphism of varieties 
\begin{align}\label{E:reduction to twisted 1}
M_{Z_G^o}^0\times M_{\widehat{G},\cL}^{\widehat d}\to M^d_{G,\cL}\times M_{\widehat{G}/D_{\widehat{G}},\cL}^{\widehat{d}} 
\end{align}
such that \begin{itemize}
    \item $D_{\widehat{G}}$ is simply-connected, and 
    \item $d\in \pi_1(G)$ admits a lift $\widehat{d} \in \pi_1(\widehat{G})\simeq \pi_1(\widehat{G}/D_{\widehat{G}})$. 
\end{itemize} 
Restricting the morphism \eqref{E:reduction to twisted 1} over a $k$-point $(E,\theta)$ in $M_{\widehat{G}/D_{\widehat{G}},\cL}^{\widehat{d}}$, we obtain a finite surjective morphism of varieties
\begin{align}\label{E:reduction to twisted 2}
M_{Z_G^o,\cL}^0\times M^{(E,\theta)}_{\widehat{G},\cL}\to M^d_{G,\cL}.
\end{align}
% By \Cref{L:isogeny-finite} and \Cref{L:reduction to D_G almost-simple}, there exists a finite surjective morphism of schemes 
% \begin{align}\label{E:reduction to twisted 1}
% M_{Z_G^o}^0\times M_{\widehat{G},\cL}^{\widehat d}\to M^d_{G,\cL}\times M_{\widehat{G}/D(\widehat{G}),\cL}^{\widehat{d}} 
% \end{align}
% where $D(\widehat{G})$ is simply-connected, $Z_G^o$, $H/D(H)$ are tori, and $\widehat{d}\in\pi_1(\widehat{G})=\pi_1(\widehat{G}/D(\widehat{G}))$ is an element mapping to $d\in\pi_1(G)$. Fix a closed point $\{(E,\theta)\}$ in $M_{\widehat{G}/D(\widehat{G}),\cL}^{\widehat{d}}$. After restricting the morphism \eqref{E:reduction to twisted 1} to the closed subvariety $M_{G,\cL}^d\times\{(E,\theta)\}\hookrightarrow M^d_{G,\cL}\times M_{H/D(H),\cL}^{\widehat{d}}$, we get a finite surjective morphism of schemes
% \begin{align}\label{E:reduction to twisted 2}
% M_{Z_G^o,\cL}^0\times M^{(E,\theta)}_{\widehat{G},\cL}\to M^d_{G,\cL}.
% \end{align}
By Cases \eqref{i:higgs rational 2} and \eqref{i:higgs rational 4}, the source of \eqref{E:reduction to twisted 2} has rational singularities. By \cite[Prop.\ 5.13]{KM1998}, the target has rational singularities too.

We conclude that $\MHiggs_{G,\cL}^{d}$ is a symplectic variety. Furthermore, the Hitchin fibration $h: \MHiggs^d_{G, \cL} \to A$ is Lagrangian %(in particular proper and equidimensional; cf. \cite{Laumon1988} or \cite{Mat1999}); 
by
\cite[Prop.\ 4.5]{Hitchin1987a} or \cite[Prop.\ 5.9 and Rem.\ 5.10]{biswas-ramanan-infinitesimal}.
\end{proof}

\begin{remark}\label{rmk:literature}
   An analogue of \Cref{prop:LagrangianfibrationHitchin}, namely the simplecticity of the singularities, holds for $G$-character varieties too. This result have been achieved for groups of type A by \cite[Thm.\ 1.1 and Cor.\ 1.2]{BS19}, for a semisimple group and genus $\geq 374$ in \cite[Thm.\ B]{AA16}, and for arbitrary reductive group with no assumption on the genus in \cite[Thm.\ G.(2)]{HSS2021} and \cite[Thm.\ A and Cor. 1.1]{Shu2024}. The simplecticity of the G-character variety implies the simplecticity of the G-Higgs moduli space $\MHiggs^d_{G,\cL}$  for $d=0$, via Simpons' isosingularity principle \cite[Thm.\ 10.6]{simpson-repnII}. %; see for instance \cite[Paragraph below Thm.\ G]{HSS2021}. 
   In the proof of the \Cref{prop:LagrangianfibrationHitchin}, we provide a direct elementary proof of  the simplecticity of the moduli space $\MHiggs^d_{G,\cL}$, independent of Simpons' isosingularity, with no assumption on the degree.  However, we do not address the rationality of the singularities of the stack $\Higgs_{G,\cL}^{ss, d}$; the latter  follows for $d=0$, in increasing degree of generality,  by \cite[Thm.\ 1.3]{Budur2021}, \cite[Thm.\ B]{AA16} and \cite[Thm.\ G.(1)]{HSS2021}.
\end{remark}

\begin{coroll}\label{cor:relativedimcanonical}
     Suppose that $\text{char}(k)=0$ and that $\cL=\omega_C$. For any given $d \in \pi_1(G)$, the Hitchin fibration $h:\MHiggs_{G,\cL}^{d} \to A$ satisfies the following
     \begin{enumerate}[(1)]
         \item $\tau_{> 2\dim(h)} Rh_*(IC_{\MHiggs_{G,\cL}^{d}}[-\dim(\MHiggs_{G,\cL}^{d})]) =0$;
         \item $R^{2\dim(h)} h_*(IC_{\MHiggs_{G,\cL}^{d}}[-\dim(\MHiggs_{G,\cL}^{d})]) = R^{2\dim(h)} h_*\bQ_{\MHiggs_{G,\cL}^{d}}$.
     \end{enumerate} 
\end{coroll}
\begin{proof}
It follows from \Cref{prop:LagrangianfibrationHitchin} and \cite[Prop. 3.3 and 3.6]{MM2024}.
\end{proof}

% \begin{defn}

% A \emph{Lagrangian fibration} $f:X \to A$ is a proper surjective morphism with connected fibers between normal varieties, such that the general fiber (equivalently any irreducible component of any closed fiber) of $f$ is a Lagrangian subvariety.
%  \end{defn}
\end{subsection}

\end{section}

\begin{section}{\texorpdfstring{$\delta$}{delta}-regularity} \label{section: delta regularity}

\begin{context}\label{cont: tt1}
    In this section, we assume that the ground field $k$ is algebraically closed and of characteristic $0$.
     We fix a reduced divisor $D \subset C$ such that $\deg(\omega_C(D))>0$.
\end{context}

In this section,  we prove \Cref{thm:deltaregularmero}, namely that in \Cref{cont: tt1},  the group scheme of symmetries of the Hitchin fibration $P^{\circ}$ (cf. \Cref{P:centr-group}(ii)) is $\delta$-regular in the sense of \Cref{defn: delta regularity}.  In \Cref{remark: delta regularity for twisted higgs}, we observe that the same argument applies in the case of isotrivial group schemes (cf. \Cref{section: isotrivial group schemes}).

\subsection{The tangent complex of the Higgs moduli stack} \label{subsection: tangent complex of the Higgs stack}

A key idea in our proof of $\delta$-regularity is the use of the following notion at the level of the moduli stack of Higgs bundles.
\begin{notn}[Tangent space of stacks] 
    For a stack $\scrX$ locally of finite type over $k$, and a geometric point $x \to \scrX$, we denote by $T_{\scrX,x}:= \cH^0(x^*\bT_{\scrX/k})$ the tangent space of the stack at $x$, given by the zeroth cohomology of the restriction of the tangent complex $\bT_{\scrX/k}$ to $x$.   By using results in \cite{olsson-cotangent-deformations}, this agrees with the usual definition of tangent space at a point in terms of morphisms from the scheme of dual numbers.
    % these \href{https://people.mpim-bonn.mpg.de/gaitsgde/grad_2009/SeminarNotes/Sept22(Dmodstack1).pdf}{notes by Gaitsgory}
\end{notn} 

In order to obtain a more canonical description of certain duality pairings appearing in tangent spaces of the moduli of Higgs bundles, for the rest of this section we adopt the following convention.
 \begin{notn}\label{rmk:dualornot}
 In this section, the Higgs stack $\Higgs_{G,\omega_C(D)}$ will be viewed as the section stack $\text{Sect}_C([\mathfrak{g}^{\vee} \otimes \omega_C(D) /G])$ for the morphism $[\fg^{\vee} \otimes \omega_C(D)/G] \to C$. 

\end{notn}

\begin{remark}\label{rmk: annoying}
Recall that we have previously viewed the Higgs stack $\Higgs_{G,\omega_C(D)}$ as a section stack $\text{Sect}_C([\mathfrak{g} \otimes \omega_C(D) /G])$ for the morphism $[\fg \otimes \omega_C(D)/G] \to C$ (see Section \ref{section: stack of higgs bundles}).
This point of view can be interchanged with the one in \Cref{rmk:dualornot} by choosing a $G$-invariant inner product on $\mathfrak{g}$ and by using it to identify $\mathfrak{g} \simeq \mathfrak{g}^{\vee}$.

 We note that the convention in \Cref{rmk:dualornot} is more natural in deformation theory  than the one in \Cref{rmk: annoying}: the tangent complex $\mathbb{T}_{\mathfrak{g}^{\vee}/G}=[\mathfrak{g} \to \mathfrak{g}^{\vee}]$ is canonically isomorphic to $\mathbb{T}^{\vee}_{\mathfrak{g}^{\vee}/G}[1]$, while the isomorphism $\mathbb{T}_{\mathfrak{g}/G} = [\mathfrak{g} \to \mathfrak{g}] \simeq [\mathfrak{g}^\vee \to \mathfrak{g}^\vee] = \mathbb{T}^{\vee}_{\mathfrak{g}/G}[1]$ depends on a choice of an identification $\mathfrak{g} \simeq \mathfrak{g}^{\vee}$.
\end{remark}

We end this subsection by recalling a standard description of the 
tangent spaces of $\Higgs_{G, \omega(D)}$. There is a naturally 
defined derived enhacement $\Higgs^{der}_{\omega(D)} = 
\text{Sect}^{der}_C([\mathfrak{g}^{\vee} \otimes \omega_C(D) /G])$ 
of the classical section stack  
$\text{Sect}_C([\mathfrak{g}^{\vee} \otimes \omega_C(D) /G])$ 
which is defined, for example, in 
\cite{halpernleistner2019mapping}. The main feature that we need 
about $\Higgs^{der}_{\omega(D)}$ is that its tangent complex 
$\mathbb{T}_{\Higgs^{der}_{\omega(D)}}$ has the following simple 
description. Let $S$ be a (classical) scheme and let $x: S \to 
\Higgs^{der}_{\omega(D)}$ be a morphism. Any such $x$ will factor through the classical substack $\text{Sect}_C([\mathfrak{g}^{\vee} 
\otimes \omega_C(D) /G])=\Higgs_{G, \omega(D)} \hookrightarrow \Higgs^{der}_{G, \omega(D)}$, and hence it corresponds to a $C$-
morphism $\sigma: C_S \to [\mathfrak{g}^{\vee} \otimes 
\omega(D)/G]$. Then, by the deformation theory developed in \cite{halpernleistner2019mapping}, the pullback $x^*
(\mathbb{T}_{\Higgs^{der}_{\omega(D)}})$ is naturally quasi-isomorphic to $R(\pi_S)_* L\sigma^* 
\mathbb{T}_{[\mathfrak{g}^{\vee} \otimes \omega(D)/G]/ \, C}$, 
where we denote by $\pi_S: C_S \to S$ the structural morphism. 

If we think of the data of $x: S \to \Higgs^{der}_{\omega(D)}$ as an $\omega(D)$-valued $G$-Higgs bundle $(E, \theta)$ on $C_S$, then, by the standard theory of (relative) tangent complexes of quotient stacks, we get the following description of the pullback $L\sigma^*\mathbb{T}_{[\mathfrak{g}^{\vee} \otimes \omega(D)/G]/ \, C}$:
\[ L\sigma^*\mathbb{T}_{[\mathfrak{g}^{\vee} \otimes \omega(D)/G]/ \, C} = \mathscr{K}^\bullet_{(E, \theta)}\coloneqq [\ad(E) \xrightarrow{\ad^*(\cdot, \theta)} \ad(E)^\vee \otimes \omega_C(D)], \]
 where the two terms of the complex  $\mathscr{K}^\bullet_{(E, \theta)}$ are placed in (cohomological) degrees $-1$ and $0$.

We have a natural identification of truncations $\tau_{\leq 0} \left(x^*\mathbb{T}_{\Higgs_G, \omega(D)}\right) = \tau_{\leq 0} \left(x^*\mathbb{T}_{\Higgs^{der}_{G, \omega(D)}}\right)$, so that we may use the complex defined above in order to compute tangents spaces of the stack $\Higgs_{G, \omega(D)}$. In particular, if the scheme $S$ above is of the form $\Spec(K)$ for some field $K$, then we get a natural identification of $K$-vector spaces $T_{\Higgs_{G, \omega(D)}, x} = H^0(C_K, \mathscr{K}^\bullet_{(E, \theta)})$.

\begin{remark}
It is known that $\Higgs^{der}_{G, \omega(D)}$ coincides with $\Higgs_{G, \omega(D)}$ either when $D>0$, or when $G$ is semisimple. We don't make use of this fact.
\end{remark}

\begin{subsection}{A general context for \texorpdfstring{$\delta$}{delta}-regularity}
    In this subsection, we describe a general context where we are able to prove $\delta$-regularity for certain smooth commutative group schemes. In Subsection \ref{subsection pf d reg}, we apply this to the case of the group of symmetries $\schP^o$ of the Hitchin fibration $h: \MHiggs_{G,\omega_C(D)}^d \to A$. 

Let us recall the definition of $\delta$-regular group scheme.

\begin{defn}[$\delta$-invariant] \label{defn: delta invariant} Let $P \to B$ be a smooth commutative group scheme of finite type over a scheme $B$. For a geometric point $b$ in $B$, we define $\delta(b)$ as the dimension of the maximal affine and connected subgroup of the neutral component $P^o_b \subset P_b$ (see \Cref{T:sCd-kclosed}). For any locally closed irreducible subscheme $Z \subseteq B$, we denote by $\delta(Z)$ the minimum value of the (upper semi-continuous, see \cite[Lm. 5.6.3.]{ngo-lemme-fondamental}) function $\delta$ on $Z$.    
\end{defn}

\begin{defn}[$\delta$-regularity] \label{defn: delta regularity}
Let $P \to B$ be a smooth commutative group scheme of finite type over a scheme $B$ which is locally of finite type over $k$.
We say that $P \to B$ is $\delta$-regular if for any 
locally closed irreducible subscheme $Z \subseteq B$ we have the inequality
\begin{equation}\label{delta reg}
\codim_{B} (Z) \geq \delta(Z),
\end{equation}
where $\delta(Z)$ is as in \Cref{defn: delta invariant}.
%$\codim Z_{\delta} \geq \delta$ for any $Z_{\delta} \coloneqq \{b \in B \,|\, \delta(b)=\delta\}$.
\end{defn}
%Let us recall recall the definition of $\delta$-regularity.

% Given any $x \in \mathrm{MHiggs}_{G, \cL}$, there exists a unique closed point $y$ in $\pi^{-1}(x)$, since $\pi$ is a good quotient. Denote by $H_y$ the isotropy of $y$. Due to the $\schP^\circ$-equivariance of $\pi$, the stabilizer $\Stab(x) \subseteq \schP^\circ_{h(x)}$ preserves the fiber $\pi^{-1}(x)$, and so it fixes the unique closed point $y$. Therefore, the action map $\act_y \colon \Stab(x) \to \mathrm{Higgs}_{G, \cL}$, given by $\act_y(j)=j \cdot y$, factors through $\Stab(x) \to BH_y \hookrightarrow \mathrm{Higgs}_{G, \cL}$. Since $T_y BH_y =0$, the differential $d\act_y \colon \Lie(\Stab(x)) \to T_{\mathrm{Higgs}_{G, \cL}, y}$ must vanish.

% Let $T_{A/A_{D}}$ be the relative (or vertical) tangent sheaf of the linear map $q \colon A \twoheadrightarrow A_{D}$.

\begin{prop}\label{prop:deltaregduality} A smooth commutative group scheme  $P \to B$ of finite type over a smooth $k$-scheme $B$ is $\delta$-regular if  the following three conditions are met:
\begin{enumerate}[(1)]
    \item\label{item: setting} there is a proper and surjective  morphism $h \colon X \to B$ from a  scheme $X$ %\Mirko{Maybe algebraic spaces: what about Borel fixed point theorem?}\andres{Yes, we know how to prove it, but if we quote the litearture, we need to find a reference for Borel fixed theorem for algebraic spaces}
    equipped with an action $\act: P\times_{B} X \to X$ over $B$; 
    \item \label{item:goodmodulispace} $\pi:\scrX \to X$ is the good moduli space of an algebraic stack $\scrX$ such that the $P$-action on $X$ lifts to $\scrX$;

    \item \label{item:duality} (Duality for the action) for any closed $k$-point $x$ of $\scrX,$ there exists a commutative diagram 
    \begin{equation}\label{eq duality act}
    \xymatrix{
	{\Lie(P_{h(\pi(x)}))} \ar[r]^-{\simeq} \ar[d]_{d\act_x}	& T^{\vee}_{B, h(\pi(x))} \ar[d]^{(h \circ \pi)^{*}}\\
    T_{\scrX, x} \ar[r]^{\simeq}  & T^{\vee}_{\scrX, x},}
    \end{equation} 
    where 
    %$\pi \colon \scrX \to X$ is the good moduli space morphism, 
    the action $\act_x \colon P_{h\circ \pi(x)} \to \scrX$ is given by $\act_x(j)=j \cdot x$, the horizontal arrows are isomorphisms,
    and $(-)^*$ denotes the codifferential.
\end{enumerate}
\end{prop} 
\begin{remark}
   In practice, the isomorphism $T_{\scrX, x} \simeq T^{\vee}_{\scrX,x}$ in (\ref{eq duality act}) is induced by a symplectic form. In fact, if furthermore $\scrX=X$ is smooth, \Cref{prop:deltaregduality} recovers \cite[Thm. 2.1]{AR2016} (which assumes in addition that $h$ has integral fibers), \cite[Section 2]{Ngo2011} 
   and \cite[Prop. 2.3.2]{dCRS2021}. However, if 
   $X$ is singular, then the arguments in loc.\ cit.\ do 
   not  go through: a symplectic form $\omega \in 
   H^0(X, \Omega^{[2]}_{X})$ does not induce an 
   isomorphism between tangent and cotangent spaces at 
   singular points. Note also that in general, in the 
   proof of \Cref{prop:deltaregduality},
   the Borel Fixed Point Theorem singles out singular 
   points of $X$, so that singularities are unavoidable. 
   On the positive side, hypotheses 
   \eqref{item:goodmodulispace} and \eqref{item:duality} 
   give the flexibility to apply the argument in \cite[Thm. 2.1]{AR2016}
   beyond the  case when $\scrX=X$ is smooth and $h$ has integral fibers. 
\end{remark}

\begin{proof}[Proof of \Cref{prop:deltaregduality}]
     We follow \cite[Prop. 2.3.2]{dCRS2021} closely. Without loss of generality (see the definition of $\delta(Z)$ in \Cref{defn: delta invariant}), we may replace $P$ with its open subgroup $P^o$ of fiberwise neutral components and assume that $P \to B$ has connected fibers. Let $Z \subset B$ be a locally closed integral subscheme. By shrinking $Z$ if necessary,  and by recalling that we are assuming that the ground field has characteristic zero, we have  the Chevalley devissage of $P_Z$ in its  affine and Abelian parts (cf.\ \Cref{T:sCd-kclosed})
    \[0 \to P^{\mathrm{aff}}_Z \to P_Z \to P^{\mathrm{ab}}_Z \to 0.\]
By assumption \eqref{item:duality}, we have the equality  $\dim(P_z) 
 = \dim(B)$ for every  $z \in Z(k)$. 
 
 It suffices to show that for a general $z\in Z(k)$  there exists a surjection
    \begin{equation}\label{eq:surjdelta}
        \Lie(P^{\mathrm{ab}}_z) \twoheadrightarrow T^{\vee}_{Z, z},
    \end{equation}
    since this implies that $\dim P^{\mathrm{ab}}_z \geq \dim Z$,  so that the desired inequality (\ref{delta reg}) follows by the ensuing chain of inequalities
\begin{equation} \label{eqn: z-inequality}
        \codim_{B} (Z) = \dim(B) - \dim(Z)  \geq \dim(P_z) - \dim(P^{\mathrm{ab}}_z) = \dim(P^{\mathrm{aff}}_z) = \delta(Z).
    \end{equation}

  Let $F_Z$ be the reduction of the closed subscheme in $h^{-1}(Z)$ fixed by the action of $P^{\mathrm{aff}}_{Z}$. Note that $h(F_{Z})=Z$ by Borel's fixed-point theorem \cite[Prop. 15.2]{borel-linear-groups}.
   %Let $\eta$ denote the generic point of $Z$. The closed subscheme \Mark{this is reduced if $X$ is reduced; but is this important?} $X_{\eta}^{P^{\mathrm{aff}}_{\eta}} \subset X_{\eta}$ 
   %of $P^{\mathrm{aff}}_{\eta}$-fixed points 
   %in $X_{\eta}$ is nonempty by the Borel Fixed-Point Theorem \cite[Prop. 15.2]{borel-linear-groups}. We denote by $F_Z \subset X_Z$ the closure of $X_{\eta}^{P^{\mathrm{aff}}_{\eta}}$ {\cb in $X_Z$} \Mark{to use generic smoothness, we need this to be reduced; but we could always take it with its reduced structure, right?}. After perhaps shrinking $Z$, we have that the closed subscheme $F_Z \subset X_Z$ is fixed by the action of $P^{\mathrm{aff}}_{Z}$.
  In  the next intermediate Claim, we find a local \'{e}tale section over $Z$ of the restriction $\pi^{-1}(F_{Z}) \to Z$  of the morphism $h \circ \pi: \scrX \to B$. 
   
   \begin{claim}\label{lem:etalesection} 
      For a general point $z \in Z(k)$, there is a morphism of stacks $f: \mathscr{Q} \to \scrX$ and a $k$-point $q \in \mathscr{Q}(k)$  mapping to $z$ and satisfying the following:
   \begin{enumerate}[(a)]
       \item  The image $y:= f(q)$ is a closed point of $\scrX$ and $\pi(f(q))=\pi(y) \in (h^{-1}(z) \cap F_{Z}) (k)$.

       \item $\mathscr{Q}$ dominates $Z$ and the differential $d(h\circ \pi \circ f): T_{\mathscr{Q},q} \to T_{Z,z}$ is an isomorphism.
   \end{enumerate} 
   \end{claim} 
   
   \noindent
{\em Proof of the {\bf Claim.}} Since the statement is \'etale local on 
  $Z$ and the formation of good moduli spaces commutes 
  with base-change, we may replace $Z$ with an \'etale over-scheme $V \to Z$ throughout the proof. By shrinking $Z$ if necessary, we may assume that the projection $F_Z \to Z$ is smooth by generic 
  smoothness. Up to base-change by an 
  \'etale morphism $V \to Z$, we may choose a section 
  $\sigma: Z \to F_Z$. Consider the fiber product
    \[
    \xymatrix{
	\mathscr{Z} \ar[rr] \ar[d]_{\pi_Z} &	& \scrX \ar[d]^{\pi} 
 \\
    Z \ar[r]^{\sigma}  & F_Z \ar[r] & X,}
    \]
  where the bottom arrow $Z\to X$ factors through the 
  section $\sigma$ followed by the closed embedding $F_Z \to X$.
   Any closed point of $\mathscr{Z}$ is sent 
   to a closed point of $\scrX$ lying over 
   $F_Z$. Therefore, to conclude, it is enough to find a closed substack 
   $\mathscr{Q} \subset \mathscr{Z}$ that 
   dominates $Z$ and satisfies (b) at any 
   closed point $q \in \mathscr{Q}(k)$. Let 
   $\eta$ denote the generic point of $Z$. 
   Since $\pi_Z$ is a good moduli space, the 
   generic fiber $\pi_Z^{-1}(\eta)$ 
   contains 
   a unique closed  point $\eta'$. 
   By taking the scheme-theoretic closure of 
   $\eta'$ in $\mathscr{Z}$ we get a closed 
   substack $\mathscr{Q} \subset 
   \mathscr{Z}$. The 
   generic fiber $\mathscr{Q}_{\eta}$ is a 
   gerbe over an algebraic space, because it is reduced and has a single point \cite[\href{https://stacks.math.columbia.edu/tag/06QK}{Tag 06QK}]{stacks-project}. Since $\mathscr{Q}_{\eta} \to \eta$ is a good moduli space, then in fact $\mathscr{Q}_{\eta}$ is actually a gerbe over $\eta$.
   In particular, since we are working in characteristic $0$,
   the relative cotangent complex 
   $\mathbb{L}_{\mathscr{Q}_{\eta}/\eta}$ 
   has Tor amplitude concentrated in 
   cohomological degree $1$, so that
   $\mathscr{Q}_{\eta} \to \eta$ is smooth. 
   By openness of smoothness, after shrinking $Z$ we can suppose that $\mathscr{Q} \to Z$ is smooth, so $\mathbb{L}_{\mathscr{Q}/Z}$ has Tor amplitude in degrees $[0,1]$. Furthermore,  after further shrinking $Z$, we may assume that $\mathbb{L}_{\mathscr{Q}/Z}$ has Tor amplitude concentrated in cohomological degree $1$. For all $k$-points $q \in \mathscr{Q}(k)$ with image $z \in Z(k)$, the exact sequence
   \[
   0=H^0(\mathbb{L}^{\vee}_{\mathscr{Q}/Z, q}) \to T_{\mathscr{Q},q} \to T_{Z,z} \to  H^1(\mathbb{L}^{\vee}_{\mathscr{Q}/Z, q})=0
   \]
   gives the desired isomorphisms of tangent spaces. The \textbf{Claim} is thus proved.

  We conclude the proof of \Cref{prop:deltaregduality}.
  Set $x:= \pi(y) \in F_Z(k) \subseteq X(k)$. Since $\pi$ is $P$-equivariant and $P^{\mathrm{aff}}_z$ fixes $x$, the action of $P^{\mathrm{aff}}_z$ preserves the fiber $\pi^{-1}(x)$, and so it fixes the unique  closed point $y$ in $\pi^{-1}(x)$. Hence, the action $\act_y \colon P^{\mathrm{aff}}_z \to \scrX$ %, given by $\act_y(j)=j \cdot y$, 
   factors through
   \[P^{\mathrm{aff}}_z \to \mathscr{G}_y \hookrightarrow \scrX,\] where $\mathscr{G}_y$ is the residual gerbe at $y$. Since $T_{\mathscr{G}_y, y}=0$, the differential $d\act_y \colon \Lie(P^{\mathrm{aff}}_z) \to T_{\scrX, y}$ must vanish, so that the infinitesimal action factors as follows: \[d\act_y \colon \Lie(P_z) \twoheadrightarrow \Lie(P^{\mathrm{ab}}_z) \to T_{\scrX, y}.\]

   Together with the assumption \eqref{item:duality},\, we obtain the following commutative diagram \[
    \xymatrix{
	\Lie(P^{\mathrm{ab}}_z) \ar[rd]& \Lie(P_z)\ar@{->>}[l] \ar[r]^{\simeq} \ar[d]_{d\act_y}	& T^{\vee}_{B, z} \ar[d]^{(h\circ \pi)^*} \ar@{->>}[r] & T^{\vee}_{Z, z} \ar[d]^{\simeq}\\
    & T_{\scrX, y} \ar[r]^{\simeq}  & T^{\vee}_{\scrX, y}\ar[r]_{f^*} & T^{\vee}_{\mathscr{Q}, q}.}
    \]
    
    This gives the required surjectivity \eqref{eq:surjdelta}.
\end{proof}

The following corollary is to be applied to the symplectic leaves of the Poisson structure associated with the Hitchin morphism.

\begin{coroll}\label{cor:dual} Assume that $P \to B$ is a smooth commutative group scheme satisfying the hypotheses \eqref{item: setting} and \eqref{item:goodmodulispace} of \Cref{prop:deltaregduality}. Suppose furthermore that:

\begin{enumerate}
   % \item[(1$^{\prime}$)] $Y \subseteq X$ is a $P_{h(Y)}$-invariant locally closed subset such that and $h \colon Y \to h(Y)$ is proper; \Mark{if $h(Y)$ is not a scheme, it is not clear to me what $P_{h(Y)}$-invariant means \ldots How to formulate (1') in a light way? I try below. Also, Do we need reducendness in the assumption of the corollary and of the proposition?}\Mirko{$h(Y)$ with its own reduced structure is a locally closed scheme {\Md Why locally closed?}. Note that the emphasis should be more on $Y$ than on $B$}

    \item[(1$^{\prime}$)] $B' \subseteq B$
    is a locally closed smooth subscheme, and
    $Y \subseteq X\times_B B'$ is a $P_{B'}$-invariant closed subscheme; 
    
    \item[(2$^{\prime}$)] (Duality for the action) for any closed $k$-point $y$ of the stack $\scrY \coloneqq \scrX \times_{X} Y$, there exists a commutative diagram  
    \[
    \xymatrix{
	\Lie(P_{h( \pi(y))}) \ar[r]^{\simeq} \ar[d]_{d\act_y}	& T^{\vee}_{B, h (\pi(y))} \ar[d]^{(h \circ \pi)^*}\\
    T_{\scrY, y} \ar[r]^{\simeq}  & T^{\vee}_{\scrY, y}.}
    \]
\end{enumerate}
Then for any $Z \subseteq B',$ we have the inequality $\codim_{B} (Z) \geq \delta(Z)$. 
\end{coroll}
\begin{proof}
    The proof of the $\delta$-regularity inequality \eqref{eqn: z-inequality} in \Cref{prop:deltaregduality} holds for any locally closed irreducible subscheme $Z\subset B'$ under the current assumptions by replacing $X$ and $\scrX$ with  $Y$ and $\scrY$, respectively.
\end{proof}
\end{subsection}

\begin{subsection}{Duality and symplectic leaves of the Higgs stack}\label{sec:symplecticleaves}

In this subsection, we include a discussion of the   foliation of the Higgs moduli stack $\Higgs_{G,\omega_{C}(D)}$ given by fixing  the  residues of  Higgs fields at the points of $D$. The goal is to show that the  hypotheses $(1')$ and $(2')$ in \Cref{cor:dual} hold for each leaf; this is achieved in \Cref{Lem:stratapreservedbyaction} and \Cref{square:delta}, respectively. This is a key ingredient in the proof of the main result of this section, namely \Cref{thm:deltaregularmero}, i.e.\ the $\delta$-regularity of the group scheme of symmetries $\schP^{o}$ for the Hitchin fibration. 
    %We limit ourselves to  develop what is needed for the  proof of $\delta$-regularity to be found in \S\ref{subsection pf d reg}. We refer the reader to {\cb REF - probably a later section} for a more complete treatment of the Poisson geometry of $\Higgs_{G,\cL}$.

    % In this subsection, we include a discussion of: the stratification of the Higgs moduli stack $\Higgs_{G,\cL}$ by conjugacy classes of residues; the properties with respect to the natural Poisson structure $\Theta$, when $\cL= \omega_{C}(D)$ for some effective divisor $D \geq 0$. 
    % The goal is to show that the {\cb duality-type hypothesis  in \Cref{cor:dual}.(2')} holds for each stratum \Mark{so that???}
    % We limit ourselves to  develop what is needed for the  proof of $\delta$-regularity to be found in \S\ref{subsection pf d reg}. We refer the reader to {\cb REF - probably a later section} for a more complete treatment of the Poisson geometry of $\Higgs_{G,\cL}$.

The following  \Cref{lemma:duality} is the (local) group-theoretic
underpinning of the  (global) duality result for the Hitchin fibration (see \Cref{square:delta}). Denote by $\chi$ the quotient map $\chi \colon \mathfrak{g}^{\vee} \twoheadrightarrow \fc^{\vee} \coloneqq \mathfrak{g}^{\vee} \sslash G = \Spec(\mathrm{Sym}^{\bullet}\mathfrak{g})^{G}$,
and let $J \to \fc^{\vee}$ be the regular centralizer group scheme (cf.\ \S \ref{SS:local-centr-gr} and Remark \ref{rmk:dualornot}). By construction, there exists a (generically injective) morphism $m \colon \chi^* J \to G \times \mathfrak{g}^{\vee}$ of smooth group schemes over $\mathfrak{g}^{\vee}$. 

 The group $G \times \bG_m$ acts on $\mathfrak{g}^{\vee}$, where the $\mathbb{G}_m$ component acts via scalar multiplication (weight $1$). Therefore, the vector bundles $T_{\fg^{\vee}}$ and $T^{\vee}_{\fg^{\vee}}$ over $\fg^{\vee}$ comes equipped with an induced $G \times \mathbb{G}_m$-equivariant structure. The pullback $\chi^* T_{\fc^{\vee}}$ under the $G\times \mathbb{G}_m$-equivariant morphism $\chi: \fg^{\vee} \to \fc^{\vee}$ is also equipped with an equivariant structure (note that here $\mathbb{G}_m$ acts on the zero fiber with weights  $(e_1,\ldots , e_r)$ as in \Cref{T:konstant}). On the other hand, recall that the group scheme $J \to \fc^{\vee}$ is naturally $\bG_m$-equivariant \cite[Prop. 3.3]{NG06}, and the group $\mathbb{G}_m$ acts on the zero fiber of $\mathrm{Lie}(J)$ with weights $(-e_1+1, \ldots, -e_r +1)$ (see \cite[Lem. A.9]{dCHM2021}, or the proof of \cite[Prop. 4.13.2]{ngo-lemme-fondamental}).

In the proposition below, the twist $(-1)$ refers to twisting of the $\bG_m$-action by the $\bG_m$-character of
weight $-1$.

\begin{lemma}[{\cite[Rem. A.7 + Lem. A.9]{dCHM2021}}]\label{lemma:duality} The $G \times \bG_m$-equivariant
maps of relative Lie algebras over $\mathfrak{g}^{\vee}$
\[dm \colon \chi^* \Lie(J)(-1) \to \mathfrak{g} \times \mathfrak{g}^{\vee} =T^{\vee}_{\mathfrak{g}^{\vee}} \quad \text{ and }\quad
d\chi \colon T_{\mathfrak{g}^{\vee}}=\mathfrak{g}^{\vee} \times \mathfrak{g}^{\vee} \to \chi^* T_{\fc^{\vee}} \simeq \fc^{\vee} \times \fg^{\vee}
%\simeq \fc
\]
are $G \times \bG_m$-equivariantly dual to each other with respect to the canonical $(G \times \bG_m)$-equivariant pairing $\langle \cdot, \cdot \rangle \colon T_{\mathfrak{g}^{\vee}} \times_{\mathfrak{g}^{\vee}} T^{\vee}_{\mathfrak{g}^{\vee}} \to k$.
\end{lemma}

The commutative group stack $\stP^\circ \to A$ of symmetries of the Hitchin fibration  acts on the stack $\Higgs_{G, \omega_C(D)}$, and its rigidification $\schP^\circ = \stP^\circ\fatslash Z_G$ acts on the rigidification 
\[
\Higgsrig_{G, \omega_C(D)} \coloneqq \Higgs_{G, \omega_C(D)} \fatslash Z_G.
\]

The inclusion $i_D: D \hookrightarrow C$ induces the restriction morphism
\[ev_D \colon \Higgs_{G, \omega_C(D)} \to [\mathfrak{g}^{\vee}_D/G_D], \qquad (E,\phi) \mapsto i^*_D(E,\phi).\]
Here $\fg^{\vee}_D$ and $G_D$ are the Weil restrictions of $\fg^{\vee} \times D$ and $G \times D$ under the structure morphism $D \to \Spec(k)$, and we are using the canonical residue isomorphism $\omega_C(D)|_{D} \xrightarrow{\sim} \cO_D$ to define $ev_D$. Note that $\fg^{\vee}_D \cong (\fg^{\vee})^{\oplus |D|}$, and that $G_D \cong G^{|D|}$ acts on $\fg^{\vee}_D$ in the natural way. Hence, we have that $[\fg^{\vee}_D/G_D] \cong [\fg^{\vee}/G]^{|D|}$. There is a commutative square
\begin{equation}\label{eq: dfgt}
\xymatrix{
\Higgs_{G, \omega_C(D)} \ar[r]^-{ev_D} \ar[d]_-{h} & [\mathfrak{g}^{\vee}_D/G_D]=[\fg^{\vee}/G]^{|D|} \ar[d]^-{h_D} \\% \ar[r]^-{=} & [\fg/G]^{|D|} \ar[d]\\
A=H^0(C, \mathfrak{c}_{\omega_C(D)}) \ar[r]^-{ev_D} & A_D\coloneqq H^0(D, \mathfrak{c}^{\vee} \otimes \mathcal{O}_{D})=(\fc^{\vee})^{|D|}. %\ar[r]^-{=} & \fc^{|D|}. 
}
\end{equation}
Note that the morphisms $ev_D$ fails to be surjective when $G$ is not semisimple.

%We denote by $R$ the image $ev_D(E, \theta)=i_{D}^*(E, \theta)$. 
\begin{defn}[Foliation by  residues]
    For any $R \in [\fg^{\vee}_D/G_D](k)$, let $\mathscr{G}_R\hookrightarrow 
    [\mathfrak{g}^{\vee}_D/G_D]$ denote the locally closed residual gerbe of 
    $R$.  The leaf   $\Higgsrig_{G,\omega_C(D)}(R)$ is 
    the $Z_G$-rigidification of the fiber product $\Higgs_{G,\omega_C(D)}
    (R):=\Higgs_{G,\omega_C(D)} \times_{[\fg^{\vee}_D/G_D]} \mathscr{G}_R$. 
    %, or equivalently the fiber product $\Higgsrig_{G,\omega_C(D)} \times_{[\fg^{\vee}_D/(G_D/Z_G)]} (\mathscr{G}_R \fatslash Z_G)$ of the $Z_G$-rigidifications. 
\end{defn}
% By construction, there is a locally closed embedding of stacks $\Higgs_{G,\omega_C(D)}(R)\hookarr\Higgs_{G,\omega_C(D)}$. After taking $Z_G$-rigidifications, we get a locally closed embedding of rigidified stacks $\cH_{G,\omega_C(D)}(R)\hookarr\cH_{G,\omega_C(D)}$.

% After taking $Z_G$-rigidifications, the locally closed embedding of stacks $\Higgs_{G,\omega_C(D)}(R)\hookarr\Higgs_{G,\omega_C(D)}$ induces a locally closed embedding of rigidified stacks $\cH_{G,\omega_C(D)}(R)\hookarr\cH_{G,\omega_C(D)}$.

By construction, we have the following Cartesian diagram of stacks
\[
\xymatrix{
\Higgs_{G,\omega_C(D)}(R)\ar[d]\ar@{^{(}->}[r]&\Higgs_{G,\omega_C(D)}\ar[d]\\
\cH_{G,\omega_C(D)}(R)\ar@{^{(}->}[r]&\cH_{G,\omega_C(D)},
}
\]
where the vertical arrows are the $Z_G$-rigidification morphisms and the horizontal ones are the locally closed embeddings induced by the residual gerbe $\mathscr{G}_R\hookrightarrow 
    [\mathfrak{g}^{\vee}_D/G_D]$.
    
\begin{lemma}\label{Lem:stratapreservedbyaction}
    For any $R \in [\fg^{\vee}_D/G_D](k)$, the locally closed substack $\Higgsrig_{G,\omega_C(D)}(R) \hookrightarrow \Higgsrig_{G,\omega_C(D)}$ is preserved by the action of the group scheme $\schP^o$.
\end{lemma}
\begin{proof}
    Let $I \to [\fg^{\vee}/G]$ denote the inertia stack of $[\fg^{\vee}/G]$ relative to $\fc^{\vee}$. There is a morphism of relatively affine commutative group $[\fg^{\vee}/G]$-stacks $J \times_{\fc^{\vee}} [\fg^{\vee}/G] \to I$ \cite[Lem. 2.1.1]{ngo-lemme-fondamental}. Since $J \to \fc^{\vee}$ is commutative, the classifying stack $BJ \to \fc^{\vee}$ is a commutative group stack over $\fc^{\vee}$. The morphism $J \times_{\fc^{\vee}} [\fg^{\vee}/G] \to I$ induces an action of $BJ$ on $[\fg^{\vee}/G]$ such that the isomorphism class of every given $k$-point of $[\fg^{\vee}/G]$ is preserved by the action of the corresponding fiber of $BJ$ over $\fc^{\vee}$ (see \Cref{rmk: BJ action}). Let $(BJ)_D \to A_D$ denote the Weil restriction of $(BJ) \times D$ under the structure morphism $D \to \Spec(k)$. Then $(BJ)_D \cong (BJ)^{|D|}$ acts on $[\fg^{\vee}_D/G_D] \cong [\fg^{\vee}/G]^{|D|}$ and preserves the isomorphism class of any given $k$-point of $[\fg^{\vee}_D/G_D]$. In particular, by the reducedness of the fibers of $BJ$, this guarantees that the residual gerbe $\mathscr{G}_{R} \hookrightarrow [\fg^{\vee}_D/G_D]$ is stable under the action of $(BJ)_D$. By construction, the commutative group stack $\stP$ is the Weil restriction of $ev^*(BJ_{\omega_C(D)}) \to C \times A$ under the morphism $C \times A \to A$, where $ev: C \times A \to \fc^{\vee}_{\omega_C(D)}$ is the natural evaluation morphism. Restriction under the inclusion $i_D: D \hookrightarrow C$ induces a commutative diagram
    \[
\xymatrix{
\stP \ar[r] \ar[d] & (BJ)_D \ar[d]\\
A \ar[r]^-{ev_D} & A_D. 
}
\]
By construction, the restriction morphism $\Higgs_{G,\omega_C(D)} \xrightarrow{ev_D} [\fg^{\vee}_D/G_D]$ is compatible with the corresponding actions of $\stP$ and $(BJ)_D$. Since the action of $(BJ)_D$ preserves the residual gerbe $\mathscr{G}_R \hookrightarrow [\fg^{\vee}_D/G_D]$, it follows that the action of $\stP$ preserves the locally closed substack $\Higgs_{G,\omega_C(D)}(R) := (ev_D)^{-1}(\mathscr{G}_R) \hookrightarrow \Higgs_{G,\omega_C(D)}$. By taking $Z_G$-rigidifications, it follows that the action of $\schP = \stP\fatslash Z_G$ preserves $\Higgsrig_{G,\omega_C(D)}(R) \hookrightarrow \Higgsrig_{G,\omega_C(D)}$. Since $\schP^o \subset \schP$ is an open subgroup scheme, it also preserves $\Higgsrig_{G,\omega_C(D)}(R)$.
\end{proof}

Let $R \in [\fg^{\vee}_{D}/G_D](k)$ with $r\coloneqq h_D(R) \in ev_D(A)$ $\subset A_D$ (cf. \eqref{eq: dfgt}), and $A(r) \coloneqq ev_D^{-1}(r) \simeq H^0(C, \mathfrak{c}_{\omega_C(D)} \otimes \mathcal{O}(-D))$.
\begin{prop}\label{square:delta} 
      For any $(E,\theta) \in \Higgs_{G,\omega_C(D)}(R)(k)$,
      we denote by the same symbol its image 
      in the rigidified stack $\cH(R):=\cH_{G,\omega_C(D)}(R)$
      and we  set $a \coloneqq h(E,\theta) $. There is a commutative diagram with vertical isomorphisms 
\[
        \xymatrix{
 T^{\vee}_{A(r), a} \ar[rr]^{h^*} \ar[d]_{\simeq} &   &T^{\vee}_{\Higgsrig(R), (E, \theta)} \ar[d]^{\simeq} \\
 \Lie(\schP^{o}_{a}) \ar[rr]^{d\act_{(E,\theta)}} & & T_{\Higgsrig(R), (E, \theta)}.
 }
\]
\end{prop}
\begin{proof}

 %The fiber of the tangent complex of the derived enhancement $\bT_{\Higgs_{G,\omega_C(D)}^{der}}|_{(E\,\theta)}$ is isomorphic in the derived category to the complex of derived global sections $R\Gamma(C,-)$ of the complex
Consider the complex of vector bundles on $C$
\[
 \mathscr{K}^\bullet_{(E, \theta)}\coloneqq [\ad(E) \xrightarrow{\ad^*(\cdot, \theta)} \ad(E)^\vee \otimes \omega_C(D)], \]
 placed in degrees $-1$ and $0$, defined in \Cref{subsection: tangent complex of the Higgs stack}. 
%
 %where $\ad^*\colon \mathfrak{g} \times \mathfrak{g}^{\vee} \to \mathfrak{g}^{\vee}$  denotes the coadjoint action of $\mathfrak{g}$  on $\mathfrak{g}^{\vee}$. In particular, \Mark{reference for tangent space being the same for stack and its rigidification} 
We have that
 \[H^0(C,\mathscr{K}^\bullet_{(E, \theta)})=T _{\Higgs_{G,  \omega_C(D)}, (E, \theta)} =T_{\Higgsrig_{G, \omega_C(D)}, (E, \theta)}.\]

 %  Recall that  the group scheme $J_{\omega_C(D)}\coloneqq J \times^{\bG_m}\omega_C(D)^{\times} \to \fc^{\vee}_{\omega_C(D)} \coloneqq \fc^{\vee} \times^{\bG_m} \omega_C(D)^{\times}$ over $\fc^{\vee}_{\omega_C(D)}$ is obtained from the group scheme
 % $J \to \fc^{\vee}$ via twisting by $\omega_C(D)$.
 Let $a: C \to \fc^{\vee}_{\omega_C(D)}$ be the section defined by $h(E, \theta)$. The fiber $\schP^{o}_a$ is the $Z_G$-rigidification of the neutral component of the stack $\stP_a$ which parametrizes $J_a \coloneqq a^* J_{\omega_C(D)}$-torsors over $C$; see \Cref{P:centr-group}. In particular, we have  $\Lie(\schP^{o}_{a})=H^1(C, \Lie(J_a))$.
 By the proof of \cite[Prop.\ 7.12]{dCHM2021}, the infinitesimal action 
  \[d\act_{(E, \theta)} \colon \Lie(\schP^{o}_{a})=H^1(C, \Lie(J_a)) \to T_{\Higgsrig_{G,  \omega_C(D)}, (E, \theta)}=H^0(C, \mathscr{K}^\bullet_{(E, \theta)})\] 
  is obtained by applying $H^0(C,-)$ to the morphism of complexes placed in degrees $-1$ and $0$
\[[\Lie(J_a) \to 0] \xrightarrow{(d m_{\omega_C(D)}, 0)} [\ad(E) \to \ad(E)^{\vee} \otimes \omega_C(D)]
= \mathscr{K}^\bullet_{(E, \theta)}.\] 
 
 Analogously, 
 the differential of the Hitchin morphism
\begin{align*}
    d(h \circ \pi)_{(E, \theta)} \colon  &T_{\Higgsrig_{G,  \omega_C(D)}, (E, \theta)} %=T_{\Higgs_{G, \omega_C(D)}, (E, \theta)}
    =H^0(C, \mathscr{K}^\bullet_{(E, \theta)}) \to  T_{A, a}=H^0(C, a^*T_{\fc^{\vee}_{\omega_C(D)}/C}) \simeq H^0(C, \fc^{\vee}_{\omega_C(D)}),
 \end{align*}
 is obtained by applying $H^0(C,-)$ to  the morphism of complexes placed in degrees $-1$ and $0$
\begin{equation}\label{eq:complex1}
 \mathscr{K}^\bullet_{(E, \theta)} =
 [\ad(E) \to \ad(E)^{\vee} \otimes \omega_C(D)] \xrightarrow{(0, d \chi_{\omega_C(D)})} [0 \to a^*T_{\fc^{\vee}_{\omega_C(D)}/C} \simeq \fc^{\vee}_{\omega_C(D)}]. 
\end{equation}

By dualizing \eqref{eq:complex1} and  by tensoring with $\omega_C[1]$, we obtain the morphism of complexes 
placed in degrees $-1$ and $0$
\[
[c^{\vee}_{\omega_C(D)}\otimes \omega_C \to 0] \xrightarrow{( (d\chi_{\omega_C(D)})^{\vee}\otimes Id_{\omega_C},0)}  [\ad(E)\otimes \omega_C(D)^{\vee} \otimes \omega_C \to \ad(E)^{\vee} \otimes \omega_C] \simeq (\mathscr{K}^\bullet_{(E, \theta)})^{\vee} \otimes \omega_C.
\]

By taking $H^0(C,-)$ and by using Serre Duality, we obtain the dual pullback morphism
\[
(h \circ \pi)^* \colon T^{\vee}_{A, a} = H^1(C, c^{\vee}_{\omega_C(D)}\otimes \omega_C) \to H^1(C, (\mathscr{K}^\bullet_{(E, \theta)})^{\vee} \otimes \omega_C)=T^{\vee}_{\Higgsrig_{G, \omega_C(D)}, (E, \theta)}.
\]

Now, the local duality of \Cref{lemma:duality} gives the following global statement: The morphism  $dm_{\omega_C(D)} \colon \Lie(J_{a}) \to E \times^{G} \mathfrak{g}=\ad(E)$ is dual up to 
tensorization with $\omega_C(D)$ to 
\[d\chi_{\omega_C(D)} \colon (E \times^{G} \mathfrak{g}^{\vee}) \otimes \omega_C(D)=\ad(E)^{\vee} \otimes \omega_C(D) \to \Lie(J_{a})^{\vee} \otimes \omega_C(D) = c_{\omega_C(D)},\]
where the last equality is \cite[Prop. 4.13.2]{ngo-lemme-fondamental}. 

 The canonical inclusion $\mathcal{O}(-D) \hookrightarrow \mathcal{O}$ induces the commutative square %\Mark{again, the tensor factor $K$ may be missing?}
%\begin{equation}\label{diagram:complexesdual} 
%\xymatrix{
%	[c^{\vee}_{\omega_C(D)}\otimes \omega_C \to 0] \ar[r]^-{((d \chi_{\omega_C(D)})^{\vee} \otimes Id_{\omega_C},0)} \ar[d] &  [\ad(E)\otimes \mathcal{O}(-D) \to \ad(E)^{\vee} \otimes \omega_C]=\mathscr{K}^\bullet_{(E, \theta)}\otimes \mathcal{O}(-D) \ar[d]	\\
% [c^{\vee}_{\omega_C(D)}\otimes \omega_C(D) \simeq \Lie(J_a)  \to 0] \ar[r]^-{(dm_{\omega_C(D)}, 0)} &  [\ad(E)\to \ad(E)^{\vee} \otimes \omega_C(D)]=\mathscr{K}^\bullet_{(E, \theta)}.
% }
%\end{equation}

\begin{equation}\label{diagram:complexesdual}
    \begin{tikzcd}[ampersand replacement=\&]
    {[c^{\vee}_{\omega_C(D)}\otimes \omega_C \to 0]} \ar[r] \ar[dd, "{((d \chi_{\omega_C(D)})^{\vee} \otimes Id_{\omega_C},0)}", labels = left] \& {[c^{\vee}_{\omega_C(D)}\otimes \omega_C(D) \simeq \Lie(J_a)  \to 0]} \ar[dd, "{(dm_{\omega_C(D)}, 0)}"]\\ \& \\
    {[\ad(E)\otimes \mathcal{O}(-D) \to \ad(E)^{\vee} \otimes \omega_C]} \ar[r] \ar[d, symbol = {=}]\& {[\ad(E)\to \ad(E)^{\vee} \otimes \omega_C(D)]} \ar[d, symbol = {=}] \\
     \mathscr{K}^\bullet_{(E, \theta)}\otimes \mathcal{O}(-D) \ar[r]\& \mathscr{K}^\bullet_{(E, \theta)}
    \end{tikzcd}
\end{equation}

     Applying $H^0(C,-)$ and using Serre duality, we obtain the following commutative square (note that we are reflecting the square along one of the diagonals)
\begin{equation}\label{eq:square}
    \xymatrix{
	T^{\vee}_{A, a} \ar[rr]^-{h^*} \ar[d]&  &T^{\vee}_{\Higgsrig_{G, \omega_C(D)}, (E, \theta)} \ar[d]^-{\Theta}	\\
 \Lie(\schP^{o}_{a}) \ar[rr]^-{d\act_{(E, \theta)}} &   &T_{\Higgsrig_{G,  \omega_C(D)}, (E, \theta)}.
 } 
\end{equation}

The bottom horizontal arrow in \eqref{diagram:complexesdual} fits in the exact triangle 
\[
\mathscr{K}^\bullet_{(E, \theta)}\otimes \mathcal{O}(-D) \to \mathscr{K}^\bullet_{(E, \theta)} \to i_{D, *}i_{D}^* \mathscr{K}^\bullet_{(E, \theta)} \xrightarrow{+1}
\]
with $i_D \colon D \hookrightarrow C$. Applying $H^0(C,-)$, we obtain the exact sequence
\begin{equation}\label{eq:evD}
    T^{\vee}_{\Higgsrig_{G, \omega_C(D)}, (E, \theta)} \xrightarrow{\Theta} T_{\Higgsrig_{G, \omega_C(D)}, (E, \theta)} = T_{\Higgs_{G, \omega_C(D)}, (E, \theta)} \xrightarrow{d \, ev_D} T_{[\mathfrak{g}_D/G_D], R}
\end{equation}

 Then the exact sequence \eqref{eq:evD} gives
\[\mathrm{Im}(\Theta)=\ker(d \, ev_{D, (E, \theta)})=T_{\Higgs(R), (E, 
\theta)} = T_{\Higgsrig(R), (E, \theta)}.\]
Since the pairing $\Theta \in \Hom(T^{\vee}_{\Higgsrig_{G, \omega_C(D)}, (E,\theta)}, 
T_{\Higgsrig_{G, \omega_C(D)}, (E,\theta)}) \simeq T^{\otimes 
2}_{\Higgsrig(R), (E,\theta)}$ is actually alternating  %This induces a pairing which is alternating; this can be seen 
by means of a similar reasoning based on a direct computation as in \cite[Lemma 7.14]{Markmann94}, or as in \cite[Prop. 3.1]{Dalakov},
%by its definition 
%using Serre duality, 
the map 
\[\Theta \colon T^{\vee}_{\Higgsrig_{G, \omega_C(D)}, (E, \theta)} \to T_{\Higgsrig(R), (E, \theta)}\] 
factors through  $T^{\vee}_{\Higgsrig(R), (E, \theta)}$. It follows that the square \eqref{eq:square} factors as follows:
\[
        \xymatrix{
	T^{\vee}_{A, a} \ar[rr]^{h^*} \ar@{->>}[d] &   &T^{\vee}_{\Higgsrig_{G, 
 \omega_C(D)}, (E, \theta)} \ar@{->>}[d] \ar[ddr]^-{\Theta} &	
 \\
 T^{\vee}_{A(r), a} \ar[rr]^{h^*} \ar[d]_{\simeq} &   &T^{\vee}_{\Higgsrig(R), (E, 
 \theta)} \ar[d]^{\simeq} & 
 \\
 \Lie(\schP^{o}_{a}) \ar[rr]^{d\act} & & T_{\Higgsrig(R), (E, \theta)} 
 %\ar[r]
 \ar@{^{(}->}[r]
 &T_{\Higgsrig_{G,  \omega_C(D)}, (E, \theta)},
 }
\]
so that the desired conclusion follows.
\end{proof}

    \begin{remark} \label{remark: poisson duality for twisted higgs}
        The same result holds with the same arguments (verbatim) if we replace $G$ with an isotrivial reductive group $\cG \to C$. More precisely, \Cref{square:delta} holds for the action of $\schP^o$ on the rigidification $\Higgsrig_{\cG, \omega_C(D)} = \Higgs_{\cG,\omega_C(D)} \fatslash \pi_*(Z_{\cG})$ as defined in \Cref{coroll: action group on rigidication of twisted higgs}.
    \end{remark}
\end{subsection}

\begin{subsection}{\texorpdfstring{$\delta$}{delta}-regularity for the group scheme \texorpdfstring{$\schP^o$}{Po}}\label{subsection pf d reg}

 Recall that the action of $\schP^\circ$ preserves the semistable locus $\Higgsrig^{ss}_{G, \omega_C(D)} = \Higgs_{G,\omega_C(D)}^{d,ss}\fatslash Z_G$, and so it acts on the good moduli space $\MHiggs_{G, \omega_C(D)}$ such that the good moduli space morphism  \[\pi \colon \Higgsrig^{ss}_{G, \omega_C(D)} \to \MHiggs_{G, \omega_C(D)}\] is $\schP^\circ$-equivariant. Recall that the connected components $ \MHiggs_{G,\omega_C(D)}^d$ of the moduli space are quasi-projective, labelled by elements of $d \in \pi_1(G)$.

\begin{thm}\label{thm:deltaregularmero}
    %If $\mathcal{L}=\omega_C(D)$ for some reduced effective divisor $D \geq 0$ such that $\deg(\omega_C(D))>0$, then 
    Suppose that the characteristic of the ground field $k$ is $0$, and let $D \subset C$ be a (possibly empty) reduced divisor such that $\deg(\omega_C(D))>0$. Then the group scheme $\schP^\circ \to A$ for $\cL= \omega_C(D)$ (as in \Cref{P:centr-group}(ii)) is $\delta$-regular.
\end{thm}
\begin{proof}

In this proof, we fix an arbitrary $d \in \pi_1(G)$.
    By \Cref{prop: properness of the hitchin fibration} the Hitchin 
    fibration $h: \MHiggs_{G,\omega_C(D)}^d \to A$ is proper, and by 
    \Cref{coroll: action on the semistable moduli space} it admits an 
    action of the group scheme $\schP^\circ \to A$. Let $Z\subseteq A$ be an integral subscheme.
Let $A'_D:= \mathrm{Im} (ev_D) \subseteq A_D.$ Choose $r \in A_D'(k)$ such that $A(r) \cap Z \neq \emptyset.$ We shall show that we can find $B' \subseteq A(r)$ and $Y\subseteq \MHiggs_{G,\omega_C(D)}^d$ as in \Cref{cor:dual} with $B' \cap Z$ dense in $Z \cap A(r)$.

% We shall show that we 
%     can find $B' \subseteq A$ and $Y\subseteq \MHiggs_{G,\omega_C(D)}^d$ as in \Cref{cor:dual} with $B' \cap Z$ dense in $Z$.

%\sout{and let $Z$ be any integral subvariety of $A(r)$}. 
Consider the fiber product $[\fg^{\vee}_D/G_D]_r := [\fg^{\vee}_D/G_D] \times_{A_D} r$. The stack $[\fg^{\vee}_D/G_D]_r$, being
a fiber of
$[\fg^{\vee}/G]^{|D|} \to (\fc^{\vee})^{|D|}$), 
has finitely many geometric points $\{R_i\}_{i \in I}$  which are all defined over $k$.  If we denote by $\Higgs_{G,\omega_C(D)}(r) \to A(r)$ the preimage of $h^{-1}(A(r)) \subset \Higgs_{G,\omega_C(D)}$, then we have a finite stratification $\Higgs_{G,\omega_C(D)}(r) = \bigsqcup_{i \in I} \Higgs_{G,\omega_C(D)}(R_i)$. The action of the group scheme $\schP^o_{A(r)} \to A(r)$ preserves each stratum $\Higgs_{G,\omega_C(D)}(R_i)$ by \Cref{Lem:stratapreservedbyaction}.

Choose an irreducible component $W$ of $Z \cap A(r).$
Choose $R_{W}$ among the $R_{i}$  such that $\Higgs_{G,\omega_C(D)} (R_W)$ dominates $W$ and is minimal with
such property, i.e.\ the closure $\overline{\Higgs_{G,\omega_C(D)}(R_W)} 
\subset \Higgs_{G,\omega_C(D)}(r)$ does not contain any other distinct 
$\Higgs_{G,\omega_C(D)}(R_i)$ dominating $W$. Consider the stack $\scrX \coloneqq (\Higgs_{G,\omega_C(D)}(r) \cap 
\Higgs_{G,\omega_C(D)}^{d,ss}) \! \fatslash \, Z_{G}$ with good moduli space $X$. % \coloneqq \MHiggs_{G,\omega_C(D)}^{d}$ 
Set also   
\[\Higgs_{G,\omega_C(D)}^{d,ss}
 (R_W):= \Higgs_{G,\omega_C(D)}(R_W) \cap \Higgs_{G,\omega_C(D)}^{d,ss}.\] 
%\qquad 
% \Higgs_{G,\omega_C(D)}^{d,ss}(r):= \Higgs_{G,\omega_C(D)}(r) \cap 
% \Higgs_{G,\omega_C(D)}^{d,ss}.\] 
By the minimality of $R_{W}$, there is  an 
open set $W^{\circ} \subset W$ such that the restriction 
\[\scrY \coloneqq 
(\Higgs_{G,\omega_C(D)}^{d,ss}(R_W)|_{W^{\circ}}) \! \fatslash \, Z_{G} \hookrightarrow \scrX|_{W^{\circ}}\] 
% \[\scrY \coloneqq 
% \Higgs_{G,\omega_C(D)}^{d,ss}(R_W)|_{W^{\circ}} \hookrightarrow \scrX|_{W^{\circ}} \coloneqq 
% \Higgs_{G,\omega_C(D)}^{d,ss}(r)|_{W^{\circ}} = \Higgs_{G,\omega_C(D)}^{d,ss}|_{W^{\circ}}\] 
is a closed immersion of stacks of finite type, all admitting actions 
of  $\schP^o_{W^{\circ}} \to {W^{\circ}}$.

By \cite[Prop. 4.7(i), Lem. 4.14]{alper-good-moduli}, it descends to a closed immersion of good moduli spaces 
% \[Y \coloneqq \MHiggs_{G,\omega_C(D)}^{d,ss}(R_W)|_{W^{\circ}} 
% \hookrightarrow X \coloneqq \MHiggs_{G,\omega_C(D)}^{d}|_{W^{\circ}} = h^{-1}(W^{\circ}).\]
\[Y \coloneqq \MHiggs_{G,\omega_C(D)}^{d,ss}(R_W)|_{W^{\circ}} 
\hookrightarrow X|_{W^{\circ}}  = h^{-1}(W^{\circ}).\]
%The morphism $Y \hookrightarrow h^{-1}(W^{\circ}) \to W^{\circ} = h(Y)$ is proper. 
We apply \Cref{cor:dual},  with $B \coloneqq A(r)$, $P\coloneqq \schP^\circ_{W^{\circ}}$, and $\scrX, \scrY, X, Y$ as above. We can do so since the hypothesis $(1')$  of the corollary is satisfied by construction and the hypothesis $(2')$ is satisfied by
 \Cref{square:delta}.  We conclude that $\codim_{A(r)}(W) \geq \delta(W)$ for any irreducible component of  $Z \cap A(r)$.

This implies that $\schP^{o}$ 
is $\delta$-regular as follows.
We have the following inequalities:
\[\codim_{A} (Z) \geq \codim_{A(r)} (W \cap A(r))\geq \delta(W \cap A(r)) \geq \delta(Z),\]
where the first one is due to the smoothness of $A$, the second one has been proved above, and the last one is due to the upper semicontinuity of the function $\delta.$ 
\end{proof}

\begin{remark} \label{remark: delta regularity for twisted higgs}
    In view of \Cref{remark: poisson duality for twisted higgs}, the same argument holds verbatim if we replace $G$ by an isotrivial group scheme $\cG$. Hence, %under the assumption that $\cL = \omega_{C}(D)$ for some reduced effective divisor $D \geq 0$ with $\deg(\cL)>0$, we get that 
    the group scheme $\schP^o=\schP^o(\omega_C(D))\to A$ defined in \S\ref{subsection: group scheme of symmetries twisted higgs} is $\delta$-regular.
\end{remark}

%%%%%%%%%%%%%%%%%%%%%%%%%%%%%%%%%%%%%%%%%%%%%%%%%%%%%%%%%%%%%%%%%%%%%%%%%%%%%%%%%%%%%%%%%

\appendix

\begin{section}{Higgs moduli space in arbitrary characteristic}
   We are not aware of references in the literature for the existence of quasiprojective adequate moduli spaces of semistable Higgs bundles for arbitrary connected reductive groups in positive characteristic. In this section, we include a short stack-theoretic argument, which relies on some mild assumptions on the characteristic of $k$. It follows similar lines as in \cite{herrero2023meromorphic}, which work with a more restricted choice of line bundle. We include the full arguments here for completeness.
    \begin{remark}
        If one is only interested in the existence of a quasiprojective moduli space (not the properness of the Hitchin fibration), then we believe that it is possible to remove the assumptions on the characteristic by modifying the GIT construction in \cite[\S 5]{glss-large-characteristic}. Since the characteristic bounds in this appendix suffice for our purposes, we have decided to provide instead the shorter stack-theoretic argument below.
    \end{remark}

    Set $G_{ab} := G/[G,G]$, which is the maximal torus quotient of $G$. The product of the quotient morphism $q: G \to G_{ab}$ and the adjoint representation $\ad: G \to \GL(\mathfrak{g})$ is a homomorphism of reductive groups $ad \times q: G \to \GL(\mathfrak{g}) \times G_{ab}$. The kernel of this morphism is the center $Z_{[G,G]}$ of the derived subgroup of $G$, and so it is a finite multiplicative group. The homomorphism factors thorugh an isogeny $G \xrightarrow{i} G/Z_{[G,G]} \xrightarrow{j} \GL(\mathfrak{g}) \times G_{ab}$.
    
    The homomorphisms induce, by extension of structure groups, morphisms of stacks
    \[ \ad_* \times q_*: \Bun_{G} \xrightarrow{i_*} \Bun_{G/Z_{[G,G]}} \xrightarrow{j_*} \Bun_{\GL(\mathfrak{g})} \times \Bun_{G_{ab}} = \Bun_{\GL(\mathfrak{g})\times G_{ab}}\]
    Similarly there are induced morphisms on Higgs moduli stacks, which we will denote by the same names
    \[ ad_*\times q_*: \Higgs_{G, \cL} \xrightarrow{i_*} \Higgs_{G/Z_{[G,G]}} \xrightarrow{j_*} \Higgs_{\GL(\mathfrak{g})\times G_{ab}, \cL} \]
    We can factor the first morphism of Higgs moduli stacks as
    \[ i_*: \Higgs_{G, \cL} \xrightarrow{f} \Higgs_{G/Z_{[G,G]}, \cL} \times_{\Bun_{G/Z_{[G,G]}}} \Bun_{G} \xrightarrow{p_1} \Higgs_{G/Z_{[G,G]}, \cL}\]
    \begin{lemma} \label{lemma: factorization adjoint morphism first}
        In the factorization above, $f$ is affine of finite type, and $p_1$ is quasifinite and proper with linearly reductive relative stabilizers. Moreover, if $\text{char}(k)=0$ or if $\text{char}(k)$ does not divide the order of the kernel $Z_{[G,G]}$, then $f$ is a closed immersion.
    \end{lemma}
    \begin{proof}
    To simplify notation, we set $H = G/Z_{[G,G]}$. We have an isogeny $i: G \to H$ with kernel $K=Z_{[G,G]}$ of multiplicative type. The factorization reads
    \[ \Higgs_{G, \cL} \xrightarrow{f} \Higgs_{H, \cL} \times_{\Bun_{H}} \Bun_{G} \xrightarrow{p_1} \Higgs_{H, \cL}  \]

    Note that both of the stacks $\Higgs_{G, \cL} \to \Bun_{G}$ and $\Higgs_{H, \cL} \times_{\Bun_{H}} \Bun_{G} \to \Bun_{G}$ are affine of finite type over $\Bun_{G}$ (the second is a base-change of the affine finite type morphism $\Higgs_{H, \cL} \to \Bun_{H}$). It follows that the morphism of $\Bun_{G}$-stacks $f:\Higgs_{G, \cL} \to \Higgs_{H, \cL} \times_{\Bun_{H}} \Bun_{G}$ is affine.

    The morphism $p_1: \Higgs_{\GL(\mathfrak{g})\times G_{ab}, \cL} \times_{\Bun_{\GL(\mathfrak{g})\times G_{ab}}} \Bun_{G} \to \Higgs_{\GL(\mathfrak{g})\times G_{ab}, \cL}$ is a base-change of the morphism $i_*: \Bun_{G} \to \Bun_{H}$, and so it suffices to show the desired properties for $i_*$. The properness and quasifiniteness follow from \cite[Lem. A.5]{gauged_theta_stratifications}. For the statement about relative stabilizers, fix a geometric $F$-point $p \in \Bun_{G}(C)$ with image $i_*(p)$. The proof of quasicompactness in \cite[Lem. A.5]{gauged_theta_stratifications} shows that the $i_*(p)$-fiber of $i_*: \Bun_{G} \to \Bun_{H}$ admits a finite schematic morphism to a stack of the form $\Bun_{K}(\widetilde{C})$, where $\widetilde{C} \to C_F$ is a smooth projective (possibly disconnected) curve over $C$. In particular the group of automorphisms of $p$ in its $i_*(p)$-fiber embed into the group of automorphisms of the corresponding point in $\Bun_{K}(\widetilde{C})$, which is isomorphic to $K^{\pi_0(\widetilde{C})}$. Since $K$ is of multiplicative type, it follows that the stablizer of $p$ relative to $i_*$ is also of multiplicative type, and hence linearly reductive.

    Finally, suppose that $\text{char}(k)=0$ or that $\text{char}(k)$ does not divide the order of $K$. Both of these hypotheses imply that the kernel $K$ is \'etale, and hence the isogeny $i: G \to H$ induces an injection of Lie algebras 
    $\text{Lie}(i): \mathfrak{g} \hookrightarrow \mathfrak{h}$. In 
    particular, for any $G$-bundle $E$, there is an inclusion of twisted adjoint bundles $\ad(E) \hookrightarrow \ad(i_*(E))$. Let $S$ be a $k$-scheme, and choose a morphism $S \to 
\Higgs_{H, \cL} \times_{\Bun_{H}} \Bun_{G}$ corresponding to a pair $(E, \phi)$ of a $G$-bundle $E$ on $C_S$ and a section $\phi: \cO_{C_S} \to \ad(i_*(E))\otimes \cL|_{C_S}$. The $S$-fiber of $f: 
\Higgs_{G, \cL} \to \Higgs_{H, \cL} \times_{\Bun_{H}} \Bun_{G}$ is the subfunctor of $S$ consisting of those morphisms $T \to S$ such that the pullback $\phi|_{C_T} : \cO_{C_T} \to \ad(i_*(E|_{C_T})) 
\otimes \cL|_{C_T}$ factors through the subbundle $\ad(E|_{C_T})\otimes \cL|_{C_T} \subset \ad(i_*(E|_{C_T})) 
\otimes \cL|_{C_T}$. This is represented by a closed subscheme of $S$, and so we conclude that $f$ is a closed immersion, as desired.
    \end{proof}

    \begin{lemma} \label{lemma: partial moduli space}
        The morphism $\ad_* \times q_*: \Higgs_{G, \cL} \to \Higgs_{\GL(\mathfrak{g}) \times G_{ab}, \cL}$ factorizes as
        \[ \ad_* \times q_*:\Higgs_{G, \cL} \xrightarrow{h} \cB \xrightarrow{g}  \Higgs_{\GL(\mathfrak{g}) \times G_{ab}, \cL}\]
        where $h$ is a good moduli space morphism and $g$ is affine and of finite type.
    \end{lemma}
    \begin{proof}
        In the factorization
        \[ ad_*\times q_*: \Higgs_{G, \cL} \xrightarrow{i_*} \Higgs_{G/Z_{[G,G]}, \cL} \xrightarrow{j_*} \Higgs_{\GL(\mathfrak{g})\times G_{ab}, \cL} \]
        $j_*$ is affine and of finite type. Therefore, to prove the lemma we can replace $\Higgs_{\GL(\mathfrak{g})\times G_{ab}, \cL}$ with $\Higgs_{G/Z_{[G,G]}, \cL}$ and prove that there is a factorization
        \[ i_*:\Higgs_{G, \cL} \xrightarrow{h} \cB \xrightarrow{g}  \Higgs_{G/Z_{[G,G]}, \cL}\]
        For simplicity we use again the notation $H := G/Z_{[G,G]}$. Consider the quasifinite, proper morphism $p_1: \Higgs_{H, \cL} \times_{\Bun_{H}} \Bun_{G} \to \Higgs_{H, \cL}$ (\Cref{lemma: factorization adjoint morphism first}). The base change of $p_1$ with any smooth scheme cover $U \to \Higgs_{H, \cL}$ is a quasifinite proper stack over $U$. By the Keel-Mori theorem \cite[\href{https://stacks.math.columbia.edu/tag/0DUT}{Tag 0DUT}]{stacks-project}, the base-change admits a coarse moduli space, which is quasifinite and proper (and hence finite) over $U$. Since the formation of the coarse space is compatible with flat-base change, this descends to a relative coarse space $\Higgs_{H, \cL} \times_{\Bun_{H}} \Bun_{G} \xrightarrow{c} \cM \to \Higgs_{H, \cL}$, where $\cM \to \Higgs_{H, \cL}$ is finite and schematic. Moreover, since the relative stabilizers of $p_1$ are linearly reductive, it follows that $\Higgs_{H, \cL} \times_{\Bun_{H}} \Bun_{G} \xrightarrow{c} \cM$ is actually a good moduli space morphism. Since the morphism $f: \Higgs_{G, \cL} \to \Higgs_{H, \cL} \times_{\Bun_{H}} \Bun_{G}$ is affine and of finite type, by \cite[Lem. 4.14+Thm. 4.16(xi)]{alper-good-moduli} we have a factorization 
        \[\Higgs_{G, \cL} \xrightarrow{h} \cB \to \cM \to \Higgs_{H, \cL}\]
        where $\cB \to \cM$ is affine and of finite type and $\Higgs_{G, \cL} \xrightarrow{h} \cB$ is a good moduli space morphism. The composition $g: \cB \to \cM \to \Higgs_{H, \cL}$ of two affine morphisms of finite type will also be affine and of finite type.
    \end{proof}
    
    \begin{lemma} \label{lemma: semistability under adjoint representation}
        Suppose that $\text{char}(k)=0$ or that it  is strictly larger than the height of the adjoint representation $\mathfrak{g}$. Then for any algebraically closed field $F \supset k$ and any $F$-point $p = (E, \phi) \in \Higgs_{G, \cL}$, we have that $p$ is semistable if and only if the associated Higgs bundle $ad_*(p)$ for the group $\GL(\mathfrak{g})$ is semistable.
    \end{lemma}
    \begin{proof}
   After replacing $k$ with $F$, we can assume without loss of generality that $p$ is defined over $k$. The homomorphism $\ad$ factors as $\ad: G \xrightarrow{u} G/Z_G \xrightarrow{\ad} \GL(\mathfrak{g})$. For any $G$-Higgs bundle $p=(E, \phi)$, there is a canonical bijection between parabolic reductions of $E$ and parabolic reductions of $u_*(E)$. This bijection respects compatibility with $\phi$ and $u_*(\phi)$ respectively, and it follows from the definition that the numerical condition for semistability of $p$ and $u_*(p)$ are identified. Hence $p$ is semistable if and only if $u_*(p)$ is semistable, and so we can replace $G$ with $G/Z_G$ and assume without loss of generality that $G$ is semisimple and $\ad$ is an embedding.

    Suppose that $\ad_*(p)$ is semistable. Then the fact that $\ad$ is an embedding implies (by a similar argument as in \cite[Prop. 5.3]{biswas-holla-hnreduction}) that $p$ is semistable. On the other hand, if $p$ is semistable, then $\ad_*(p)$ is semistable because $G$ is semisimple and the assumptions on the characteristic of $k$ (cf. \cite[Thm. 7.5]{balaji-parameswaran-tensor} and \cite[Thm. 8.17]{balaji-parameswaran-tensor}).
    \end{proof}

    \begin{thm} \label{thm: adequate moduli spaces for higgs}
        Suppose that $\text{char}(k)=0$ or is strictly larger than the height of the adjoint representation $\mathfrak{g}$. Then,
        \begin{enumerate}[(a)]
            \item The set of semistable geometric points of $\Higgs_{G,\cL}$ are exactly the geometric points of an open substack $\Higgs_{G, \cL}^{ss}$.
            \item For any $d \in \pi_1(G)$, the open and closed substack $\Higgs^{d,ss}_{G, \cL} \subset \Higgs_{G, \cL}^{ss}$ of semistable Higgs bundles of degree $d$ admits an adequate moduli space $\MHiggs_{G, \cL}^{d}$ which is quasiprojective over $k$. 
            % \item If $\text{char}(k)=0$ or $\text{char}(k)$ does not divide the order of $Z_{[G,G]}$, then the induced Hitchin morphism $\MHiggs_{G, \cL}^{d} \to A$ is projective.
        \end{enumerate}
    \end{thm}
    \begin{proof}
    If $G = \GL_n$ for some positive integer $n$, then the theorem is true by the GIT constructions of the moduli space of $\Lambda$-modules \cite{Simpson-repnI, langer-moduli-lie-algebroids}. We use the morphism $\ad_* \times q_*$ to extend the result to arbitrary $G$.
    
\noindent \textit{Proof of (a):} By the result for the general linear group $\GL(\mathfrak{g})$, the semistable geometric points of $\Higgs_{\GL(\mathfrak{g}), \cL}$ are exactly the points of an open substack $\Higgs_{\GL(\mathfrak{g}), \cL}^{ss} \subset \Higgs_{\GL(\mathfrak{g}), \cL}$. By \Cref{lemma: semistability under adjoint representation}, it follows that the semistable geometric points of $\Higgs_{G, \cL}$ are the points of the open preimage $(\ad_*\times q_*)^{-1}(\Higgs_{\GL(\mathfrak{g}), \cL}^{ss}) \subset \Higgs_{G, \cL}$.

 \medskip

 \noindent \textit{Proof of (b):} Consider the factorization from \Cref{lemma: partial moduli space}:
 \[ \ad_* \times q_*:\Higgs_{G, \cL} \xrightarrow{h} \cB \xrightarrow{g}  \Higgs_{\GL(\mathfrak{g}) \times G_{ab}, \cL}\]
 We have an open substack $\Higgs_{\GL(\mathfrak{g}), \cL}^{ss} \times \Higgs_{G_{ab}, \cL} \subset \Higgs_{\GL(\mathfrak{g}) \times G_{ab}, \cL}$, with open preimage $\mathcal{U} := g^{-1}(\Higgs_{\GL(\mathfrak{g}), \cL}^{ss} \times \Higgs_{G_{ab}, \cL})$. By \Cref{lemma: semistability under adjoint representation}, we have 
 \[h^{-1}(\mathcal{U}) = (\ad_*\times q_*)^{-1}(\Higgs_{\GL(\mathfrak{g}), \cL}^{ss} \times \Higgs_{G_{ab}, \cL}) = \Higgs_{G, \cL}^{ss} \subset \Higgs_{G, \cL}.\]
 Hence we get a factorization
 \[ \Higgs_{G, \cL}^{ss} \xrightarrow{h} \mathcal{U} \xrightarrow{g} \Higgs_{\GL(\mathfrak{g}), \cL}^{ss} \times \Higgs_{G_{ab}, \cL}\]
 where $h$ is a good moduli space morphism and $g$ is affine and of finite type.

Note that $G_{ab}$ is a split torus, since $G$ is split. In particular $G_{ab} \cong \GL_1^r$ for some $r$, and we have $\Higgs_{G_{ab}, \cL} \cong \left(\Higgs_{\GL_1, \cL}\right)^r$. The homomorphism $\ad\times q: G \to \GL(\mathfrak{g}) \times \GL_1^r$ induces a morphism of fundamental groups $\pi_1(G) \to \pi_1(\GL(\mathfrak{g})) \times \pi_1((\GL_1)^r) = \mathbb{Z}^{r+1}$. Let $d'= (d_0, d_1, \ldots, d_r)$ denote the image of $d$ under this homomorphism. For Higss bundles of rank $1$ we have $\Higgs_{\GL_1, \cL}^{ss} = \Higgs_{\GL_1, \cL}$, since there are no proper parabolic subgroups. It follows that
\[ \Higgs_{\GL(\mathfrak{g}) \times G_{ab}, \cL}^{ss} \times \Higgs_{G_{ab}, \cL}^{ss} = \Higgs_{\GL(\mathfrak{g}), \cL}^{ss} \times \left( \Higgs_{\GL_1 , \cL}^{ss}\right)^r = \Higgs_{\GL(\mathfrak{g}), \cL}^{ss} \times \left( \Higgs_{\GL_1 , \cL}\right)^r\]
Furthermore, the result for general linear groups implies that the open and closed substack 
\[\Higgs_{\GL(\mathfrak{g})\times G_{ab}, \cL}^{d',ss} = \Higgs_{\GL(\mathfrak{g}), \cL}^{d_0, ss} \times \prod_{i=1}^r \Higgs_{\GL_1 , \cL}^{d_i, ss}\]
admits a quasiprojective adequate moduli space $\MHiggs_{\GL(\mathfrak{g})\times G_{ab}, \cL}^{d'}$ over $k$. Let us denote $\mathcal{U}^{d'}:= g^{-1}(\Higgs_{\GL(\mathfrak{g})\times G_{ab}, \cL}^{ss, d'})$, which is an open and closed substack of $\mathcal{U}$. We have a chain or morphisms
\[ \Higgs_{G, \cL}^{d, ss} \hookrightarrow h^{-1}(\mathcal{U}^{d'}) \xrightarrow{h} \mathcal{U}^{d'} \xrightarrow{g} \Higgs_{\GL(\mathfrak{g})\times G_{ab}, \cL}^{d',ss}\]
The left-most inclusion $\Higgs_{G, \cL}^{d, ss} \hookrightarrow h^{-1}(\mathcal{U}^{d'})$ is an open and closed immersion. Therefore, the image $\mathcal{U}^d := h(\Higgs_{G, \cL}^{d, ss})$ is an open and closed substack of $\mathcal{U}^{d'}$, and the induced morphism $h: \Higgs_{G, \cL}^{d, ss} \to \mathcal{U}^d$ is a good moduli space morphism. We conclude that we have a factorization
\[\Higgs_{G, \cL}^{d, ss} \xrightarrow{h} \mathcal{U}^d \xrightarrow{g}  \Higgs_{\GL(\mathfrak{g})\times G_{ab}, \cL}^{d',ss}\]
where $h$ is a good moduli space morphism and $g$ is affine of finite type. Since $\Higgs_{\GL(\mathfrak{g})\times G_{ab}, \cL}^{d',ss}$ admits an adequate moduli space $\Higgs_{\GL(\mathfrak{g})\times G_{ab}, \cL}^{d',ss} \to \MHiggs_{\GL(\mathfrak{g})\times G_{ab}, \cL}^{d'}$, by \cite[Lem. 5.2.11+Thm. 6.3.3]{alper_adequate} it follows that $\mathcal{U}^{d}$ admits an adequate moduli space $U$ that is affine and of finite type over the quasiprojective scheme $\MHiggs_{\GL(\mathfrak{g})\times G_{ab}, \cL}^{d'}$. In particular $U$ is quasiprojective over $k$. The composition $\Higgs_{G, \cL}^{d, ss} \xrightarrow{h} \mathcal{U}^d \to U$ of two adequate moduli space morphisms is an adequate moduli space morphism, and hence $\MHiggs_{G, \cL}^{d} := U$ is an adequate moduli space for $\Higgs_{G, \cL}^{d, ss}$.

%  \medskip 

%  \noindent \textit{Proof of (c):}
%      Suppose that $\text{char}(k)=0$ or that $\text{char}(k)$ does not divide the order of $Z_{[G,G]}$. Then have shown in part (b) that $\Higgs_{G, \cL}^{d, ss}$ admits a quasiprojective adequate moduli space $\MHiggs_{G, \cL}^{d}$ such that the naturally induced morphism $g: \MHiggs_{G, \cL}^{d} \to \MHiggs_{\GL(\mathfrak{g})\times G_{ab}, \cL}^{d'}$ if finite. This fits into a commutative diagram, where the vertical morphism are the Hitchin fibrations:
%      \begin{figure}[H]
% \centering
% \begin{tikzcd}
%   \MHiggs_{G, \cL}^{d} \ar[d, "h_G"] \ar[r, "g"]  & \MHiggs_{\GL(\mathfrak{g})\times G_{ab}, \cL}^{d'} \ar[d, "h_{\GL(\mathfrak{g}) \times G_{ab}}"] \\  
%   A_{G} \ar[r, "i_{A}"] & A_{\GL(\mathfrak{g}) \times G_{ab}}
% \end{tikzcd}
% \end{figure}
% By the properness of the Hitchin fibration for general linear groups $\GL_n$ (\andres{add reference}), it follows that $h_{\GL(\mathfrak{g}) \times G_{ab}}$ is proper. Since $g$ is finite, the composition $h_{\GL(\mathfrak{g}) \times G_{ab}} \circ g = \iota_{A} \circ h_G$ is proper. Notice that $\iota_{A}$ is separated, because it is a morphism of affine $k$-schemes, and therefore it follows from right cancellation that the Hitchin morphism $h_G$ is proper.
    \end{proof}
    
\end{section}

\end{subsection}
\end{section}

\begin{section}{Ancillary results on the moduli of Higgs bundles}

In this appendix, we collect several auxiliary results, used in the proof of some of the main results above, concerning the functoriality of the construction of Higgs stacks or moduli spaces as the structural group varies.

\begin{lemma}\label{L:Higgs products} Suppose that the commutative diagram of reductive groups 
\[
\xymatrix{
H\ar[d]\ar[r]& G_1\ar[d]^{f}\\
G_2\ar[r]& G
}
\]
is Cartesian, and $f$ is faithfully flat. Then, the induced morphism of stacks
\[
\xymatrix{
\Higgs_{H,\cL}\ar[d]\ar[r]& \Higgs_{G_1,\cL}\ar[d]^{f_*}\\
\Higgs_{G_2,\cL}\ar[r]& \Higgs_{G,\cL}
}
\]
is Cartesian. %Moreover, if $f$ is central, then the restriction $f_*:\Higgs^{d_1}_{G_1,\cL}\to \Higgs^{d}_{G,\cL}$ is fppf, for any $d_1\in\pi_1(G_1)$ and $d:=\pi_1(f)\in\pi_1(G)$.
\end{lemma}

\begin{proof}The statement about the underlying principal bundles follows by arguing as in the proof of \cite[Lem. 2.2.1]{biswas-hoffmann-line-bundles}. The statement about Higgs fields follows from the fact the induced diagram of Lie algebras \[
\xymatrix{
\Lie(H)\ar[d]\ar[r]& \Lie(G_1)\ar[d]\\
\Lie(G_2)\ar[r]& \Lie(G)
}
\]
is Cartesian.
%\Mirko{@Note for myself: I do not see the implication} If $\ker(f)=\mu$ is a commutative linear algebraic group, then there is a Cartesian diagram of stacks
%\[
%\xymatrix{
%\Higgs_{G_1,\cL}^{d_1}\ar[r]^{f_*}\ar[d]&\Higgs^d_{G,\cL}\ar[d]\\
%H^0(C,\Lie(\mu)\otimes\cL)\times \Bun^{d_1}_{G_1}\ar[r]& \Bun^d_G,
%}
%\]
%where the lower horizontal arrow sends a pair $(\sigma,E)$ to the $G$-bundle $E\times^{G_1}G$. Since the lower arrow is a fppf morphism \cite[Lemma 2.2.2]{biswas-hoffmann-line-bundles}, the base-change $f_*$ is fppf.
\end{proof}

\begin{lemma}\label{L:higgs isogeny fppf}Let $f:G_1\to G_2$ be a faithfully-flat homomorphism of reductive groups over $k$ with central kernel. Fix $d_1\in\pi_1(G_1)$, and let $d_2$ be the image of $d_1$ via the morphism $\pi_1(G_1)\to \pi_1(G_2)$. Then, the induced morphism of stacks $f_*:\Higgs_{G_1,\cL}^{d_1}\to \Higgs_{G_2,\cL}^{d_2}$ is fppf.
\end{lemma}

\begin{proof}Set $\mu:=\ker(f)$. By \Cref{R:decom-pi}, there exists a $G$-equivariant splitting of Lie algebras  $s:\Lie(G_1)\cong\Lie(D_{G_1})\oplus\Lie(Z_{G_1})$, where $D_{G_1} = [G_1, G_1] \subset G_1$ denotes the derived subgroup. Since $f$ is faithfully-flat, $\mu$ is contained in the center of $G_1$, and $G_1$ acts trivially on $\Lie(Z_{G_1})$, the splitting $s$ induces a $G$-equivariant splitting of Lie algebras $s':\Lie(G_1)\cong\Lie(G_2)\oplus\Lie(\mu)$. 
In particular, for any $G_1$-bundle $E$, the surjective morphism of Lie algebras $\Lie(G_1)\to \Lie(\mu)$, induced by the splitting $s'$, gives a morphism of linear spaces $$H^0(C,\ad(E)\otimes\cL)=H^0(C,E(\Lie(G_1))\otimes\cL)\xrightarrow{p_E}H^0(C,E(\Lie(\mu))\otimes\cL)=H^0(C,\Lie(\mu)\otimes\cL),$$
where the last equality comes from the fact that $G_1$ acts trivially on $\Lie(\mu)$. Then, we have a Cartesian diagram of stacks
\[
\xymatrix{
\Higgs_{G_1,\cL}^{d_1}\ar[r]^{f_*}\ar[d]^{(\zeta,\pi_{G_1})}&\Higgs^d_{G_2,\cL}\ar[d]^{\pi_{G_2}}\\
H^0(C,\Lie(\mu)\otimes\cL)\times \Bun^{d_1}_{G_1}\ar[r]& \Bun^{d_2}_{G_2},
}
\]
where $\pi_{G_i}(E,\theta)=E$, $\zeta(E,\theta)= p_E(\theta)\in H^0(C,\Lie(\mu)\otimes\cL)$, and the lower horizontal arrow sends a pair $(\sigma,E)$ to the $G_2$-bundle $E\times^{G_1}G_2$. Since the lower arrow is a fppf morphism \cite[Lem. 2.2.2]{biswas-hoffmann-line-bundles}, the base-change $f_*$ is fppf.
\end{proof}

\begin{coroll}\label{R:G to G/Z Higgs}
Let $Z_G^o$ be the maximal central torus of $G$, with quotient map $\pi:G\to G/Z_G^o$. Fix $d\in\pi_1(G)$ and let $e$ be the image of $d$ via the natural morphism of fundamental groups $ \pi_1(G)\to\pi_1(G/Z_G^o)$. 
The induced morphism of stacks $\pi_*:\Higgs_{G,\cL}^d\to\Higgs^e_{G/Z^o_G,\cL}$ is smooth and its geometric fibers are non-canonically isomorphic to the moduli stack $\Higgs^0_{Z_G^o,\cL}\cong H^0(C,\Lie(Z_G)\otimes\cL)\times \Bun_{Z_G^o}^0$ of $Z_G^o$-Higgs bundles of degree $0\in\pi_1(Z_G^o)$.    
\end{coroll}

\begin{proof} Consider the Cartesian diagram 
\[
\xymatrix{
Z_G^o\times G\ar[d]_-{pr_G}\ar[r]^-{m}& G\ar[d]^-{\pi}\\
G\ar[r]^-{\pi}& G/Z_G^o
}
\]
where $pr_G$ is the projection on the second factor, $m$ is the multiplication morphism $(z,g)\mapsto zg$. Observe that the quotient $\pi$ is fppf and central. The statement follows then from  \Cref{L:Higgs products}, \Cref{L:higgs isogeny fppf} and the smoothness of $\Higgs^0_{Z_G^o,\cL}$. 
\end{proof}

\begin{context}For the rest of this section, we always assume that either $\text{char}(k) =0$, or $\text{char}(k)$ is strictly larger than the height of the adjoint representation. In this way, the stack $\Higgs_{G,\cL}^{d,ss}$ of semistable $G$-Higgs bundles admits an adequate moduli space $M_{G,\cL}$ \Cref{thm: adequate moduli spaces for higgs}.
\end{context}

\begin{lemma}\label{L:induced morphism semistable} Let $f:G_1\to G_2$ be a morphism of reductive groups over $k$ and let $f_*:\Higgs_{G_1,\cL}\to \Higgs_{G_2,\cL}$ be the induced morphism of stacks. 

\begin{enumerate}[(i)]
    \item\label{L:induced morphism semistable 1} Assume that $f$ is central and that the induced morphism of adjoint groups $G_1/Z_{G_1}\to G_2/Z_{G_2}$ is an isomorphism (e.g. $f$ is faitfully-flat). Then, a pair $(E,\theta)$ in $\Higgs_{G_1,\cL}$ is (semi)stable if and only if its image $f_*(E,\theta)$ in $\Higgs_{G_2,\cL}$ is (semi)stable.
    \item\label{L:induced morphism semistable 2} Assume that $G_2$ is a torus. Then, for any pair $(E,\theta)$ in $\Higgs_{G_1,\cL}$, its image $f_*((E,\theta))$ in $\Higgs_{G_2,\cL}$ is stable.
\end{enumerate}
 In particular, if either \eqref{L:induced morphism semistable 1} or \eqref{L:induced morphism semistable 2} hold, then the morphism $f$ induces a morphism $f_*:\MHiggs_{G_1,\cL}\to \MHiggs_{G_2,\cL}$ of adequate moduli spaces.
\end{lemma}

\begin{proof}Point \eqref{L:induced morphism semistable 2} follows from the fact that the a torus does not have proper parabolic subgroups. We now focus on Point \eqref{L:induced morphism semistable 1}. In the proof of  \Cref{lemma: semistability under adjoint representation}, we have seen that a $G$-Higgs bundle is (semi)stable if and only if the induced $(G/Z_G)$-Higgs bundle is (semi)stable. Hence, it is enough to check the statement for the induced morphism $G_1/Z_{G_1}\to G_2/Z_{G_2}$. By assumptions, the latter morphism is an isomorphism and, so, we have the statement.
\end{proof}

\begin{lemma}\label{L:isogeny-finite}Let $f:G_1\to G_2$ be an isogeny of reductive groups over $k$ with \'etale kernel. Fix $d_1\in\pi_1(G_1)$, and let $d_2$ be the image of $d_1$ via the morphism $\pi_1(H)\to\pi_1(G)$. Then, the induced morphism $\MHiggs_{G_1,\cL}^{d_1}\to \MHiggs_{G_2,\cL}^{d_2}$ of adequate moduli spaces (cf. \Cref{L:induced morphism semistable}\eqref{L:induced morphism semistable 1}) is finite and surjective.
\end{lemma}

\begin{proof}By \Cref{L:higgs isogeny fppf}, the induced morphism of stacks $f_*:\Higgs_{G_1, \cL}^{d_1} \to \Higgs_{G_2, \cL}^{d_2}$ is fppf, and so in particular it is surjective. \Cref{L:induced morphism semistable}\eqref{L:induced morphism semistable 1} implies that $(f_*)^{-1}(\Higgs_{G_2,\cL}^{d_2, ss}) = \Higgs_{G_1, \cL}^{d_1, ss}$, and hence the induced morphism $f_*: \Higgs_{G_1, \cL}^{d_1, ss} \to \Higgs_{G_2, \cL}^{d_2, ss}$ is also surjective. It follows that the morphism of adequate moduli spaces $\MHiggs_{G_1,\cL}^{d_1} \to \MHiggs_{G_2, \cL}^{d_2}$ is also surjective.

We are left to show that the morphism of adequate moduli spaces is finite. We may factor the morphism $f_*: \Higgs_{G_1,\cL} \xrightarrow{h} \Higgs_{G_2, \cL} \times_{\Bun_{G_2}} \Bun_{G_1} \xrightarrow{p_1} \Higgs_{G_2, \cL}$ similarly as in \Cref{lemma: factorization adjoint morphism first}. Using the assumption that the kernel of the isogeny $f$ is \'etale, the argument in \Cref{lemma: factorization adjoint morphism first} applies verbatim to show that $h$ is a closed immersion and $p_1$ is quasifinite and proper, and hence $f_*$ is quasifinite and proper. Using $(f_*)^{-1}(\Higgs_{G_2,\cL}^{d_2, ss}) = \Higgs_{G_1, \cL}^{d_1, ss}$, we get that $f_*: \Higgs_{G_1, \cL}^{d_1, ss} \to \Higgs_{G_2, \cL}^{d_2, ss}$ is quasifinite and proper. Similarly as in the proof of \Cref{lemma: partial moduli space}, we may define a factorization $f_*: \Higgs_{G_1, \cL}^{d_1, ss} \to \cM \to \Higgs_{G_2, \cL}^{d_2, ss}$, where $ \Higgs_{G_1, \cL}^{d_1, ss} \to \cM$ is a good moduli space morphism and $\cM \to \Higgs_{G_2, \cL}^{d_2, ss}$ is schematic and finite. By \cite[Lem. 5.2.11+Thm. 6.3.3]{alper_adequate}, the stack
 $\cM$ admits an adequate moduli $M$ space such that the induced morphism $M \to \MHiggs_{G_2, \cL}^{d_2}$ is finite. Similarly as in the argument of \Cref{lemma: partial moduli space}, we have that the good moduli space morphism $\Higgs_{G_1, \cL}^{d_1, ss} \to \cM$ induces an isomorphism of adequate moduli spaces $\MHiggs_{G_1, \cL}^{d_1} \xrightarrow{\sim} M$, and hence we conclude that $\MHiggs_{G_1, \cL}^{d_1} \to \MHiggs_{G_2, \cL}^{d_2}$ is finite.
\end{proof}

\begin{defn}\label{N:twisted higgs} Let $G$ be a reductive group over $k$, with derived subgroup $D_G$ and abelianization $G_{ab} \coloneqq G/D_G$. By \Cref{L:induced morphism semistable}\eqref{L:induced morphism semistable 2}, the quotient $\nu \colon G\to G_{ab}$ induces a morphism \[\nu_* \colon \MHiggs_{G,\cL}\to \MHiggs_{G_{ab},\cL}.\] For any $k$-point $(E,\theta)$ in $\MHiggs_{G_{ab},\cL}$, we denote by $\MHiggs_{G,\cL}^{(E,\theta)}$ the fiber over $(E,\theta)$ of the morphism $\nu_*$.
%$\MHiggs_{G,\cL}\to \MHiggs_{G/D_G,\cL}$ induced by the quotient morphism $G\to G/D_G$.  
For any $d\in\pi_1(G)$, we set $\MHiggs_{G,\cL}^{d,(E,\theta)}:=M_{G,\cL}^{(E,\theta)}\cap \MHiggs_{G,\cL}^d$. 

Similarly, we define $\Higgs_{G, \cL}^{d, (E, \theta)}$ to be the fiber over $(E, \theta)$ of the morphism $\nu_*: \Higgs_{G, \cL}^d \to \Higgs_{G_{ab},\cL}$ and set $\Higgs_{G, \cL}^{d, (E, \theta),ss} := \Higgs_{G, \cL}^{d, (E, \theta)}\times_{\Higgs_{G, \cL}^d} \Higgs_{G, \cL}^{d,ss}$.
\end{defn}

\begin{prop} \label{prop: moduli space commutes with base change on abelianization}
    With notation as above, let $X$ be a scheme, and choose a morphism $X \to \Higgs_{G_{ab}, \cL}$. Then, the induced morphism of fiber products $X \times_{\Higgs_{G_{ab}, \cL}} \Higgs_{G, \cL}^{d, ss} \to X \times_{\MHiggs_{G_{ab}, \cL}} \MHiggs_{G, \cL}^d$ is an adequate moduli space.
\end{prop}
\begin{proof}
    The assertion of the proposition can be checked fpqc locally on $X$. We may extend scalars and assume that $k$ is algebraically closed. Consider the Cartesian diagram of reductive groups
\begin{equation}\label{E:Cartesian product groups locally}
    \xymatrix{
    G\times Z_G^o\ar[d]_-{(\nu,\mathrm{Id}_{Z_G^o})}\ar[r]^-m& G\ar[d]^-{\nu}\\
    G_{ab}\times Z_G^o\ar[r]^-{n}& G_{ab},
    }
\end{equation}
where $m(g, z)=gz$ and $n(t,z)=t\nu(z)$.

Fix a $k$-point $(E,\theta)$ in $\Higgs_{G_{ab},\cL}$. By \Cref{L:Higgs products}, the induced diagram of stacks
\[
    \xymatrix{
    \Higgs_{G,\cL}^{d,(E,\theta)}\times \Higgs^0_{Z_G^o,\cL}\ar[r]\ar[d]&\Higgs^d_{G,\cL}\times \Higgs^0_{Z_G^o,\cL}\ar[d]_-{(\nu_*,\mathrm{Id})}\ar[r]^-{m_*}& \Higgs^d_{G,\cL}\ar[d]^{\nu_*}\\
    \{(E,\theta)\}\times \Higgs_{Z_G^o,\cL}^0\ar[r]&\Higgs_{G_{ab},\cL} \times \Higgs_{Z_G^o,\cL}\ar[r]& \Higgs_{G_{ab},\cL}
    }
\]
is Cartesian. Moreover, since $m$ is faithfully-flat and with central kernel, we have $(m_*)^{-1}(\Higgs_{G,\cL}^{d,ss})=\Higgs^{d,ss}_{G,\cL}\times \Higgs^0_{Z_G^o,\cL}$, cf. \Cref{L:induced morphism semistable}\eqref{L:induced morphism semistable 1}. If we denote by $\nu_*(d)$ the image of $d$ in $\pi_1(G_{ab})$, then we may replace $\Higgs_{G_{ab}, \cL}$ with open and closed substack $\Higgs_{G_{ab}, \cL}^{\nu_*(d)}$. Hence, we get the Cartesian diagram of stacks
\begin{equation}\label{E:Cartesian product higgs locally}
    \xymatrix{
    \Higgs_{G,\cL}^{d,(E,\theta),ss}\times \Higgs^0_{Z_G^o,\cL}\ar[r]\ar[d]& \Higgs^{d,ss}_{G,\cL}\ar[d]^{\nu_*}\\
    \{(E,\theta)\}\times \Higgs^0_{Z_G^o,\cL}\ar[r] & \Higgs_{G_{ab},\cL}^{\nu_*(d)},
    }
\end{equation}
where the horizontal arrows are fppf and surjective. 
Under our assumption that $k$ is algebraically closed, there is a $k$-point $c \in C(k)$ which induces a splitting of the good moduli space morphism $\Bun_{\mathbb{G}_m} \to \Pic(C)$. Note that $\Bun_{Z_G^{\circ}} \cong \left(\Bun_{\mathbb{G}_m}\right)^r$ for some nonnegative integer $r$, and hence we also get an induced splitting of the good moduli space morphism $\Bun_{Z^{\circ}_G} \to \Pic(C)^r$. The Higgs stack $\Higgs_{Z_G^{\circ}, \cL} \to \Bun_{Z_G^{\circ}}$ is isomorphic to the total space of a vector bundle on $\Bun_{Z_G^{\circ}}$ which descends to the moduli space $\Pic^r(C)$. If follows that the splitting of $\Bun_{Z^{\circ}_G} \to \Pic(C)^r$ also induces a splitting of $\Higgs_{Z_G^{\circ}, \cL} \to \MHiggs_{Z_G^{\circ}, \cL}$, which is surjective and fppf. Applying a further base-change of diagram \eqref{E:Cartesian product higgs locally} under the section $\MHiggs_{Z_G^{\circ}, \cL}^{\circ} \to \Higgs_{Z_G^{\circ}, \cL}^0$, we get a Cartesian diagram where the horizontal arrows are surjective and fppf 
\begin{equation} \label{E:Cartesian product higgs locally 2}
    \xymatrix{
    \Higgs_{G,\cL}^{d,(E,\theta),ss}\times \MHiggs^0_{Z_G^o,\cL}\ar[r]\ar[d]& \Higgs^{d,ss}_{G,\cL}\ar[d]^{\nu_*}\\
    \{(E, \theta)\} \times \MHiggs_{Z_G^o,\cL}^0\ar[r]& \Higgs_{G_{ab},\cL}^{\nu_*(d)}
    }
\end{equation}
We may assume without loss of generality that $X \to \Higgs_{G_{ab}, \cL}$ factors through $\Higgs_{G_{ab}, \cL}^{\nu_*(d)}$. Furthermore, up to passing to an fppf cover of $X$, we may assume that $X \to \Higgs_{G_{ab}, \cL}^{\nu_*(d)}$ factors through $\{(E, \theta)\} \times \MHiggs_{Z_G^o,\cL}^0$. We get a Cartesian diagram
\begin{equation} \label{E:Cartesian product higgs locally 3}
    \xymatrix{
    X \times_{\Higgs_{G_{ab}, \cL}} \Higgs^{d, ss}_{G, \cL} \ar[r] \ar[d] & \Higgs_{G,\cL}^{d,(E,\theta),ss}\times \MHiggs^0_{Z_G^o,\cL}\ar[d]\\
    X \ar[r] & \{(E, \theta)\} \times \MHiggs_{Z_G^o,\cL}^0
    }
\end{equation} 
Now, in view of diagrams \eqref{E:Cartesian product higgs locally 2} and \eqref{E:Cartesian product higgs locally 3}, the proposition is reduced to showing the following.

\noindent \textbf{Claim.} Both of the following statements (a) and (b) hold.
\begin{enumerate}[(a)]
    \item The induced morphism $\Higgs_{G, \cL}^{d, (E, \theta), ss} \to \MHiggs_{G, \cL}^{(E, \theta)}$ is an adequate moduli space morphism. Equivalently, the following diagram of schemes (obtained by taking moduli spaces of diagram \eqref{E:Cartesian product higgs locally 2}) is Cartesian
    \[
\xymatrix{
\MHiggs_{G,\cL}^{d,(E, \theta)}\times \MHiggs_{Z_G^o,\cL}^0\ar[r]\ar[d]&\MHiggs^d_{G,\cL}\ar[d]\\
\MHiggs_{Z_G^o,\cL}^0\ar[r]&\MHiggs_{G_{ab},\cL}.
}
\]
    \item The formation of the adequate moduli space of $\Higgs_{G,\cL}^{d,(E,\theta),ss}\times \MHiggs^0_{Z_G^o,\cL}$ commutes with arbitrary base-change on $\{(E, \theta)\} \times \MHiggs_{Z_G^o,\cL}^0$. In particular, the square obtained by taking adequate moduli spaces in diagram \eqref{E:Cartesian product higgs locally 3} remains Cartesian.
\end{enumerate}

%For part (a), it is enough to show the first sentence (the induced morphism $\Higgs_{G, \cL}^{d, (E, \theta), ss} \to \MHiggs_{G, \cL}^{(E, \theta)}$ is an adequate moduli space). 
Part (a) has been proven in the context of the moduli of $G$-bundles in \cite[Lem. 5.1]{distinguish_algebraicspaces_schemes}. We note that the same argument can be adapted to the case of $G$-Higgs bundles, we omit the details.

Part (b) follows because the morphism  $\Higgs_{G,\cL}^{d,(E,\theta),ss}\times \MHiggs^0_{Z_G^o,\cL}\to \{(E, \theta)\} \times \MHiggs_{Z_G^o,\cL}^0$ is a base-change of $\Higgs_{G,\cL}^{d,(E,\theta),ss} \to \Spec(k)$, and the formation of the adequate moduli space of $\Higgs_{G,\cL}^{d,(E,\theta),ss}$ commutes with arbitrary (flat) base-change on $\Spec(k)$.
\end{proof}

Note that when $D_G$ is simply-connected, we have $\pi_1(G)=\pi_1(G_{ab})$, so in that case there exists a unique $d\in\pi_1(G)$ such that $\MHiggs_{G,\cL}^{(E,\theta)}=\MHiggs_{G,\cL}^{d,(E,\theta)}$. \Cref{L:twisted locally trivial} states that the morphism $\MHiggs_{G,\cL}\to \MHiggs_{G_{ab},\cL}$ is isotrivial, so in particular the isomorphism type of $M_{G,\cL}^{d,(E,\theta)}$ is independent of $(E, \theta)$.
\begin{lemma} \label{L:twisted locally trivial}
    Let $(E, \theta)$ be a $k$-point of $\Higgs_{G_{ab}, \cL}$. The induced morphism $\Higgs_{G, \cL}^{d, (E, \theta),ss} \to \MHiggs_{G, \cL}^{d, (E, \theta)}$ is an adequate moduli space. Furthermore, there is a Cartesian diagram of schemes
\[
\xymatrix{
\MHiggs_{G,\cL}^{d,(E, \theta)}\times \MHiggs_{Z_G^o,\cL}^0\ar[r]\ar[d]&\MHiggs^d_{G,\cL}\ar[d]\\
\MHiggs_{Z_G^o,\cL}^0\ar[r]&\MHiggs_{G_{ab},\cL}.
}
\]
\end{lemma}
\begin{proof}
    This is part (a) in the claim shown within the proof of \Cref{prop: moduli space commutes with base change on abelianization}. 
\end{proof}

\begin{notn}\label{N:D(G) simply-connected} When $D_G$ is simply-connected, the derived subgroup $D_G$ decomposes into a product $H_1\times \cdots \times H_s$, where each group $H_i$ is simply-connected and almost-simple. Write $G_i = G / \prod_{i\neq j}H_j$. %We denote by $G_i$ the quotient of $G$ by the normal subgroup $\prod_{i\neq j}H_j$.
Note that the morphism $\pi_i:G\to G_i$ induces an isomorphism of tori $G_{ab}=G/D_G\cong G_i/D_{G_i})$.
\end{notn}

\begin{lemma}\label{L:twisted split} With notation as in \ref{N:twisted higgs} and \ref{N:D(G) simply-connected}, there exists an isomorphism of moduli spaces
$
\MHiggs_{G,\cL}^{(E, \theta)}\cong \MHiggs_{G_1,\cL}^{(E, \theta)}\times\cdots\times \MHiggs_{G_s,\cL}^{(E, \theta)}.
$
\end{lemma}

\begin{proof}Consider the Cartesian diagram of reductive groups
\begin{equation}\label{E:Cartesian product groups}
    \xymatrix{
    G\ar[d]\ar@{^{(}->}[rr]^-{(\pi_1,\ldots,\pi_s)}&& G_1\times\cdots\times G_s\ar[d]\\
    G_{ab}\ar@{^{(}->}[rr]^{\Delta}&& (G_{ab})^{s}
    }
\end{equation}
where the vertical arrows are the abelianization morphisms and $\Delta$ is the diagonal embedding. By \Cref{L:Higgs products}, the induced diagram of stacks
\begin{equation}\label{E:Cartesian product higgs}
    \xymatrix{
    \Higgs_{G,\cL}\ar[d]\ar[r]& \Higgs_{G_1,\cL}\times\cdots\times \Higgs_{G_s,\cL}\ar[d]\\
    \Higgs_{G_{ab},\cL}\ar[r]& (\Higgs_{G_{ab},\cL})^s
    }
\end{equation}
is Cartesian. The closed embedding $G\hookrightarrow G_1\times\cdots \times G_s$ induces an isomorphism of adjoint groups $G/Z_G\cong (G_1/Z_{G_1})\times\cdots \times (G_s/Z_{G_s})$. By \Cref{L:induced morphism semistable}\eqref{L:induced morphism semistable 1}, the inverse image of the semistable locus of $\Higgs_{G_1,\cL}\times\cdots\times \Higgs_{G_s,\cL}$ via the top arrow in \eqref{E:Cartesian product higgs} coincides with the semistable locus of $\Higgs_{G,\cL}$, and hence we get a Cartesian diagram
\begin{equation}\label{E:Cartesian product higgs 2}
    \xymatrix{
    \Higgs_{G,\cL}^{ss}\ar[d]\ar[r]& \Higgs_{G_1,\cL}^{ss}\times\cdots\times \Higgs_{G_s,\cL}^{ss}\ar[d]\\
    \Higgs_{G_{ab},\cL}\ar[r]& (\Higgs_{G_{ab},\cL})^s
    }
\end{equation}
To conclude the proof of the lemma, it suffices to show that diagram \eqref{E:Cartesian product higgs 2} remains Cartesian after taking adequate moduli spaces. For this, we may extend scalars and assume that $k$ is algebraically closed. As in the proof of \Cref{prop: moduli space commutes with base change on abelianization}, the choice of a $k$-point $c \in C(k)$ induces a splitting $\MHiggs_{G_{ab}, \cL} \to \Higgs_{G_{ab}, \cL}$ of the adquate moduli space morphism. Applying a further base-change of diagram \eqref{E:Cartesian product higgs 2} by $\MHiggs_{G_{ab}, \cL} \to \Higgs_{G_{ab}, \cL}$, we obtain the following Cartesian diagram of stacks
\begin{equation}\label{E:Cartesian product higgs 3}
    \xymatrix{
     \MHiggs_{G_{ab}, \cL} \times_{\Higgs_{G_{ab}, \cL}} \Higgs_{G,\cL}^{ss} \ar[d]\ar[r]& \Higgs_{G_1,\cL}^{ss}\times\cdots\times \Higgs_{G_s,\cL}^{ss}\ar[d]\\
    \MHiggs_{G_{ab},\cL}\ar[r]& (\Higgs_{G_{ab},\cL})^s
    }
\end{equation}
By \Cref{prop: moduli space commutes with base change on abelianization}, the diagram \eqref{E:Cartesian product higgs 3} remains Cartesian after taking moduli spaces. Since the adequate moduli space of $\MHiggs_{G_{ab}, \cL} \times_{\Higgs_{G_{ab}, \cL}} \Higgs_{G,\cL}^{ss}$ is equal to $\MHiggs_{G, \cL}$ (by 
\Cref{prop: moduli space commutes with base change on abelianization}), we conclude that the following diagram is Cartesian, as desired
\[
    \xymatrix{
     \MHiggs_{G, \cL} \ar[d]\ar[r]& \MHiggs_{G_1,\cL}\times\cdots\times \MHiggs_{G_s,\cL}\ar[d]\\
    \MHiggs_{G_{ab},\cL}\ar[r]& (\MHiggs_{G_{ab},\cL})^s.
    }
\]
\end{proof}

\begin{lemma}[Biswas--Hoffmann's trick]\label{L:reduction to D(G) almost-simple}Let $G$ be a reductive group over a field $k$ and $d\in\pi_1(G)$. Then, there exists a morphism $\pi: \widehat{G}\to G$ of reductive group schemes such that the followings hold
\begin{enumerate}[(i)]
    \item the derived subgroup $D_{\widehat{G}}$ of $\widehat{G}$ is simply-connected;
    \item the element $d\in\pi_1(G)$ is contained in the image of the induced morphism of fundamental groups $\pi_1(\widehat{G})\to \pi_1(G)$;
\end{enumerate}
\end{lemma}

\begin{proof}The first two points are the content of \cite[Lem. 5.3.2]{biswas-hoffmann-line-bundles}. The third point follows from the construction of $\widehat{G}$.
\end{proof}
\end{section}

\bibliographystyle{alpha}
\footnotesize{\bibliography{decomposition-g-higgs.bib}}

\end{document}